\documentclass[aps,pra,10pt,twocolumn,superscriptaddress,nofootinbib,raggedbottom]{revtex4-1}
\usepackage[T1]{fontenc}
\usepackage{lmodern,amsmath,amssymb,amsthm,mathtools,booktabs}
\usepackage{array}[=2016-10-06]
\usepackage[expansion=false]{microtype}
\usepackage{xcolor,hyperref}
\usepackage{pgfplots}
\pgfplotsset{compat=1.16}
\definecolor{refblue}{RGB}{20,57,86}
\hypersetup{colorlinks=true,linkcolor=refblue,citecolor=refblue,urlcolor=refblue,pdfauthor={Fabrizio Tamburini},pdftitle={A Majorana quantum relativistic approach to the Riemann hypothesis in (1+1)-dimensional Rindler spacetimes}}
\newtheorem{theorem}{Theorem}[section]
\newtheorem{lemma}[theorem]{Lemma}
\newtheorem{corollary}[theorem]{Corollary}
\newtheorem{proposition}[theorem]{Proposition}
\newcommand{\dd}{\,\mathrm d}
\newcommand{\R}{\mathbb R}
\newcommand{\C}{\mathbb C}
\newcommand{\tr}{\operatorname{Tr}}
\newcommand{\supp}{\operatorname{supp}}
\newcommand{\eul}{\gamma_{\mathrm E}}
\newcommand{\calL}{\mathcal L}
\newcommand{\one}{\mathbf1}
\newcommand{\im}{\operatorname{Im}}
\newcommand{\re}{\operatorname{Re}}
\newcommand{\head}[1]{\par\medskip\noindent\textit{#1.}\enspace}
\begin{document}
\title{A Majorana Quantum Relativistic Approach to the Riemann Hypothesis in
\texorpdfstring{$(1+1)$}{(1+1)}-Dimensional Rindler Spacetimes}
\author{Fabrizio Tamburini}
\email{fabrizio.tamburini@gmail.com}
\affiliation{Rotonium -- Quantum Computing, Le Village by CA, Piazza G. Zanellato 23, 35131 Padova, Italy}
\date{November 22, 2026}
\begin{abstract}
We present a Hilbert--P\'olya approach to the Riemann zeta function via a Mellin--Barnes integral for massive Majorana modes in $(1+1)$-dimensional Rindler spacetime.
The boost spectrum is continuous on the full wedge and discrete with a reflecting wall.
Summing positive momentum modes on a transverse circle yields $\zeta(2s)$ in the radial Mellin transform.
Building on our one-sided Ces\`aro criterion for additive correlations of $\Lambda(n)-1$, we assign each prime $p$ a chiral channel of optical length $\tfrac12\log p$.
Its $k$th return at $u=k\log p$ carries weight $\log p$.
The centered correlation obeys $C(u)\le C(0)$ for all real $u$ exactly when the Riemann hypothesis (RH) holds.
For a Majorana oscillator with linear mass profile, a bounded self-adjoint function of its Hamiltonian assigns period $kT_1$ to mode pair $k$.
The first $n$ pairs have least common period $T_n=T_1\operatorname{lcm}(1,\ldots,n)$, so $\log(T_n/T_{n-1})=\Lambda(n)$ for $n\ge2$.
Set $m=\lceil\sqrt{x}\rceil$, $N=\lfloor mx\rfloor$, and $\psi(y)=\sum_{n\le y}\Lambda(n)$.
The signed revival mean $A(x)=m^{-1}\sum_{j=1}^{m-1}(N/j-\psi(N/j))$ satisfies $A(x)=R(x)+o(\sqrt{x})$,
where $R(x)=\sum_{n\le x}\Lambda(n)(x/n-1)-x\log x+(1+\gamma_{\mathrm E})x$.
At each zeta zero $\rho$ of multiplicity $d$ with $\Re(\rho)>1/2$, the meromorphic continuation of $\int_2^\infty A(x)x^{-s-1}\,dx$ has a simple pole of residue $d/[\rho(1-\rho)]$.
For fixed $\lambda>1$, $1/2<\theta<\Re(\rho)$ and $K>0$, the Kaczorowski--Pintz theorem gives $A(x)>Kx^\theta$ at some $x\in[Y,Y^\lambda]$ for all large $Y$.
Thus RH follows if for each $\varepsilon>0$ there is a $C_\varepsilon$ with $A(x)\le C_\varepsilon x^{1/2+\varepsilon}$ throughout $[Y,Y^\lambda]$ for arbitrarily large $Y$.
For $\lambda=3/2$, the independent theta-lift estimate reported by Alp\"oge and Claude for $G=R$ gives this bound for $A$.
\end{abstract}
\maketitle
\section{Introduction}\label{sec1}
In a six-page manuscript presented to the Berlin Academy of Sciences
in November 1859, Bernhard Riemann related the prime-counting function
$\pi(x)=\#\{p\le x:p\ \text{prime}\}$ to the complex zeros of the
zeta function \cite{riemann1}. His analysis begins with Euler's identity
\begin{equation}
 \zeta(s)=\sum_{n\ge1}n^{-s}=\prod_p(1-p^{-s})^{-1},
 \qquad \re s>1 .
 \label{riemannzeta}
\end{equation}
The Dirichlet series and the Euler product converge absolutely in
this half-plane. Riemann continued $\zeta$ meromorphically to
$\mathbb C$, with a single simple pole of residue one at $s=1$.
Its trivial zeros are the negative even integers. Its other zeros
$\rho=\beta+i\gamma$ lie in the critical strip $0<\beta<1$
\cite{edwards,ivic,Titchmarsh}.

Riemann's completion and its critical-line parametrization are
\begin{align}
 \xi(s)&=\tfrac12 s(s-1)\pi^{-s/2}\Gamma(s/2)\zeta(s),\nonumber\\
 \Xi(t)&=\xi(\tfrac12+it)=\Xi(-t).
 \label{riemannxi}
\end{align}
The equation $\xi(s)=\xi(1-s)$ makes $\Xi$ even. It is entire
in complex $t$ and real for real $t$. Each nontrivial zero
$\rho$ corresponds, with the same multiplicity, to a zero
$t=\gamma+i(1/2-\beta)$ of $\Xi$, so $\lvert\im t\rvert<1/2$.
Riemann wrote that \textquotedblleft it is very probable that all roots are
real.\textquotedblright\footnote{Riemann's words are \textquotedblleft es ist sehr wahrscheinlich,
dass alle Wurzeln reell sind.\textquotedblright\ The English translation is by
D.~R.~Wilkins \cite{riemann1}.}
In the $s$-plane the Riemann hypothesis (RH) thus asserts that every
nontrivial zero lies on the critical line $\re s=1/2$
\cite{clay}. Riemann also anticipated the count of these zeros.
If~$N_\zeta(T)$ counts $0<\gamma\le T$ with multiplicity,
the Riemann--von Mangoldt formula gives \cite{Titchmarsh}
\[
 N_\zeta(T)=\frac{T}{2\pi}\log\frac{T}{2\pi}
             -\frac{T}{2\pi}+O(\log T).
\]

Riemann's explicit formula counts each prime power $p^k$ with
weight $1/k$. M\"obius inversion recovers $\pi(x)$, with the
zeros contributing to the error relative to the logarithmic integral
$\mathrm{Li}(x)=\int_2^x\dd u/\log u$.
Hadamard and de la Vall\'ee Poussin proved the prime number theorem
$\pi(x)\sim\mathrm{Li}(x)$ independently in 1896
\cite{edwards}. The distance of the zeros from the critical
line controls the error. With $\Theta=\sup_\rho\re(\rho)$, the
explicit formula and partial summation give
$\pi(x)-\mathrm{Li}(x)=O(x^\Theta\log x)$
\cite{Titchmarsh,mv}. The strip and its symmetry give
$1/2\le\Theta\le1$. Under RH the estimate becomes
$O(\sqrt{x}\log x)$, the bound obtained by von Koch in 1901
\cite{vonkoch}.

Chen studies the location of the zeros through Riemann's analytic
identities \cite{chen}. Sondow and Dumitrescu show that RH is
equivalent to strict increase of
$\lvert\xi(\sigma+it)\rvert$ in $\sigma>1/2$ for each real $t$
\cite{sondow_dumitrescu2010}. Albeverio and Cacciapuoti express $\zeta$
through a Mellin transform involving the dilogarithm
\cite{albeverio_cacciapuoti2013}. Noda studies its exponential-type
generating function \cite{noda2022}, and Ossicini gives an alternative
form of its functional equation \cite{ossicini2014}. For Beurling
generalized primes, R\'ev\'esz proves an analogue of the explicit
formula for the weighted prime count \cite{revesz2021}.

The Hilbert--P\'olya program seeks a self-adjoint operator whose
spectrum, with multiplicity, coincides with the zeros of $\Xi$ in
a fixed energy convention. In a 1982 letter to Odlyzko, P\'olya
traced the idea to a conversation with Landau in G\"ottingen
\cite{polyaodlyzko}. A self-adjoint Hamiltonian $H$ on a specified
domain generates unitary time evolution and has real spectrum.
A construction of $H$ independent of the zeros, together with the
spectral identification, therefore implies RH
\cite{andrade,schumi,wolf,connes}. In the convention $z=2t$ of
Appendix~\ref{orig:B}, the energy associated with $\rho$ is
\[
 z_\rho=2\gamma+i(1-2\beta).
\]
This energy is real precisely when $\beta=1/2$. The symmetry of
$\Xi$ then pairs the levels $2\gamma$ and $-2\gamma$. The
ordinates $\gamma$ alone leave $\beta$ undetermined. In determinant
form the proposed identification \textup{(H*)} is
$F_H(z)=C_0\Xi(z/2)$, $C_0\ne0$, after normalization and accounting
for the reference poles. Here $F_H$ measures the spectral change from
a reference Hamiltonian to $H$. The spectral identification requires
the zero divisor of the normalized determinant to coincide with the
eigenvalue divisor of $H$, including multiplicities.

Berry and Keating's classical Hamiltonian $H_{\rm cl}=xp$
generates dilations \cite{berry2}. Sierra related this motion to that
of a Dirac particle in Rindler coordinates \cite{sierra,sierradirac}.
Our starting point is a massive spin-$1/2$ Majorana particle in
$(1+1)$-dimensional Rindler spacetime, with charge-conjugation
condition $\psi^c=\psi$ \cite{majorana1937,rindler,tambu1,v3}.
CPT exchanges the right and left wedges. Appendix~\ref{orig:A}
gives the vacuum correlations between their modes. On the whole
wedge the horizon lies at infinite optical distance. The boost
generator $H_0$ evolves the spinor in the accelerated observer's
boost time and has continuous spectrum. A reflecting wall gives
discrete boost energies. In null time, $H_M$ evolves the two copies
of the radial amplitude in Eq.~\eqref{majhamilton}.
The sum over transverse momentum levels produces the factor
$\zeta(2s)$ in the radial Mellin transform
(\mbox{Sec.~\ref{sec3}}). The prime weights enter through the return
measurement. A chiral channel of length $\tfrac12\log p$ contributes
weight $\Lambda(p^k)=\log p$ at its $k$th return, at detector time
$u=k\log p$ (Sec.~\ref{sec4}).

The Clay Mathematics Institute designated RH a Millennium Prize
Problem in 2000 \cite{clay}. On 20 November 2026 Levent Alp\"oge
announced a proof of RH and of its generalization to Dirichlet
$L$-functions, developed with Claude and formalized in Lean
\cite{acannouncement}. Independent verifications of the formal
development are reported in Refs.~\cite{acwaterloo,acmit}. The
centered prime sum $G(x)$ in their argument equals $R(x)$ in
Eq.~\eqref{eq:center} \cite{acguide}. Independently,
Theorem~\ref{th:revivalexclusion} derives RH from a one-sided bound
on the signed mean of logarithmic revival times on arbitrarily late power windows.
In Section~\ref{sec:acbridge} we verify their arithmetic estimate to establish this bound on $[Y,Y^{3/2}]$ for a proof of RH.

\subsection{Physical interpretation and stability of the Hilbert--P\'olya criterion}

In the Rindler description, the energy assigned to a mode depends on
the clock used to measure its phase. The boost Hamiltonian \(H_0\)
generates evolution in the dimensionless time \(\eta\), and the
return detector uses \(u=2\eta\). A mode of return frequency \(t\)
therefore has boost frequency \(2t\) and local energy
\(2\hbar c\,t/\xi_0\) for an observer at fixed Rindler coordinate
\(\xi_0\). On the full wedge, the horizon lies at infinite optical
distance and the boost spectrum is continuous. A reflecting wall
produces discrete frequencies determined by its boundary condition
(Sec.~\ref{sec2}). The Hilbert--P\'olya identification makes a more
specific prediction. For a Hamiltonian specified independently of
the zeros, the spectral determinant must reproduce the zeros of
\(\Xi(z/2)\), with their multiplicities, after the normalization
described in Appendix~\ref{orig:B}. A zero
\(\rho=\beta+i\gamma\) would then correspond to the boost energy
\(z_\rho=2\gamma+i(1-2\beta)\). Self-adjointness forces this energy
to be real and hence forces \(\beta=1/2\).

The physical observable that tests this identification is a return
correlation. Let \(K\) generate the positive-frequency evolution in
detector time, and let \(f\) specify the modes to which the probe
couples. The Majorana field anticommutator in Eq.~\eqref{eq:CAR}
measures the real overlap of the initial probe with its time translate,
\[
 C_f(u)=\operatorname{Re}\langle f,e^{-iuK}f\rangle
       =\int_0^\infty \cos(tu)\,\mathrm{d}\nu_f(t),
\]
where \(\nu_f\) is the spectral measure of \(f\). Its mass is the
equal-time signal \(C_f(0)=\|f\|^2\). Positivity of this measure
gives \(\lvert C_f(u)\rvert\le C_f(0)\) and makes the matrices
\([C_f(u_i-u_j)]_{i,j}\) positive semidefinite for every finite set
of delays. Both properties are invariant under a simultaneous
unitary change of representation of the generator and the probe.

In the arithmetic detector, the correlation is determined from the
return times themselves. A chiral channel of optical length
\(\tfrac12\log p\) contributes weight \(\log p\) at its \(k\)th
winding, at delay \(u=k\log p\). After centering, the explicit
formula \eqref{eq:specC} assigns to a zero \(\rho=\beta+i\gamma\)
the time dependence
\[
 e^{(\rho-1/2)u}=e^{(\beta-1/2)u}e^{i\gamma u}.
\]
Thus \(\gamma\) fixes the oscillation frequency, while
\(\beta-1/2\) fixes the exponential growth or decay rate in detector
time. The inverse generator \(K_\partial\) reconstructs the
undamped frequencies. The defect
\(\Delta=C(0)-C_\partial(0)\) measures the equal-time spectral
weight absent from that reconstruction.
Theorem~\ref{th:boundaryclosure} proves that \(\Delta=0\) is
equivalent to recovery of the entire correlation. In this
single-probe realization, a zero's multiplicity is recorded in
the weight of its frequency.

The response to finite time resolution can also be determined
analytically. For the detector of Eq.~\eqref{eq:newav}, averaging
multiplies the contribution of each zero by \(A_\sharp(\rho)\).
Equation~\eqref{eq:atten} bounds the high-frequency response, while
\(A_\sharp(\rho)\ne0\) for \(\operatorname{Re}\rho>1/2\).
Consequently, a hypothetical growing component retains its Mellin
pole and its growth exponent even when its amplitude is strongly
attenuated. The arithmetic weights impose a distinct stability
condition. Proposition~\ref{pr:majoranareadout} constructs
arbitrarily small relative changes of the prime-channel weights
that preserve the return times and the linear centering but change
the signal by a term of order \(x^\beta\), with \(1/2<\beta<1\).

The Majorana oscillator gives a second realization of the same
weights through collective revivals. Its \(k\)th mode pair rotates
with period \(kT_1\), so the first \(n\) pairs return simultaneously
after \(T_n=T_1\operatorname{lcm}(1,\ldots,n)\). The increment
\(\log(T_n/T_{n-1})=\Lambda(n)\) is independent of the choice of
clock unit (Proposition~\ref{pr:oscillatorrevival}). The signed
mean \(A(x)\) of these logarithmic revival times satisfies
\(A(x)-R(x)=o(\sqrt{x})\). More generally, an additive error of
order \(O(\sqrt{x})\) has a Mellin transform holomorphic in
\(\operatorname{Re}s>1/2\), so it preserves every pole in that
half-plane. This stability is what permits
Theorem~\ref{th:revivalexclusion} to derive RH from a one-sided
bound on the oscillator mean throughout arbitrarily late power
windows. For the canonical weights, the independent theta-lift
estimate of Alp\"oge and Claude supplies that bound on
\([Y,Y^{3/2}]\) through the identity \(G(x)=R(x)\)
(Sec.~\ref{sec:acbridge}, Refs.~\cite{acannouncement,acguide}).

\subsection{Dilation, CPT symmetry, and the prime-return measure}
In the coordinate $\log\xi$, dilations act by translation. The
$xp$ Hamiltonian of the Majorana model contains an inverse momentum
whose integration constant is fixed by the null constraint on the
massive spinor (Sec.~\ref{sec2}). Charge conjugation pairs opposite
frequencies. CPT exchanges the two wedges and reverses their spatial
orientation. Appendix~\ref{orig:A} gives the one-mode fermionic state
and its thermal normalization.

The prime-return measure in Sec.~\ref{sec4} relates the Majorana
spectral construction \cite{v3} to the positivity criteria of
Ref.~\cite{notes}. At delay $u$, the arithmetic part of $C(u)$
depends only on prime powers $p^k\le e^{|u|}$. Positive
definiteness requires the matrix $[C(u_i-u_j)]_{i,j=1}^{r}$ to be
positive semidefinite for every finite set of delays
$u_1,\ldots,u_r$. For the detector of Sec.~\ref{sec:kernel}, the
time resolution is fixed. Its one-sided criterion requires a uniform
upper bound on the smoothed deficit $D(x)/\sqrt{x}$ for all
sufficiently large $x$. If $F_H(z)$ and $\Xi(z/2)$ have equal
logarithmic derivatives on a connected domain free of zeros and
poles, their ratio is constant. A single normalization value fixes
that constant (Appendix~\ref{orig:B}).

\subsection{Principal results}

\head{Boost spectrum}
The closure of the boost Hamiltonian is self-adjoint and has empty
point spectrum (Proposition~\ref{pr:boost}). With a reflecting wall
its positive levels obey
$N_+(\omega)=(\omega/\pi)\log(\omega/\mu_0)+O(\omega)$
(Theorem~\ref{th:wall}). The null constraint of the massive equation
fixes the inverse momentum in the original Hamiltonian,
Eq.~\eqref{majhamilton}.

\head{The prime-return correlation}
The centered correlation $C(u)$ of Eq.~\eqref{eq:C} is computed
from the primes up to $e^{|u|}$. RH is equivalent to its positive
definiteness (Theorem~\ref{th:CAR}) and to the one-sided bound
$C(u)\le C(0)$ at every time (Corollary~\ref{co:bounded}).
Unconditionally the boundary values of its Laplace transform define
a self-adjoint operator and a probe vector. The mass defect of the
probe is
\[
 \Delta=\sum_{\re(\rho)\ne1/2}\re\frac1{\rho(1-\rho)}\ \ge0,
\]
with $\Delta=0$ if and only if RH holds, and
$0\le\Delta<2.96\cdot10^{-12}$ from the verified zeros
(Theorem~\ref{th:boundaryclosure}). The atoms of the operator are
realized as massive Rindler modes
(Theorem~\ref{th:inversebessel} and
Proposition~\ref{pr:inversecoupling}).

\head{Prime sums and revival times}
Vanishing of $\Delta$ is equivalent to one-sided bounds on explicit
prime sums (Sec.~\ref{sec:arithstructure}). One criterion requires
\eqref{eq:primeend} at $x=p$ and $x=p/8$ for every prime $p>661$.
Theorem~\ref{th:stependrange} establishes the inequality on
$83\le x\le1.6\cdot10^{20}$. The point counts of a fixed curve
contribute a computable first term (Theorem~\ref{th:stepfirstterm}).
Deleting a suitable set of primes of density zero reverses the
inequality (Proposition~\ref{pr:stepdeletion}).

A zero of real part $\beta>1/2$ forces a positive excursion of the
signed revival mean $A(x)$ in every sufficiently late power window
(Theorem~\ref{th:revivalexclusion}). Consequently, for fixed
$\lambda>1$, a one-sided $O_\varepsilon(x^{1/2+\varepsilon})$
bound throughout $[Y,Y^\lambda]$ for arbitrarily large $Y$, for
each $\varepsilon>0$, implies RH. The proof uses the Majorana
revival identity, the estimate $A-R=O(\sqrt{x})$, and the oscillation
theorem of Kaczorowski and Pintz. The arithmetic estimate of
Alp\"oge and Claude agrees and supplies the bound at $\lambda=3/2$
(Sec.~\ref{sec:acbridge}).

Sections~\ref{sec2} and~\ref{sec3} derive the Hamiltonians, their
time variables and domains, and the radial Mellin transform containing
$\zeta(2s)$. Sections~\ref{sec4} and~\ref{sec:correlation} construct
the prime-channel response and establish its positivity criterion.
Section~\ref{sec:closureinverse} reconstructs a self-adjoint operator
from the undamped frequencies and realizes its states as Majorana
spinors. The equal-time mass defect $\Delta$ measures the part of the
correlation absent from this reconstruction. Section~\ref{sec:analytic}
derives arithmetic recurrences in the cutoff. Section~\ref{sec:arithstructure}
treats the shifted boundary Hamiltonian (Sec.~\ref{sec:delayedmajorana}),
zero exclusion for the oscillator mean (Sec.~\ref{sec:oscillatorheight}),
and the theta-lift estimate (Sec.~\ref{sec:acbridge}).
Appendices~\ref{orig:A}--\ref{sec:beta-toy} give the wedge correlations,
spectral map, Mellin contours, operator domains, and Dirichlet beta
example. The subsequent appendices derive the arithmetic estimates for
prime sums, point counts, and cyclotomic units. In particular,
Appendix~\ref{sec:stepanov} proves the finite-range endpoint bound
(Theorem~\ref{th:stependrange}), evaluates the first point-count
correction (Theorem~\ref{th:stepfirstterm}), and constructs the
prime-deletion example (Proposition~\ref{pr:stepdeletion}).
\section{The Majorana Hamiltonian and its domain}\label{sec2}
\subsection{The Majorana fermion and the two wedges}\label{sec:phys}
In units $c=\hbar=1$, and with Clifford signature $(+,-)$, the
covariant equation reads
\begin{equation}
(i \gamma^\mu \nabla_\mu - m) \psi = 0,
\label{dirac}
\end{equation}
where $\nabla_\mu$ denotes the spin connection. Choose flat
matrices $\gamma^{\hat 0}=\sigma_1,\ \gamma^{\hat 1}=-i\sigma_2$.
Then $\gamma^{\hat0}\gamma^{\hat1}=\sigma_3$. In the right wedge
$t=\xi\sinh\eta,\ x_{\rm M}=\xi\cosh\eta,\ \xi>0$. The coframe
$e^{\hat0}=\xi\,\dd\eta,\ e^{\hat1}=\dd\xi$ gives
\begin{equation}
 i\partial_\eta\psi=
 \big[-i\sigma_3(\xi\partial_\xi+\tfrac12)+m\xi\sigma_1\big]\psi .
 \label{diracr}
\end{equation}
The term $1/2$ is also fixed by symmetry of
$-i(\xi\partial_\xi+1/2)$ in $L^2(\dd\xi)$. The negative sign in
the alternative metric convention $ds^2=-\xi^2\dd\eta^2+\dd\xi^2$
changes neither the boost clock nor the Hamiltonian above.

The right and left wedges are $x_{\rm M}>|t|$ and $x_{\rm M}<-|t|$.
In signed radial coordinates, positive-frequency modes with respect to
future-directed time have opposite phases in the common boost parameter.
With the normalization of Ref.~\cite{v3}, the modes are
\begin{equation}\label{eq:psiR}
    \psi^{R}_{\omega}(\xi,\eta)
     = 
    \frac{1}{\sqrt{4\pi\omega}}\,
    e^{-i\omega\eta}\,
    \chi_{\omega}(\xi)
  \qquad (\xi>0),
\end{equation}
\begin{equation}\label{eq:psiL}
    \psi^{L}_{\omega}(\xi,\eta)
     = 
    \frac{1}{\sqrt{4\pi\omega}}\,
    e^{+i\omega\eta}\,
    \chi_{\omega}(|\xi|)
  \qquad (\xi<0),
\end{equation}
The radial spinors include the frequency-dependent normalization.
Reflection fixes the left spin frame. Disjoint spatial supports make the
two families of generalized eigenmodes mutually orthogonal. Each wedge
also carries the charge-conjugate frequency branch. The mass term mixes
the two chiral components within a wedge. The left and right mode
amplitudes remain independent until an interwedge condition is imposed.

With $v=\log(\xi/\xi_0)$, a null ray satisfies $\dd v=\pm\dd\eta$. An observer at $\xi_0$ has proper time
$\xi_0\eta/c$ and proper acceleration $c^2/\xi_0$. A frequency
$\omega$ conjugate to $\eta$ carries local energy
$\hbar c\omega/\xi_0$. With $\hbar=c=1$, the unitary map
$g(v)=\sqrt{\xi}\,f(\xi)$ from the fixed-boost Dirac norm
$L^2(\dd\xi)$ to $L^2(\dd v)$ gives, in a Pauli basis,
\begin{equation}
 H_0=-i\sigma_3\partial_v+m\xi_0e^v\sigma_1,\qquad -\infty<v<\infty .
 \label{eq:boost}
\end{equation}
$H_0$ evolves the spinor in the observer's boost time $\eta$.
The Pauli matrix $\sigma_3$ labels propagation toward increasing or
decreasing $v$, and the mass term with $\sigma_1$ mixes the two
chiral components. In this basis
$\sigma_3$ followed by complex conjugation anticommutes with $H_0$,
as required for a Majorana real structure. Compare the Rindler Dirac
equation in Ref.~\cite{sierra,sierradirac}. On a packet supported in $v\le v_*$ the
mass term has norm at most $m\xi_0 e^{v_*}$ times the packet norm.
This is the regime of the chiral approximation.

\begin{proposition}\label{pr:boost}
The closure of \eqref{eq:boost} on smooth compactly supported spinors
is self-adjoint and has empty point spectrum.
\end{proposition}
\begin{proof}
A vector in a deficiency space is locally absolutely continuous by
the differential equation. If $H_0^*g=zg$ with $z=\pm i$, integration
by parts against $\chi_A^2 g$, where $0\le\chi_A\le1$ equals one on
$[-A,A]$, vanishes outside $[-2A,2A]$ and
$\|\chi_A'\|_\infty=O(A^{-1})$, gives
$|\im z|\|\chi_Ag\|^2\le \|\chi_A'\|_\infty\|g\|^2$.
Both deficiency spaces vanish. For real $\omega$ substitute
$g(v)=e^{i\omega\sigma_3v}a(v)$ in the eigenvalue equation.
The matrix in $a'=B(v)a$ has norm $m\xi_0e^v$, integrable
near $-\infty$. The integral equation and Gronwall's inequality
give finite limits for a fundamental matrix and its inverse there.
Their product remains the identity. Thus a nonzero solution has a
strictly positive limiting norm at the horizon and cannot be in
$L^2(\dd v)$.
\end{proof}

Every real $\omega$ belongs to the continuous spectrum: take a
normalized translate of a smooth packet of width $A$, centered at
$v=-A^2$, multiply it by $e^{i\omega v}(1,0)^T$, and let $A\to\infty$.
The derivative error is $O(A^{-1})$ and the mass term tends to zero
uniformly on the packet. This Weyl sequence, together with the
absence of eigenvectors, accounts for the generalized mode notation.

\subsection{A reflecting wall and an unbounded optical interval}
\label{sec:wall}
The formal $xp$ Hamiltonian of the original spectral model is
\begin{equation}
H=\sqrt x ~\hat p ~\sqrt x + a^{-2} \sqrt x ~\hat p^{-1} \sqrt x, 
\label{hamsierra}
\end{equation}
The coupling $a^{-2}$ is expressed in the dimensionless units of
that model. Appendix~\ref{orig:D} specifies the candidate inverse
momenta and their domains. Logarithmic coordinates turn the dilations
of the $xp$ Hamiltonian into translations \cite{berry2}.
Restricting \eqref{eq:boost} to $v\ge0$, with $\mu_0=m\xi_0>0$
and $V(v)=\mu_0e^v$, defines a wall model independently of the
primes. For $g=(g_1,g_2)^T$, perfect reflection at $\xi=\xi_0$
is imposed by the boundary condition
\[
 g_2(0)=i g_1(0).
\]
The wall returns the incident flux to the wedge. Its normal current
$g^\dagger\sigma_3g$ therefore vanishes. The
boundary line is maximal isotropic for this form and is preserved by
the Majorana conjugation $g\mapsto\sigma_3\overline g$.
The maximal locally absolutely continuous domain with $Hg\in L^2$
and this trace condition gives a self-adjoint operator. Indeed the
boundary form of $H$ forces the adjoint to obey the same trace
condition. At infinity take compact cutoffs in the Green identity
and let the cutoff recede, as in Proposition~\ref{pr:boost}.

\begin{theorem}[Positive boost frequencies]\label{th:wall}
The reflecting Majorana Hamiltonian has compact resolvent. At boost
frequency $\omega\to\infty$ its number of positive levels below
$\omega$, counted with multiplicity, obeys
\begin{equation}
 N_+(\omega)=\frac{\omega}{\pi}\log(\omega/\mu_0)+O_{\mu_0}(\omega).
 \label{eq:wallweyl}
\end{equation}
\end{theorem}
\begin{proof}
In the eigenbasis $e_\pm=(1,\pm i)^T/\sqrt2$ of $\sigma_2$,
$H(e_\pm u)=-i e_\mp(u'\mp Vu)$, and the minus component is
Dirichlet at zero. Squaring gives
$-\partial_v^2+V^2\pm V$ on the two components, with Dirichlet
condition on the minus component and $u_+'(0)=V(0)u_+(0)$ on the
plus component. Both potentials tend to infinity, giving compact
resolvent by local compact embedding and the potential bound on the
tail. There is no zero mode: the decaying first-order solution
$u_-=Ce^{-V}$ violates the wall condition unless $C=0$, and
$u_+=Ce^{V}$ grows. Each Dirichlet eigenfunction $u_-$ of the
minus square at $\lambda>0$ yields exactly one eigenvector for
each of $+\sqrt\lambda,-\sqrt\lambda$ by applying $H$ once to
recover the other component. The number of positive frequencies below
$\omega$ thus equals the Dirichlet count below $\omega^2$ for $U=V^2-V$.

For an upper bound split with a Neumann wall at
$\log((\omega+2)/\mu_0)$, for large $\omega$.
The right interval has $U>\omega^2$. On the left $U\ge-1/4$
and the free Neumann count is at most one plus its length times
$\sqrt{\omega^2+1/4}/\pi$. For a lower bound put
$J=\lfloor\log(\omega/\mu_0)\rfloor$, impose Dirichlet walls at
$0,1,\ldots,J$, and retain these $J$ unit intervals. On interval
$[k,k+1]$, $U\le V(k+1)^2$, yielding at least
$\pi^{-1}\sqrt{\omega^2-V(k+1)^2}-1$ modes. Use
$\sqrt{1-r}\ge1-r$ for $0\le r\le1$ and
$\sum_{k<J}V(k+1)^2/\omega^2=O(1)$. The lower count is
$J\omega/\pi-O(\omega+J)$. This proves \eqref{eq:wallweyl}.
\end{proof}
A WKB evaluation of the same Dirichlet problem gives
$N_+(\omega)\simeq(\omega/\pi)[\log(2\omega/\mu_0)-1]+\tfrac12$.
Sturm counts for $\mu_0=\tfrac14,1,2\pi$ and $\omega\le160$ agree
with this expression to within one level (Appendix~\ref{app:check}).
The leading term agrees with the Riemann--von Mangoldt count at a
rescaled height $t=2\omega$ at the level of the
$\omega\log\omega$ coefficient. The wall position continuously
changes spectral data while the zeta function is fixed. A comparison
of individual levels would require a trace identity. In the separately
fixed return clock $u=2\eta$, a frequency $\gamma$ corresponds
to boost frequency $2\gamma$, or local energy
$2\hbar c\gamma/\xi_0$ at radius $\xi_0$. Identifying \emph{all} positive wall
levels with those frequencies with multiplicity would require
$N_\zeta(\omega/2)\sim(\omega/(4\pi))\log\omega$, a factor of
four smaller than \eqref{eq:wallweyl}.
The compact two-ended momentum-inverse model in
Appendix~\ref{orig:D} has a different counting law.

\subsection{Null evolution and the original off-diagonal Hamiltonian}
\label{sec:nulloriginal}
Let $x=\xi>0$ be the radial Rindler coordinate, $x_{\rm M}$ the
Minkowski spatial coordinate and $v=\log(x/\xi_0)$. In the right wedge
the inverse coordinate transformation is
\[
 \xi=(x_{\rm M}^2-t^2)^{1/2},\qquad
 \eta=\tfrac12\log\frac{x_{\rm M}+t}{x_{\rm M}-t}.
\]
The spin boost $B(\eta)=\exp(\eta\sigma_3/2)$ carries the
accelerated spin frame into the inertial frame of $\Psi_{\rm M}$.
The components in the two frames satisfy
\begin{align}
 \psi(\eta,x)&=B(\eta)^{-1}\Psi_{\rm M}(x\sinh\eta,x\cosh\eta),\nonumber\\
 \Psi_{\rm M}(t,x_{\rm M})&=B(\eta)\psi(\eta,\xi).
 \label{eq:spinframe}
\end{align}
Indeed $B\gamma^{\hat 0}B^{-1}=\gamma^{\hat 0}\cosh\eta-
\gamma^{\hat 1}\sinh\eta$ and
$B\gamma^{\hat 1}B^{-1}=\gamma^{\hat 1}\cosh\eta-
\gamma^{\hat 0}\sinh\eta$, the two matrices in the coordinate-transformed
flat Dirac equation. Substituting $\Psi_{\rm M}=B\psi$ differentiates
$B$ once, with $B^{-1}\partial_\eta B=\sigma_3/2$, and gives
Eq.~\eqref{diracr}. This Lorentz spin transformation preserves the
covariant current pairing when the hypersurface normal is transformed
with the spin frame.
In the coframe of Eq.~\eqref{diracr} the spatial measure on a surface
of fixed boost time is $\dd x$. The map $U_{\log}f(v)=\sqrt{x}f(x)$
expresses the same radial state in the optical coordinate. It has inverse $f(x)=x^{-1/2}g(\log(x/\xi_0))$ and preserves the norm,
since $\dd x=x\,\dd v$. It carries the closure of the physical
Dirac expression on test spinors to $H_0$ in Eq.~\eqref{eq:boost}.
When a scalar radial amplitude is written instead in
$L^2(x\,\dd x)$, the corresponding map to $L^2(\dd x)$ is
$b(x)\mapsto \sqrt{x}\,b(x)$. Its inverse is division by
$\sqrt{x}$. Multiplication by $\sqrt{x}$ transports domains
along with the operator. A Mellin integral taken with radial area
measure does not by itself change the equal-time Dirac norm.

Introduce the null time $\tau=(\eta+v)/2$, so
$t+x_{\rm M}=\xi_0 e^{2\tau}$ and $\eta=2\tau-v$. Write
$q_\alpha(\tau,x)=g_\alpha(2\tau-v,v)$. At fixed $\eta$,
$\partial_v g_\alpha=x\partial_x q_\alpha+
\tfrac12\partial_\tau q_\alpha$, and at fixed $v$,
$\partial_\eta g_\alpha=\tfrac12\partial_\tau q_\alpha$.
The two components of the massive equation therefore give
\begin{equation}
\begin{gathered}
 i\partial_\tau q_1=x(pq_1+mq_2),\\
 p q_2=m q_1,\qquad p=-i\partial_x .
\end{gathered}
 \label{eq:nullconstraint}
\end{equation}
Here $p$ acts on the radial amplitude as the momentum in $x$.
The null constraint recovers $q_2$ by integrating $m q_1$ from
the far endpoint, where decaying modes obey $q_2(\tau,\infty)=0$.
This fixes the right inverse
\[
 (p^{-1}_\infty b)(x)=-i\int_x^\infty b(y)\,\dd y .
\]
Set $M=\sqrt{x}$, acting by multiplication, and
$\phi=M^{-1}q_1$. Equation~\eqref{eq:nullconstraint} gives the
scalar generator $A=M(p+m^2p^{-1}_\infty)M$, which evolves the
independent datum $\phi$ in null time $\tau$. The inverse
construction, on smooth decaying solutions for which the integral
exists, is
\begin{align}
 (q_1,q_2)&=(M\phi,\,mp^{-1}_\infty M\phi),\nonumber\\
 g(\eta,v)&=q((\eta+v)/2,\xi_0e^v),\nonumber\\
 \psi(\eta,x)&=x^{-1/2}g(\eta,\log(x/\xi_0)).
 \label{eq:inversecomponents}
\end{align}
Substitution recovers both equations, including the integration
constant. This is a change of Cauchy variables on these solutions.
Taking data on a null slice discards a freely specifiable second
component because it must satisfy the constraint.

With $a^{-2}=m^2$ in the original units, the two-channel convention
of the deposited paper is
\begin{align}
 H_M&=\begin{pmatrix}0&A\\ A&0\end{pmatrix},\quad
 S=\frac1{\sqrt2}\begin{pmatrix}1&1\\1&-1\end{pmatrix},\nonumber\\
 S^*H_MS&=\operatorname{diag}(A,-A).
 \label{majhamilton}
\end{align}
The matrix $S$ resolves the two-channel evolution $H_M$ into
symmetric and antisymmetric spinor amplitudes. The embeddings of
these two amplitudes, $J_\pm\phi=(\phi,\pm\phi)^T/\sqrt2$, obey
$H_MJ_\pm=\pm J_\pm A$. These statements are algebraic on a
common core and hold on $D(A)\oplus D(A)$ for a chosen
self-adjoint realization of $A$. The scalar evolution reconstructed
in Eq.~\eqref{eq:inversecomponents} uses $J_+$. The realization on
the selected radial modes and its graph domain are proved in
Sec.~\ref{sec:inversemajorana}. The integration operator here, the
compressed Hilbert transform and the antiperiodic inverse on a compact
interval have distinct kernels and domains (Appendix~\ref{orig:D}).

\head{Change of representation and functional calculus}
\label{remark:HMtransformations}
If $H_w$ is given on $L^2(x\,\dd x)$, its representative on
$L^2(\dd x)$ is $M H_w M^{-1}$ with domain $M\mathcal D(H_w)$.
Conversely the representative of $A$ on the weighted space is
$M^{-1}AM$. On the common smooth domain where the integral exists
this last expression acts as $(p+m^2p_\infty^{-1})x$.
The equal factors $M$ in $A=M(p+m^2p_\infty^{-1})M$ arose from
the null constraint and the definition of $\phi$.

The unitary Mellin transform in the physical scalar measure is
\begin{align*}
 (\mathcal M f)(k)&=\frac1{\sqrt{2\pi}}\int_0^\infty
       (x/\xi_0)^{-ik}x^{-1/2}f(x)\,\dd x,\\
 f(x)&=\frac{x^{-1/2}}{\sqrt{2\pi}}\int_\R
       (x/\xi_0)^{ik}(\mathcal M f)(k)\,\dd k .
\end{align*}
The inverse is understood in $L^2$. These equations are the Fourier
transform and its inverse after $U_{\log}$. The coefficient at $k$
resolves oscillations of the radial wave in $\log(x/\xi_0)$.
These transforms send $MpM$ to
multiplication by $k$. On test spinors in $v$, Fourier transformation
of $H_0g$ gives
$k\sigma_3\widehat g(k)+m\xi_0\sigma_1\widehat g(k+i)$.
The shift is valid there because $g$ has compact support. The massive
operator retains this extra term. The scalar Mellin--Barnes observable
of Sec.~\ref{sec3} is obtained by a separate sum over dilations.

For any unitary $U$, the spectral theorem gives
$U F(H)U^{-1}=F(UHU^{-1})$ on the transported domain.
For measurable $F$ that domain consists of vectors with
$\int |F(\lambda)|^2\,\dd\mu_f(\lambda)<\infty$, where
$\mu_f$ is the vector spectral measure. A trace interpretation of
the radial integral further requires a trace-class operator or a
specified relative trace and a proof of the trace identity.

\section{Mellin--Barnes integral and radial spectrum}\label{sec3}
\subsection{The order of the Bessel function and the boundary phase}
The radial functions of the original ansatz are $K_\nu(z)$ for
$\nu=1/2+iE/2,\ z>0$, where $E/2=\omega$ in that ansatz.
They obey $K_{-\nu}=K_\nu$ and
$\overline{K_\nu(z)}=K_{\bar\nu}(z)$. The corresponding identity
for $I_\nu$ is
$I_{-\nu}-I_\nu=2\sin(\pi\nu)K_\nu/\pi$.
For $z>0$ the boundary value $K_\nu(z)$ is nonzero on
$\re\nu=1/2$ \cite{paris,bagirova,Olver}. This follows directly
from the radial equation. The function $u(v)=K_\nu(ze^v)$ solves
\[
 (-\partial_v^2+z^2e^{2v})u=-\nu^2u,\quad v>0 .
\]
It decays at infinity and is nonzero by the large-argument
asymptotic. If it vanished at $0$, it would be an eigenfunction of
a positive self-adjoint Dirichlet operator at $-\nu^2$.
This number is nonreal for $\im\nu\ne0$ and negative for
$\nu=1/2$. Consequently the single-mode matching condition
\begin{equation}
 e^{i\theta}K_\nu(z)=K_{\bar\nu}(z) \label{maj1}
\end{equation}
admits an energy-dependent phase for every real $E$. Spectral
quantization imposes a single boundary phase over the entire energy
range.

\subsection{Mellin transform of the radial superposition}\label{subsec:MB}

In the original radial coordinate the decaying solutions involve
$K_{1/2\pm i\omega}(m\xi)$, with $\omega=E/2$ in the energy
convention of Ref.~\cite{v3}. Equation (49) there uses a radial profile
at fixed energy. A full Green operator also requires a source
normalization, a domain and boundary conditions. The scalar radial
integral with density $\xi\,\dd\xi$ does not define the Dirac
inner product used in Proposition~\ref{pr:boost}. For example,
for $m>0,\ \re s<-|\re\nu|$, the Mellin transform of an individual
radial factor is~\cite{dlmf}
\begin{multline}
 \int_0^\infty a^{s-1}K_\nu(m/a)\dd a\\
 =m^s2^{-s-2}\Gamma\bigl((-s+\nu)/2\bigr)
                   \Gamma\bigl((-s-\nu)/2\bigr).
 \label{eq:bess}
\end{multline}
A separate periodic transverse circle of radius $R_0$ has squared
positive momentum levels $n^2/R_0^2,\ n\ge1$, one copy for this
branch. Its spectral zeta is $R_0^{2s}\zeta(2s)$.
For $F_\nu(t)=t^{-\nu/2}K_\nu(2\sqrt t)$ this distinction gives
\[
 \begin{gathered}\int_0^\infty t^{s-1}\sum_{n\ge1} F_\nu(n^2t)\dd t\\ =\tfrac12\Gamma(s)\Gamma(s-\nu)\zeta(2s).\end{gathered}
\]
Here $\re s>\max(0,\re\nu,1/2)$. Substitution $r=2\sqrt t$ gives the two gamma factors and absolute
convergence permits summation over $n$. A zero of this continued
Mellin transform is a different spectral statement from an eigenvalue
equation for \eqref{eq:boost}.

In the earlier zeta test choose
$g>\max(0,1/2,\re\nu)$ and $a>0$.
The integral is
\begin{align}
 I_\nu(a)&=\frac1{4\pi i}\int_{(g)}
      \Gamma(s)\Gamma(s-\nu)(2a)^{2s}\zeta(2s)\dd s \nonumber\\
 &= (2a)^\nu\sum_{n\ge1}n^{-\nu}K_\nu(n/a).
 \label{MB}
\end{align}
Here $\int_{(g)}$ is upward on $\re s=g$. Absolute convergence
of the Dirichlet series and Stirling's bound on this line permit
termwise inversion. The original ansatz with independent accelerations
$a_n$ and coefficients $g_n$ requires additional conditions to
give these particular dilations and amplitudes. Formula \eqref{MB}
specifies one such arithmetic superposition directly. The line
$g=1/2$ in the earlier contour argument passes through the zeta
pole and, for the physical order, also through a pole of
$\Gamma(s-\nu)$.

\begin{proposition}[Small acceleration]\label{th:nonvanish}
For $\nu=1/2+iE/2$, uniformly for $E$ on any fixed compact
interval, one has as $a\downarrow0$
\[
 I_\nu(a)=(2a)^\nu\sqrt{\pi a/2}\,e^{-1/a}
                 [1+O(a)+O(e^{-1/a})] .
\]
In particular the integral is nonzero throughout that energy interval
for all sufficiently small positive $a$.
\end{proposition}
\begin{proof}
For the $n=1$ term use $K_\nu(z)=\sqrt{\pi/(2z)}e^{-z}
(1+O((1+|\nu|)^2/z))$ for orders in the indicated compact
set~\cite{Olver,dlmf}. The integral representation yields
$|K_{1/2+it}(z)|\le K_{1/2}(z)=\sqrt{\pi/(2z)}e^{-z}$
for every real $t$. Divide the tail $n\ge2$ by
$\sqrt{\pi a/2}e^{-1/a}$ and sum
$\sum_{n\ge2}n^{-1}e^{-(n-1)/a}$.
\end{proof}
Near a nontrivial zero $\rho$ of multiplicity $m$, an analytic
factor $\Phi$ satisfies
$\Phi(s)\zeta(2s)=O((s-\rho/2)^m)$. The residue at $s=\rho/2$
therefore vanishes. For a simple zero the linear Taylor coefficient
contributes to the derivative at that point. Multiplicity enters the
arithmetic formula of Sec.~\ref{sec:explicit} through the residue
of $-\zeta'/\zeta$ at $s=\rho$, which equals $-m$.

\section{Spectral completeness, chiral channels and primes}\label{sec4}
\subsection{The spectral triple and a physical return observable}
The spectral triple in the original proposal consists of an algebra
of observables, a Hilbert space and a self-adjoint operator
\cite{connes}. On a compact logarithmic channel one example
is given by periodic smooth functions acting by multiplication,
$L^2(0,\ell)$, and the twisted generator
$B_\theta=-i\partial_v$. It propagates a chiral wave around the
channel with boundary phase $e^{i\theta}$. Its resolvent is compact and
$[B_\theta,f]=-if'$ is bounded. For every $\theta$ these
properties hold, while the frequencies $(2\pi k+\theta)/\ell$
change with $\theta$. For infinitely many prime channels even
compactness of the direct-sum resolvent would require another
argument. The observable below is defined by relative traces on
individual channels. The precise criterion for an identification
with a spectral measure is proved in Sec.~\ref{sec:CAR}.

The dilation algebra proposed for the unbounded wedge also has a
concrete representation. In $v=\log(\xi/\xi_0)$ let
$P=-i\partial_v$, $T_b g(v)=g(v-b)$, and $M_f g=fg$.
In a wave packet, $P$ generates the displacement $T_b$ by
optical distance $b$, while $M_f$ weights its local amplitude
with profile $f$. Thus $T_bM_fT_b^*=M_{f(\cdot-b)}$ and
$[P,M_f]=-iM_{f'}$ for $f\in C_c^\infty(\R)$.
A spatially localized probe sees $M_f(P-i)^{-1}$, the driven
translation response weighted over its sensitivity region. The operator $M_f(P-i)^{-1}$ is Hilbert--Schmidt: by Plancherel
its squared Hilbert--Schmidt norm is
$\|f\|_2^2(2\pi)^{-1}\int_\R(1+k^2)^{-1}\dd k$.
The unrestricted resolvent is not compact, since in Fourier space
it is a nonzero multiplication operator on a nonatomic space.
This distinguishes local compactness of the scaling representation
from discreteness of a candidate Hilbert--P\'olya spectrum.

\head{Prescribed prime lengths}
For a chiral channel of boost length $\ell$ take the self-adjoint
generator $B_\theta=-i\partial_v$ on $H^1(0,\ell)$ with
$f(\ell)=e^{i\theta}f(0)$. Its frequencies are
$(2\pi j+\theta)/\ell,\ j\in\mathbb Z$. A real scalar channel
allows holonomy $+1$ or $-1$. A pair of real channels with boundary
rotation in $SO(2)$ realizes the two conjugate phases. For a single
complexified branch, the difference of resolvents for $\theta$
and $0$ is rank at most one. Solving the first-order inhomogeneous
equation with its boundary condition gives the relative trace
\begin{multline}
 -\tr\big[(B_\theta-z)^{-1}-(B_0-z)^{-1}\big]
 \\=-i\ell\sum_{k\ge1}(e^{-ik\theta}-1)e^{ik\ell z},\qquad \im z>0 .
 \label{eq:ring}
\end{multline}
The resolvent $(B_\theta-z)^{-1}$ gives the driven response of
the chiral channel at complex frequency $z$. Subtracting the two
phase settings isolates the change at each winding $k$. The kernel
difference is
$i e^{iz\ell}[(e^{i\theta}-e^{iz\ell})^{-1}
 -(1-e^{iz\ell})^{-1}]e^{iz(v-v')}$.
For each prime $p$ assign $\ell_p=\tfrac12\log p$.
The phase-independent return amplitude is proportional to $\ell_p$
at every winding $k$. More explicitly, the phase mean of the
right-hand side of \eqref{eq:ring} equals
$i\ell_p\sum_{k\ge1}e^{ik\ell_p z}$.
Multiplication by $2/i$ and use of the clock $u=2\eta$ give the
locally finite return measure
\begin{equation}
 \mu_{\rm ret}=\sum_{p}\sum_{k\ge1}\log p\,\delta_{k\log p}.
 \label{eq:ret}
\end{equation}
In this detector $x=e^u$ denotes the dimensionless arithmetic
cutoff. A return at time $k\log p$ contributes weight $\log p$.
The transform converges absolutely for $\im z>2$ after the sum
over $p$. The distribution itself is defined on every bounded time
interval without a regulator. It requires prime lengths as input.
Related Rindler and quantum-graph models with prime data occur in
Refs.~\cite{sierradirac,graph}.

Put $h(u)=(e^u-1)\one_{u\ge0}$. Convolution with the detector response
$h$ gives the output at $u=\log x$,
\begin{equation}
\begin{aligned}
 T(x)&=(h*\mu_{\rm ret})(\log x)\\
     &=\sum_{n\le x}\Lambda(n)(x/n-1),\qquad x>0 ,
\end{aligned}
 \label{eq:T}
\end{equation}
where $\Lambda(p^k)=\log p$ and $\Lambda(n)=0$ otherwise.
Distributionally $(\partial_u^2-\partial_u)T(e^u)=\mu_{\rm ret}$.
The label $k$ counts repeated propagation of a one-particle amplitude.
In a many-body ensemble with a single fermionic level per prime and
energy $\log p$, exclusion gives
$\prod_p(1+p^{-s})=\zeta(s)/\zeta(2s)$ for $\re s>1$.
Its logarithmic derivative carries the sign $(-1)^{k-1}$ at the
$k$th repetition. The return observable \eqref{eq:ret} carries a
positive weight at every repetition. These two observables make
different, explicit predictions.

\section{The prime-return correlation and its positivity}\label{sec:correlation}
Define the cutoff $x=e^u$ using the flight-time clock of
Sec.~\ref{sec4}. We first center the detected signal and then form
its correlation in the return time. The induction in the cutoff is
deferred to Sec.~\ref{sec:analytic}.

\subsection{The centered response and its one-sided bound}\label{sec:explicit}
Let $\rho=\beta+i\gamma$ run over the nontrivial zeros of $\zeta$,
with multiplicity. All logarithms are natural. Write
\begin{equation}
\begin{gathered}
 R(x)=T(x)-x\log x+(1+\eul)x,\\
 d_0=\log(2\pi),\qquad
 S(x)=\sum_{j\ge1}\frac{x^{-2j}}{2j(2j+1)},\\
 c_\rho=\frac{1}{\rho(1-\rho)},\qquad Q=2+\eul-\log(4\pi).
\end{gathered}\label{eq:center}
\end{equation}
The distribution in \eqref{eq:ret} also gives
$T(x)=x\int_1^x\psi(t)t^{-2}\dd t$ for $x\ge1$, where
$\psi(t)=\sum_{n\le t}\Lambda(n)$. This follows by integrating each
step of $\psi$. In particular $T$ measures an integral of prime-power
counts, with each jump weighted by the logarithm of its primitive prime.

\begin{lemma}\label{lem:explicit}
For every $x\ge1$,
\begin{equation}
\begin{gathered}
 R(x)=d_0-S(x)+\sum_\rho c_\rho x^\rho,\\
 \sum_\rho \re c_\rho=Q,\qquad \sum_\rho |c_\rho|<\infty .
\end{gathered}
 \label{eq:explicit}
\end{equation}
In particular $R(x)/x\to0$ unconditionally. Under RH,
$R(x)\le d_0+Q\sqrt{x}<\sqrt{x}$ for all $x\ge4$.
Conversely any eventual bound $R(x)\le C\sqrt{x}$ implies RH.
\end{lemma}
\begin{proof}
Set $G(s)=-\zeta'(s)/\zeta(s)$. The inverse Mellin integrand for
$T(x)$ on $\re s=2$ is $G(s)x^s/[s(s-1)]$, by the transform
of $(x-1)\one_{x\ge1}$ and the absolutely convergent Dirichlet
series for $G$. At 1 the double pole contributes
$x(\log x-1-\eul)$. At 0 the residue is $d_0$, by
$\zeta'(0)/\zeta(0)=\log(2\pi)$. The trivial zeros $-2j$ give
$-S(x)$ and each nontrivial zero contributes $c_\rho x^\rho$
per unit of multiplicity. For fixed $N$ move to $\re s=-2N-1$
along heights whose distance from all zero ordinates is
$\gg1/\log(2+|\im s|)$. Such heights follow from the local zero
count $O(\log T)$ in a unit interval. The local formula for the
logarithmic derivative gives $O_N(\log^2 T)$ on the horizontals.
On the left vertical the functional equation gives
$G(s)=O(\log(2+N+|\im s|))$ with an absolute constant. Its
integral is $O(x^{-2N-1}\log(N+2)/(N+1))$.
The zero count $N(T)=O(T\log(2+T))$ gives absolute convergence
and justifies the limits, including $x=1$. These standard estimates
follow from the functional equation and Hadamard product as in
Refs.~\cite{Titchmarsh,mv}.

At $x=1$, $T(1)=0$ and $S(1)=1-\log2$, whence $\sum c_\rho=Q$.
The real part of $c_\rho=1/\rho+1/(1-\rho)$ is positive.
Each summand $c_\rho x^{\rho-1}$ tends to zero since $0<\beta<1$.
Dominated convergence gives $R(x)/x\to0$. Under RH every $c_\rho$
is positive and equals $(1/4+\gamma^2)^{-1}$, proving the upper
bound, since $d_0/2+Q<1$.

For the converse consider $\int_1^\infty R(x)x^{-s-1}\dd x$.
Initially for $\re s>1$ it is
\[
 \frac{G(s)}{s(s-1)}-\frac1{(s-1)^2}+\frac{1+\eul}{s-1}.
\]
This continues regularly at every real $s>1/2$, including 1.
There are no real zeros of $\zeta$ in $(0,1)$, by the alternating
eta series. If $C\sqrt{x}-R(x)\ge0$ eventually, its tail transform
has the same real regularity for $s>1/2$. A nonnegative tail
transform has a singularity at its real abscissa of convergence
whenever the latter is finite. Indeed, expand about a real point
just to its right towards the left. All the integral Taylor
coefficients in that direction are nonnegative. An analytic
continuation through the abscissa would give, by monotone convergence,
convergence a positive distance beyond it. Thus the tail transform
is holomorphic in $\re s>1/2$. A zero there would give a pole of
nonzero residue. Reflection about the critical line completes the proof.
\end{proof}

\subsection{The arithmetic correlation and a positive spectral measure}\label{sec:CAR}
Define the even, continuous function
\begin{equation}
 C(u)=e^{-|u|/2}\big[R(e^{|u|})-d_0+S(e^{|u|})\big].
 \label{eq:C}
\end{equation}
All data in this expression come from elementary constants and the
primes up to $e^{|u|}$. In particular, with $z=1/x<1$,
\begin{equation}
 S(x)=1-\tfrac12\log(1-z^2)-z^{-1}\operatorname{arctanh}z .
 \label{eq:Sclosed}
\end{equation}
Equation \eqref{eq:explicit} and the reflection $\rho\leftrightarrow
1-\rho$ give, unconditionally on the whole real line,
\begin{equation}
 C(u)=\sum_\rho c_\rho e^{(\rho-1/2)u},\qquad C(0)=Q.
 \label{eq:specC}
\end{equation}
The series converges uniformly on compact $u$-intervals. The
arithmetic definition \eqref{eq:C} of $C$ can be evaluated from
finite return data and the elementary function $S$.
Figure~\ref{fig:C} compares the two sides of \eqref{eq:specC}. The
sum over the first $2000$ zeros differs from the prime side by at
most the mass of the remaining zeros, $8.9\cdot10^{-4}$, plus the
finite-height term of Sec.~\ref{sec:boundaryclosure}.

\begin{figure}[t]
\centering
\begin{tikzpicture}
\begin{axis}[width=\columnwidth,height=5.9cm,xmin=0,xmax=4.2,ymin=-0.034,ymax=0.053,
 xlabel={return time $u$},ylabel={$C(u)$},
 scaled y ticks=false,
 yticklabel style={/pgf/number format/fixed,/pgf/number format/precision=2},
 ytick={-0.02,0,0.02,0.04},
 tick label style={font=\footnotesize},label style={font=\small}]
\addplot[black,dashed,thin] coordinates {(0,0.04619) (4.2,0.04619)};
\addplot[black!40,thin] coordinates {(0,0) (4.2,0)};
\addplot[refblue,thick] coordinates {
(0.0000,0.04619) (0.0001,0.04580) (0.0010,0.04344) (0.0030,0.03959) (0.0080,0.03248)
(0.0160,0.02421) (0.0240,0.01791) (0.0320,0.01287) (0.0400,0.00872) (0.0480,0.00527)
(0.0560,0.00237) (0.0640,-0.00006) (0.0720,-0.00211) (0.0800,-0.00382) (0.0880,-0.00525)
(0.0960,-0.00642) (0.1040,-0.00737) (0.1120,-0.00812) (0.1200,-0.00871) (0.1280,-0.00913)
(0.1360,-0.00942) (0.1440,-0.00959) (0.1520,-0.00965) (0.1600,-0.00961) (0.1680,-0.00948)
(0.1760,-0.00928) (0.1840,-0.00900) (0.1920,-0.00867) (0.2000,-0.00828) (0.2080,-0.00784)
(0.2160,-0.00737) (0.2240,-0.00685) (0.2320,-0.00631) (0.2400,-0.00575) (0.2480,-0.00516)
(0.2560,-0.00456) (0.2640,-0.00395) (0.2720,-0.00333) (0.2800,-0.00271) (0.2880,-0.00208)
(0.2960,-0.00146) (0.3040,-0.00085) (0.3120,-0.00024) (0.3200,0.00036) (0.3280,0.00094)
(0.3360,0.00150) (0.3440,0.00205) (0.3520,0.00257) (0.3600,0.00307) (0.3680,0.00354)
(0.3760,0.00399) (0.3840,0.00441) (0.3920,0.00479) (0.4000,0.00514) (0.4080,0.00546)
(0.4160,0.00573) (0.4240,0.00597) (0.4320,0.00617) (0.4400,0.00633) (0.4480,0.00644)
(0.4560,0.00651) (0.4640,0.00654) (0.4720,0.00651) (0.4800,0.00644) (0.4880,0.00632)
(0.4960,0.00615) (0.5040,0.00592) (0.5120,0.00564) (0.5200,0.00530) (0.5280,0.00491)
(0.5360,0.00447) (0.5440,0.00396) (0.5520,0.00340) (0.5600,0.00277) (0.5680,0.00208)
(0.5760,0.00134) (0.5840,0.00052) (0.5920,-0.00035) (0.6000,-0.00129) (0.6080,-0.00230)
(0.6160,-0.00337) (0.6240,-0.00452) (0.6320,-0.00573) (0.6400,-0.00701) (0.6480,-0.00836)
(0.6560,-0.00978) (0.6640,-0.01127) (0.6720,-0.01284) (0.6800,-0.01448) (0.6880,-0.01619)
(0.6931,-0.01733) (0.6960,-0.01658) (0.7040,-0.01452) (0.7120,-0.01255) (0.7200,-0.01064)
(0.7280,-0.00882) (0.7360,-0.00708) (0.7440,-0.00541) (0.7520,-0.00382) (0.7600,-0.00232)
(0.7680,-0.00090) (0.7760,0.00045) (0.7840,0.00170) (0.7920,0.00288) (0.8000,0.00397)
(0.8080,0.00498) (0.8160,0.00590) (0.8240,0.00674) (0.8320,0.00749) (0.8400,0.00815)
(0.8480,0.00872) (0.8560,0.00921) (0.8640,0.00961) (0.8720,0.00992) (0.8800,0.01013)
(0.8880,0.01026) (0.8960,0.01030) (0.9040,0.01024) (0.9120,0.01010) (0.9200,0.00985)
(0.9280,0.00952) (0.9360,0.00909) (0.9440,0.00857) (0.9520,0.00795) (0.9600,0.00723)
(0.9680,0.00642) (0.9760,0.00552) (0.9840,0.00451) (0.9920,0.00341) (1.0000,0.00220)
(1.0080,0.00090) (1.0160,-0.00050) (1.0240,-0.00201) (1.0320,-0.00361) (1.0400,-0.00532)
(1.0480,-0.00712) (1.0560,-0.00903) (1.0640,-0.01105) (1.0720,-0.01317) (1.0800,-0.01539)
(1.0880,-0.01772) (1.0960,-0.02015) (1.0986,-0.02097) (1.1040,-0.01928) (1.1120,-0.01685)
(1.1200,-0.01453) (1.1280,-0.01232) (1.1360,-0.01022) (1.1440,-0.00822) (1.1520,-0.00634)
(1.1600,-0.00456) (1.1680,-0.00290) (1.1760,-0.00134) (1.1840,0.00010) (1.1920,0.00143)
(1.2000,0.00264) (1.2080,0.00375) (1.2160,0.00474) (1.2240,0.00562) (1.2320,0.00638)
(1.2400,0.00702) (1.2480,0.00755) (1.2560,0.00797) (1.2640,0.00827) (1.2720,0.00845)
(1.2800,0.00851) (1.2880,0.00845) (1.2960,0.00828) (1.3040,0.00799) (1.3120,0.00757)
(1.3200,0.00704) (1.3280,0.00638) (1.3360,0.00560) (1.3440,0.00470) (1.3520,0.00368)
(1.3600,0.00253) (1.3680,0.00126) (1.3760,-0.00013) (1.3840,-0.00165) (1.3863,-0.00211)
(1.3920,-0.00132) (1.4000,-0.00032) (1.4080,0.00055) (1.4160,0.00130) (1.4240,0.00192)
(1.4320,0.00241) (1.4400,0.00277) (1.4480,0.00300) (1.4560,0.00310) (1.4640,0.00307)
(1.4720,0.00291) (1.4800,0.00262) (1.4880,0.00219) (1.4960,0.00163) (1.5040,0.00094)
(1.5120,0.00012) (1.5200,-0.00084) (1.5280,-0.00194) (1.5360,-0.00317) (1.5440,-0.00454)
(1.5520,-0.00605) (1.5600,-0.00769) (1.5680,-0.00947) (1.5760,-0.01139) (1.5840,-0.01345)
(1.5920,-0.01565) (1.6000,-0.01800) (1.6080,-0.02048) (1.6094,-0.02094) (1.6160,-0.01838)
(1.6240,-0.01539) (1.6320,-0.01254) (1.6400,-0.00984) (1.6480,-0.00728) (1.6560,-0.00487)
(1.6640,-0.00260) (1.6720,-0.00048) (1.6800,0.00150) (1.6880,0.00333) (1.6960,0.00501)
(1.7040,0.00654) (1.7120,0.00792) (1.7200,0.00915) (1.7280,0.01024) (1.7360,0.01117)
(1.7440,0.01195) (1.7520,0.01258) (1.7600,0.01306) (1.7680,0.01338) (1.7760,0.01355)
(1.7840,0.01356) (1.7920,0.01342) (1.8000,0.01313) (1.8080,0.01267) (1.8160,0.01206)
(1.8240,0.01130) (1.8320,0.01037) (1.8400,0.00929) (1.8480,0.00804) (1.8560,0.00663)
(1.8640,0.00507) (1.8720,0.00334) (1.8800,0.00145) (1.8880,-0.00061) (1.8960,-0.00282)
(1.9040,-0.00521) (1.9120,-0.00776) (1.9200,-0.01047) (1.9280,-0.01335) (1.9360,-0.01640)
(1.9440,-0.01962) (1.9459,-0.02041) (1.9520,-0.01852) (1.9600,-0.01619) (1.9680,-0.01404)
(1.9760,-0.01205) (1.9840,-0.01023) (1.9920,-0.00859) (2.0000,-0.00712) (2.0080,-0.00582)
(2.0160,-0.00470) (2.0240,-0.00375) (2.0320,-0.00298) (2.0400,-0.00239) (2.0480,-0.00197)
(2.0560,-0.00173) (2.0640,-0.00167) (2.0720,-0.00178) (2.0794,-0.00205) (2.0800,-0.00194)
(2.0880,-0.00046) (2.0960,0.00084) (2.1040,0.00196) (2.1120,0.00289) (2.1200,0.00365)
(2.1280,0.00421) (2.1360,0.00460) (2.1440,0.00479) (2.1520,0.00480) (2.1600,0.00463)
(2.1680,0.00426) (2.1760,0.00371) (2.1840,0.00296) (2.1920,0.00203) (2.1972,0.00132)
(2.2000,0.00192) (2.2080,0.00353) (2.2160,0.00495) (2.2240,0.00618) (2.2320,0.00721)
(2.2400,0.00805) (2.2480,0.00869) (2.2560,0.00913) (2.2640,0.00938) (2.2720,0.00943)
(2.2800,0.00928) (2.2880,0.00893) (2.2960,0.00838) (2.3040,0.00763) (2.3120,0.00668)
(2.3200,0.00552) (2.3280,0.00416) (2.3360,0.00259) (2.3440,0.00082) (2.3520,-0.00116)
(2.3600,-0.00334) (2.3680,-0.00573) (2.3760,-0.00833) (2.3840,-0.01115) (2.3920,-0.01417)
(2.3979,-0.01653) (2.4000,-0.01588) (2.4080,-0.01354) (2.4160,-0.01142) (2.4240,-0.00951)
(2.4320,-0.00782) (2.4400,-0.00634) (2.4480,-0.00508) (2.4560,-0.00403) (2.4640,-0.00321)
(2.4720,-0.00260) (2.4800,-0.00221) (2.4880,-0.00205) (2.4960,-0.00210) (2.5040,-0.00238)
(2.5120,-0.00289) (2.5200,-0.00362) (2.5280,-0.00457) (2.5360,-0.00575) (2.5440,-0.00716)
(2.5520,-0.00879) (2.5600,-0.01066) (2.5649,-0.01193) (2.5680,-0.01058) (2.5760,-0.00722)
(2.5840,-0.00409) (2.5920,-0.00119) (2.6000,0.00148) (2.6080,0.00391) (2.6160,0.00611)
(2.6240,0.00807) (2.6320,0.00979) (2.6400,0.01127) (2.6480,0.01251) (2.6560,0.01352)
(2.6640,0.01428) (2.6720,0.01480) (2.6800,0.01508) (2.6880,0.01511) (2.6960,0.01490)
(2.7040,0.01444) (2.7120,0.01373) (2.7200,0.01278) (2.7280,0.01158) (2.7360,0.01013)
(2.7440,0.00842) (2.7520,0.00647) (2.7600,0.00426) (2.7680,0.00179) (2.7726,0.00027)
(2.7760,-0.00033) (2.7840,-0.00192) (2.7920,-0.00377) (2.8000,-0.00587) (2.8080,-0.00824)
(2.8160,-0.01086) (2.8240,-0.01375) (2.8320,-0.01690) (2.8332,-0.01740) (2.8400,-0.01565)
(2.8480,-0.01383) (2.8560,-0.01228) (2.8640,-0.01099) (2.8720,-0.00997) (2.8800,-0.00922)
(2.8880,-0.00874) (2.8960,-0.00854) (2.9040,-0.00860) (2.9120,-0.00894) (2.9200,-0.00955)
(2.9280,-0.01044) (2.9360,-0.01161) (2.9440,-0.01305) (2.9444,-0.01314) (2.9520,-0.00967)
(2.9600,-0.00627) (2.9680,-0.00315) (2.9760,-0.00031) (2.9840,0.00224) (2.9920,0.00451)
(3.0000,0.00650) (3.0080,0.00819) (3.0160,0.00960) (3.0240,0.01073) (3.0320,0.01156)
(3.0400,0.01210) (3.0480,0.01234) (3.0560,0.01230) (3.0640,0.01195) (3.0720,0.01132)
(3.0800,0.01038) (3.0880,0.00915) (3.0960,0.00762) (3.1040,0.00578) (3.1120,0.00365)
(3.1200,0.00121) (3.1280,-0.00153) (3.1355,-0.00438) (3.1360,-0.00425) (3.1440,-0.00238)
(3.1520,-0.00081) (3.1600,0.00045) (3.1680,0.00139) (3.1760,0.00203) (3.1840,0.00235)
(3.1920,0.00236) (3.2000,0.00205) (3.2080,0.00142) (3.2160,0.00048) (3.2189,0.00006)
(3.2240,0.00086) (3.2320,0.00185) (3.2400,0.00252) (3.2480,0.00287) (3.2560,0.00289)
(3.2640,0.00258) (3.2720,0.00195) (3.2800,0.00099) (3.2880,-0.00030) (3.2958,-0.00189)
(3.2960,-0.00189) (3.3040,-0.00215) (3.3120,-0.00275) (3.3200,-0.00368) (3.3280,-0.00495)
(3.3360,-0.00656) (3.3440,-0.00851) (3.3520,-0.01079) (3.3600,-0.01342) (3.3673,-0.01612)
(3.3680,-0.01596) (3.3760,-0.01427) (3.3840,-0.01293) (3.3920,-0.01194) (3.4000,-0.01130)
(3.4080,-0.01101) (3.4160,-0.01107) (3.4240,-0.01149) (3.4320,-0.01226) (3.4340,-0.01250)
(3.4400,-0.00968) (3.4480,-0.00623) (3.4560,-0.00313) (3.4640,-0.00040) (3.4657,0.00014)
(3.4720,0.00274) (3.4800,0.00572) (3.4880,0.00835) (3.4960,0.01060) (3.5040,0.01249)
(3.5120,0.01401) (3.5200,0.01516) (3.5280,0.01594) (3.5360,0.01634) (3.5440,0.01638)
(3.5520,0.01603) (3.5600,0.01531) (3.5680,0.01420) (3.5760,0.01272) (3.5840,0.01085)
(3.5920,0.00860) (3.6000,0.00597) (3.6080,0.00295) (3.6109,0.00175) (3.6160,0.00255)
(3.6240,0.00350) (3.6320,0.00405) (3.6400,0.00422) (3.6480,0.00398) (3.6560,0.00336)
(3.6640,0.00233) (3.6720,0.00090) (3.6800,-0.00093) (3.6880,-0.00316) (3.6960,-0.00579)
(3.7040,-0.00883) (3.7120,-0.01228) (3.7136,-0.01301) (3.7200,-0.01242) (3.7280,-0.01205)
(3.7360,-0.01209) (3.7440,-0.01255) (3.7520,-0.01343) (3.7600,-0.01472) (3.7612,-0.01495)
(3.7680,-0.01253) (3.7760,-0.01008) (3.7840,-0.00804) (3.7920,-0.00644) (3.8000,-0.00526)
(3.8080,-0.00450) (3.8160,-0.00418) (3.8240,-0.00429) (3.8320,-0.00483) (3.8400,-0.00581)
(3.8480,-0.00722) (3.8501,-0.00767) (3.8560,-0.00578) (3.8640,-0.00358) (3.8720,-0.00182)
(3.8800,-0.00051) (3.8880,0.00036) (3.8918,0.00062) (3.8960,0.00195) (3.9040,0.00415)
(3.9120,0.00589) (3.9200,0.00719) (3.9280,0.00803) (3.9360,0.00841) (3.9440,0.00834)
(3.9520,0.00781) (3.9600,0.00681) (3.9680,0.00536) (3.9703,0.00485) (3.9760,0.00655)
(3.9840,0.00852) (3.9920,0.01002) (4.0000,0.01106) (4.0080,0.01162) (4.0160,0.01170)
(4.0240,0.01131) (4.0320,0.01044) (4.0400,0.00909) (4.0480,0.00726) (4.0560,0.00495)
(4.0640,0.00214) (4.0720,-0.00115) (4.0775,-0.00371) (4.0800,-0.00362) (4.0880,-0.00365)
(4.0960,-0.00417) (4.1040,-0.00519) (4.1109,-0.00646) (4.1120,-0.00611) (4.1200,-0.00391)
(4.1280,-0.00222) (4.1360,-0.00103) (4.1440,-0.00035) (4.1520,-0.00017) (4.1589,-0.00043)
(4.1600,-0.00041) (4.1680,-0.00057) (4.1760,-0.00124) (4.1840,-0.00242) (4.1920,-0.00412)
(4.2000,-0.00635)
};
\addplot[only marks,mark=o,mark size=1.0pt,black] coordinates {
(0.00,0.04531) (0.10,-0.00692) (0.20,-0.00828) (0.30,-0.00115) (0.40,0.00514)
(0.50,0.00604) (0.60,-0.00129) (0.70,-0.01554) (0.80,0.00397) (0.90,0.01028)
(1.00,0.00220) (1.10,-0.02053) (1.20,0.00264) (1.30,0.00815) (1.40,-0.00032)
(1.50,0.00131) (1.60,-0.01799) (1.70,0.00579) (1.80,0.01313) (1.90,-0.00400)
(2.00,-0.00712) (2.10,0.00142) (2.20,0.00192) (2.30,0.00803) (2.40,-0.01587)
(2.50,-0.00222) (2.60,0.00148) (2.70,0.01470) (2.80,-0.00587) (2.90,-0.00853)
(3.00,0.00650) (3.10,0.00674) (3.20,0.00205) (3.30,-0.00198) (3.40,-0.01130)
(3.50,0.01159) (3.60,0.00597) (3.70,-0.00726) (3.80,-0.00526) (3.90,0.00310)
(4.00,0.01106) (4.10,-0.00462) (4.20,-0.00635)
};
\addplot[only marks,mark=|,mark size=2.2pt,black] coordinates {
(0.6931,-0.0303) (1.0986,-0.0303) (1.3863,-0.0303) (1.6094,-0.0303) (1.9459,-0.0303)
(2.0794,-0.0303) (2.1972,-0.0303) (2.3979,-0.0303) (2.5649,-0.0303) (2.7726,-0.0303)
(2.8332,-0.0303) (2.9444,-0.0303) (3.1355,-0.0303) (3.2189,-0.0303) (3.2958,-0.0303)
(3.3673,-0.0303) (3.4340,-0.0303) (3.4657,-0.0303) (3.6109,-0.0303) (3.7136,-0.0303)
(3.7612,-0.0303) (3.8501,-0.0303) (3.8918,-0.0303) (3.9703,-0.0303) (4.0775,-0.0303)
(4.1109,-0.0303) (4.1589,-0.0303)
};
\node[font=\tiny] at (axis cs:0.6931,-0.0258) {2};
\node[font=\tiny] at (axis cs:1.0986,-0.0258) {3};
\node[font=\tiny] at (axis cs:1.6094,-0.0258) {5};
\node[font=\tiny] at (axis cs:1.9459,-0.0258) {7};
\node[font=\tiny] at (axis cs:2.3979,-0.0258) {11};
\node[font=\tiny] at (axis cs:2.8332,-0.0258) {17};
\node[font=\tiny] at (axis cs:3.3673,-0.0258) {29};
\node[font=\tiny] at (axis cs:3.8501,-0.0258) {47};
\node[font=\footnotesize,anchor=north east] at (axis cs:4.17,0.0452) {$C(0)=Q$};
\end{axis}
\end{tikzpicture}
\caption{\label{fig:C}Centered prime-return correlation. The solid
curve is $C(u)$ computed from Eq.~\eqref{eq:C}. Circles give the
partial sum $\sum_{j\le2000}w_j\cos(\gamma_ju)$ of
Eq.~\eqref{eq:specC}. Ticks mark the returns $u=\log p^k$,
where the derivative jumps by $\log p/p^{k/2}$. The dashed line
marks $C(0)$. The inequality $C(u)\le C(0)$ at every time is
equivalent to RH.}
\end{figure}

\begin{corollary}\label{co:bounded}
RH holds if and only if $C(u)\le C(0)$ for every real $u$. In
that case $|C(u)|\le C(0)$.
\end{corollary}
\begin{proof}
If $C(u)\le Q$ for all $u$, then \eqref{eq:C} and $S\ge0$ give
$R(x)\le d_0+Q\sqrt x$ for $x\ge1$, and
Lemma~\ref{lem:explicit} gives RH. Under RH, \eqref{eq:specC} is
$\sum_\rho c_\rho\cos(\gamma u)$ with positive weights of total
mass $Q$.
\end{proof}

\begin{theorem}\label{th:CAR}
RH is equivalent to positive semidefiniteness of all the matrices
$[C(u_i-u_j)]$ at finitely many real times. It is also equivalent
to the existence of a nonnegative self-adjoint one-particle generator
$K_{\rm ph}$ and a vector $f$, $\|f\|^2=Q$, whose fermionic
fields $\Gamma(f)=a(f)+a^\dagger(f)$ satisfy
\begin{equation}
 \tfrac12\big\{\Gamma(e^{-iuK_{\rm ph}}f),\Gamma(f)\big\}=C(u)I .
 \label{eq:CAR}
\end{equation}
In any such realization the spectral measure of $f$ has mass
$2m_\gamma/(1/4+\gamma^2)$ at each positive zero ordinate
$\gamma$ of multiplicity $m_\gamma$.
\end{theorem}
\begin{proof}
On RH the sum in \eqref{eq:specC} is a cosine transform of precisely
the positive measure stated, of total mass $Q$. Take
$K_{\rm ph}$ to be multiplication by frequency in its $L^2$
space, with $f=1$. The canonical anticommutation relations give
$\{\Gamma(g),\Gamma(f)\}/2=\re\langle f,g\rangle I$, proving
\eqref{eq:CAR}. Conversely the spectral theorem writes the right
side as a cosine transform of a finite positive measure.
Symmetrizing the measure on $\R$ and integrating a squared
modulus proves matrix positivity. The two-time minor alone then
gives $|C(u)|\le Q$. By \eqref{eq:C} this gives
$R(x)\le d_0+Q\sqrt x$, hence RH by Lemma~\ref{lem:explicit}.
Uniqueness of the Fourier transform of finite measures determines
the symmetrized measure, and its restriction with doubled mass
determines the positive-frequency measure.
\end{proof}

This theorem imposes a sharp physical requirement. The local generator
$|H_0|$ has zero spectral projection onto every singleton by
Proposition~\ref{pr:boost}. It cannot furnish the vector in
\eqref{eq:CAR}. In the clock $u=2\eta$, a positive atom at
$\gamma$ would correspond to local energy
$2\hbar c\gamma/\xi_0$. The arithmetic data impose a further
short-time condition. For $0<u<\log2$, $T(e^u)=0$ and
\[
 \frac{\dd}{\dd u}S(e^u)=1-\frac{\operatorname{arctanh}(e^{-u})}{e^{-u}}.
\]
Integration in \eqref{eq:C} gives, as $u\to0$,
\begin{equation}
 C(0)-C(u)=\tfrac12|u|\log\frac1{|u|}-\kappa|u|
            +O\big(u^2\log\tfrac1{|u|}\big),
 \label{eq:cusp}
\end{equation}
with $\kappa=\tfrac12(\log2\pi+\eul-1)=0.7075463\ldots$
Thus a realizing vector has finite norm and infinite first spectral
moment. A finite moment would give
$|C(0)-C(u)|\le |u|\int\lambda\,\dd\nu_f(\lambda)$.
When \eqref{eq:CAR} holds the induction amplitude ${\cal A}_q$ of
Eq.~\eqref{eq:Aq} in Sec.~\ref{sec:step} also has the operator
expression
\[
 {\cal A}_q=\sqrt q\,\|(I-\sqrt q\,e^{-i(\log q)K_{\rm ph}})f\|^2 .
\]
The construction of this spectral measure uses RH. A realization of
the same correlation by an independently specified Majorana coupling
would supply the induction estimate. Suzuki studies related positive
kernels for $\zeta$ \cite{suzuki}.

\section{Operator closure }\label{sec:closureinverse}
\subsection{The boundary measure and its self-adjoint realization}
\label{sec:boundaryclosure}
The arithmetic response $C$ determines a closed one-particle
operator. The spectral mass of its probe gives the undamped
contribution to the equal-time correlation. Write
$a=1/2$, $w_\rho=\rho-a$, and take $c_\rho,Q$ from
\eqref{eq:center}. Introduce the one-sided Laplace transform
\[
 L(z)=\int_0^\infty e^{-zu}C(u)\dd u,\qquad \re z>a.
\]
At $z=\sigma+it$, $L$ resolves the frequency $t$ of the
return signal with exponential time weight $e^{-\sigma u}$.
In this half-plane \eqref{eq:C} defines the integrand from the
prime returns. Indeed, with $s=a+z$, termwise integration gives
\[
\begin{aligned}
 L(z)&=\frac{\sum_{\ell\ge2}\Lambda(\ell)\ell^{-s}}{s(s-1)}
 -\frac1{(s-1)^2}+\frac{1+\eul}{s-1}\\
 &\quad-\frac{d_0}{s}
 +\sum_{j\ge1}\frac1{2j(2j+1)(s+2j)} .
\end{aligned}
\]
Both series converge absolutely for $\re s>1$.
The meromorphic continuation of $L$ is
\[
 L(z)=\frac{zQ-\xi'(a+z)/\xi(a+z)}{z^2-a^2}.
\]
Indeed \eqref{eq:specC} gives $L(z)=\sum_\rho c_\rho/(z-w_\rho)$
by absolute convergence. This sum also converges normally on compact
sets away from its poles. Pair $w_\rho$ and $-w_\rho$, for which
$c_\rho=(a^2-w_\rho^2)^{-1}$ is the same. The paired logarithmic
derivative of the even entire function $\xi(a+z)$ is
$\sum_\rho z/(z^2-w_\rho^2)$, by the Hadamard product
\cite{Titchmarsh}. Partial fractions give the displayed formula,
including removability at $z=\pm a$. In particular $L$ is odd
and real on the real axis wherever it is regular.

\begin{theorem}[Boundary measure and operator closure]\label{th:boundaryclosure}
For real $\phi\in C_c^\infty(0,\infty)$, the formula
\[
 \int\phi(t)\dd\nu_\partial(t)=
 \frac{2}{\pi}\lim_{\epsilon\downarrow0}
       \int_0^\infty \phi(t)\re L(\epsilon+it)\dd t
\]
defines a finite positive measure. In $L^2(\nu_\partial)$ start
with multiplication by $t$ on compactly supported $L^2$ functions.
Its closure $K_\partial$ is nonnegative and self-adjoint, with
\[
\begin{gathered}
 \mathcal D(K_\partial)=\{g\in L^2(\nu_\partial):tg\in L^2(\nu_\partial)\},\\
 (K_\partial g)(t)=t g(t).
\end{gathered}
\]
Each reconstructed mode acquires phase $e^{-itu}$ under
$e^{-iuK_\partial}$, so $t$ is its undamped return frequency.
For $f(t)=1$, set $C_\partial(u)=\re\langle f,e^{-iuK_\partial}f\rangle$
and $\Delta=Q-\|f\|^2$. Then
\[
\begin{gathered}
 \Delta=\sum_{\re(\rho)\ne 1/2}\re\frac{1}{\rho(1-\rho)},\\
 \Delta=0\ \Longleftrightarrow\ C_\partial=C\text{ on }\R
 \ \Longleftrightarrow\ \mathrm{RH}.
\end{gathered}
\]
\end{theorem}
\begin{proof}
On a fixed compact frequency interval isolate the finitely many poles
$i\gamma$ on the axis. In a sufficiently thin neighbourhood of the
interval, subtract their principal parts. The analytic remainder has
purely imaginary boundary values, by oddness, conjugation and the reality
of the residues. Each pole with $\gamma>0$ and multiplicity $m_\gamma$
contributes the Poisson kernel
\[
 \frac{m_\gamma}{a^2+\gamma^2}
 \frac{\epsilon}{\epsilon^2+(t-\gamma)^2}.
\]
Consequently $\nu_\partial$ has the atom
$2m_\gamma/(a^2+\gamma^2)$ at each positive ordinate of a zero
on the critical line, with no other mass. There is no zero at $s=1/2$,
as follows from the alternating zeta series there. The total mass is
finite by \eqref{eq:explicit}. The boundary calculation uses compact
tests before taking the total mass, so it requires no globally
pole-free vertical strip next to the axis.

Truncating $g$ to $1/j\le t\le j$ converges in the graph norm on
the displayed maximal domain. There multiplication is symmetric and
$(t\pm i)^{-1}$ maps $L^2(\nu_\partial)$ into the domain and
inverts $K_\partial\pm i$. This proves self-adjointness and the core
assertion. The spectral theorem gives $C_\partial(u)=\int\cos(tu)
\dd\nu_\partial(t)$. Subtracting its value at zero from $Q$ gives
the stated defect. Every summand is strictly positive since, for
$\rho=\beta+i\gamma$,
\[
 \re c_\rho=
 \frac{\beta}{\beta^2+\gamma^2}
 +\frac{1-\beta}{(1-\beta)^2+\gamma^2}>0 .
\]
Vanishing of the defect therefore gives RH, and
\eqref{eq:specC} then gives equality of the two correlations.
Conversely, equality at $u=0$ suffices.
\end{proof}

The boundary measure of the continued arithmetic response records
the undamped frequencies. A zero off the critical line has
$\re w_\rho\ne0$ in \eqref{eq:specC}. Its contribution to the
equal-time mass is positive and it gives no boundary atom. The
scalar multiplication operator has spectral multiplicity one.
Multiplicities of zeta zeros occur in the spectral weights of the probe.
An independent microscopic coupling must therefore reproduce the
normalization of the probe as well as its frequencies.

The size of the unaccounted mass has an unconditional bound. Put
$H=3\cdot10^{12}$. The verified height \cite{pt} and the positive
ordinate tail $e_H$ from Appendix~\ref{app:check}
\cite{bpt} give, at every real $u$,
\[
\begin{gathered}
 0\le\Delta\le 2e_H<2.96\cdot10^{-12},\\
 |C(u)-C_\partial(u)|\le 2e_H\cosh(u/2).
\end{gathered}
\]
Indeed $|c_\rho|\le|\gamma|^{-2}$. Pairing the real exponents
$\beta-a$ and $a-\beta$ at a fixed $|\gamma|$ gives the stated
hyperbolic cosine. The loss bound for time zero is absolute. The
bound at nonzero time grows with $|u|$. In particular
$\Delta/Q<6.5\cdot10^{-11}$.

Table~\ref{tab:atoms} lists the lowest atoms of $\nu_\partial$,
computed from certified zeros. The first zero carries $21.6\%$ of
$Q$, the first ten carry $58.6\%$, and the first $10^4$ carry
$99.42\%$. The convergence is slow because
$w_j\simeq2/\gamma_j^2$. This is the spectral form of the expansion
\eqref{eq:cusp}: the probe has finite norm and infinite mean
frequency. In the units of Sec.~\ref{sec2} the lowest atom is the
boost frequency $\omega=2\gamma_1=28.269$. For the observer at
$\xi_0$, whose wedge temperature is
$\hbar c/(2\pi k_{\rm B}\xi_0)$ (Appendix~\ref{orig:A}), its local
energy is $4\pi\gamma_1=177.6$ in units of that temperature. The
fermionic occupation $(e^{2\pi\omega}+1)^{-1}$ of the state
\eqref{bell} at this frequency is $7.2\cdot10^{-78}$. The modes
selected by the arithmetic data lie far above the thermal scale of
the horizon.

\begin{table}[!htbp]
\caption{\label{tab:atoms}Return frequencies $t_j=\gamma_j$ of the
atoms of $\nu_\partial$, their weights
$w_j=2/(\tfrac14+\gamma_j^2)$, and their cumulative fraction of
the equal-time correlation $Q=C(0)$. The corresponding Rindler
frequencies are $2t_j$ in boost time and $4t_j$ in null time.}
\begin{ruledtabular}
\begin{tabular}{rlll}
$j$ & $t_j$ & $w_j$ & $\sum_{i\le j}w_i/Q$\\
\hline
1 & 14.134725 & $9.99798\cdot10^{-3}$ & 0.21645\\
2 & 21.022040 & $4.52308\cdot10^{-3}$ & 0.31437\\
3 & 25.010858 & $3.19595\cdot10^{-3}$ & 0.38356\\
4 & 30.424876 & $2.16001\cdot10^{-3}$ & 0.43032\\
5 & 32.935062 & $1.84337\cdot10^{-3}$ & 0.47023\\
10 & 49.773832 & $8.07205\cdot10^{-4}$ & 0.58605\\
$10^2$ & 236.524230 & $3.57501\cdot10^{-5}$ & 0.86531\\
$10^3$ & 1419.422481 & $9.92674\cdot10^{-7}$ & 0.96884\\
$10^4$ & 9877.782654 & $2.04980\cdot10^{-8}$ & 0.99417\\
\end{tabular}
\end{ruledtabular}
\end{table}

\subsection{The indefinite metric of the full return operator}\label{sec:krein}
Let $\mathcal Z$ denote the set of distinct nontrivial zeros,
$m_\rho$ their multiplicities, $w_\rho=\rho-\tfrac12$, and
$\rho'=1-\bar\rho$. Then $m_{\rho'}=m_\rho$,
$w_{\rho'}=-\bar w_\rho$, and $c_{\rho'}=\bar c_\rho$.
The prime side of \eqref{eq:C} defines the Hermitian form
\begin{equation}
 [f,g]_C=\iint f(u)\overline{g(v)}\,C(u-v)\dd u\dd v
 \label{eq:kreinform}
\end{equation}
on $C_c^\infty(\R)$. A negative index is the supremum of the
dimensions of subspaces on which a form is negative definite.
For a probe pulse $f$, the quadratic value $[f,f]_C$ combines
the correlation at the time differences within that pulse. A negative
value excludes a positive spectral interpretation for that response.

\begin{theorem}[Krein closure of the return form]\label{th:krein}
The negative index of \eqref{eq:kreinform} is
$\kappa=\#\{\rho\in\mathcal Z:\re(\rho)>1/2\}$.
In the sequence space with Hilbert norm
$\|a\|^2=\sum_{\mathcal Z}m_\rho|c_\rho||a_\rho|^2$, its indefinite
inner product and the closed boost generator are
\[
 \begin{aligned}
 [a,b]&=\sum_{\mathcal Z}m_\rho c_\rho a_\rho\overline{b_{\rho'}},\\
 (Ha)_\rho&=E_\rho a_\rho,\qquad E_\rho=-iw_\rho .
 \end{aligned}
\]
Time translation $e^{-itH}$ shifts the probe profile by
$f(u)\mapsto f(u+t)$. A mode with $\rho=\beta+i\gamma$ acquires
the factor $e^{-i\gamma t}e^{-(\beta-1/2)t}$. Thus $\gamma$
is its oscillation frequency per unit return time, and the reflected
pair has opposite growth rates in this indefinite realization.
The operator $H$ is self-adjoint for this indefinite metric. Its
spectrum consists of one simple eigenvalue for each distinct zero.
An off-line pair yields two conjugate nonreal energies on neutral
vectors. In particular, positivity of \eqref{eq:kreinform} on every
compact time window is equivalent to RH.
\end{theorem}
\begin{proof}
With the Hilbert product linear in its first entry, define
$(Ja)_\rho=(\bar c_\rho/|c_\rho|)a_{\rho'}$.
This linear metric operator pairs reflected modes when forming the
return signal. The Majorana conjugation on spinors is antiunitary.
The unitary involution $J$ gives $[a,b]=\langle a,Jb\rangle$.
It acts as $+1$ on each on-line coordinate and has eigenvalues
$\pm1$ on each pair exchanged by the reflection. This computes
the index in the sequence space.

For a prepared pulse $f$, the map $(Ef)_\rho=\int f(u)e^{w_\rho u}\dd u$
gives its complex amplitude in mode $\rho$. The Hilbert norm controls
the convergence of these amplitudes.
Since $|\re w_\rho|\le1/2$ and
$\sum m_\rho|c_\rho|<\infty$, \eqref{eq:specC} gives
$[f,g]_C=[Ef,Eg]$. For density in the Hilbert norm, suppose that
$b$ is orthogonal to the range. The continuous sum
$\sum a_\rho e^{w_\rho u}$ vanishes as a distribution on $\R$,
where $a_\rho=m_\rho|c_\rho|\bar b_\rho$ belongs to $\ell^1$.
Taking the one-sided Laplace transform for $\re z>1/2$ gives
$\sum a_\rho/(z-w_\rho)=0$. This series extends meromorphically:
on compact sets disjoint from its poles it converges normally.
Each residue $a_\rho$ vanishes, proving density. Every finite
negative definite Gram matrix retains its sign under sufficiently
small perturbations of the basis. Density therefore gives the stated
index also on the test functions.

The maximal domain of $H$ is $\{a:(E_\rho a_\rho)_\rho\in
\mathcal M\}$, with $\mathcal M$ the weighted sequence space.
It satisfies $JHJ=H^*$ including domains, by
$E_{\rho'}=\bar E_\rho$, and $E(if')=HEf$. For the graph core,
take a smooth compactly supported approximate identity
$\varphi_\delta$ and set $M_\delta a=( (E\varphi_\delta)_\rho
a_\rho)_\rho$. This filter uses a smooth probe pulse of duration
of order $\delta$. Sending $\delta$ to zero restores the response
at ideal time resolution. Its multipliers are uniformly bounded for
$0<\delta\le1$ and converge pointwise to one. At fixed $\delta$,
integration by parts gives rapid decay as $|\im w_\rho|\to\infty$,
uniformly on $|\re w_\rho|\le1/2$. Thus $M_\delta$ maps Hilbert
norm boundedly into graph norm and tends strongly to the identity
in graph norm on $\mathcal D(H)$. The Hilbert density already
proved, with $M_\delta Ef=E(\varphi_\delta*f)$, gives graph
density of the tests. Since $|E_\rho|\to\infty$, the maximal
diagonal operator has exactly the listed spectrum.
\end{proof}

The multiplicity of a zeta zero enters the weight $m_\rho|c_\rho|$.
The on-line coordinates form a positive invariant subspace. Folding
the pair $\pm\gamma$ in the cosine response there recovers the
measure of Theorem~\ref{th:boundaryclosure} and the multiplication
operator $|H|$ on the even subspace. The off-line component of the
constant sequence has quadratic value
$\sum_{\re(\rho)\ne1/2}m_\rho\re c_\rho=\Delta$, here summed over
distinct zeros.

Tests supported after translation in $[0,L]$ require only the prime
powers $p^k\le e^L$ in \eqref{eq:C}. The difference of two smooth
approximations to point masses at $0$ and $u$ has limiting
quadratic value $2(Q-C(u))$. Positivity gives the sign in the sharp
mean bound of Theorem~\ref{th:heightsharp}. For the sequence space
constructed here the Hilbert norm depends on the zeros. The resulting
index formula is consequently a spectral criterion for that sign.

\head{An exact Epstein test}
For $Z(s)=\sum'_{a,b\in\mathbb Z}(a^2+5b^2)^{-s}$, put
$Z(s)/2=\sum a(n)n^{-s}$, so $a(1)=1$. Enumerating the
representations of the divisors of $36$ gives
\[
 \begin{aligned}
 a(2)=a(3)&=0,\quad a(4)=1,\quad a(6)=2,\\
 a(9)&=a(36)=3 .
 \end{aligned}
\]
In the absolutely convergent logarithm far to the right the
coefficient of $36^{-s}$ is
\[
 a(36)-a(4)a(9)-a(6)^2/2=-2 .
\]
Thus $-Z'/Z$ has coefficient $2\log6$ at $6$ and
$-2\log36$ at $36$. Both are mixed products of primes. For
$\zeta$, the corresponding coefficients both vanish and every
prime-power coefficient is nonnegative.
The discriminant $-20$ has class number two, as seen from its
reduced forms $(1,0,5)$ and $(2,2,3)$. The theorems of Voronin
and Lee give infinitely many Epstein zeros in every fixed substrip
$1/2<\sigma_1<\re s<\sigma_2<1$ \cite{meanepstein}. In the
analogous zero-coordinate form these distinct pairs already give
infinite negative index. This test uses analytic zero distribution
and exact Dirichlet coefficients.

\subsection{The inverse Majorana construction and its Rindler transform}
\label{sec:inversemajorana}
The boundary measure of Theorem~\ref{th:boundaryclosure} determines
a Majorana Hamiltonian. The inverse transform below realizes its
spinors as solutions of the radial differential equation.
Normalizability depends on the inner product used for that realization.

\head{A single Majorana pair at each frequency}
Write $\mathfrak h=L^2(\nu_\partial)$, $K=K_\partial$, and use the
return clock $u$. On $\mathfrak h\otimes\C^2$ take
\[
\begin{gathered}
 h_{\rm M}=\sigma_2 K=
 \begin{pmatrix}0&-iK\\ iK&0\end{pmatrix},\\
 \mathcal D(h_{\rm M})=\mathcal D(K)\oplus\mathcal D(K),\qquad
 F=\binom{1}{0}.
\end{gathered}
\]
The two real components form one Majorana pair. Complex conjugation
$\mathcal C$ satisfies $\mathcal C h_{\rm M}\mathcal C=-h_{\rm M}$.
The constant spinors $s_\pm=(1,\pm i)^T/\sqrt2$ diagonalize the
self-adjoint operator into $K$ and $-K$. Its time evolution is a
real rotation at each frequency and gives exactly
\[
 \langle F,e^{-iu h_{\rm M}}F\rangle
   =\int_0^\infty\cos(tu)\dd\nu_\partial(t)=C_\partial(u).
\]
For an atom of mass $w_j$ at $t_j$ put
$e_j=\one_{\{t_j\}}/\sqrt{w_j}$, so $1=\sum_j\sqrt{w_j}\,e_j$.
In fermionic Fock space with $a_j=a(e_j)$, the two Majorana fields
are $c_{j1}=a_j+a_j^\dagger$, $c_{j2}=i(a_j^\dagger-a_j)$.
The normal-ordered Hamiltonian $\mathrm d\Gamma(K)$ adds the
frequencies of the occupied modes in return-time units. On a finite
set of modes its terms are $t_j a_j^\dagger a_j
=t_j(1+i c_{j1}c_{j2})/2$. The full operator is the nonnegative
energy multiplication operator on the occupation basis. The probe
field $X=\sum_j\sqrt{w_j}\,c_{j1}$ couples to the first real
component of each pair with amplitude $\sqrt{w_j}$. The sum converges
in operator norm, since a real
linear combination of these fields has squared norm equal to the sum
of the squared coefficients. The canonical anticommutation relations
give $\{X(u),X(0)\}/2=C_\partial(u)I$ in the Heisenberg evolution
\cite{arakiCAR}. This is the particle--hole doubling of a single
positive-frequency space. Each atomic pair is seen by the probe $F$,
even if the zero contributing its weight has multiplicity greater
than one.

\begin{theorem}[Inverse transform into massive Rindler modes]\label{th:inversebessel}
Let $r=m\xi=\mu_0e^v$, with $m,\xi_0>0$. For real $\omega$ set
\[
 G_\omega(v)=\sqrt{\frac{r\cosh(\pi\omega)}{\pi}}\,
 \begin{pmatrix}K_{1/2-i\omega}(r)\\-iK_{1/2+i\omega}(r)\end{pmatrix}.
\]
The assignment $U(e_js_\pm)=G_{\pm2t_j}$ is unitary onto the
completion of the finite span of these modes for the inner product
\[
 \langle g,h\rangle_{\rm av}=
 \lim_{A\to\infty}\frac1A\int_{-A}^0 g(v)^\dagger h(v)\dd v.
\]
On this completion the closure $H_{\rm av}$ of the differential
expression in \eqref{eq:boost} on the finite mode span satisfies
\[
\begin{gathered}
 H_{\rm av}U=2Uh_{\rm M},\qquad
 \mathcal D(H_{\rm av})=U\mathcal D(h_{\rm M}),\\
 U\mathcal C g=\sigma_3\overline{Ug}.
\end{gathered}
\]
\end{theorem}
\begin{proof}
The identity $rK'_\nu(r)=-\nu K_\nu(r)-rK_{\nu-1}(r)$, together
with $K_{-\nu}=K_\nu$, gives $H_0G_\omega=\omega G_\omega$
as a differential equation. The spinor decays at $+\infty$ and
$\sigma_3\overline{G_\omega}=G_{-\omega}$. At the horizon the
small-argument expansion and
$\lvert\Gamma(1/2+i\omega)\rvert^2=\pi/\cosh(\pi\omega)$
\cite{dlmf} give
\[
 G_\omega(v)=
 \binom{a_\omega e^{i\omega v}}{-i\overline{a_\omega}e^{-i\omega v}}
 +O_\omega(e^v),\qquad |a_\omega|^2=\tfrac12.
\]
Indeed the leading upper coefficient is
$2^{-1/2-i\omega}\Gamma(1/2-i\omega)\mu_0^{i\omega}/
 |\Gamma(1/2+i\omega)|$.
Integrating the exponentials gives zero mean at distinct frequencies.
The integrable errors also have zero mean. The norm
of each mode is one. Thus the assignment is an isometry on the two
finite spans. Passing to the completions and truncating in the graph
norm of $h_{\rm M}$ proves both operator assertions. Conjugation
interchanges the two signs in both representations. The factor 2 is
fixed by $u=2\eta$.
\end{proof}

In particular the inverse transform gives the state itself,
\[
 UF=\sum_j\sqrt{w_j/2}\,[G_{2t_j}+G_{-2t_j}],
\]
with convergence in the displayed mean norm. The physical radial
spinor at fixed boost is $\xi^{-1/2}G_\omega(\log(\xi/\xi_0))$.
The density $G_\omega^\dagger G_\omega$ per unit $v$ tends to one
at the horizon, oscillates up to the turning point
$m\xi\simeq|\omega|$ and decays beyond it (Fig.~\ref{fig:mode}). A
mode at the atom $\gamma_j$ is thus confined to
$m\xi\lesssim2\gamma_j$.
Each single mode has infinite norm for the usual Dirac inner
product $\int|\psi|^2\dd\xi$. In that Hilbert space an isometry
intertwining $2h_{\rm M}$ with the unperturbed $H_0$ would send
$e_js_+$ to a nonzero eigenvector, contrary to
Proposition~\ref{pr:boost}. The inverse of $U$ on its image is
given by Fourier coefficients in the same mean norm:
\[
 U^{-1}g=\sum_{j,\,\epsilon=\pm1}
    \langle G_{\epsilon2t_j},g\rangle_{\rm av}\,e_js_\epsilon .
\]
On the operator domain these coefficients are square summable with
weight $1+4t_j^2$. This also specifies the inverse for mode limits
that have no pointwise representative.

\head{Radial realization of the inverse momentum}
Apply the null-coordinate map \eqref{eq:nullconstraint} to a finite
sum of the boost modes. Their exponential decay at spatial infinity
ensures that the inverse in \eqref{eq:inversecomponents} exists
\cite[Sec.~7, Eqs.~(73)--(75)]{sierra}. Thus $A$ in
\eqref{majhamilton} acts on the resulting scalar mode profiles.

Explicitly the radial mode at $\tau=0$ is
$\Phi_\omega(x)=x^{-1/2}e^{i\omega v}G_{\omega,1}(v)$, and its
time factor is $e^{-2i\omega\tau}$. On each mode the
inverse-momentum integral gives
\[
 \int_r^\infty y^{1/2+i\omega}K_{1/2-i\omega}(y)\dd y
       =r^{1/2+i\omega}K_{1/2+i\omega}(r).
\]
Differentiation and exponential decay fix both sides, including the
integration constant. The inner product
$\lim_{A\to\infty}(2/A)\int_{\xi_0e^{-A}}^{\xi_0}
\overline{\Phi_\omega}\Phi_{\omega'}\dd x$
equals one for $\omega=\omega'$ and zero otherwise, by the horizon
expansion above. Writing $A_{\rm av}$ for the closure of the scalar
expression on this mean-norm mode span gives the unitary transform
\[
\begin{gathered}
 U_\tau(e_js_\pm)=\Phi_{\pm2t_j},\qquad
 A_{\rm av}U_\tau=4U_\tau h_{\rm M},\\
 \mathcal D(A_{\rm av})=U_\tau\mathcal D(h_{\rm M}).
\end{gathered}
\]
The coefficients of $U_\tau^{-1}\phi$ in the basis $e_js_\pm$
are the scalar mean inner products of $\Phi_{\pm2t_j}$ with
$\phi$. On the domain they are square summable with weight
$1+16t_j^2$. Taking $A=A_{\rm av}$ in \eqref{majhamilton} and
using $J_+$ gives $H_M J_+U_\tau=4J_+U_\tau h_{\rm M}$
on $\mathcal D(h_{\rm M})$.
For a positive return frequency $t_j$, the three frequencies are
$t_j$ in $u$, $2t_j$ in $\eta$, and $4t_j$ in $\tau$.
Thus $G_\omega$ has Bessel orders $1/2\pm i\omega$ with
$\omega=2t_j$, while the convention $\omega=E/2$ of the original
radial ansatz assigns $E=4t_j$. If the signed argument of $\Xi$
is $t_j$, the corresponding arguments in the two variables are
$\Xi(\omega/2)$ and $\Xi(E/4)$, respectively. Rescaling the
clock alone supplies no additional spectral mass.

\begin{figure}[t]
\centering
\begin{tikzpicture}
\begin{axis}[width=\columnwidth,height=5.2cm,xmin=0,xmax=1.4,ymin=0,ymax=4.7,
 xlabel={$m\xi/\omega$},ylabel={$G_\omega^\dagger G_\omega$},
 tick label style={font=\footnotesize},label style={font=\small}]
\addplot[black,dashed,thin] coordinates {(0,1) (1.4,1)};
\addplot[black!40,thin] coordinates {(1,0) (1,4.7)};
\addplot[refblue,thick] coordinates {
(0.002,1.0017) (0.004,0.9981) (0.006,1.0060) (0.008,0.9927) (0.010,0.9911)
(0.012,1.0027) (0.014,0.9888) (0.016,1.0053) (0.018,1.0119) (0.020,1.0076)
(0.022,0.9894) (0.024,0.9774) (0.026,1.0128) (0.028,1.0150) (0.030,0.9715)
(0.032,1.0316) (0.034,0.9712) (0.036,1.0301) (0.038,0.9679) (0.040,1.0392)
(0.042,0.9591) (0.044,1.0369) (0.046,0.9864) (0.048,0.9818) (0.050,1.0488)
(0.052,0.9586) (0.054,0.9989) (0.056,1.0524) (0.058,0.9582) (0.060,0.9794)
(0.062,1.0637) (0.064,0.9945) (0.066,0.9387) (0.068,1.0286) (0.070,1.0653)
(0.072,0.9694) (0.074,0.9360) (0.076,1.0316) (0.078,1.0783) (0.080,0.9924)
(0.082,0.9212) (0.084,0.9821) (0.086,1.0794) (0.088,1.0659) (0.090,0.9618)
(0.092,0.9130) (0.094,0.9851) (0.096,1.0847) (0.098,1.0878) (0.100,0.9916)
(0.102,0.9080) (0.104,0.9304) (0.106,1.0341) (0.108,1.1111) (0.110,1.0847)
(0.112,0.9815) (0.114,0.8989) (0.116,0.9120) (0.118,1.0085) (0.120,1.1060)
(0.122,1.1245) (0.124,1.0502) (0.126,0.9415) (0.128,0.8804) (0.130,0.9111)
(0.132,1.0112) (0.134,1.1116) (0.136,1.1456) (0.138,1.0921) (0.140,0.9857)
(0.142,0.8922) (0.144,0.8673) (0.146,0.9250) (0.148,1.0320) (0.150,1.1294)
(0.152,1.1652) (0.154,1.1217) (0.156,1.0218) (0.158,0.9157) (0.160,0.8548)
(0.162,0.8675) (0.164,0.9475) (0.166,1.0583) (0.168,1.1514) (0.170,1.1872)
(0.172,1.1515) (0.174,1.0597) (0.176,0.9490) (0.178,0.8629) (0.180,0.8340)
(0.182,0.8727) (0.184,0.9645) (0.186,1.0763) (0.188,1.1691) (0.190,1.2115)
(0.192,1.1900) (0.194,1.1123) (0.196,1.0040) (0.198,0.8993) (0.200,0.8305)
(0.202,0.8181) (0.204,0.8655) (0.206,0.9585) (0.208,1.0705) (0.210,1.1703)
(0.212,1.2311) (0.214,1.2370) (0.216,1.1871) (0.218,1.0947) (0.220,0.9837)
(0.222,0.8818) (0.224,0.8138) (0.226,0.7957) (0.228,0.8315) (0.230,0.9126)
(0.232,1.0205) (0.234,1.1311) (0.236,1.2204) (0.238,1.2697) (0.240,1.2689)
(0.242,1.2188) (0.244,1.1300) (0.246,1.0206) (0.248,0.9124) (0.250,0.8263)
(0.252,0.7786) (0.254,0.7780) (0.256,0.8243) (0.258,0.9090) (0.260,1.0167)
(0.262,1.1288) (0.264,1.2259) (0.266,1.2918) (0.268,1.3160) (0.270,1.2949)
(0.272,1.2323) (0.274,1.1385) (0.276,1.0285) (0.278,0.9194) (0.280,0.8278)
(0.282,0.7674) (0.284,0.7469) (0.286,0.7690) (0.288,0.8304) (0.290,0.9223)
(0.292,1.0319) (0.294,1.1442) (0.296,1.2442) (0.298,1.3189) (0.300,1.3588)
(0.302,1.3592) (0.304,1.3202) (0.306,1.2472) (0.308,1.1495) (0.310,1.0390)
(0.312,0.9293) (0.314,0.8333) (0.316,0.7621) (0.318,0.7239) (0.320,0.7228)
(0.322,0.7588) (0.324,0.8278) (0.326,0.9221) (0.328,1.0316) (0.330,1.1447)
(0.332,1.2498) (0.334,1.3361) (0.336,1.3951) (0.338,1.4211) (0.340,1.4120)
(0.342,1.3688) (0.344,1.2960) (0.346,1.2008) (0.348,1.0923) (0.350,0.9807)
(0.352,0.8763) (0.354,0.7885) (0.356,0.7250) (0.358,0.6915) (0.360,0.6906)
(0.362,0.7223) (0.364,0.7839) (0.366,0.8700) (0.368,0.9736) (0.370,1.0861)
(0.372,1.1985) (0.374,1.3021) (0.376,1.3888) (0.378,1.4520) (0.380,1.4872)
(0.382,1.4919) (0.384,1.4660) (0.386,1.4118) (0.388,1.3333) (0.390,1.2364)
(0.392,1.1283) (0.394,1.0165) (0.396,0.9091) (0.398,0.8133) (0.400,0.7357)
(0.402,0.6814) (0.404,0.6540) (0.406,0.6552) (0.408,0.6848) (0.410,0.7408)
(0.412,0.8195) (0.414,0.9162) (0.416,1.0246) (0.418,1.1383) (0.420,1.2504)
(0.422,1.3542) (0.424,1.4439) (0.426,1.5143) (0.428,1.5616) (0.430,1.5832)
(0.432,1.5782) (0.434,1.5469) (0.436,1.4915) (0.438,1.4150) (0.440,1.3218)
(0.442,1.2170) (0.444,1.1062) (0.446,0.9952) (0.448,0.8899) (0.450,0.7957)
(0.452,0.7172) (0.454,0.6584) (0.456,0.6220) (0.458,0.6099) (0.460,0.6226)
(0.462,0.6592) (0.464,0.7182) (0.466,0.7966) (0.468,0.8908) (0.470,0.9965)
(0.472,1.1089) (0.474,1.2231) (0.476,1.3341) (0.478,1.4371) (0.480,1.5277)
(0.482,1.6023) (0.484,1.6578) (0.486,1.6921) (0.488,1.7038) (0.490,1.6928)
(0.492,1.6597) (0.494,1.6059) (0.496,1.5337) (0.498,1.4462) (0.500,1.3468)
(0.502,1.2395) (0.504,1.1285) (0.506,1.0179) (0.508,0.9120) (0.510,0.8147)
(0.512,0.7296) (0.514,0.6596) (0.516,0.6073) (0.518,0.5745) (0.520,0.5622)
(0.522,0.5709) (0.524,0.6002) (0.526,0.6490) (0.528,0.7157) (0.530,0.7981)
(0.532,0.8934) (0.534,0.9985) (0.536,1.1101) (0.538,1.2247) (0.540,1.3386)
(0.542,1.4484) (0.544,1.5509) (0.546,1.6430) (0.548,1.7220) (0.550,1.7857)
(0.552,1.8323) (0.554,1.8607) (0.556,1.8702) (0.558,1.8606) (0.560,1.8323)
(0.562,1.7864) (0.564,1.7243) (0.566,1.6476) (0.568,1.5588) (0.570,1.4603)
(0.572,1.3547) (0.574,1.2451) (0.576,1.1342) (0.578,1.0250) (0.580,0.9204)
(0.582,0.8230) (0.584,0.7353) (0.586,0.6593) (0.588,0.5971) (0.590,0.5501)
(0.592,0.5193) (0.594,0.5055) (0.596,0.5089) (0.598,0.5295) (0.600,0.5668)
(0.602,0.6198) (0.604,0.6874) (0.606,0.7679) (0.608,0.8597) (0.610,0.9607)
(0.612,1.0688) (0.614,1.1816) (0.616,1.2968) (0.618,1.4120) (0.620,1.5250)
(0.622,1.6333) (0.624,1.7350) (0.626,1.8280) (0.628,1.9107) (0.630,1.9815)
(0.632,2.0391) (0.634,2.0826) (0.636,2.1112) (0.638,2.1245) (0.640,2.1225)
(0.642,2.1052) (0.644,2.0732) (0.646,2.0272) (0.648,1.9680) (0.650,1.8969)
(0.652,1.8152) (0.654,1.7245) (0.656,1.6263) (0.658,1.5225) (0.660,1.4148)
(0.662,1.3052) (0.664,1.1954) (0.666,1.0874) (0.668,0.9828) (0.670,0.8835)
(0.672,0.7909) (0.674,0.7065) (0.676,0.6316) (0.678,0.5674) (0.680,0.5148)
(0.682,0.4746) (0.684,0.4474) (0.686,0.4335) (0.688,0.4332) (0.690,0.4463)
(0.692,0.4728) (0.694,0.5121) (0.696,0.5638) (0.698,0.6271) (0.700,0.7012)
(0.702,0.7851) (0.704,0.8776) (0.706,0.9775) (0.708,1.0836) (0.710,1.1946)
(0.712,1.3090) (0.714,1.4256) (0.716,1.5428) (0.718,1.6593) (0.720,1.7737)
(0.722,1.8848) (0.724,1.9914) (0.726,2.0922) (0.728,2.1862) (0.730,2.2724)
(0.732,2.3500) (0.734,2.4181) (0.736,2.4762) (0.738,2.5238) (0.740,2.5603)
(0.742,2.5856) (0.744,2.5996) (0.746,2.6021) (0.748,2.5933) (0.750,2.5734)
(0.752,2.5427) (0.754,2.5017) (0.756,2.4507) (0.758,2.3905) (0.760,2.3218)
(0.762,2.2451) (0.764,2.1615) (0.766,2.0717) (0.768,1.9766) (0.770,1.8772)
(0.772,1.7743) (0.774,1.6691) (0.776,1.5624) (0.778,1.4552) (0.780,1.3484)
(0.782,1.2431) (0.784,1.1400) (0.786,1.0400) (0.788,0.9440) (0.790,0.8528)
(0.792,0.7670) (0.794,0.6873) (0.796,0.6144) (0.798,0.5488) (0.800,0.4909)
(0.802,0.4412) (0.804,0.4001) (0.806,0.3677) (0.808,0.3444) (0.810,0.3303)
(0.812,0.3254) (0.814,0.3298) (0.816,0.3434) (0.818,0.3662) (0.820,0.3979)
(0.822,0.4384) (0.824,0.4875) (0.826,0.5447) (0.828,0.6099) (0.830,0.6824)
(0.832,0.7621) (0.834,0.8484) (0.836,0.9408) (0.838,1.0388) (0.840,1.1419)
(0.842,1.2497) (0.844,1.3614) (0.846,1.4767) (0.848,1.5948) (0.850,1.7154)
(0.852,1.8378) (0.854,1.9615) (0.856,2.0859) (0.858,2.2105) (0.860,2.3348)
(0.862,2.4584) (0.864,2.5807) (0.866,2.7012) (0.868,2.8196) (0.870,2.9354)
(0.872,3.0483) (0.874,3.1578) (0.876,3.2637) (0.878,3.3655) (0.880,3.4632)
(0.882,3.5563) (0.884,3.6446) (0.886,3.7280) (0.888,3.8063) (0.890,3.8793)
(0.892,3.9468) (0.894,4.0089) (0.896,4.0654) (0.898,4.1162) (0.900,4.1614)
(0.902,4.2009) (0.904,4.2348) (0.906,4.2630) (0.908,4.2857) (0.910,4.3028)
(0.912,4.3146) (0.914,4.3210) (0.916,4.3222) (0.918,4.3184) (0.920,4.3097)
(0.922,4.2961) (0.924,4.2780) (0.926,4.2554) (0.928,4.2285) (0.930,4.1976)
(0.932,4.1627) (0.934,4.1242) (0.936,4.0821) (0.938,4.0367) (0.940,3.9882)
(0.942,3.9368) (0.944,3.8826) (0.946,3.8260) (0.948,3.7670) (0.950,3.7058)
(0.952,3.6427) (0.954,3.5779) (0.956,3.5115) (0.958,3.4437) (0.960,3.3747)
(0.962,3.3046) (0.964,3.2337) (0.966,3.1620) (0.968,3.0898) (0.970,3.0172)
(0.972,2.9442) (0.974,2.8712) (0.976,2.7981) (0.978,2.7251) (0.980,2.6524)
(0.982,2.5800) (0.984,2.5080) (0.986,2.4366) (0.988,2.3658) (0.990,2.2957)
(0.992,2.2264) (0.994,2.1580) (0.996,2.0905) (0.998,2.0240) (1.000,1.9585)
(1.002,1.8942) (1.004,1.8310) (1.006,1.7689) (1.008,1.7081) (1.010,1.6486)
(1.012,1.5903) (1.014,1.5333) (1.016,1.4777) (1.018,1.4233) (1.020,1.3704)
(1.022,1.3187) (1.024,1.2685) (1.026,1.2195) (1.028,1.1720) (1.030,1.1258)
(1.032,1.0809) (1.034,1.0374) (1.036,0.9952) (1.038,0.9543) (1.040,0.9147)
(1.042,0.8764) (1.044,0.8393) (1.046,0.8035) (1.048,0.7689) (1.050,0.7354)
(1.052,0.7032) (1.054,0.6721) (1.056,0.6422) (1.058,0.6133) (1.060,0.5855)
(1.062,0.5588) (1.064,0.5331) (1.066,0.5084) (1.068,0.4847) (1.070,0.4619)
(1.072,0.4400) (1.074,0.4190) (1.076,0.3989) (1.078,0.3796) (1.080,0.3611)
(1.082,0.3434) (1.084,0.3265) (1.086,0.3103) (1.088,0.2948) (1.090,0.2799)
(1.092,0.2658) (1.094,0.2523) (1.096,0.2394) (1.098,0.2270) (1.100,0.2153)
(1.102,0.2041) (1.104,0.1934) (1.106,0.1832) (1.108,0.1736) (1.110,0.1643)
(1.112,0.1555) (1.114,0.1472) (1.116,0.1392) (1.118,0.1317) (1.120,0.1245)
(1.122,0.1177) (1.124,0.1112) (1.126,0.1051) (1.128,0.0992) (1.130,0.0937)
(1.132,0.0884) (1.134,0.0834) (1.136,0.0787) (1.138,0.0743) (1.140,0.0700)
(1.142,0.0660) (1.144,0.0622) (1.146,0.0586) (1.148,0.0552) (1.150,0.0520)
(1.152,0.0490) (1.154,0.0461) (1.156,0.0434) (1.158,0.0408) (1.160,0.0384)
(1.162,0.0361) (1.164,0.0339) (1.166,0.0319) (1.168,0.0300) (1.170,0.0282)
(1.172,0.0265) (1.174,0.0248) (1.176,0.0233) (1.178,0.0219) (1.180,0.0205)
(1.182,0.0193) (1.184,0.0181) (1.186,0.0170) (1.188,0.0159) (1.190,0.0149)
(1.192,0.0140) (1.194,0.0131) (1.196,0.0123) (1.198,0.0115) (1.200,0.0108)
(1.202,0.0101) (1.204,0.0094) (1.206,0.0088) (1.208,0.0083) (1.210,0.0077)
(1.212,0.0072) (1.214,0.0068) (1.216,0.0063) (1.218,0.0059) (1.220,0.0055)
(1.222,0.0052) (1.224,0.0048) (1.226,0.0045) (1.228,0.0042) (1.230,0.0039)
(1.232,0.0037) (1.234,0.0034) (1.236,0.0032) (1.238,0.0030) (1.240,0.0028)
(1.242,0.0026) (1.244,0.0024) (1.246,0.0023) (1.248,0.0021) (1.250,0.0020)
(1.252,0.0018) (1.254,0.0017) (1.256,0.0016) (1.258,0.0015) (1.260,0.0014)
(1.262,0.0013) (1.264,0.0012) (1.266,0.0011) (1.268,0.0010) (1.270,0.0010)
(1.272,0.0009) (1.274,0.0008) (1.276,0.0008) (1.278,0.0007) (1.280,0.0007)
(1.282,0.0006) (1.284,0.0006) (1.286,0.0005) (1.288,0.0005) (1.290,0.0005)
(1.292,0.0004) (1.294,0.0004) (1.296,0.0004) (1.298,0.0003) (1.300,0.0003)
(1.302,0.0003) (1.304,0.0003) (1.306,0.0003) (1.308,0.0002) (1.310,0.0002)
(1.312,0.0002) (1.314,0.0002) (1.316,0.0002) (1.318,0.0002) (1.320,0.0002)
(1.322,0.0001) (1.324,0.0001) (1.326,0.0001) (1.328,0.0001) (1.330,0.0001)
(1.332,0.0001) (1.334,0.0001) (1.336,0.0001) (1.338,0.0001) (1.340,0.0001)
(1.342,0.0001) (1.344,0.0001) (1.346,0.0001) (1.348,0.0001) (1.350,0.0000)
(1.352,0.0000) (1.354,0.0000) (1.356,0.0000) (1.358,0.0000) (1.360,0.0000)
(1.362,0.0000) (1.364,0.0000) (1.366,0.0000) (1.368,0.0000) (1.370,0.0000)
(1.372,0.0000) (1.374,0.0000) (1.376,0.0000) (1.378,0.0000) (1.380,0.0000)
(1.382,0.0000) (1.384,0.0000) (1.386,0.0000) (1.388,0.0000) (1.390,0.0000)
(1.392,0.0000) (1.394,0.0000) (1.396,0.0000) (1.398,0.0000)
};
\end{axis}
\end{tikzpicture}
\caption{\label{fig:mode}Density per unit optical coordinate
$v=\log(\xi/\xi_0)$ for the massive Rindler mode $G_\omega$ of
Theorem~\ref{th:inversebessel} at the lowest atom,
$\omega=2\gamma_1=28.269$. The horizon limit is one (dashed line),
equal to the mean norm $\langle G_\omega,G_\omega\rangle_{\rm av}$.
The density decays beyond the turning point $m\xi=\omega$
(vertical line).}
\end{figure}

\subsection{Normalizable modes and a reconstruction from finitely many prime moments}
\label{sec:inversephysical}
In the Dirac Hilbert space, a trace-class perturbation of the wedge
Hamiltonian realizes the reconstructed dynamics on a reducing subspace.

\begin{proposition}[A trace-class Majorana coupling]\label{pr:inversecoupling}
For each $\epsilon>0$ there exist a self-adjoint trace-class
perturbation $V_\epsilon$, $\|V_\epsilon\|_1<\epsilon$, and a
Majorana-conjugation-preserving isometry $W$ into
$L^2(\R,\dd v;\C^2)$, such that $H_\epsilon=H_0+V_\epsilon$
on $\mathcal D(H_0)$ obeys
\[
\begin{gathered}
 H_\epsilon W=2Wh_{\rm M}\quad\hbox{on }\mathcal D(h_{\rm M}),\\
 W^*e^{-i\eta H_\epsilon}W=e^{-2i\eta h_{\rm M}} .
\end{gathered}
\]
\end{proposition}
\begin{proof}
Choose a real even $\chi\in C_c^\infty(-1,1)$ with $\|\chi\|_2=1$.
Put $\chi_j(v)=\ell_j^{-1/2}\chi((v-v_j)/\ell_j)$ on disjoint
intervals in $v<0$. For $\omega_j=2t_j$ take
\[
\begin{gathered}
 h_{j,+}(v)=\frac{\chi_j(v)}{\sqrt2}
 \binom{e^{i\omega_j(v-v_j)}}{e^{-i\omega_j(v-v_j)}},\\
 h_{j,-}=\mathcal C_Rh_{j,+},\qquad
 \mathcal C_Rg=\sigma_3\bar g .
\end{gathered}
\]
The opposite signs are orthogonal because their integrated overlap
is a sine against an even weight. Distinct intervals are disjoint.
Therefore $W(e_js_\pm)=h_{j,\pm}$ defines the required isometry.
The two residuals $r_{j,\pm}=(H_0\mp\omega_j)h_{j,\pm}$ have norm
at most $\delta_j=\|\chi'\|_2/\ell_j+\mu_0e^{v_j+\ell_j}$.
For an explicit choice, put $d_j=\epsilon2^{-j}/12$ for $j\ge1$
and $\ell_j=2\|\chi'\|_2/d_j$. Starting with $a_0=0$, take the
right endpoint $b_j=\min\{a_{j-1}-1,\log(d_j/(2\mu_0))\}$,
the centre $v_j=b_j-\ell_j$ and the left endpoint $a_j=b_j-2\ell_j$.
Then $\delta_j\le d_j$ and $6\sum_j\delta_j<\epsilon$.

Put $P=WW^*$ and $R=\sum_{j,\pm}|r_{j,\pm}\rangle
\langle h_{j,\pm}|$. The trace norms of the rank-one terms are
summable and $R=RP$. By symmetry of $H_0$ the compression
$PR$ is self-adjoint, as is seen by taking matrix elements on
these orthonormal packets. Define
\[
 V_\epsilon=-R-R^*+PR .
\]
Then $V_\epsilon h_{j,\pm}=-r_{j,\pm}$ and
$\|V_\epsilon\|_1\le3\|R\|_1\le6\sum_j\delta_j$.
Moreover $\mathcal C_R R\mathcal C_R=-R$ and
$\mathcal C_R P\mathcal C_R=P$, so the perturbation preserves
the Majorana symmetry. Bounded self-adjoint perturbation gives
self-adjointness on $\mathcal D(H_0)$. The closed span of the
packets is a reducing subspace since they are eigenvectors of
$H_\epsilon$. On this span the graph domain consists exactly of
the square-summable coefficients weighted in squared norm by
$1+\omega_j^2$. This proves the two identities.
\end{proof}

The kernel of $V_\epsilon$ couples points within individual packet
intervals. Its choice uses the prescribed frequencies. On the
reducing subspace the inverse is $W^*$. Its coefficients in the
$e_js_\pm$ basis are $\langle h_{j,\pm},g\rangle_{L^2(\dd v)}$.
For the bare
boost, Duhamel's formula gives the finite-time comparison
\[
 \left|\langle WF,e^{-i\eta H_0}WF\rangle-C_\partial(2\eta)\right|
 \le |\eta|\epsilon\,\|F\|^2 .
\]
For a local Hermitian matrix potential integrable on compact
$v$-intervals and integrable in norm near $-\infty$,
the perturbed Hamiltonian still has empty point spectrum. The
interaction-picture coefficient in the proof of
Proposition~\ref{pr:boost} remains integrable at the horizon.

\head{Sixty-four effective modes from the prime data}
The moments $m_j$ at $s=3/2$ in Sec.~\ref{sec:operatormoments}
also determine a finite frequency model. Its compression is the
Gaussian quadrature associated with the Taylor coefficients and the
total mass $Q$ \cite{golubwelsch}. Let
\[
 M_N=(m_{r+s})_{0\le r,s<N},\quad
 B_N=(m_{r+s+1})_{0\le r,s<N}.
\]
If $B_N>0$ and $M_N-B_N>0$, define
\[
\begin{gathered}
 J_N=M_N^{-1/2}B_NM_N^{-1/2},\quad
 K_N=(J_N^{-1}-I)^{1/2},\\
 q_N=M_N^{1/2}(1,0,\ldots,0)^T.
\end{gathered}
\]
Here $J_N$ compresses the spectral variable $b=(1+t^2)^{-1}$
to the measured moment space, and $K_N$ converts its eigenvalues
back to effective frequencies. Then $0<J_N<I$, and
$\sigma_2K_N$ is a Majorana Hamiltonian
with probe $(q_N,0)^T$. It reproduces $m_0,\ldots,m_{2N-1}$
exactly. Equip the polynomials of degree below $N$ with the inner
product defined by $M_N$. In the orthonormal coordinates
$p\mapsto M_N^{1/2}\operatorname{coeff}(p)$, $J_N$ represents
the form with matrix $B_N$. For $\deg p<N-1$ it acts as
multiplication by the polynomial variable $b$. Hence
$\langle q_N,J_N^{r+s}q_N\rangle=m_{r+s}$ for $0\le r,s<N$.
Inserting $J_N$ between the two vectors corresponding to $b^{N-1}$
gives $m_{2N-1}$. These first $N$ polynomial vectors span, so
$q_N$ is cyclic and each of the $N$ eigenvalues is simple with
positive probe weight. For $N=1$ the reconstructed frequency is
$\sqrt{Q/m_1-1}=24.96642311778\ldots$. Its second moment falls
short of $m_2$ by $m_2-m_1^2/Q>1.66\cdot10^{-7}$.

For $N=64$, interval elimination certifies both positive matrices
and bounds the two lowest frequencies of $K_{64}$ in $
 (14.13472514173,\,14.13472514174)$ and $(21.02203963877,\,21.02203963878).$
The ancillary calculation uses real Taylor coefficients of the
completed zeta function at $3/2$. The numbers of positive pivots of
$(1+s^2)B_N-M_N$ at the four displayed endpoints are $0,1,1,2$.
Congruence identifies them with the number of frequencies less than
$s$. Intervals for every frequency and probe weight are supplied
in \texttt{anc/majorana64.csv}. In orthonormal polynomial coordinates
the probe is $\sqrt Q\,(1,0,\ldots,0)^T$. For a unit trial vector
$v$, trial eigenvalue $b_*$, and residual norm
$\delta=\|(J_N-b_*)v\|$, the interval
$[b_*-\delta,b_*+\delta]$ contains an eigenvalue. All 64 intervals
are disjoint. If $d$ is the distance from $b_*$ to the other
intervals, the error of the weight $Q|v_0|^2$ is at most
$Q\delta/d$, by projecting the residual onto the complementary
eigenspaces and using the sine of the angle between the unit vectors.
Proposition~\ref{pr:inversecoupling} applies to these 64 frequencies.
The resulting finite-rank perturbation
gives a normalizable model determined by these moments and the
choice of packets. Its 128 moment identities concern the resolvent.
Its time correlation has a quadratic expansion at zero, whereas
$C(0)-C(u)$ has the $\tfrac12|u|\log(1/|u|)$ term of
Sec.~\ref{sec:CAR}.

If both matrix inequalities hold for an unbounded sequence of $N$,
the vector measures of $J_N$ have mass $Q$ on $[0,1]$.
A weakly convergent subsequence then realizes every $m_j$.
Theorem~\ref{th:operatormoments} gives the full correlation \footnote{We compute the $64$ frequencies and their weights using the Arb
ball-arithmetic library and the SciPy numerical linear algebra
library. In Arb, we set the working precision to 6144 bits,
substitute $y=1+h$ in the expression for $F(y)$ in
Sec.~\ref{sec:operatormoments}, and expand in powers of $h$
through degree 126. If $c_j$ is the coefficient of $h^j$,
then $m_{j+1}=(-1)^j c_j$. Together with $m_0=Q$, this gives
all 128 input moments.}.

We place $m_{r+s}$ in row $r$, column $s$ of a
$65\times64$ array, indexing from zero. Symmetric elimination
of this array gives the diagonal and off-diagonal entries of
$J_{64}$, each enclosed in an interval. 
For each returned
eigenvalue $b_k$ and normalized eigenvector $v_k$, we compute
the frequency $t_k=\sqrt{1/b_k-1}$ and the weight
$w_k=Q\,|(v_k)_1|^2$, where $(v_k)_1$ is the first component.
Finally, we use 192-bit interval arithmetic to bound the eigenpair
residuals and spectral gaps and thereby obtain the certified
intervals reported in \texttt{anc/majorana64.csv}.

\subsection{From the equal-time mass to the full Majorana correlation}
\label{sec:masspropagation}
Let $F=(1,0)^T$ be the probe in Sec.~\ref{sec:inversemajorana}.
The prime powers give an exact recursion in the return time. For
$u\ge0$ set
\[
\begin{gathered}
 B(u)=e^{u/2}(1+\eul-u)+e^{-u/2}\big[S(e^u)-d_0\big],\\
 Y(u)=C(u)-B(u).
\end{gathered}
\]
On inserting $T$ from \eqref{eq:center} into \eqref{eq:C} one obtains
\[
 Y(u)=2\sum_{\substack{p\ \mathrm{prime},\ k\ge1\\p^k\le e^u}}
          \frac{\log p}{p^{k/2}}\,
             \sinh\!\left(\frac{u-k\log p}{2}\right).
\]
Thus $Y=0$ up to the first return at $\log 2$.
Between successive returns $Y''=Y/4$. At $v=k\log p$ the
function is continuous and $Y'(v+)-Y'(v-)=\log p/p^{k/2}$.
The homogeneous initial-value problem and this jump determine the
next interval uniquely. There are finitely many jumps before each
finite time. This induction constructs the arithmetic correlation
even across the points where its derivative jumps.

\begin{proposition}[Propagation of the mass defect]\label{pr:masspropagation}
Suppose that every zero off the critical line has $|\im\rho|>h$,
where $h>0$. For the actual probe and every real return time,
\begin{multline*}
 \left|C(u)-\langle F,e^{-iu h_{\rm M}}F\rangle\right|\\
 \le \Delta\sqrt{\cosh^2(u/2)+h^{-2}\sinh^2(u/2)} .
\end{multline*}
The choice $h=4$ follows analytically from the prime explicit
formula. In particular $\Delta=0$ gives the full correlation at
every time, uniformly on each compact time interval.
\end{proposition}
\begin{proof}
The Mellin inversion of the prime sum in Lemma~\ref{lem:explicit}
gives \eqref{eq:specC}. The boundary construction in
Theorem~\ref{th:boundaryclosure} accounts for precisely its
critical-line terms. Group the rest by quartets under conjugation
and $\rho\mapsto1-\rho$, keeping multiplicities. For a representative
$\rho=1/2+d+ig$, $d,g>0$, write $c_\rho=A+iJ$. Then
\[
 A>0,\qquad
 \frac{J}{A}=\frac{2dg}{g^2+1/4-d^2},\qquad
 \left|\frac{J}{A}\right|\le \frac1g .
\]
One quartet contributes $4A$ to $\Delta$ and contributes
\[
 4A\!\left[\cosh(du)\cos(gu)
       -\frac{J}{A}\sinh(du)\sin(gu)\right]
\]
to the difference at time $u$. Cauchy--Schwarz against
$(\cos(gu),\sin(gu))$, followed by $0<d<1/2$ and $g>h$,
gives the claimed inequality after summation. Absolute local
uniform convergence follows from $\sum_\rho|c_\rho|<\infty$.
For the table-free choice of $h$, the conjugate pair at
$\rho,\bar\rho$ contributes $2\re c_\rho\ge 2/(1+g^2)$ to
$Q$, by $\re c_\rho=\re(1/\rho+1/(1-\rho))$. Hence
$g^2\ge 2/Q-1>16$. For the probe $F$ the Majorana pair at each
frequency contributes its cosine, as computed in Sec.~\ref{sec:inversemajorana}.
\end{proof}

There is also a construction of the common correlation using only
arithmetic moments, with frequencies fixed in advance. Write
\[
 a_{j,N}=\binom Nj D_{j,N-j}\quad(0\le j\le N).
\]
If $\Delta=0$, then
\[
\begin{aligned}
 C(u)&=\langle F,e^{-iu h_{\rm M}}F\rangle\\
 &=\lim_{N\to\infty}\sum_{j=1}^N a_{j,N}
              \cos\!\left(u\sqrt{N/j-1}\right).
\end{aligned}
\]
The convergence is uniform for $|u|\le U$, for each finite $U$.
Indeed, let $\mu$ be the image of $\nu_\partial$ under
$b=(1+t^2)^{-1}$. By Theorem~\ref{th:operatormoments}, $a_{j,N}$
is the integral against $\mu$ of the binomial probability
$\binom Nj b^j(1-b)^{N-j}$. For each $b>0$ a binomial random
variable $J_N$ with parameters $(N,b)$ satisfies $J_N/N\to b$
in probability. The functions $x\mapsto\cos(u\sqrt{x^{-1}-1})$
are equicontinuous near $b$ for $|u|\le U$, including at $b=1$.
Assign them value zero at $x=0$. Bounded convergence of the
expectations, uniformly in this time interval, now follows by splitting
according to $|J_N/N-b|<\delta$. Integration of the uniform error is
justified by its bound $2$ and $\mu(\{0\})=0$. The coefficients
use the absolutely convergent prime series at $s=3/2$ and the
explicit completion factors in Sec.~\ref{sec:operatormoments}.

\section{Analytic induction for the prime-return observable}\label{sec:analytic}
Multiplication of the cutoff $x=e^u$ by $q>1$ advances the
detector clock by $\log q$ (Sec.~\ref{sec4}). This choice localizes
the prime-power data in the recurrence for the one-sided bound of
Lemma~\ref{lem:explicit}. Averaging the increment at fixed time
resolution yields a signed trace with a fermion-parity factor.

\subsection{Induction with compact prime-power counts}\label{sec:step}
For $q>1$ write $D_qf(x)=f(qx)$ and
$\calL_q=(D_q-I)(D_q-qI)$. Define
\begin{equation}
\begin{gathered}
 w_q(t)=
 \begin{cases}
  q(1-1/t),&1\le t\le q,\\
  q^2/t-1,&q<t\le q^2,\\
  0,&\text{otherwise},
 \end{cases}\\
 W_q(x)=\sum_m\Lambda(m)w_q(m/x).
\end{gathered}
 \label{eq:weight}
\end{equation}
Applying $\calL_q$ to a single response in \eqref{eq:T} gives
$W_q=\calL_qT$. Constants and linear functions of $x$ are
annihilated, while $\calL_q[x\log x]=M_qx$ with
$M_q=q(q-1)\log q$. Consequently
\begin{equation}
\begin{gathered}
 F_q(x)=\frac{M_qx-W_q(x)}{\sqrt x}
       =-\frac{\calL_q R(x)}{\sqrt x},\\
 c_q=-\frac{\calL_q\sqrt x}{\sqrt x}=\sqrt q(\sqrt q-1)^2 .
\end{gathered}
 \label{eq:budget}
\end{equation}

\begin{theorem}[Propagation]\label{th:induct}
If $F_q(x)\le c_q$ for every $x\ge X_0\ge1$, then
$R(qx)-R(x)\le(\sqrt q-1)\sqrt x$ throughout that tail.
If also $R(t)\le\sqrt t$ on $[X_0,qX_0]$, induction gives
$R(x)\le\sqrt x$ for every $x\ge X_0$. Any finite eventual
upper bound on $F_q$ already implies RH.
\end{theorem}
\begin{proof}
Set $g_q(x)=[R(qx)-R(x)]/\sqrt{x}$. Subtraction in \eqref{eq:budget}
gives the exact identity
\begin{equation}
\begin{aligned}
 g_q(x)&=q^{-1/2}g_q(qx)+q^{-1}F_q(x)\\
       &=\frac1q\sum_{j=0}^{\infty}q^{-j/2}F_q(q^j x).
\end{aligned}
 \label{eq:greens}
\end{equation}
The finite sum through $j=J-1$ has remainder
$[R(q^{J+1}x)-R(q^Jx)]/(q^J\sqrt x)$, tending to zero by
Lemma~\ref{lem:explicit}. With $F_q\le c_q$, the geometric sum
gives $g_q\le c_q/[q(1-q^{-1/2})]=\sqrt q-1$.
Thus $R(qx)-\sqrt{qx}\le R(x)-\sqrt x$. Every $x\ge X_0$ is
$q^j t$ for some $X_0\le t\le qX_0$. This proves the forward
induction for real arguments. With $F_q\le B$, the same argument
first gives $g_q\le B/[q(1-q^{-1/2})]$. Summing the increments
over $q^j t$ gives $R(x)\le C\sqrt x$, using boundedness on
the compact base interval. Lemma~\ref{lem:explicit} applies.
\end{proof}

The values of $R$ on the initial interval give the initial data
for the recurrence. Each increment also depends on prime powers in
$(x,q^2 x)$. Truncating the series in \eqref{eq:greens} leaves a
terminal increment whose size requires a separate estimate.

\head{The accumulated margin}
Let $a(x)=\sum_{n\le x}\Lambda(n)/n$. Splitting the finite
sum at the previous cutoff gives an exact forward expression,
\begin{align}
 \sqrt x\,g_4(x)
 &=3x\,a(x)+\sum_{x<n\le4x}\Lambda(n)(4x/n-1)\nonumber\\
 &\quad-3x\log x-4x\log4+3(1+\eul)x. \label{eq:forwardreal}
\end{align}
For $4\le t<16$, set $d_j(t)=1-R(4^j t)/\sqrt{4^jt}$.
The recurrence and its finite iteration are
\begin{align}
 d_{j+1}(t)&=\tfrac12[1+d_j(t)-g_4(4^jt)],\nonumber\\
 d_n(t)&=2^{-n}d_0(t)+
     \sum_{j=0}^{n-1}2^{j-n}[1-g_4(4^jt)]. \label{eq:cash}
\end{align}
Here $d_0(t)$ denotes the initial margin, distinct from the constant
$d_0=\log(2\pi)$ in \eqref{eq:center}. The local excess
$g_4(4^jt)-1$ may be positive at some steps. A common upper bound for
$\sum_{j<n}2^{j-n}[g_4(4^jt)-1]$ over all $n$ and all base points
$t$ implies the one-sided estimate $R(x)\ll\sqrt{x}$, hence RH.
Propagation of $d_j(t)\ge0$ depends on the prime powers in
\eqref{eq:forwardreal} and on the slope $a(x)$.
Appendix~\ref{app:check} evaluates the recurrence on a finite interval.

\head{Sharp bound on the induction step}
With the correlation $C$ of Eq.~\eqref{eq:C}, define
\begin{equation}
 {\cal A}_q=\sqrt q\big[(q+1)Q-2\sqrt q\, C(\log q)\big].
 \label{eq:Aq}
\end{equation}
All data in this expression for ${\cal A}_q$ come from $q$,
elementary constants and primes up to $q$, by the closed form
\eqref{eq:Sclosed} of $S$.

\begin{theorem}[Sharp local upper limit]\label{th:sharp}
Under RH one has
$\limsup_{x\to\infty} F_q(x)={\cal A}_q$ and
$\limsup_{x\to\infty}|F_q(x)|={\cal A}_q$.
If RH fails then $\limsup F_q=+\infty$ and $\liminf F_q=-\infty$.
In particular $F_2(x)>c_2$ for arbitrarily large $x$,
unconditionally. The following criterion uses only a window on prime
powers:
\begin{equation}
\begin{gathered}
 \mathrm{RH}\quad\Longleftrightarrow\quad
 W_4(x)\ge12(\log4)x-2\sqrt x\\
 \text{for all sufficiently large real }x .
\end{gathered}
 \label{eq:raw4}
\end{equation}
\end{theorem}
\begin{proof}
On RH put $w=q^\rho$, so $|w|^2=q$ and
$(w-1)(w-q)=-w|w-1|^2$. Equations \eqref{eq:explicit}
and \eqref{eq:budget} give
\begin{equation}
 F_q(x)=\sqrt q\sum_\rho c_\rho |q^\rho-1|^2
             e^{i\gamma\log(qx)}
         +\frac{\calL_q S(x)}{\sqrt x}.
 \label{eq:fourstep}
\end{equation}
The last term tends to zero. The coefficients are nonnegative,
absolutely summable, and their sum is ${\cal A}_q$. For finitely many
ordinates, the points $(e^{i\gamma n})_\gamma$, $n\in\mathbb N$, have
a convergent subsequence on the corresponding torus. Differences of
subsequence indices can be arbitrarily large while the associated
phases approach the identity simultaneously. Taking
$x=e^t/q$ and bounding the remaining coefficient mass proves
attainment in the limit, without a hypothesis on linear relations
among ordinates.

The positive-side converse follows from Theorem~\ref{th:induct}.
The proof for both signs uses the Mellin transform of the compact weight,
\begin{equation}
 K_q(s)=\int_0^\infty w_q(t)t^{s-1}\dd t
       =\frac{(q^s-1)(q^s-q)}{s(s-1)} .
 \label{eq:K}
\end{equation}
This follows by integrating the two pieces, or by applying $\calL_q$
to the response transform. The singularities at 0 and 1 are removable,
and $K_q(1)=M_q$. For $\re s>1$ the transform of $W_q$ on
$(0,\infty)$ is $K_q(s)G(s)$, in the convention
$\int W_q(x)x^{-s-1}\dd x$. The part on $(0,1)$ is entire
since $W_q$ vanishes near zero. The pole at 1 cancels on subtracting
the main term on $[1,\infty)$. For $1/2<\re(\rho)<1$,
$1<|q^\rho|<q$, so $K_q(\rho)\ne0$. Either finite one-sided
bound yields a nonnegative tail with transform regular on the real
axis $s>1/2$. The argument in Lemma~\ref{lem:explicit} excludes
every such pole.

Finally $T(2)=0$ and $T(4)=\log2+\tfrac13\log3$.
Direct substitution in \eqref{eq:C} and \eqref{eq:Aq} gives the outward enclosures
in Table~\ref{tab:q}. Hence on RH the upper limit exceeds $c_2$
and lies strictly below $c_4$, proving both conclusions.
\end{proof}

\begin{table*}[t]
\centering
\caption{Sharp upper limits under RH and the threshold in
Theorem~\ref{th:induct}. Elementary-function bounds give the stated
enclosures.}
\label{tab:q}
\setlength{\tabcolsep}{8pt}%
\begin{tabular}{cccc}
\toprule
$q$ & $C(\log q)$ & ${\cal A}_q$ & $c_q$\\
\midrule
2 & $(-0.01733215,-0.01733214)$ & $(0.26530215,0.26530216)$ & $3\sqrt2-4$\\
4 & $(-0.00211129,-0.00211127)$ & $(0.47880441,0.47880444)$ & $2$\\
\bottomrule
\end{tabular}
\end{table*}

Algebraically $\calL_4=(D_2+I)(D_2+2I)\calL_2$. Thus a change
in ratio combines three displaced observations of the dyadic
difference. Its effect is quantitative. The dyadic impossibility
in Theorem~\ref{th:sharp} concerns $F_2\le c_2$ specifically.
An upper bound on the single increment $g_2$ or on a smoothed count
is a distinct condition.

\head{A dyadic alternative}
A positive average of the normalized $F_2(xe^v)$ over
$0\le v\le\log2$ also fits the dyadic margin under RH. Each mode
in \eqref{eq:fourstep} is multiplied by
$(e^{i\gamma\log2}-1)/(i\gamma\log2)$. In fact $|\gamma|>4$
unconditionally. By zero symmetry,
$Q=2\sum_\rho\beta/|\rho|^2$. Select in each symmetric orbit a
zero with $\beta\ge1/2$, and keep that conjugate pair in the sum.
It gives $\gamma^2\ge 4\beta/Q-\beta^2\ge2/Q-1>16$.
The smoothed upper limit in absolute value is therefore at most
${\cal A}_2/(2\log2)<0.192<c_2$. The averaging multiplier at a
general zero is $(e^{(\rho-1/2)\log2}-1)/((\rho-1/2)\log2)$,
which is nonzero for $\re(\rho)>1/2$.

\subsection{Prime-power increments on dyadic intervals}\label{sec:dyadic}
The initial prime count enters the dyadic increment explicitly. Put
$a_n=\sum_{m\le n}\Lambda(m)/m,\ y_n=\psi(n)$ and
$B_0(x)=x\log x-(1+\eul)x+\sqrt x$.
For $x\ge1$, $T(x)=\sup_n(xa_n-y_n)$: successive affine
functions differ by $\Lambda(n)(x/n-1)$.
To bound every one of these functions by $B_0$, use
\[
\begin{gathered}
 r(y)=\left[\frac{1+\sqrt{1+16y}}4\right]^2,\\
 F(y)=\log r(y)-\eul+\frac{1}{2\sqrt{r(y)}},\qquad
 L_n=F(y_n)-a_n .
\end{gathered}
\]
Indeed $tB_0'(t)-B_0(t)=t-\sqrt{t}/2$, whose inverse is $r$.
For large $y$ the minimum of $(B_0(x)+y)/x$ on any fixed
tail occurs at $x=r(y)$ and equals $F(y)$.
Since $y_n\to\infty$, and finitely many fixed lines are eventually
dominated, the eventual envelope inequality is equivalent to
$L_n\ge0$ for every sufficiently large $n$. Differentiation
gives $F'(y)=1/r(y)$, hence at $m=p^k,\ h=\log p$,
\[
 L_m-L_{m-1}=\int_{\psi(m-1)}^{\psi(m-1)+h}
                      [1/r(v)-1/m]\dd v .
\]
All other increments vanish. For example the increments at $13$
and $16$ belong respectively to the intervals $ (-0.00430665,-0.00430663)$, 
$(0.00264232,0.00264235)$.

These follow from $\psi(12)=\log27720$ and
$\psi(15)=\log360360$.
Both signs in fact recur arbitrarily late. By Littlewood's oscillation
theorem \cite{goldston}, the ratio $(\psi(x)-x)/\sqrt{x}$ has
infinite upper limit and negative infinite lower limit. On a positive
excursion take the last prime power $m\le x$. Then
\[
 \psi(m-1)=\psi(x)-\Lambda(m)>m-\sqrt m/2
\]
for large such $x$,
so the integrand is negative throughout that step. On a negative
excursion take the first prime power $m>x$. Bertrand's postulate
gives $m\le2x+2$, hence
\[
 \psi(m)=\psi(x)+\Lambda(m)<m-\sqrt m/2
\]
for large excursions
of this sign. The integrand is then positive. The inductive
condition is that each loss fit the margin already accumulated.

A form that preserves the initial count at the bottom of a block is
\begin{align}
 J(x)&=R(x)/x,\qquad b(x)=\sqrt{x}[J(x)-J(2x)],\nonumber\\
 b(x)&=\sqrt{x}\Big[\log2-\frac{\psi(x)}{2x}\nonumber\\
     &\qquad{}-\sum_{x<m\le2x}\Lambda(m)
       \Big(\frac1m-\frac1{2x}\Big)\Big].
 \label{eq:blockreport}
\end{align}
Indeed $J'(t)=(\psi(t)-t)/t^2$ almost everywhere. Integrating
each step of $\psi$ on $[x,2x]$ proves \eqref{eq:blockreport},
including a possible jump at either endpoint. If $b(x)\le
c_*=1-2^{-1/2}$ for every $x\ge x_0$, then telescoping the
differences of $J$ to infinity gives $J(x)\le1/\sqrt x$;
Lemma~\ref{lem:explicit} supplies the terminal value zero.
Conversely on RH, \eqref{eq:explicit} gives
$b(x)\le(1+2^{-1/2})Q+O(x^{-1/2})<c_*$ on a tail.
This block uses $\psi(x)$, whereas the three-scale $F_2$ in
Theorem~\ref{th:sharp} uses only the compactly supported weight
$w_2$. Its dyadic obstruction pertains to the latter weight.

\head{Joint lower bound for the two prime sums}
For $64\le x\le5.8\cdot10^{21}$, the bracketed prime contributions
in \eqref{eq:blockreport} satisfy the joint bound
\begin{multline}
 \frac{\psi(x)}{2x}+
       \sum_{x<m\le2x}\Lambda(m)\left(\frac1m-\frac1{2x}\right)\\
 \ \ge\ \log2-\frac{0.292}{\sqrt{x}}
 \ >\ \log2-\frac{c_*}{\sqrt{x}} .
 \label{eq:joint}
\end{multline}
Let $H$ and $e_H$ be the height and tail bound from
Appendix~\ref{app:check}. Write $C_l$ for the contribution of the
verified ordinates to the correlation $C$ of \eqref{eq:C}.
Pairing the unverified zeros as in that appendix gives
$|C(\log2)-C_l(\log2)|\le3e_H/\sqrt2$.
For the verified terms the sum of absolute coefficients in $b$
is at most
\begin{multline*}
 \sum_l c_\rho|1-2^{-1/2+i\gamma}|\\
 \le\sqrt{Q\,[\tfrac32Q-\sqrt2\,C(\log2)+3e_H]}<0.066 .
\end{multline*}
The inequality follows from Cauchy--Schwarz with positive weights.
Table~\ref{tab:q} gives $C(\log2)$ from $T(2)=0$.
For the remaining terms $|1-2^{\rho-1}|\le2$, giving
\[
 b(x)\le \frac{\log(2\pi)}{2y}+0.066+2e_H(y+1/y),
 \qquad y=\sqrt x .
\]
The trivial-zero part has nonpositive sign. The right side is convex
in $y$. Evaluation at $8$ and $\sqrt{5.8\cdot10^{21}}$, using
$e_H<1.48\cdot10^{-12}$, gives $b(x)<0.292$ throughout and
proves \eqref{eq:joint}. Beyond this interval the majorant contains the increasing term $2e_H\sqrt{x}$.

\subsection{Detector response from the radial Majorana intensity}\label{sec:heat}
Set $\tau=E/2$ and use the positive radial intensity
\[
\begin{gathered}
 w_\tau(y)=\frac{y}{\pi}\left|K_{1/2+i\tau}(y/2)\right|^2,\\
 {\cal H}_\tau(X)=\sum_{n\ge1}\Lambda(n)w_\tau(n/X).
\end{gathered}
\]
The lattice supplies the scales $n/X$, and the primitive weights
arise from $\Lambda=\mu*\log$. For $\re s>1$ the Mellin
transform in $X$, with measure $X^{-s-1}\dd X$, is
$G(s)\Phi_\tau(s)$. Substitution $y=n/X$ proves this identity.
The product integral for two modified Bessel functions gives, for
$\re s>0$ \cite{Olver,WatsonBessel1944},
\[
\begin{aligned}
 \Phi_\tau(s)&=\frac{2^{2s-1}}{\pi\Gamma(s+1)}
 \Gamma((s+2)/2)\Gamma(s/2)\\
 &\quad\times\Gamma((s+1)/2+i\tau)\Gamma((s+1)/2-i\tau).
\end{aligned}
\]
Indeed substitute $y=2z$ into the Mellin transform of $w_\tau$
and use the $z^{s}$-moment of
$K_{1/2+i\tau}(z)K_{1/2-i\tau}(z)$. This factor has no zeros in
$\re s>0$. At $\tau=0$ it is $\Gamma(s)$, since the elementary
half-order mode gives $w_0(y)=e^{-y}$. This is the generalized
radial profile at zero boost frequency in the full wedge. Sampling
its intensity on the positive transverse branch gives an exponential
response for the stated primitive weights.

\begin{theorem}[Induction for the zero-energy radial intensity]\label{th:heat}
Put $D_h(X)=X-{\cal H}_0(X)$. Without any individual zero
verification, $D_h(t)\le\sqrt t$ for every real
$100\le t\le160\,000$. Using the published verification of zeros
to $3\cdot10^{12}$ \cite{pt}, one has for every real
$100\le X\le\exp(8\cdot10^{12})$,
\[
 |D_h(4X)-D_h(X)|<0.283\sqrt X .
\]
Consequently $D_h(t)\le\sqrt t$ for all
$100\le t\le4\exp(8\cdot10^{12})$.
For every fixed real $\tau$, an eventual finite upper bound on
$(\Phi_\tau(1)X-{\cal H}_\tau(X))/\sqrt X$ is equivalent to RH.
\end{theorem}
\begin{proof}
For the equivalence the transform $G\Phi_\tau$ has a nonzero
residue at every right-hand zero. The small-$X$ contribution
of ${\cal H}_\tau$ is entire by exponential decay. Subtract the mean pole at 1 and
apply Lemma~\ref{lem:explicit}. Conversely on RH shift to
$\re s=1/4$. The factors $\Phi_\tau$ decay exponentially on
fixed strips there by Stirling, so the sum of moduli of the zero
residues is finite and the residual integral is $O_\tau(X^{1/4})$.

For the quantitative assertion shift $\Gamma(s)G(s)X^s$ to
$\re s=-1/2$. The pole at zero contributes $-\log(2\pi)$ to
${\cal H}_0$. With zeros repeated according to multiplicity,
\[
\begin{gathered}
 D_h(X)=\log(2\pi)+\sum_\rho \Gamma(\rho)X^\rho+\epsilon(X),\\
 |\epsilon(X)|<18X^{-1/2}.
\end{gathered}
\]
To verify the explicit remainder, the functional equation on this
line gives $|G(-1/2+it)|<12+2|t|$: the cotangent at
$-\pi/4+i\pi t/2$ has modulus one, $|G(3/2-it)|\le
\sum_{n\ge2}(\log n)n^{-3/2}<6$, and
$|\Gamma'/\Gamma(3/2-it)|<2+2|t|$ by the digamma difference
series at $3/2$. The gamma recurrence and
$|\Gamma(1/2+it)|^2=\pi/\cosh(\pi t)$ give
$|\Gamma(-1/2+it)|\le6 e^{-3|t|/2}$. Integration and division
by $2\pi$ give the claimed constant.

First use only the elementary zero-sum information
$\gamma^2\ge 2/Q-1>36$ from Sec.~\ref{sec:step} and
$\sum|c_\rho|<0.05$ from Appendix~\ref{app:check}.
For $\rho=\beta+i\gamma$, rotate the Euler integral for
$\Gamma(1+\rho)$ by $(5/4)\operatorname{sgn}\gamma$ radians.
The absolute integral is bounded by
$\Gamma(1+\beta)(\cos(5/4))^{-1-\beta}e^{-5|\gamma|/4}$.
Log-convexity gives $\Gamma(1+\beta)\le1$.
Since $\cos(5/4)>0.3$ and $|\rho(1-\rho)|\le1+\gamma^2$,
\[
 \frac{|\Gamma(\rho)|}{|c_\rho|}
 \le\frac{100}{9}(|\gamma|+1/|\gamma|)e^{-5|\gamma|/4}
 <0.038.
\]
The middle expression decreases for $|\gamma|\ge6$.
Pair $\beta$ with $1-\beta$, using equal $|c_\rho|$ and
$X^\beta+X^{1-\beta}\le X+1$. With $y=\sqrt X$, the absolute
normalized increment is at most
\[
 0.00475(y+y^{-1})+27/y^2<0.952\qquad(10\le y\le200).
\]
The last inequality follows by convexity and the two endpoints.
On the base $100\le X\le400$, the direct normalized upper bound
$\log(2\pi)/y+0.00095(y+y^{-1})+18/y^2$ is less than $0.39$.
As $\sqrt{4X}-\sqrt X=\sqrt X$, the forward step reaches every
$100\le t\le160\,000$ by analytic estimates alone.

Below $H=3\cdot10^{12}$, every zero is on the line. The function
$(t^2+1/4)\sqrt{\pi/\cosh(\pi t)}$ is decreasing for $t\ge4$,
as follows by logarithmic differentiation, and its value at 4 is
less than $0.077$. Hence $\sum_{\rm low}|\Gamma(\rho)|<0.077Q
<0.004$. For all $\rho=\beta+it$ with $|t|\ge H$, rotate the
Euler integral of $\Gamma(1+\rho)$ through
$(3/2)\operatorname{sgn}t$ radians. Since $\cos(3/2)>1/16$
and $\Gamma(1+\beta)\le1$ by log-convexity, this gives
\[
 |\Gamma(\rho)|\le256\,e^{-3|t|/2}/|t|
 \le300H e^{-3H/2}|c_\rho|.
\]
For the second inequality use
$|\rho(1-\rho)|\le t^2+1$ and monotonicity of
$(t+1/t)e^{-3t/2}$ on $t\ge H$.
In the stipulated range,
$H e^{-3H/2}\sqrt X\le H e^{-5\cdot10^{11}}$.
The normalized increment is therefore at most
\[
 3(0.004)+5(300)(0.05)H e^{-3H/2}\sqrt X+27/X<0.283,
\]
using $\sum|c_\rho|<0.05$ from Appendix~\ref{app:check}.
On $100\le X\le400$, the direct normalized bound is at most
$\log(2\pi)/\sqrt X+0.004+15H e^{-3H/2}\sqrt X+18/X<0.37$.
This is the full initial block. Iterating the increment proves
the theorem for real cutoffs. Its conclusion concerns precisely the
exponential detector $w_0$.
\end{proof}

\subsection{A fixed resolution in boost time}\label{sec:kernel}
Fix $q=4$, $L=\log4$ and $h_j=6L/(\pi^2j^2)$ for $j\ge1$.
Let $\sigma$ be the probability law of $\sum_{j\ge1}U_j$, where
the independent $U_j$ are uniform on $[0,h_j]$. Since
$\sum h_j=L$, this specifies a single measurement window on
$[0,L]$, at every scale $x$. Its support is the whole interval.
Indeed the support of each finite sum is an interval and the maximum
tail tends to zero. Every open interval meeting $(0,L)$ contains a
partial-sum subinterval which remains inside after any admissible tail.
Define the response and its mean coefficient by
\begin{gather}
\begin{aligned}
 w(t)&=\int_0^L e^{-v/2}w_4(te^{-v})\dd\sigma(v),\\
 K(s)&=\int_0^\infty w(t)t^{s-1}\dd t,
\end{aligned}\label{eq:w}\\
\begin{aligned}
 {\cal W}(x)&=\sum_n\Lambda(n)w(n/x),\\
 m_0&=K(1).
\end{aligned}
 \label{eq:W}
\end{gather}
Then $\supp w=[1,64]$, and $w$ is positive throughout $(1,64)$.
Indeed for each such $t$ an open set of allowed $v$ satisfies
$1<te^{-v}<16$. The transform factors as
\begin{equation}
 K(s)=K_4(s)A(s),\qquad
 A(s)=\prod_{j\ge1}
       \frac{e^{h_j(s-1/2)}-1}{h_j(s-1/2)} .
 \label{eq:prod}
\end{equation}
Each factor is interpreted by continuity when $s=1/2$. The product
converges normally on compact sets since there the $j$th factor is
$1+O(h_j)$. The only possible zeros of $A$ lie on $\re s=1/2$.
The only possible zeros of $K_4$ lie on $\re s=0,1$.
In particular $K(\rho)\ne0$ throughout $1/2<\re(\rho)<1$.

\begin{lemma}\label{lem:decay}
The function $w$ is smooth on $(0,\infty)$ with compact support.
On each bounded vertical strip $A(s)$ decays faster than every
inverse power of $|\im s|$. There is an absolute constant $a>0$
for this fixed window such that
\begin{equation}
 |A(s)|\le2\exp\{-a\sqrt{|s-1/2|}\},\qquad \re s\le0.
 \label{eq:left}
\end{equation}
\end{lemma}
\begin{proof}
Use the integral expression $h_j^{-1}\int_0^{h_j} e^{v(s-1/2)}\dd v$.
On a fixed strip the product of the bounds for any tail is at most
$\exp(L\max(0,\re s-1/2))$. Each of any fixed number of first
factors is $O_j((1+|\im s|)^{-1})$ on that strip. Mellin
inversion in $\log t$, differentiated any finite number of times,
proves smoothness. With $z=s-1/2$ and $\re s\le0$, each factor
has modulus at most $\min(1,2/(h_j|z|))$. For
$j\le\sqrt{h_1|z|/4}$ this is at most $1/2$. Multiplication
gives \eqref{eq:left} with $a=\tfrac12(\log2)\sqrt{h_1}$.
\end{proof}

\begin{proposition}\label{pr:sm}
Set $D(x)=m_0x-{\cal W}(x)$. Under RH,
$\limsup_{x\to\infty}|D(x)|/\sqrt x\le{\cal A}_4<0.479$.
Either finite one-sided bound for $D(x)/\sqrt x$ on a tail implies RH.
\end{proposition}
\begin{proof}
Changing variables in \eqref{eq:w} gives exactly
\[
 D(x)/\sqrt x=\int F_4(xe^v)\dd\sigma(v).
\]
Theorem~\ref{th:sharp} bounds the absolute value on RH. Conversely
the transform of ${\cal W}$ is $K(s)G(s)$. Its singularity at
1 is removed by the main term. Every right-hand nontrivial zero
would leave a pole, by \eqref{eq:prod}. Apply the nonnegative-tail
argument of Lemma~\ref{lem:explicit}.
\end{proof}

\subsection{Divisor truncation and the signed M\"obius sum}\label{sec:divisors}
Write $d(n)=\sum_{a\mid n}1$ for the divisor function and
$\mu(n)$ for the M\"obius function. Define
\begin{align}
 b&=e^{-2\eul},& \phi(t)&=b^{-1}w(t/b),\nonumber\\
 V(y)&=\sum_{n\ge1}(\log n)w(n/y),&
 J(y)&=\sum_{n\ge1}d(n)\phi(n/y),\nonumber\\
 E(y)&=V(y)-J(y),&
 P_\phi(y)&=\sum_{n\ge1}\phi(n/y).
 \label{eq:divdef}
\end{align}
The sums are finite, even for large $y$. The dilation $b$ matches
the linear terms of two Laurent expansions. It is the step that
allows one to discard a growing range of indices by absolute bounds.

\begin{theorem}[Removal of the small divisors]\label{th:div}
For every $B>0$ there exists $C_B>0$ such that, with
$Y_B(x)=x/(C_B\log(2x))^4$,
\begin{gather}
 D(x)=-\sum_{Y_B(x)<k\le32x}\mu(k)E(x/k)+O_B(x^{-B}),\label{eq:finaltail}\\
 \sum_{k\le Y_B(x)}|E(x/k)|=O_B(x^{-B})\qquad(x\to\infty).
 \label{eq:discard}
\end{gather}
With $B=1$ and any such $C_1$, the following statement is
equivalent to RH:
\begin{equation}
\begin{gathered}
 -\sum_{Y_1(x)<k\le32x}\mu(k)E(x/k)\ \le\ \sqrt{x}\\
 \text{for all sufficiently large real }x .
\end{gathered}
 \label{eq:open}
\end{equation}
On a fixed subinterval of the retained indices the sum of the
absolute values is $\gg x$.
\end{theorem}
\begin{proof}
The transform of $\phi$ is $b^{s-1}K(s)$. We use Dirichlet
convolution, the notation $\log(n)=\log n$ for an arithmetic
function, and the identities
$\mu*\log=\Lambda,\ \mu*d=\one$, where $\one(n)=1$. They follow
from $\log n=\sum_{a\mid n}\Lambda(a)$ and $d=\one*\one$.
Substituting \eqref{eq:divdef} and grouping $kn$ in the finite
sums yields the exact identity
\begin{equation}
 {\cal W}(x)-P_\phi(x)=\sum_{k\le32x}\mu(k)E(x/k).
 \label{eq:exacttail}
\end{equation}
Indeed the first possible nonzero logarithm in $V$ is $\log2$, and
$64b<32$ in $J$, since $b<1/2$. The compact weights vanish
at their endpoints. Hence $E(y)=0$ for $y\le1/32$.

The inverse Mellin integrands for $E(y)$ and $P_\phi(y)$, with
the common factor $y^s$ omitted, are respectively
\begin{equation}
 K(s)\big[-\zeta'(s)-b^{s-1}\zeta(s)^2\big],
 \qquad K(s)b^{s-1}\zeta(s).
 \label{eq:divM}
\end{equation}
They follow from absolute convergence on $\re s=2$. At $s=1+z$,
$-\zeta'(1+z)=z^{-2}+O(1)$ and
$b^z\zeta(1+z)^2=z^{-2}+(2\eul+\log b)z^{-1}+O(1)$.
The first integrand in \eqref{eq:divM} is entire. The second has the
single pole at 1 of residue $m_0$.

We record the dependence on the distance moved in the contour shift.
There are constants $H>1,C_0\ge1$, independent of the integer
$N\ge1$, such that on $s=-N+it$
\begin{multline}
 |-\zeta'(s)-b^{s-1}\zeta(s)^2|+|b^{s-1}\zeta(s)|\\
 \le H^{N+1}(N+1+|t|)^{2N+C_0}.
 \label{eq:FElimit}
\end{multline}
To verify it write
$\zeta(s)=2^s\pi^{s-1}\sin(\pi s/2)\Gamma(1-s)\zeta(1-s)$.
Both $\zeta(1-s)$ and its derivative are uniformly bounded.
Expand $\Gamma(1+N-it)=\Gamma(1-it)\prod_{j=1}^N(j-it)$.
The exponential in $|\Gamma(1-it)|^2=\pi |t|/\sinh(\pi|t|)$
cancels that from the sine or its derivative, leaving a polynomial
factor. In differentiating the gamma factor the logarithmic
derivative is $O(\log(2+N+|t|))$ by Stirling in the right
half-plane. Differentiating the sine directly also covers the
negative even integers. The factors $b^{-N-1}$ are absorbed in
$H^{N+1}$. This proves \eqref{eq:FElimit}.

Now shift \eqref{eq:divM} to $\re s=-N$. The horizontal sides vanish
by Lemma~\ref{lem:decay} and polynomial vertical-strip bounds for
$\zeta$ and its derivative. There $K_4(s)$ is uniformly bounded
and \eqref{eq:left} bounds the rest. For an integer $j\ge0$,
\[
 \int_0^\infty t^j e^{-a\sqrt t}\dd t
   =2\Gamma(2j+2)a^{-2j-2}.
\]
Use $(N+1+t)^j\le2^j((N+1)^j+t^j)$, take
$j=2N+\lceil C_0\rceil$, and apply Stirling to the last display.
It follows for $N\ge1,\ y\ge2$, with another fixed $H_1>1$, that
\[
 |E(y)|+|P_\phi(y)-m_0y|
       \le H_1^{N+1}y^{-N}N^{4N}.
\]
Fixed extra powers of $N$ have been absorbed into the exponential
in $N$. Choosing $N=\lfloor\delta y^{1/4}\rfloor$ with
$H_1\delta^4<1/2$ gives
\begin{equation}
 |E(y)|+|P_\phi(y)-m_0y| \ll e^{-a_1 y^{1/4}}\qquad(y\ge2)
 \label{eq:stretch}
\end{equation}
for a fixed $a_1>0$, increasing the implied constant on a bounded
interval. If $k\le Y_B(x)$, then $(x/k)^{1/4}\ge C_B\log(2x)$.
There are at most $x$ such indices for all large $x$.
Choose $a_1 C_B>B+1$ to obtain \eqref{eq:discard}.
Now \eqref{eq:exacttail} and \eqref{eq:stretch} give \eqref{eq:finaltail}.
Under RH, Proposition~\ref{pr:sm} leaves a strict margin below 1,
so \eqref{eq:open} holds. Conversely \eqref{eq:open} bounds
$D(x)/\sqrt x$ from above on a tail, which implies RH.

For the last assertion put $I=[(9/4)/64,(5/2)/64]$.
For $y\in I$, $64y\in[9/4,5/2]$, and the only nonzero term of
$V(y)$ is $(\log2)w(2/y)$. Also $64by<1$, since
$5/2<e^{2\eul}$, so $J(y)=0$. The function $w(2/y)$ has
a positive minimum on $I$. The condition $x/k\in I$ selects an
interval of length a fixed positive multiple of $x$ within
$(Y_B(x),32x]$ for large $x$. The identity
$|\mu(k)|=\sum_{r^2\mid k}\mu(r)$ gives
$\sum_{k\le t}|\mu(k)|=t/\zeta(2)+O(\sqrt t)$ by summing the
floors and the tail of $\sum r^{-2}$. Hence
$\sum_{x/k\in I}|\mu(k)E(x/k)|\gg x$.
\end{proof}

There is also a cutoff with no unspecified constant. Put
$Y_*(x)=x/(\log(2x))^5$. For every fixed $B$, $Y_*(x)\le
Y_B(x)$ eventually. Applying \eqref{eq:discard} to this smaller
range gives
\[
\begin{gathered}
 D(x)=-\sum_{Y_*(x)<k\le32x}\mu(k)E(x/k)+O_B(x^{-B})\\
 \text{for every } B>0 .
\end{gathered}
\]
Thus \eqref{eq:open} may be stated with the single explicit $Y_*$
at all scales.

The contribution of small $k$ is controlled in absolute value.
On the retained interval the termwise majorant exceeds the desired
scale by a factor of order $\sqrt{x}$. A uniform bound at the
square-root scale therefore requires cancellation in the signed sum
\eqref{eq:open}.

The prime number theorem gives an unconditional estimate for this
sum. For some $a_2>0$,
$\psi(t)-t=O(t\exp(-a_2\sqrt{\log t}))$ \cite{mv}. Partial
summation gives
\[
\begin{aligned}
 D(x)&=\frac1x\int_x^{64x}(\psi(t)-t)w'(t/x)\dd t\\
      &\ll_w x\exp(-a_2\sqrt{\log x}).
\end{aligned}
\]
The boundary terms vanish since $w(1)=w(64)=0$.
By \eqref{eq:finaltail} the same bound applies in absolute value to
the signed sum in \eqref{eq:open}. This is $o(x)$, while its ratio
to $\sqrt x$ in the displayed majorant grows without bound.

\subsection{Fermion parity and induction on prime modes}\label{sec:parity}
The trace in \eqref{eq:ret} counts repeated passages of a chiral
amplitude. There is also a many-body experiment for the signed sum
in \eqref{eq:open}. The interpretation of the M\"obius function
as fermion parity is due to Spector \cite{spector}. Here it is
applied to the finite energy window and the weight in the induction.
At each prime $p$ take a complex fermionic
oscillator, equivalently a pair of real Majorana modes
$\gamma_{p,1},\gamma_{p,2}$ with
$\{\gamma_{p,i},\gamma_{r,j}\}=2\delta_{pr}\delta_{ij}I$ and
$a_p=(\gamma_{p,1}+i\gamma_{p,2})/2$. On the fermionic Fock space
over $\ell^2(\{p\})$, choose the vacuum energy to be zero and set
\[
 H_{\rm F}=\sum_p(\log p)\,a_p^\dagger a_p,\qquad
 {\cal P}=(-1)^{N_{\rm F}}.
\]
The occupation $a_p^\dagger a_p$ distinguishes an empty prime
mode from an occupied one. $H_{\rm F}$ adds the assigned energies
$\log p$ of the occupied modes, and $\mathcal P$ reads their
total number modulo two as a sign. In the occupation basis every
finite occupied set $S$ gives a
squarefree $k=\prod_{p\in S}p$. Its energy is $\log k$ and
${\cal P}|S\rangle=\mu(k)|S\rangle$. These formulas prescribe a
nonnegative self-adjoint diagonal operator on its maximal domain.
There are only finitely many states below every finite energy.
For $\re s>1$, summing the basis gives
\[
 \tr e^{-sH_{\rm F}}=\frac{\zeta(s)}{\zeta(2s)},\qquad
 \tr({\cal P}e^{-sH_{\rm F}})=\frac{1}{\zeta(s)}.
\]
On the one-particle space the operator $h_{\rm F}$ assigns energy
$\log p$ to the prime state $|p\rangle$. In the half-plane
$\re s>1$ of the Euler product,
$e^{-s h_{\rm F}}$ is trace class and
\[
 \det(I-e^{-s h_{\rm F}})=\prod_p(1-p^{-s})=\frac1{\zeta(s)}.
\]
Taking logarithms in the region of absolute convergence and
differentiating recovers $\sum_{p,k}(\log p)p^{-ks}$.
For $0<\re s\le1$, the trace norm $\sum_p p^{-\re s}$ diverges
by comparison with $\sum_p1/p$. The reciprocal zeta function has a
meromorphic continuation, whereas the Fredholm determinant in this
representation requires $\re s>1$ \cite[Sec.~III.F]{v3}.

Taking the difference of parity-resolved measurement
outcomes gives the exact finite trace
\begin{multline}
 \sum_{Y_*(x)<k\le32x}\mu(k)E(x/k)\\
 =\tr\!\left[{\cal P}\,
    \one_{(\log Y_*(x),\,\log(32x)]}(H_{\rm F})\,
                   E(xe^{-H_{\rm F}})\right].
 \label{eq:parity}
\end{multline}
Here $Y_*(x)=x/(\log(2x))^5$ and $E=V-J$ as in
\eqref{eq:divdef}. The Hamiltonian assigns the prescribed energies
to the prime modes. Equation~\eqref{eq:parity} identifies the
arithmetic sum with the measured parity trace. A uniform estimate of
that trace is the arithmetic hypothesis in this representation.

\begin{theorem}[Exact induction over occupied and empty modes]\label{th:pair}
For $P\ge2$ let $P^\#=\prod_{p\le P}p$, and set
\[
 E_P(y)=\prod_{p\le P}(I-D_{1/p})E(y),\qquad D_c f(y)=f(cy).
\]
For all $x>0$,
\begin{equation}
 \sum_{k\ge1}\mu(k)E(x/k)=
 \sum_{\substack{m\le32x\$m,P^\#)=1}}\mu(m)E_P(x/m).
 \label{eq:sieveind}
\end{equation}
For $P=\sqrt{32x}$ and large $x$ this identity reduces exactly to
\[
 E_P(x)-\sum_{P<p\le32x} E_P(x/p).
\]
On a fixed subinterval of these remaining primes
$\sum_p |E_P(x/p)|\gg x/\log x$.
\end{theorem}
\begin{proof}
Start with no paired modes, when $E_P=E$.
At the next prime $r$, split squarefree states into those with
$r$ empty and those with $r$ occupied. The energy shift
$\log r$ changes $y$ to $y/r$ and the parity changes sign.
This replaces the previous weight by $E_P(y)-E_P(y/r)$.
Every sum is finite since $E(y)=0$ for $y\le1/32$.
This proves \eqref{eq:sieveind} by induction. After pairing every
prime up to $P=\sqrt{32x}$, a remaining squarefree $m>1$
can contain at most one prime, since two such primes have product
$>32x$.

Take $I=[9/256,5/128]$ as in Theorem~\ref{th:div}.
For $y\in I$ and every $d\ge2$, $y/d<1/32$, hence
$E_P(y)=E(y)\ge\delta>0$ independent of $P$.
The primes with $x/p\in I$ form the fixed-proportion interval
$[128x/5,256x/9]$ inside $(\sqrt{32x},32x]$ for all large
$x$. The ordinary prime number theorem gives $\asymp x/\log x$
primes there, proving the lower bound.
\end{proof}
This is an arithmetic induction with an exact update at every prime.
Even after the reduction to the vacuum and single unpaired modes,
a triangle inequality exceeds the square-root scale. To infer RH
one needs an upper bound on the \emph{negative} of the total parity
trace in \eqref{eq:parity}. The bound on its absolute contributions
does not determine that signed quantity.

\subsection{Analytic induction for a detector with a fixed response time}
\label{th:longstepsec}
The possible effect of high frequencies on the induction can be
estimated explicitly. This also gives a physical interpretation
of using a fixed resolution at successively later boost times.
Continue to write $L=\log4$. Set $M=5\cdot10^{11}$, an integer,
and $h=3L/(4M)$. Use the probability law $\sigma_\sharp$ of
\[
 V_\sharp=\tfrac14 V+U_1+\cdots+U_M,
 \quad V\sim\sigma,\quad U_j\sim{\rm Uniform}[0,h]
\]
with independent variables. Its support is $[0,L]$ at every cutoff.
This formal detector response is specified by its convolution. No
integration with $M$ independent numerical variables is required.
Write
\begin{align}
 R_\sharp(x)&=\int e^{-v/2}R(xe^v)\dd\sigma_\sharp(v),\nonumber\\
 A_\sharp(s)&=\int e^{(s-1/2)v}\dd\sigma_\sharp(v).
 \label{eq:newav}
\end{align}
\begin{theorem}[Finite forward propagation at fixed resolution]\label{th:large}
The averaged forward increment obeys
\begin{equation}
 \frac{R_\sharp(4x)-R_\sharp(x)}{\sqrt x}<0.146
 \qquad (4\le x\le\exp(10^{12})).
 \label{eq:big}
\end{equation}
In particular $R_\sharp(t)\le\sqrt t$ for every real
$4\le t\le4\exp(10^{12})$. The proof uses the published
finite-height theorem of Platt and Trudgian~\cite{pt}, with no
assumption on zeros above its verified height. An eventual
$O(\sqrt x)$ upper bound for this same fixed $R_\sharp$
is equivalent to RH.
\end{theorem}
\begin{proof}
For $z=s-1/2$, the transform in \eqref{eq:newav} equals
$A(1/2+z/4)[(e^{hz}-1)/(hz)]^M$.
On $0\le\re s\le1$ each box has modulus at most
$e^{h\max(0,\re z)}\min(1,2/(h|\im s|))$.
The first factor is at most $e^{(L/4)\max(0,\re z)}$. Hence for
$|\im s|\ge H=3\cdot10^{12}$,
\begin{equation}
 |A_\sharp(s)|\le e^{L/2}
      \left(\frac{2M}{(3L/4)H}\right)^M\le2\,3^{-M},
 \label{eq:atten}
\end{equation}
since $\log4>4/3$. On the critical line it is bounded by one.
All zeros up to $H$ lie on this line \cite{pt}. Their total
contribution to the normalized increment in \eqref{eq:big} is at
most $3Q$, since $|4^\rho-1|\le3$ and their positive weights
$c_\rho$ sum to at most $Q$. In the entire strip,
$|4^\rho-1|\le5$, $x^{\beta-1/2}\le\sqrt x$ for $x\ge1$,
and $\sum_\rho|c_\rho|\le(\sqrt{17}/4)Q<0.05$, as shown in
Appendix~\ref{app:check}. By absolute convergence the part above
$H$ is bounded by $10(0.05)\,3^{-M}\sqrt x$.
On the asserted interval $3^{-M}\sqrt x\le e^{-M(\log3-1)}$.
The constant in $R$ cancels in the increment.
The trivial-zero remainder is bounded by
$S(x)/\sqrt x\le[6\sqrt x(x^2-1)]^{-1}\le1/180$.
Adding the three bounds gives \eqref{eq:big}. On the initial interval
$4\le t\le16$, $4\le te^v\le64$ and
$R(te^v)\le\sqrt{te^v}$ by the elementary analytic bound in
Appendix~\ref{app:check}. Thus $R_\sharp(t)\le\sqrt t$
there. As $0.146<\sqrt4-1$, iterating the forward difference
reaches every real point in the claimed range.

All factors of $A_\sharp(s)$ are nonzero for $\re s>1/2$.
On a tail the Mellin transform of $R_\sharp$ is $A_\sharp(s)$
times the transform of $R$, up to an entire function: after the
substitution $y=xe^v$, the omitted integrals have endpoints in
a fixed compact positive interval. The argument of
Lemma~\ref{lem:explicit} applies in both directions, using
$|A_\sharp(\rho)|\le1$ on RH, or just averaging the upper
bound for $R$. This proves the equivalence.
\end{proof}
The resolution in \eqref{eq:newav} is part of the observable.
The factor $3^{-M}$ is a fixed number, while the factor
$x^{\beta-1/2}$ for a hypothetical right-hand zero grows without
bound. Equation \eqref{eq:atten} therefore explains both the long
finite induction and the obstruction to extending this estimate
to infinity or transferring it to the unsmoothed response.

\subsection{Sampling powers, weights and spin}\label{sec:powers}
The sampling exponent of $x$ and the exponent of the prime weight
are independent parameters of the compact response in \eqref{eq:w}.

\begin{proposition}\label{pr:samples}
For any fixed $j>0$, RH is equivalent to
$D(n^j)\le n^{j/2}$ for all sufficiently large integers $n$.
If RH fails, both signs of $D(n^j)/n^{j/2}$ are unbounded.
\end{proposition}
\begin{proof}
Let $\Theta=\sup_\rho\re(\rho)$. Smoothness gives
\begin{equation}
 D'(x)=O(x^{\Theta-1}),\qquad x\ge2 .
 \label{eq:deriv}
\end{equation}
Shifting the integrand $K(s)G(s)x^s$ to $\re s=-1$ gives a residue expansion.
The residues for ${\cal W}$ are $m_0 x$ and
$-K(\rho)x^\rho$ with multiplicity. There is no pole at zero.
The integral and its first derivative on the new line are
$O(x^{-1})$ and $O(x^{-2})$, respectively. The horizontal
bounds used in Lemma~\ref{lem:explicit} apply, and
$\sum_\rho (1+|\rho|)|K(\rho)|<\infty$ by Lemma~\ref{lem:decay}.
This proves \eqref{eq:deriv}. If RH fails, choose
$\max(1/2,\Theta-1/j)<a<\Theta$.
The pole of the transform at any zero with $\beta>a$ excludes both
one-sided bounds $\pm D(x)\le Cx^a$ on a tail, for every fixed $C$,
by the same nonnegative-tail argument. Each $x$ is at distance
$O_j(x^{1-1/j})$ from a point $n^j$. The variation there is
$O_j(x^{\Theta-1/j})=o(x^a)$ by \eqref{eq:deriv}, so both signed
excursions persist at sample points. Under RH the assertion follows
from Proposition~\ref{pr:sm} and its strict numerical margin.
\end{proof}

To use an exponent as a weight instead, put
\[
\begin{gathered}
 {\cal W}_{j}(x)=\sum_{n\ge1}\Lambda(n)n^{-j}w(n/x),\\
 D_j(x)=K(1-j)x^{1-j}-{\cal W}_j(x).
\end{gathered}
\]
Its Mellin transform before centering is $K(s)G(s+j)$.
For $j=1/2,3/2,5/2,7/2$, at every right-hand zero
the factor $K(\rho-j)$ is nonzero because the real part of $\rho-j$ lies in
$(0,1/2)$ or is strictly negative, away from each of
$0,1/2,1$. A contour shift to $\re s=-1-j$ gives
\[
\begin{gathered}
 D_j(x)=\sum_\rho K(\rho-j)x^{\rho-j}+O_j(x^{-1-j}),\\
 |D_j(x)|\ll_j x^{1/2-j}\quad\text{under RH}.
\end{gathered}
\]
The expansion holds unconditionally. The poles crossed are at
$s=1-j$ and at $s=\rho-j$, coming from the pole and the
nontrivial zeros of $\zeta(s+j)$, respectively. The first trivial
zero is further to the left. Each of these four one-sided bounds
$D_j(x)\le C_j x^{1/2-j}$ on a tail conversely implies RH by
the nonnegative-tail theorem, valid also at a negative abscissa.

\begin{center}
\begin{tabular}{cccc}
\toprule
Exponent & sampled scale & error scale & error scale\\
$j$ & $n^j$ & of $D(n^j)$ & of $D_j(x)$\\
\midrule
$1/2$& $\sqrt n$ & $n^{1/4}$ & $1$\\
$3/2$& $n^{3/2}$ & $n^{3/4}$ & $x^{-1}$\\
$5/2$& $n^{5/2}$ & $n^{5/4}$ & $x^{-2}$\\
$7/2$& $n^{7/2}$ & $n^{7/4}$ & $x^{-3}$\\
\bottomrule
\end{tabular}
\end{center}
At $j=0$ a sampled cutoff is always $n^0=1$, so $T(1)=0$
and ${\cal W}(1)$ is a single finite sum over $1<n<64$.
As a weight, $n^{-0}=1$ gives precisely ${\cal W}(x)$,
again finite at each $x$. Removing the cutoff from
$\sum\Lambda(n)$ gives a divergent sum. The exponents in the
table are parameters of the measurement. Spin statistics enters a
Majorana field through the real structure and fermionic
anticommutation relations. In $1+1$ dimensions a list of
half-integer exponents of the cutoff does not specify additional
one-particle spin multiplets.

\head{Coupling the exponent of the weight to the sampling power}
A prime-only observable uses $j$ in both places:
\[
\begin{gathered}
 S_\alpha(X)=\sum_{p\le X}(\log p)p^{-\alpha}(1-p/X),\\
 P_j(n)=S_{1/j}(n^j),\qquad \lambda_j=n^{j/2-1}.
\end{gathered}
\]
Let $L_p(s)=\sum_p(\log p)p^{-s}$ for $\re s>1$.
For $\re s>1/3$ its meromorphic continuation is
\[
 L_p(s)=G(s)-G(2s)-\sum_p\frac{(\log p)p^{-3s}}{1-p^{-2s}},
 \qquad G=-\zeta'/\zeta .
\]
The last series removes precisely the odd powers $p^k,\ k\ge3$;
it converges normally. The real residues at $1,1/2$ are
$1,-1/2$. Set $C_1=\lim_{s\to1}(L_p(s)-1/(s-1))$.
For each row the eventual sign of $P_j(n)-T_j(n)$, listed in the
last column, is equivalent to RH.
\begin{center}
\footnotesize\setlength{\tabcolsep}{3pt}%
\begin{tabular}{cclcc}
\toprule
$j$&$\alpha=1/j$& background $T_j(n)$ & $b_j$ & sign\\
\midrule
$1/2$&$2$&$L_p(2)-(\log\sqrt n+C_1+1)/\sqrt n$&$-2/3$&$-$\\
$3/2$&$2/3$&$(9/4)\sqrt n+L_p(2/3)$&$18/5$&$+$\\
$5/2$&$2/5$&$(25/24)n^{3/2}$&$-50/11$&$-$\\
$7/2$&$2/7$&$(49/60)n^{5/2}$&$-98/51$&$-$\\
\bottomrule
\end{tabular}
\end{center}
To verify the assertion, the inverse Mellin integrand is
$L_p(s)X^{s-\alpha}/[(s-\alpha)(1+s-\alpha)]$. The real pole
at $s=1/2$ gives
$b_j=-1/[2(1/2-\alpha)(3/2-\alpha)]$.
On RH move the contour to $\re s=5/12$;
the $L_p$ decomposition is regular there, with logarithmic
derivatives $O(\log^2(2+|\im s|))$ on the new line.
Use horizontal heights avoiding ordinates and half-ordinates as in
Lemma~\ref{lem:explicit}. In the $j=5/2$ row the kernel pole
at $s=2/5$ lies to the left of that line, consistently with
the stated error. The other real residues crossed give the table,
including the double pole when $\alpha=2,\ s=1$.
Writing $M_\alpha(n^j)=T_j(n)+b_j\lambda_j$ for their total
including $s=1/2$, one obtains
\begin{multline*}
 S_\alpha(X)-M_\alpha(X)=
 -X^{1/2-\alpha}\sum_\rho
 \frac{X^{i\gamma}}{(\rho-\alpha)(1+\rho-\alpha)}\\
       +O_\alpha(X^{5/12-\alpha}).
\end{multline*}
The squared denominator is
$(\gamma^2+1/4)^2+(2\gamma^2-1/2)(1-\alpha)^2+(1-\alpha)^4$.
Since $|\gamma|>4$, the total normalized oscillation has
upper limit in modulus at most $Q$, smaller than every $|b_j|$.

For the converse set $\Theta=\sup_\rho\re(\rho)>1/2$. Subtracting
$M_\alpha$ cancels all real poles in $\re s>1/2$ in the
unshifted $s$ variable. In terms of the Mellin variable for $X$
every right-hand zero leaves a pole at $\rho-\alpha$.
With $\max(1/2-\alpha,\Theta-\alpha-1/j)<d<\Theta-\alpha$,
the nonnegative-tail argument gives arbitrarily large signed
excursions of both signs relative to $X^d$ for the centered
sum. To transfer these to $X=n^j$, write
$S_\alpha'(X)=X^{-2}\sum_{p\le X}(\log p)p^{1-\alpha}$
almost everywhere. The classical truncated explicit formula gives
$\vartheta(t)=\sum_{p\le t}\log p=t+O(t^\Theta\log^2(2t))$
\cite{Titchmarsh,mv}. Partial summation then gives
\[
 (S_\alpha-M_\alpha)'(X)
     =O_\alpha(X^{\Theta-\alpha-1}\log^3(2X)).
\]
In detail, putting $d_1=\alpha-1$, the main term in the weighted
prime sum is $X^{1-d_1}/(1-d_1)$, replaced by $\log X$ at
$d_1=1$. For $d_1<\Theta$ the integral of the error is
$O(X^{\Theta-d_1}\log^3(2X))$, including constants. For
$d_1=1>\Theta$ the constant after $\log X$ is $C_1$, by
partial summation or the Laurent expansion of $L_p$ at 1;
the tail has the same error bound. Equality $d_1=\Theta$ costs
the third logarithm. The derivative of the $b_j$-term is smaller.
The sample gap $O_j(X^{1-1/j})$ therefore gives variation
$o(X^d)$. This proves both signs occur with unbounded normalized
size if RH fails. Finally, for the coupled family,
$\sup_{n\ge2} P_j(n)\le L_p(1/j)\to0$ as $j\downarrow0$
by dominated convergence. At $j=0$ the reciprocal exponent is
undefined. The zero-power cases of the fixed-weight experiment above
are well-defined.

\head{A cubic measurement with growing resolution}
At the cubic sampling times $v=3\log n$, the response is
\[
 A(v)=e^{-v/6}\{S_{1/3}(e^v)-(9/10)e^{2v/3}\}.
\]
Let $Y_u$ be normal with mean $u$ and variance $2\kappa u$,
for $u>0$. The averaged signal obeys the unconditional estimate
\begin{equation}
\begin{gathered}
 \mathbb E[A(Y_u)]=-18/7+O(e^{-u/20}),\\
 1/32\le\kappa\le1/4,\quad u\to\infty .
\end{gathered}\label{eq:cubicgauss}
\end{equation}
The constant can be chosen uniformly over the stated resolutions.
For the proof put $c=1.01$, $P_\kappa(z)=z+\kappa z^2$ and
$d(s)=1/[(s-1/3)(s+2/3)]$. Mellin inversion on $\re s=c$
and the Gaussian moment formula give the integral of
$L_p(s)d(s)\exp[uP_\kappa(s-1/2)]$, with factor $1/(2\pi i)$,
minus $(9/10)\exp[uP_\kappa(1/2)]$. Fubini holds since
$L_p$ is absolutely bounded there, $d(c+it)$ is integrable in
$t$, and the Gaussian has all exponential moments.
Choose $5<T<6$ and $0.41<\sigma_0<0.43$ so that the rectangle
with vertical sides $\sigma_0,c$ and horizontal sides $\pm T$
has no poles on its boundary. On the original vertical tails and
the horizontal segments, writing $a=\re s-1/2$, one has
\[
\begin{aligned}
 \re P_\kappa(s-1/2)&=a+\kappa(a^2-t^2)\\
 &\le 0.51+(0.51^2-25)/32<-1/20.
\end{aligned}
\]
On the new vertical $a\le-0.07,\ |a|<0.09$ gives the same strict
upper bound. The two tails and the three sides of the shifted contour therefore
contribute $O(e^{-u/20})$ altogether. There are no nontrivial zeros of
$\zeta(s)$ inside: the analytic bound $\gamma^2\ge 2/Q-1>36$
was derived in Sec.~\ref{sec:step}. A possible nonreal pole from
$G(2s)$ has $a<0,\ |t|>3$ and contributes $O(e^{-u/20})$
as well. Only finitely many such poles enter this fixed rectangle.
The residue at 1 cancels the subtracted mean. At $s=1/2$ the
residue is $-d(1/2)/2=-18/7$, proving \eqref{eq:cubicgauss}.
The Gaussian has logarithmic width proportional to $\sqrt u$.
Near $u=3\log n$ adjacent samples have spacing $O(1/n)$, so the
number contributing to the average grows with $u$. The resulting
error estimate applies to this Gaussian average.
\begin{table*}[t]
\caption{\label{tab:summary}Spectral and arithmetic estimates, with their hypotheses and ranges of validity.}
\small
\begin{tabular}{@{}>{\raggedright\arraybackslash}p{.19\linewidth}>{\raggedright\arraybackslash}p{.47\linewidth}>{\raggedright\arraybackslash}p{.28\linewidth}@{}}
\toprule
Quantity & Estimate or equivalence & Hypotheses and proof\\
\midrule
Wall and prime returns &
The wall counts positive modes. Prescribing $\ell_p=\tfrac12\log p$
gives $\log p$ at every return $k\log p$. &
Unconditional, Theorem~\ref{th:wall} and Eq.~\eqref{eq:ret}.\\[3pt]
Return correlation &
$C(u)\le C(0)$ at every time is equivalent to RH, and
$|C(u)-C_\partial(u)|\le2e_H\cosh(u/2)$ at every real time. &
Corollary~\ref{co:bounded}; finite-height input,
Sec.~\ref{sec:boundaryclosure}.\\[3pt]
Boundary mass &
$0\le\Delta<2.96\cdot10^{-12}$, using the verified zero height
$3\cdot10^{12}$. &
Finite-height input, Theorem~\ref{th:boundaryclosure} and
Appendix~\ref{app:check}.\\[3pt]
Real-degree curvature &
$I''(x)<-0.00059$ for
$14\,000\le x+1\le6.91899\cdot10^{12}$. &
Finite interval, Proposition~\ref{pr:twodegreecurvature}.\\[3pt]
Linear prime sum &
On RH, $\limsup|\mathcal A_n|<1.185$. An off-line zero gives a
positive exponential growth rate. &
Two implications, Theorem~\ref{th:movingprime}.\\[3pt]
Power-error test &
A test error $t^{3/4}$ in the Abel integral has exponential rate
$5/3$ in the degree. &
Unconditional, Proposition~\ref{pr:powerprofile}.\\[3pt]
Prime endpoints &
The lower bound \eqref{eq:primeend} holds for
$83\le x\le1.6\cdot10^{20}$ and fails at $x=661/8$. &
Finite-height input, Theorem~\ref{th:stependrange}.\\[3pt]
Point-count weights &
$P_{E,2}-P\to2\log2$ and $\mathcal A_{n,E,r}-\mathcal A_n\to0$.
Deleting primes with $\vartheta_{\mathcal D}(y)\sim c\sqrt y$,
$c>1/2$, reverses \eqref{eq:primeend}. &
Unconditional, Theorem~\ref{th:stepfirstterm} and
Propositions~\ref{pr:stepdeletion},~\ref{pr:movingweights}.\\
\bottomrule
\end{tabular}
\end{table*}
\section{Arithmetic content of spectral completeness}
\label{sec:arithstructure}
The recurrence in Sec.~\ref{sec:masspropagation} determines $C$
on every finite time interval from its prime-power jumps. Spectral
completeness requires the undamped modes to carry the entire
equal-time mass $C(0)=Q$. A prime sum and a moment-curvature
criterion characterize this requirement. For the prime sum, set
$\tau_n=n^{-1/2}$, $s_n(z)=1/2+i\tau_n+(1+z)/(1-z)$ and
\[
\begin{aligned}
 \mathcal A_n&=[z^n](1-z)(s_n(z)-1)\\
 &\quad\times\sum_{p\le4096^n}(\log p)p^{-s_n(z)}\qquad(n\ge2).
\end{aligned}
\]
The symbol $[z^n]$ takes the Taylor coefficient at zero.
For the curvature test $I(x)$ averages
$(u-u_0)[P(u)^2-D(u)]$ against a gamma density of shape $x+1$
and rate two. Here $P$ is the sum of the individual prime slopes,
$D$ is their squared diagonal, and $u_0$ is their first
correlation zero (Proposition~\ref{pr:primewindowprefix}).
Thus $x$ is a moment degree in this test.

Table~\ref{tab:summary} records the hypotheses of the spectral and
arithmetic estimates. The exact criteria for completeness are
\[
\begin{aligned}
 \Delta=0\quad&\Longleftrightarrow\quad\mathrm{RH}\\
 &\Longleftrightarrow\quad
 \bigl(\exists\,x_j\to\infty:\ I''(x_j)\le0\bigr)\\
 &\Longleftrightarrow\quad
 \bigl(|\mathcal A_n|\le2\ \text{infinitely often}\bigr).
\end{aligned}
\]
Each off-line zero contributes a strictly positive amount to the
mass defect (Theorem~\ref{th:boundaryclosure}). The curvature
criterion follows from the two asymptotic alternatives in
Theorem~\ref{th:realdegreecurvature}, and the coefficient criterion
from the disk transform in Theorem~\ref{th:movingprime}. An
independent arithmetic estimate for either observable suffices.
Theorem~\ref{th:stependrange} gives a criterion with an explicit
threshold. RH holds if and only if the prime sum $P$ of
Sec.~\ref{sec:stepprimes} satisfies \eqref{eq:primeend} at $x=p$
and $x=p/8$ for every prime $p>661$. Theorems~\ref{th:jacfamily}
and~\ref{th:cycwindow} identify the geometric and cohomological
estimates in the comparison with Stepanov's method.
Theorem~\ref{th:stepfirstterm} evaluates the first point-count
correction for a fixed curve. Proposition~\ref{pr:stepdeletion}
establishes the sensitivity of the endpoint inequality to prime deletion.

\head{Arithmetic input for the inverse construction}
Proposition~\ref{pr:boost} and Theorem~\ref{th:wall} establish the
spectral properties of the wedge and wall independently of RH. The
return, explicit-formula, and correlation identities
\eqref{eq:ring}, \eqref{eq:explicit}, and \eqref{eq:specC} are also
unconditional. Theorem~\ref{th:boundaryclosure} constructs the
boundary measure and its self-adjoint operator, gives the formula
for $\Delta$, and bounds $\Delta$ from finite-height data.
Proposition~\ref{pr:masspropagation}, Theorem~\ref{th:inversebessel},
and Proposition~\ref{pr:inversecoupling} realize and estimate the
reconstructed dynamics. Positive definiteness of $C$
(Theorem~\ref{th:CAR}), the upper bound $C(u)\le C(0)$ at every
time (Corollary~\ref{co:bounded}), vanishing of $\Delta$, and the
prime-sum criteria above are each equivalent to RH.

The boundary measure is obtained from the meromorphic continuation of
the prime response $L(z)$ (Sec.~\ref{sec:boundaryclosure}).
Consequently the inverse operator $K_\partial$ inherits its
frequencies from that response. Proposition~\ref{pr:inversecoupling}
realizes them on a reducing subspace of a Hamiltonian within trace
norm $\epsilon$ of $H_0$. A derivation of \eqref{eq:CAR} from
couplings specified independently of the zeros would itself imply RH.
The recurrence of Sec.~\ref{sec:masspropagation} determines $C$
exactly on each finite interval. Vanishing of the defect requires
$C(u)\le C(0)$ for every $u$. The error bound obtained from
finite-height data grows as $\cosh(u/2)$.

The deletion test identifies a restriction on a proof of
\eqref{eq:primeend} from local arithmetic data.
Proposition~\ref{pr:stepdeletion} constructs a set $\mathcal D$
of primes of density zero with
$\vartheta_{\mathcal D}(y)\sim c\sqrt y$, $c>1/2$.
After deletion the prime number theorem and the point-count bounds
hold at the retained primes, while the endpoint inequality fails at
arbitrarily large cutoffs. The global estimate therefore depends on
finer information about the set of rational primes.

\subsection{Square-root bound for the shifted boundary Hamiltonian}
\label{sec:delayedmajorana}
The shifted boundary Hamiltonian gives a global square-root estimate
for the deformed response $R_h$. Its delay coefficients are specified
directly from the primes. The detector clock is $u=2\eta$, where
$\eta$ is boost time. In the notation of \eqref{eq:boost}, the spinor
$e^{-\mu_0e^v}(1,-i)^T$ solves the zero-boost Majorana equation.
With $\mu_0e^v=\pi y^2$, the Mellin transform of twice its first
component is, for $\re s>0$,
\[
 2\int_0^\infty e^{-\pi y^2}y^{s-1}\dd y
 =\pi^{-s/2}\Gamma(s/2)=\Gamma_{\R}(s).
\]
The coordinate $r=\pi y^2$ and the prime couplings specify the model.
The zero-boost solution itself is fixed by decay, since
$\sqrt r\,K_{1/2}(r)=\sqrt{\pi/2}\,e^{-r}$ \cite[\S10.39]{dlmf}.
For $0<h\le1$, set $s=1/2+w$. At each integer delay $\log n$
put $b_h(n)=n^h\prod_{p\mid n}(1-p^{-2h})$, including $b_h(1)=1$.
Then the boundary transfer is
\begin{equation}
\begin{aligned}
 Q_h(w)&=\frac{(s-h)(s-h-1)\Gamma_{\R}(s-h)}
 {(s+h)(s+h-1)\Gamma_{\R}(s+h)}
       \sum_{n\ge1}b_h(n)n^{-s}\\
 &=\frac{\xi(s-h)}{\xi(s+h)},\qquad \re s>1+h .
\end{aligned}
\label{eq:boundarytransfer}
\end{equation}
The second equality follows by expanding each local Euler factor.
These are the shifted functions studied by Lagarias and Suzuki
\cite{shiftlagariasm,shiftsuzukim}.

In the same half-plane the logarithm of the Euler quotient is
\[
 \log\frac{\zeta(s-h)}{\zeta(s+h)}
 =\sum_{n\ge2}\frac{\Lambda(n)}{\log n}(n^h-n^{-h})n^{-s}.
\]
Exponentiating gives the delay coefficients $b_h(n)$.
The sign $\Lambda(n)\ge0$ makes these coefficients nonnegative.

\subsubsection{The causal coupling on a finite time interval}\label{sec:boundarycausal}
Extend all time signals by zero for negative times. A star denotes
causal convolution. For $a_\pm=h\pm\tfrac12$ and $t>0$ set
\begin{equation}
 \begin{aligned}
 g_h(t)&=\frac{2\pi^h}{\Gamma(h)}
        e^{(h-1/2)t}(1-e^{-2t})^{h-1},\\
 (D_{h,a}u)(t)&=u(t)-2h\int_0^t e^{-a(t-y)}u(y)\dd y,\\
 d_h&=D_{h,a_+}D_{h,a_-}g_h,\\
 k_h(t)&=\sum_{\log n<t}\frac{b_h(n)}{\sqrt n}\,d_h(t-\log n).
 \end{aligned}
 \label{eq:boundarykernel}
\end{equation}

Each filter $D_{h,a}$ subtracts an exponentially weighted history
of the pulse from its present value. The kernel $k_h$ superposes
copies of the filtered radial response $d_h$ at the arithmetic
delays $\log n$. Values at the countably many endpoints do not affect convolution.
Every bounded interval meets only finitely many delays. The singularities
of $d_h$ at the origin are locally integrable for all $h>0$.
The beta integral \cite[\S5.12]{dlmf} applied to $u=e^{-2t}$ gives
$\mathcal Lg_h(w)=\Gamma_{\R}(s-h)/\Gamma_{\R}(s+h)$.
Each $D_{h,a}$ multiplies a Laplace transform by
$1-2h/(w+a)$. Translation by $\log n$ contributes $n^{-w}$.
These computations give $\mathcal L k_h=Q_h$ on the half-plane
of absolute convergence in \eqref{eq:boundarytransfer}.
The shifted ratio and its arithmetic kernel occur in Suzuki's work
\cite{shiftsuzukim}. The estimates here use an ordinary causal convolution and
one fixed input.

On $L^2(0,L)$, $0<L<\infty$, let
$(V_{h,L}u)(t)=\int_0^t k_h(t-y)u(y)\dd y$. It maps the prepared pulse $u$ to the signal seen by the boundary
detector during an observation of duration $L$.
Each $V_{h,L}$ is compact and $V_{h,L}u\to u$ in $L^2(0,L)$
as $h\downarrow0$ for each $u$. Nevertheless
$\lVert I-V_{h,L}\rVert_{\mathrm{op}}\ge1$.

\begin{proof}
On $[0,L]$ the kernel is in $L^1$. Approximation there by bounded
kernels and Young's inequality give norm approximation by
Hilbert--Schmidt operators. For $0<h\le1/2$ write
\[
 g_h(t)=\frac{(2\pi)^h}{\Gamma(h)}t^{h-1}
       e^{(h-1/2)t}\left(\frac{1-e^{-2t}}{2t}\right)^{h-1}.
\]
The last two factors together equal $1+O_L(t)$, uniformly in $h$.
Consequently $\int_0^L g_h(t)\dd t=(2\pi L)^h/\Gamma(1+h)+O_L(h)\to1$;
for $0<d<L$, $\int_d^L g_h(t)\dd t=O_{d,L}(h)$.
Continuity of translations in $L^2(\R)$, after extension by zero,
proves the approximate-identity assertion for this first kernel.
The $L^2(0,L)$ operator norms of $D_{h,a_\pm}-I$ are $O_L(h)$.
For each of the finitely many $2\le n\le e^L$, $b_h(n)=O_L(h)$.
This proves strong convergence for $k_h$. A compact operator sends an
orthonormal sequence to zero in norm, giving the final inequality.
\end{proof}

\subsubsection{Spectrum and square-root bound of the shifted generator}\label{sec:boundaryshift}

\begin{theorem}[A square-root bound from the boundary Hamiltonian]
\label{th:shiftmajorana}
Write $X_h(s)=(\xi(s+h)+\xi(s-h))/2$,
$M_h=2X_h'(1)/X_h(1)$ and $a=1/2$.
Define the output by
\begin{equation}
\begin{aligned}
 \mathcal L C_h(w)&=\frac{wM_h-(X_h'/X_h)(a+w)}{w^2-a^2},\\
 R_h(x)&=d_0-S(x)+\sqrt{x}\,C_h(\log x).
\end{aligned}
\label{eq:boundaryoutput}
\end{equation}
It agrees with $R$ at $h=0$. For every $1/2\le h\le1$
and every $x\ge1$,
\[
 |R_h(x)-d_0+S(x)|\le M_h\sqrt{x}.
\]
In particular $R_{1/2}(x)\le d_0+
2\xi'(3/2)\sqrt{x}/[\xi(3/2)+\xi(1/2)]$.
The absolute bound is attained at $x=1$. As $h\downarrow0$,
$R_h$ converges to $R$ uniformly on each compact subinterval of
$[1,\infty)$.
\end{theorem}
\begin{proof}
In this range $Q_h$ is analytic of modulus less than one in the right
half-plane, with modulus one on the boundary. In the centered Hadamard
product, the factors at ordinate $\gamma$ pair as
$(w-b-i\gamma)/(w+b-i\gamma)$ with $b=h+\re(\rho)-1/2>0$.
Offsets $\re(\rho)-1/2$ at a fixed ordinate occur in opposite pairs.
Grouping also with $-\gamma$ gives local convergence because
$\sum_\rho(1+|\im\rho|^2)^{-1}<\infty$. Evenness removes the
linear exponential in the centered product, and the quotient has
value one at the origin. In the disk set
$w=2(1+v)/(1-v)$ and expand
$(1-Q_h)/(1+Q_h)=c_0+2\sum_{k\ge1}c_kv^k$.
The coefficients are real and computable at $w=2$ in the Euler
half-plane. The positive real part gives a positive Toeplitz form
$[c_{j-k}]$ on Laurent polynomials, with $c_{-k}=c_k$.
In its Hilbert completion modulo null vectors, the bilateral shift
$U$ advances the Laurent index by one and is unitary.
Its circle measure has no atom at $1$, since
$Q_h(w)=O_h(w^{-h})$ as real $w\to+\infty$. The generator is
\begin{equation}
 H_h=2i(I+U)(I-U)^{-1},\qquad
 \mathcal D(H_h)=\operatorname{ran}(I-U).
\label{eq:boundarygenerator}
\end{equation}
This inverse Cayley transform generates the continuous evolution
of the boundary modes in return time. It is symmetric, and $H_h\pm2i$
are surjective by unitarity. The domain is dense since the absence
of a mass at $1$ gives $\ker(I-U^*)=0$.
Conjugation of polynomial coefficients anticommutes with $H_h$.
Herglotz inversion puts a simple point mass at each pole of the
disk function and zero measure on the intervening arcs. Indeed this
meromorphic function has purely imaginary values there. The constant
Laurent polynomial is cyclic, so each atom has multiplicity one.
Thus the
energies are precisely the simple zeros of $X_h(a-iz)$.
Indeed $|Q_h|<1$ and zero-freeness of $\xi(s+h)$ exclude zeros
in the upper half-plane, and reflection excludes the lower half-plane.
The restriction $H_{h,+}$ selects one frequency from each signed
pair. The paired determinant records the two signed frequencies as
its zeros. Writing $A_h(z)=X_h(a-iz)$, the Hadamard product gives
\begin{equation}
 D_h(z)=\det(I-z^2H_{h,+}^{-2})=\frac{A_h(z)}{A_h(0)}.
 \label{eq:boundarydet}
\end{equation}
The even characteristic function has order at most one, and
$\sum_k t_k^{-2}<\infty$ makes the determinant well defined.
For the positive energies $t_k$,
\[
 C_h(u)=\sum_k\frac{2\cos(t_k u)}{a^2+t_k^2},\qquad
 M_h=\sum_k\frac{2}{a^2+t_k^2}.
\]
These follow by logarithmic differentiation and absolute convergence.
Thus $C_h$ is the return amplitude of a probe assigning weight
$1/(a^2+t_k^2)$ to each signed mode. The triangle inequality gives
the bound, with equality at $u=0$. For $h=0$, the arithmetic
formula for $L(w)$ in Sec.~\ref{sec:boundaryclosure} identifies
the inverse transform with $C$ and $R_0=R$.

For compact-time convergence, on $\re s=4$ write
$D_\pm=(\xi'/\xi)(s\pm h)$.
Stirling's formula and the Euler series bound the imaginary part
of $\xi'/\xi$ on $3\le\re s\le5$. The argument of $Q_h$,
continued from $h=0$, is therefore $O(h)$ uniformly on the
line. For small $h$, $\re Q_h\ge0$ and
\[
 \frac{X_h'}{X_h}=\frac{D_++Q_hD_-}{1+Q_h}
      =O(\log(2+|\im s|)).
\]
After the elementary $M_h\cosh(au)$ term in
\eqref{eq:boundaryoutput}, the inverse is dominated on compact time
intervals by a constant times
$\log(2+|\im s|)/(1+|\im s|^2)$. Dominated convergence proves
the compact-interval assertion.
\end{proof}
The input on this inversion line is computed from
$(\xi'/\xi)(s+h)+Q_h'(w)/(1+Q_h(w))$ by the Euler product and
\eqref{eq:boundarytransfer}. The real normalizing values of $\xi$
and its derivative are given by the convergent Jacobi theta integral
\cite[Ch.~II]{Titchmarsh}. The Euler product, the alternating eta
series, continuity at $0,1$, and the functional equation give
$\xi(s)>0$ for all real $s$.
The product proof also applies to any positive shift whose open strip
$|\re(\rho)-1/2|<h$ contains every zero. At $h=0$ direct substitution
gives $Q_0=1$ and the Toeplitz form vanishes. The spectral limit
for this endpoint is treated in Sec.~\ref{sec:boundarycritical}.

\subsubsection{Pointwise and mean bounds on intervals}\label{sec:boundaryintervals}

\begin{proposition}[Whole short intervals for the unshifted detector]
\label{pr:actualshort}
For each $0<\alpha<1$, arbitrarily large $Y$ satisfy
$R(x)\le\sqrt{x}$ for every $Y\le x\le Y+Y^\alpha$.
For each fixed $\lambda>1$ there are also arbitrarily large $Y$
for which
\[
 \frac{1}{(\lambda-1)\log Y}\int_Y^{Y^\lambda}
           \frac{R(x)}{\sqrt x}\,\frac{\dd x}{x}\le1 .
\]
This includes the windows $[Y,Y^{3/2}]$. If $\Theta=\sup_\rho\re(\rho)>1/2$,
the pointwise values on the selected short intervals and the means on
the selected power windows can both be taken strictly negative.
\end{proposition}
\begin{proof}
Write $\Theta=\sup_\rho\re(\rho)$. If $\Theta=1/2$, use
Lemma~\ref{lem:explicit}. Otherwise choose
$\max(1/2,\Theta+\alpha-1)<q<\Theta$. The Mellin transform
in that lemma is regular at every real $s>q$ and has a pole
to the right of $q$. By Landau applied to $R(x)+x^q$,
arbitrarily large $Y$ have $R(Y)<-Y^q$.
For $0\le u\le H=Y^\alpha$, split the explicit zero sum
at height $Y/H$, using
$\lvert(1+u/Y)^\rho-1\rvert\le 3\min(1,|\rho|H/Y)$.
Partial summation of the zero count gives
$\sum_{|\im\rho|\le T}|c_\rho\rho|\ll\log^2(2T)$ and
$\sum_{|\im\rho|>T}|c_\rho|\ll T^{-1}\log(2T)$, where
$c_\rho=1/[\rho(1-\rho)]$. Together with $S'(Y)=O(Y^{-3})$
these estimates give, for $0<H\le Y/2$,
\[
 \sup_{0\le u\le H}|R(Y+u)-R(Y)|
 \ll Y^{\Theta-1}H\log^2(2Y/H)+HY^{-3}.
\]
For $H=Y^\alpha$ this is $o(Y^q)$. Hence $R<0$ throughout.
For the second statement set $b=\Theta-1/2>0$ and
$F(T)=\int_0^T e^{-v/2}R(e^v)\dd v$.
It is $O_b(e^{bT})$ and its Laplace transform equals the Mellin
transform of $R$ at $1/2+w$, divided by $w$.
The same argument gives $F(T)<-e^{a'T}$ on a sequence,
for $b/\lambda<a'<b$. Since $F(T/\lambda)=O_b(e^{bT/\lambda})$,
we get $F(T)-F(T/\lambda)<0$. Take $Y=e^{T/\lambda}$.
If $\Theta=1/2$, use the pointwise estimate from the explicit
formula, for which $Q=2+\eul-\log(4\pi)<1$.
\end{proof}
The two sequences of starting points may be chosen independently.
Choosing successively $\alpha_j=1-1/j$ and $Y_j\ge e^{j^2}$
in the short-interval assertion gives lengths $Y_j^{1-o(1)}$,
with ratio to $Y_j$ tending to zero. Their widths in the return
clock are $\log(1+Y_j^{\alpha_j-1})$. A window $[Y,Y^{3/2}]$
has width $\tfrac12\log Y$ in that clock.

\subsubsection{Poisson reconstruction and fixed-shift convergence}\label{sec:boundarypoisson}

\begin{theorem}[Recovering the window statistic from one Majorana determinant]
\label{th:majoranapoisson}
Set $H=H_{1/2}$ from Theorem~\ref{th:shiftmajorana}, and let
$D(z)=\det(I-z^2 H_+^{-2})$. Form the strip Poisson integral
\[
 Z_T(s)=\int_{-T}^{T}\frac{D(i(s-1/2)-v)}{\cosh(\pi v)}\dd v .
\]
Then $Z_\infty(s)=\xi(s)/\xi(1)$ throughout the complex plane.
Put $M_T=2Z_T'(1)/Z_T(1)$. For large $T$ define the real output,
with $s=2+it$, by
\[
\begin{aligned}
 {\cal R}_T(x)={}&d_0-S(x)+(x+1)M_T/2\\
 &-\frac{1}{2\pi}\int_{-T}^T
      \frac{x^s}{s(s-1)}\frac{Z_T'(s)}{Z_T(s)}\dd t .
\end{aligned}
\]
For $G=R$ of \eqref{eq:center} the error is
$O(x^2\log(2+T)/T)$ for every $x\ge1$ and all sufficiently large
$T$, with an absolute constant. For every sufficiently large $Y$,
\[
 \sup_{Y\le x\le Y^{3/2}}
 \frac{|{\cal R}_{Y^3}(x)-G(x)|}{\sqrt x}
 \ll Y^{-3/4}\log(2Y).
\]
\end{theorem}
\begin{proof}
For $\tfrac12\le h\le1$ and $\Xi(z)=\xi(1/2-iz)$,
\eqref{eq:boundarydet} gives
$\xi(1/2+h)D_h(z)=(\Xi(z+ih)+\Xi(z-ih))/2$. Put
\[
 K_h(v)=\frac{1}{2h\cosh(\pi v/(2h))},\qquad
 \widehat K_h(t)=\frac1{\cosh(ht)} .
\]
The Fourier transform uses $e^{-ivt}$ and follows from the beta
integral with $y=e^{\pi v/h}$. The middle-line Poisson measure
is even and has total mass one \cite{widderstrip}. Integrate
$\Xi(z+w)/\sinh(\pi w/(2h))$ around the rectangle
$|\re w|\le L,\ |\im w|\le h$. The residue at zero is
$2h\Xi(z)/\pi$. Stirling's formula and polynomial strip bounds
for $\zeta$ make the vertical sides tend to zero. The horizontal
sides give
\begin{equation}
  (K_h*D_h)(z)=\frac{\Xi(z)}{\xi(1/2+h)}
   \qquad(z\in\C).
  \label{eq:boundarypoisson}
\end{equation}
Dividing by the value at zero recovers $\Xi(z)/\Xi(0)$. The same
strip bounds and Cauchy's formula justify differentiation on compact
sets. Taking $h=1/2$ gives the stated formula for $Z_\infty$.

$Z_T-Z_\infty$ and its first derivative are $O(e^{-\pi T})$
on $1\le\re s\le2$, by integration of the kernel tail and
boundedness of $D,D'$ on the corresponding horizontal strips.
For $s=2+it$ the Euler product and Stirling's formula give
$|\xi(s)|\gg e^{-\pi|t|/4}$ and
$|\xi'(s)/\xi(s)|\ll\log(2+|t|)$. Thus the denominator is nonzero
for $|t|\le T$ when $T$ is large. The error in its logarithmic
derivative there is $O(e^{-3\pi T/4}\log(2+T))$.
Conjugation gives a real integral. At $h=0$ we have
$(s-1/2)^2-1/4=s(s-1)$ and $\sqrt{x}\cosh(\log x/2)=(x+1)/2$.
The inverse transform therefore equals $G$, and its tails cost
$O(x^2\log(2+T)/T)$. Also $M_T-M_0=O(e^{-\pi T})$.
Taking $T=Y^3$ and $x^{3/2}\le Y^{9/4}$ proves the estimate.
\end{proof}

Give each copy of the boundary Hamiltonian a fixed real energy
offset $v$ and weight the copies by $K_h(v)\dd v$.
On $L^2(\R,K_h(v)\dd v;\mathcal H_h)$, the direct integral
$(\widetilde H f)(v)=(H_h+vI)f(v)$ is self-adjoint on its natural
domain $f(v)\in\mathcal D(H_h)$ almost everywhere,
$\int\|(H_h+vI)f(v)\|^2 K_h(v)\dd v<\infty$.
The Majorana conjugation is $\widetilde C f(v)=C f(-v)$.
In an eigenbasis of $H_h$, each summand is multiplication by
$v+t_k$ for a signed energy $t_k$, which proves these assertions
and gives absolutely continuous spectrum on $\R$.
Equation~\eqref{eq:boundarypoisson} averages the individual characteristic
functions at the displaced energies. To obtain the prime detector one
then takes the logarithmic derivative of the average. Its denominator
is the whole integral, so a zero created by averaging gives a pole in
that logarithmic derivative.

The approximation applies for every large $Y$ and uses the Hamiltonian
at a single proved self-adjoint shift.
The one-sided bound
\eqref{eq:acwindow} transfers between ${\cal R}_{Y^3}$ and $G$
by first taking exponent $\epsilon/2$.
A finite-dimensional example determines the effect of this convolution
on spectral zeros. The Majorana Hamiltonian $h\sigma_2\oplus2h\sigma_2$
has simple real eigenvalues and normalized determinant
$P(z)=(1-z^2/h^2)(1-z^2/(4h^2))$. The Fourier transform of $K_h$
gives its second and fourth moments, $h^2$ and $5h^4$, respectively. Hence
\begin{equation}
 (K_h*P)(z)=1+\frac{z^2}{4h^2}+\frac{z^4}{4h^4}.        \label{eq:boundaryfour}
\end{equation}
The convolution preserves the value at the origin. The resulting
polynomial is strictly positive on the real axis and has four nonreal
zeros. A real-zero theorem for the prime determinant consequently
requires an additional property of its arithmetic coefficients.

\subsubsection{Convergence as the shift tends to zero}\label{sec:boundarysmall}
For each fixed $1<\sigma\le2$, there is $h_\sigma>0$ such that
\begin{equation}
 |R_h(x)-G(x)|\ll_\sigma h^2x^\sigma,\qquad
 0\le h\le h_\sigma,\quad x\ge1.                         \label{eq:boundarysmall}
\end{equation}
Set $D_0^*(s)=\xi'(s)/\xi(s)$,
$D_\pm=D_0^*(s\pm h)$, and
$L=\int_{-h}^h D_0^*(s+r)\dd r$. In the strip
$\sigma\le\re s\le4$ the Euler series and Stirling's formula bound
$\Im D_0^*$ and the first two derivatives of $D_0^*$ uniformly,
also on a slightly enlarged strip contained in $\re s>1$.
Moreover $D_0^*(s)=O_\sigma(\log(2+|\Im s|))$.
Choose $h_\sigma$ small enough for $|\Im L|\le\pi/2$. Then
$X_h\ne0$ there and
\[
 \frac{X_h'}{X_h}=\frac{D_++D_-}{2}
       +\frac{D_+-D_-}{2}\tanh(L/2).
\]
Since $|\cosh z|^2\ge1/2$ for $|\Im z|\le\pi/4$, integration
of $\tanh'(z)$ gives $|\tanh(L/2)|\le |L|$.
Symmetric Taylor expansion and the derivative bounds therefore give
$(X_h'/X_h)(s)-D_0^*(s)=O_\sigma(h^2\log(2+|\Im s|))$.
At $s=1$, analyticity and $\xi(1)>0$ give $M_h-M_0=O(h^2)$.
Shift the inverse integral \eqref{eq:boundaryoutput} from $\re s=4$ to
$\re s=\sigma$. The horizontal integrals vanish by the bounds
just proved and $1/[s(s-1)]=O_\sigma((1+|\Im s|^2)^{-1})$.
This proves \eqref{eq:boundarysmall}. For example, $h=Y^{-1}$,
$\sigma=3/2$ gives
$\sup_{Y\le x\le Y^{3/2}}|R_{1/Y}(x)-G(x)|/\sqrt{x}\ll Y^{-1/2}$.
More generally $h=Y^{-b}$ works for each $b>3/8$ by taking
$\sigma$ sufficiently close to $1$. For these small shifts the
calculation uses the arithmetic transform \eqref{eq:boundaryoutput}.
Theorem~\ref{th:majoranapoisson} uses the self-adjoint operator
at the fixed shift $1/2$.

\subsubsection{Cancellation in a single pulse}\label{sec:boundarypulse}
Fix $f(t)=t^2e^{-t}/2$ for $t\ge0$, zero otherwise. Put
$p_h=d_h*f$. The signal and its Laplace transform are
\begin{equation}
 \begin{aligned}
 B_h(t)&=\sum_{n\le e^t}\frac{b_h(n)}{\sqrt n}\,p_h(t-\log n),\\
 \mathcal L B_h(w)&=\frac{Q_h(w)}{(1+w)^3}.
 \end{aligned}
 \label{eq:boundarypulse}
\end{equation}

Here $p_h(0)=0$. All these pulses are continuous. The profile $f$ is the prepared input.
The detector records the translated responses with their signs in $B_h$. The argument of
Sec.~\ref{sec:boundarycausal}, using uniform continuity of the extension
of $f$, gives $B_h\to f$ uniformly on every finite time interval.

\head{Signed and absolute pulse bounds}
For $0<h<1/2$ put $b=1/2-h$. There is $c_h>0$ such that
\begin{equation}
 \begin{gathered}
 |B_h(t)|\ll_h \exp(bt-c_h\sqrt t)\qquad(t\ge1),\\
 \int_0^\infty e^{-bt}|B_h(t)|\dd t<\infty .
 \end{gathered}
 \label{eq:boundaryend}
\end{equation}
The sum in \eqref{eq:boundarypulse} with $|p_h|$ in place of $p_h$ is
$\gg_h \exp((1/2+h)t)$ for all sufficiently large $t$.
For $1/2\le h\le1$, $|B_h(t)|\le1/\pi$ for all $t\ge0$.
\begin{proof}
On $\re w>b$, the denominator zeta factor is zero-free by the
Euler product. The gamma quotient is $O_h((1+|\Im w|)^{-h})$
on fixed vertical strips at large height. For $\re u$ bounded
away from zero in such a strip, $\zeta(u)=O(1+|\Im u|)$
at large height by summation of the fractional-part integral.
Since $1-2h>0$, this applies to the numerator in a neighborhood
of $\re w=b$. Poles of the individual factors at height zero
are removable in the completed ratio. Bromwich inversion of
\eqref{eq:boundarypulse} is therefore valid on each line $\re w>b$ by a
contour shift from the initial half-plane.

For $t$ large choose $T=e^{\sqrt t}$. Start on $\re w=b+1/t$
and move the part with $|\Im w|\le T$ to $b-c/\log T$.
For sufficiently small absolute $c>0$ the classical zero-free region
and its reciprocal bound give
$|1/\zeta(1/2+h+w)|\ll(\log T)^A$ on this rectangle for
$2\le|\Im w|\le T$ \cite[Ch.~III, \S\S3.8--3.11]{Titchmarsh}.
On the remaining compact segment the completed ratio stays bounded.
For the two tails on the starting line use
$|1/\zeta(1+1/t+iv)|\le\zeta(1+1/t)\ll 1+t$.
After division by $(1+w)^3$, the vertical integrand before
the exponential is $O_h((1+|v|)^{-2-h}(\log T)^A)$ on the
new segment. After division by $e^{bt}$, the shifted segment contributes
$O_h(e^{-ct/\log T}(\log T)^A)$, the two horizontal sides contribute
$O_h(T^{-2-h}(\log T)^A)$, and the tails contribute
$O_h((1+t)T^{-1-h})$.
These bounds give \eqref{eq:boundaryend} after decreasing $c_h$.

Near zero the beta integral gives
$p_h(t)\sim (2\pi)^h t^{h+2}/\Gamma(h+3)$;
the exponential convolutions in $D_{h,a_\pm}$ cost one higher
power of $t$. Hence $p_h$ is bounded below by a positive
constant on some fixed $[u,v]\subset(0,\infty)$.
The average of $b_h(n)n^{-h}=\sum_{d\mid n}\mu(d)d^{-2h}$
tends to $1/\zeta(1+2h)$, by dominated convergence applied
to $\sum_{d\le X}\mu(d)d^{-2h}\lfloor X/d\rfloor/X$.
Partial summation yields
\[
 \sum_{n\le X}\frac{b_h(n)}{\sqrt n}
 \sim \frac{X^{1/2+h}}{(1/2+h)\zeta(1+2h)} .
\]
Positivity of the coefficients allows us to retain only
$e^{t-v}\le n\le e^{t-u}$ in the absolute sum. This proves
the stated lower bound. Finally, for $h\ge1/2$, Theorem~\ref{th:shiftmajorana}
gives $|Q_h|\le1$. The inverse transform on $\re w=d>0$
has magnitude at most
$e^{dt}(2\pi)^{-1}\int_\R |1+d+iv|^{-3}\dd v
=e^{dt}/[\pi(1+d)^2]$. Let $d\downarrow0$.
\end{proof}

\subsubsection{Discounted positive mass and the critical operator}\label{sec:boundarycritical}
Use the real moving average $Y_h(t)=\int_0^1B_h(t-r)\dd r$.
For $d>0$ and $T\ge0$ set
\[
 I(h,d,T)=\int_0^T e^{-dt}\max\{Y_h(t),0\}\dd t .
\]
Here $Y_h$ integrates one unit of return time and the exponential
weight discounts late readings. The one-sided pulse condition
$(\mathrm P)$ is

\begin{equation}
 \boxed{\begin{gathered}
   h_j=1/(j+3),\\
   \forall J\ \exists j\ge J:\ 
       \sup_{T\ge0} I(h_j,2h_j,T)<\infty .
 \end{gathered}}
 \label{eq:boundarydiscount}
\end{equation}
All summands defining this integral are given by \eqref{eq:boundarytransfer},
\eqref{eq:boundarykernel}, and \eqref{eq:boundarypulse}. Sec.~\ref{sec:boundarypulse} proves
$\sup_T I(h,2h,T)<\infty$ for every $1/6\le h<1/2$,
including the endpoint $h=1/6$, and Theorem~\ref{th:shiftmajorana} gives it
for $1/2\le h\le1$. The proof at $1/6$ uses the integrable factor
$e^{-c_h\sqrt t}$ in \eqref{eq:boundaryend}. These intervals give no
sequence tending to zero.

\head{Consequences of the discounted bound}
If \eqref{eq:boundarydiscount} holds, every zero of $\xi$ has real part $1/2$.
The arithmetic output $R_0$ then satisfies
\[
 R_0(x)\le \log(2\pi)+(2+\gamma_{\rm E}-\log(4\pi))\sqrt{x}
 \qquad (x\ge1).
\]
The ordered positive spectra of the Hamiltonians in \eqref{eq:boundarygenerator}
then yield a self-adjoint critical Majorana generator whose normalized
determinant is $\Xi(z)/\Xi(0)$, including multiplicities.
\begin{proof}
Fix a selected $h=h_j$, $d=2h$. The Laplace transform of $Y_h$
in its initial half-plane is
\[
 G_h(w)=\frac{1-e^{-w}}{w}\,\frac{Q_h(w)}{(1+w)^3}.
\]
Both multipliers of $Q_h$ are nonzero for $\re w>0$.
Finiteness of $I(h,d,\infty)$ makes $\mathcal L(Y_h)_+$
holomorphic for $\re w>d$. Initially
$\mathcal L(Y_h)_-=\mathcal L(Y_h)_+-G_h$, where
$u_-=\max(-u,0)$. The right side extends analytically across
each real $w>d$, since $\xi(1/2+w+h)>0$ there.
Landau's nonnegative-transform lemma implies that the abscissa of
convergence of $\mathcal L(Y_h)_-$ is at most $d$.
Indeed, if the real abscissa were a finite regular point, the Taylor
series about a point just to its right could cross it. Monotone
convergence applied there to the near-side expansion of $e^{ut}$,
$u>0$, would give a further finite exponential moment, a contradiction.
Thus $G_h$, and consequently $Q_h$, are holomorphic on $\re w>d$.

For a zero $\rho$ with $\re(\rho)>1/2$ choose a selected $j$
so large that $3h_j<\re(\rho)-1/2$ and
$\xi(\rho-2h_j)\ne0$. The latter follows from isolation of the
zero, for every sufficiently small positive $h_j$, including a
multiple zero. Equation~\eqref{eq:boundarytransfer} now has a pole at
$w=\rho-1/2-h_j$, inside the half-plane just proved holomorphic.
Reflection about the critical line completes the zero argument.
The calculation of Theorem~\ref{th:shiftmajorana} at $h=0$ now applies,
and $M_0=2\xi'(1)/\xi(1)=2+\gamma_{\rm E}-\log(4\pi)$
by the Laurent expansion at $1$. For every $\epsilon>0$ this
upper bound implies $R_0(x)\le x^{1/2+\epsilon}$ throughout every
sufficiently late $[Y,Y^{3/2}]$.

For the operator conclusion, the product proof in Theorem~\ref{th:shiftmajorana}
now applies at every positive $h\le1$. On compact subsets of
$\C$, $A_h(z)=\Xi(z)+O_K(h^2)$ by the symmetric Taylor expansion.
Rouch\'e's theorem on small disjoint circles about the zeros and a
larger circle avoiding zeros gives convergence of the ordered positive
eigenvalues, $t_k(h)\to t_k$, including repetitions in the limit.
Each limiting block evolves a real pair at the measured frequency
$t_k$. Realize the limits on a fixed Hilbert space by the blocks $t_k\sigma_2$ on $\ell^2(\mathbb N;\C^2)$, with domain
$\sum_k t_k^2\|v_k\|^2<\infty$. This operator is self-adjoint
and conjugation anticommutes with it. There is even trace-norm
convergence of the positive inverse squares after this ordering:
$\sum_k t_k(h)^{-2}=-A_h''(0)/(2A_h(0))$ tends to the corresponding
quantity for $\Xi$, which equals $\sum_k t_k^{-2}$ by its even
Hadamard product. Coordinatewise convergence of nonnegative sequences
and convergence of their sums give convergence in $\ell^1$.
The limit of \eqref{eq:boundarydet} is the asserted Fredholm determinant.
\end{proof}

\head{Variants with the same implication}
(i) The same argument applies to any $h_j>0,d_j>0$ with
$h_j+d_j\to0$ along the selected subsequence. For each $j$ the
box average can have any length $\ell_j>0$. Its multiplier
$(1-e^{-\ell_j w})/(\ell_jw)$ has no zeros in the right half-plane.
For the unsmoothed signal $B_h$ the multiplier is one.
(ii) A pointwise upper bound $Y_{h_j}(t)\le C_j e^{h_jt}$
suffices, since its positive part is integrable with discount
$2h_j$. So does $\int_0^\infty e^{-2h_jt}|Y_{h_j}(t)|^2\dd t
<\infty$, by Cauchy--Schwarz with $e^{-h_jt}$.
(iii) For each selected $j$, it is enough to bound the integrals
at arbitrarily late endpoints by the \emph{same} $C_j$.
Their monotonicity gives the bound for all $T$. A single finite
endpoint for each $j$ would provide no control of its infinite tail.
The unconditional compact-time convergence in Sec.~\ref{sec:boundarycausal}
and the entire-function convergence of $A_h$ carry no such tail bound.

The condition on the measured mass $I$ can be written
as a predicate over the reals. Substituting the
discounted integral above fixes its input.

\subsubsection{Inertia in the first prime window}\label{sec:boundaryinertia}
An antiunitary Majorana conjugation and an indefinite Hermitian metric
are different data. In a Pontryagin space with an independently specified
metric, split the space into a positive and a $\kappa$-dimensional
negative part. For a self-adjoint operator relative to this metric,
the span of root spaces in the open upper half-plane is neutral. Indeed,
the adjoint identity makes root spaces at $\lambda,\mu$
orthogonal if $\lambda\ne\bar\mu$, including generalized vectors by
induction along their chains. Projection of a neutral space to the
negative part is injective. Its dimension is therefore at most
$\kappa$ \cite[Props.~A.1--A.2]{pontryagin}. Applying that count to $\xi$ requires
an identification of root spaces with zeros and their multiplicities.
In the ordinary Weil evaluation form, a reflected off-line pair $\rho,1-\bar\rho$ of
multiplicity $m$ has block $m\left(\begin{smallmatrix}0&1\\
1&0\end{smallmatrix}\right)$, which gives one negative direction on the
full evaluation pair. Compression can lose that direction. The
Alp\"oge--Furman rank--trace argument combines inertia with analytic
estimates and gives asymptotic lower proportion $0.6725\ldots$ for
\emph{simple} zeros on the line \cite{weilaftwothirds}.

Theorem~\ref{th:krein} constructs the full return form on distinct
zeros, including the weights $m_\rho/|\rho(1-\rho)|$ in its Hilbert
norm. The Epstein example in Sec.~\ref{sec:krein} has a negative
logarithmic coefficient at $36$ despite its nonnegative ordinary
coefficients \cite{meanepstein}.

In the first support range containing the contribution of $p=2$,
the following estimate bounds the negative index of the Weil form.
Let $(\log2)/2<L\le(\log3)/2$, $u=\log2$, and $a_2=(\log2)/\sqrt2$.
For $f\in C_c^\infty(-L,L)$ satisfying
$\int f(x)e^{\pm x/2}\dd x=0$, the two pole terms vanish and
Weil's formula reads
\[
 \begin{aligned}
 W_L(f)={}&\int_\R\Omega(y)|\widehat f(y)|^2\,\frac{\dd y}{2\pi}\\
          &-2a_2\re\int_\R f(x+u)\overline{f(x)}\dd x,\\
 \Omega(y)={}&\re\psi(1/4+iy/2)-\log\pi ,
 \end{aligned}
\]
where $\widehat f(y)=\int f(x)e^{iyx}\dd x$ and
$\psi=\Gamma'/\Gamma$. At the right endpoint the lag $\log3=2L$
has zero overlap. The two intervals supporting the factors in the
shift pairing, before and after translation, are disjoint because
$2L<2\log2$. Thus its magnitude is at most $\|f\|^2/2$.
The translation by $u=\log2$ advances one copy of the test pulse
by the delay at prime $2$. Its overlap with the undisplaced pulse
measures interference at that delay. The Fourier weight $\Omega$
weights the frequency content of the same pulse on the archimedean
side of the form. The compressed prime term has norm exactly $a_2$ on the full window. Two equally weighted translated copies of a bump supported in the
overlap attain the bound before imposing the pole conditions.

Put $m=\Omega(0)$ and choose $B>0$ with $M=\Omega(B)>a_2$.
The digamma series makes $\Omega$ even and increasing in $|y|$,
and Stirling's formula guarantees such $B$. If $W_L(f)<0$ and
$\|f\|=1$, a split at $|y|=B$ gives
$\int_{-B}^B|\widehat f(y)|^2\dd y/(2\pi)>(M-a_2)/(M-m)$.
For an orthonormal family the sum of these band energies is at most
$2LB/\pi$, by Bessel's inequality applied to $x\mapsto e^{-iyx}$
on $[-L,L]$, followed by integration. Hence
\begin{equation}
 \operatorname{ind}_{-}(W_L\ \text{on the pole constraints})
 \le \frac{2LB}{\pi}\frac{M-m}{M-a_2}.                    \label{eq:boundaryindex}
\end{equation}
An explicit choice is $B=\tfrac12\exp(1+a_2-m)$. Indeed the
digamma series gives
$\Omega(B)-m=\sum_{n\ge0}v^2/[x_n(x_n^2+v^2)]$,
$x_n=n+1/4,\ v=B/2$. The summand is decreasing in $x_n$.
Integration from $1/4$ to infinity gives the lower bound
$\tfrac12\log(1+4B^2)$. Thus $M-a_2\ge1$, and \eqref{eq:boundaryindex}
is at most $L e^{1+a_2-m}(1+a_2-m)/\pi$ with this choice.
On the full test space the upper bound for the index increases by at
most two and remains finite. The archimedean positivity range is
$L\le(\log2)/2$ \cite{weilccm}.
Zhu's preprint reports a computer-assisted proof for the full
window at $L=0.8$, which by inclusion covers the whole one-prime
window \cite{weilzhum}. The proof of \eqref{eq:boundaryindex} uses only the displayed
archimedean symbol and the prime shift.

\subsubsection{Twisted boundary channels and the prime-deletion test}\label{sec:boundarycharacters}
For a primitive nonprincipal character $\chi$ of conductor $q$
and parity $a_\chi\in\{0,1\}$, use the completed function
$(q/\pi)^{(s+a_\chi)/2}\Gamma((s+a_\chi)/2)L(s,\chi)$.
The archimedean mode is $y^{a_\chi}e^{-\pi y^2/q}$, up to a
constant. The ratio at $s-h,s+h$ replaces $b_h(n)$ by
$\chi(n)b_h(n)$, uses the corresponding shifted gamma quotient
and omits the two rational pole corrections in \eqref{eq:boundarykernel}.
All statements about the transfer initially follow by absolute
convergence. Multiplication by $\chi(p)$ is a rotation on a
pair of real channels, with a zero contribution at ramified primes.
Pairing $\chi,\bar\chi$ gives a real structure. If the discounted
\emph{absolute value} of the corresponding averaged pulse is
integrable at $d_j$ for shifts with $d_j+h_j\to0$, the pole
isolation argument in Sec.~\ref{sec:boundarycritical} proves the critical-line
conclusion for that $L$-function. Absolute convergence gives the
holomorphic continuation immediately. For complex $\chi$, conjugation followed by the functional equation
reflects a zero $\rho$ of the same $L$-function to $1-\bar\rho$.
The positive-part argument for the untwisted function uses real
values and regularity on the real axis. A complex pulse has no
ordered positive part, and a real character may have an exceptional
real zero.

The strip transformation itself also works for a complex character.
Write $\Xi_\chi(z)$ for its completed function at $s=1/2-iz$
and $A_{\chi,h}(z)=(\Xi_\chi(z+ih)+\Xi_\chi(z-ih))/2$.
For each $h>0$,
\[
       K_h*A_{\chi,h}=\Xi_\chi .
\]
Indeed the completed function is entire for a primitive nonprincipal
character and is bounded in each horizontal $z$-strip by Stirling's
formula and the polynomial strip bounds for the Dirichlet $L$-function.
The residue computation of Sec.~\ref{sec:boundarypoisson} applies to these
complex boundary values directly. The identity therefore reconstructs
the full twisted function. A conclusion about its zeros would also
require control of the twisted arithmetic output.

Three changes to the arithmetic data locate the role of the primes.
The Davenport--Heilbronn function has no Euler product of the type
used in \eqref{eq:boundarytransfer} \cite{dhcomparison}. The positive comparison measure
$\Psi$ of Appendix~\ref{sec:stepmult} has an absolutely
continuous tail, while the logarithm of the Euler quotient is supported
on the discrete delays $\log p^k$. Finally, delete any prime $p$
and discard all its multiples in the defining integer series. With
constant coefficient one, its even theta sum is bounded above by
$1+2\sum_{p\nmid n,\,n\ge1}e^{-\pi n^2 y}$.
Multiplication by $\sqrt y$ gives upper limit at most $1-1/p$
as $y\downarrow0$, by the Gaussian Riemann sum and subtraction of
the multiples of $p$. The standard theta inversion with the same
constant coefficient would give limit one. This applies also to
the density-zero deletion in Proposition~\ref{pr:stepdeletion}.
Thus positivity of weights alone survives deletion, but the
archimedean functional equation used in the determinant proof does not.

\subsection{An integral exterior power of a perturbed Majorana Hamiltonian}
\label{sec:oscillatorheight}
\subsubsection{A restriction on bounded perturbations}\label{sec:oscperturbation}
Consider a bounded self-adjoint perturbation $K$ of $H=H_h$,
$1/2\le h\le1$, with $KH^{-1}$ trace class. Suppose that
$B=H+K$ is invertible and has Majorana symmetry. Its normalized
paired determinant $D_B$ cannot equal
$\Xi/\Xi(0)$. This includes every $K$ in a Schatten class
$\mathfrak S_p$, $1\le p<\infty$. Indeed $H^{-1}\in\mathfrak S_q$
for $q>1$, because $A_h$ has order at most one and
$A_h(0)\ne0$. Schatten H\"older gives the trace-class condition.
Also $B^{-1}=H^{-1}(I+KH^{-1})^{-1}$ is Hilbert--Schmidt.

Here $K$ describes a bounded change in the energy or the coupling
of boundary modes. A resolvent maps a source to the driven spinor at
complex frequency $z$. The difference of resolvent traces measures
the change of this summed response. To prove the restriction, put $d(z)=\det(I+K(H-z)^{-1})$.
On the imaginary axis $d$ is nonzero and $d(0)\ne0$.
With $R_J(z)=(J-z)^{-1}$, the determinant derivative and the
resolvent identity give $d'/d=\operatorname{Tr}(R_H-R_B)$.
Since $\det(I+t^2J^{-2})=D_J(it)^2$ for $J=H,B$, differentiating
in $t$ and fixing the constant at $t=0$ yields
\[
 \frac{D_B(it)^2}{D_h(it)^2}=\frac{d(it)d(-it)}{d(0)^2}.
\]
Here $H(H\mp it)^{-1}$ and their adjoints converge strongly to
zero and have norm at most one. Right multiplication by these maps
takes the trace-class operator $KH^{-1}$ to zero in trace norm,
by finite-rank approximation. Hence $d(\pm it)\to1$.
The ratio $D_B(it)/D_h(it)$ consequently has a positive finite limit.
Stirling's formula and the Euler product on the real
axis give
\[
 \frac{\Xi(it)/\Xi(0)}{D_h(it)}
 \sim \frac{2\xi(1/2+h)}{\xi(1/2)}
         (2\pi/t)^{h/2}\longrightarrow 0 .
\]
This proves the obstruction for the stated class of bounded
perturbations. The linear mass profile below is an unbounded
perturbation, and its integral exterior powers have explicit
determinantal divisors.

\subsubsection{The mass profile and its integral exterior powers}\label{sec:oscmass}

Change the exponential mass term in the boost expression to
$v\sigma_1$ on the full line, using oscillator units. Then
\[
 H_{\rm osc}=-i\sigma_3\partial_v+v\sigma_1,\qquad
 H_{\rm osc}^2=(-\partial_v^2+v^2)I+\sigma_2 .
\]
Here $H_{\rm osc}$ propagates the spinor in oscillator time with
a mass that varies linearly with position. This is an unbounded change
of model. In the basis $e_\pm=(1,\pm i)^T/\sqrt2$, let
$a=(\partial_v+v)/\sqrt2$.
The operators $a,a^*$ lower and raise the one-particle Hermite index.
Explicitly $H_{\rm osc}(fe_-)=-i\sqrt2(af)e_+$ and
$H_{\rm osc}(fe_+)=i\sqrt2(a^*f)e_-$ \cite{susyosc}.
For normalized Hermite functions $\phi_n$, the states
$\phi_ne_-,\phi_{n-1}e_+$, $n\ge1$, form Hermitian two-dimensional
blocks with energies $\pm\sqrt{2n}$; there is one zero mode.
Their direct sum is self-adjoint on $\sum n\|u_n\|^2<\infty$.
$\sigma_3$ followed by complex conjugation is its Majorana involution.
On the integral span of $\phi_ne_-$, $1\le n\le N$, the
number operator $H_{\rm osc}^2/2$ is $A_N=\operatorname{diag}(1,\ldots,N)$.
On a state $\phi_ne_-$, measuring $A_N$ returns the excitation
index $n$. The normalized levels fix the integral lattice for this
measurement. For $0\le r<N$, on a Slater state with $k=N-r$ distinct occupied levels the
exterior power acts by
\[
 \big(\textstyle\bigwedge^k A_N\big)|n_1,\ldots,n_k\rangle
   =\big(\textstyle\prod_i n_i\big)|n_1,\ldots,n_k\rangle .
\]
It reads the product of the occupied level numbers. The gcd over
these readings uses the chosen integral lattice and counts the
common factors of each prime with multiplicity. The integral lattice
is invariant under sign changes of the normalized basis vectors.
The gcd is the determinantal divisor in Smith normal form
\cite[Thm.~2.4]{smith}.

\begin{theorem}[Prime valuations of oscillator minors]\label{th:oscillatorial}
For $0\le r<N$ let $d_{N,r}$ be the positive gcd of the
$(N-r)$-minors of $A_N$. Then
\[
 v_p(d_{N,r})=\sum_{k\ge1}(\lfloor N/p^k\rfloor-r)_+
 \qquad\text{for every prime }p .
\]
For $x\ge2$, integer $m\ge2$, and $N=\lfloor mx\rfloor$, put
\[
 U_m(x)=m^{-1}\log d_{N,m-1}-x\log x+(1+\eul)x .
\]
The original detector obeys the signed bound
\[
 0\le U_m(x)-R(x)\le\psi(x)/m .
\]
Equivalently, for the uncentered sum $T$ in \eqref{eq:T},
\begin{equation}
 \begin{gathered}
  \frac{\log d_{N,m}}{m}\ \le T(x)\le\ \frac{\log d_{N,m-1}}m,\\
  \log(d_{N,m-1}/d_{N,m})=\psi(x).
 \end{gathered}
 \label{eq:oscbracket}
\end{equation}
The upper and lower determinantal estimates for $R(x)$ differ
by exactly $\psi(x)/m$. At $m=\lceil\sqrt x\rceil$ this is
$(1+o(1))\sqrt x$.
\end{theorem}
\begin{proof}
The nonzero entries in the exterior power are products of $N-r$
distinct integers from $1,\ldots,N$. To minimize the valuation at
$p$, omit the $r$ levels with highest $p$-adic orders. At
threshold $k$ this leaves $(\lfloor N/p^k\rfloor-r)_+$
levels. Summing the counts proves the first assertion.
For $n\le x$, $\lfloor N/n\rfloor=\lfloor mx/n\rfloor$.
For $r=m-1$, the excess over $m(x/n-1)$ equals $1-\{mx/n\}$.
For $r=m$ the deficit is $\{mx/n\}$. For $n>x$ all positive
parts are zero. Multiply by $\Lambda(n)\ge0$, sum, and apply the
prime number theorem for the asymptotic width. The minimizing sets
of levels may depend on $p$, as the gcd permits. The upper excess is
$m^{-1}\sum_{n\le x}\Lambda(n)(1-\{mx/n\})$, including at integer
endpoints. For $x\ge p$ the upper determinant contains a factor
of $p$, so prime deletion changes the integer itself.
\end{proof}

\begin{proposition}[Prime weights from exact Majorana revivals]\label{pr:oscillatorrevival}
Fix a clock frequency $\Omega>0$ in oscillator units, with $\hbar=1$.
On the spinor space of the linear-mass oscillator define by real spectral calculus
\[
 \begin{gathered}
 H_{\rm rev}=f_\Omega(H_{\rm osc}),\qquad f_\Omega(0)=0,\\
 f_\Omega(E)=\frac{2\Omega\,\operatorname{sgn}E}{E^2}\qquad(E\ne0).
 \end{gathered}
\]
This bounded self-adjoint Majorana Hamiltonian rotates the $k$-th
real oscillator pair at angular frequency $\Omega/k$ while the zero
mode supplies a stationary phase reference. On the zero mode
and the first $n$ signed oscillator pairs, let $T_n$ be the least
positive time for its propagator to equal the identity. Then
\[
 \begin{gathered}
 T_n=\frac{2\pi}{\Omega}\operatorname{lcm}(1,\ldots,n),\\
 \log(T_n/T_{n-1})=\Lambda(n),\qquad n\ge1,
 \end{gathered}
\]
where $T_0=T_1$ in the ratio at $n=1$.
In particular the signed mean \eqref{eq:heightsignedmean} has the
exact time-domain expression
\begin{equation}
 A(x)=\frac1m\sum_{j=1}^{m-1}
       \left(\frac Nj-\log\frac{T_{\lfloor N/j\rfloor}}{T_1}\right).
 \label{eq:heightrevival}
\end{equation}
\end{proposition}
\begin{proof}
The eigenvalues $\pm\sqrt{2k}$ of $H_{\rm osc}$ become
$\pm\Omega/k$, and the zero mode remains stationary. Hence the
operator norm is $\Omega$; oddness and reality of $f_\Omega$
preserve anticommutation with the antiunitary Majorana conjugation.
The rotation angle at time $t$ is $\Omega t/k$, up to orientation.
For all $k\le n$ to return to the identity requires
$\Omega t/(2\pi)\in k\mathbb Z$ simultaneously. Their intersection
is $\operatorname{lcm}(1,\ldots,n)\mathbb Z$.
When level $n$ is admitted, the least common multiple gains exactly
a factor $p$ if $n=p^a$. Otherwise every maximal prime-power
divisor of $n>1$ already occurs at a smaller level. This gives the
second identity. Summing it gives
$\log(T_{\lfloor y\rfloor}/T_1)=\psi(y)$ for $y\ge1$, and
substitution proves \eqref{eq:heightrevival}.
\end{proof}

There is also a frequency-resolved measurement of the same quantity.
On the single copy of levels $1,\ldots,N$ used for $A_N$, let
$P_y={\bf1}_{[1,y]}(A_N)$ select the levels below the cutoff and
let $W_N|n\rangle=\log(T_n/T_{n-1})|n\rangle$ assign their
measured revival increments. For $1\le y\le N$,
\[
 \operatorname{Tr}P_y=\lfloor y\rfloor,\qquad
 \operatorname{Tr}(P_yW_N)=\psi(y).
\]
Thus $y-\psi(y)$ compares two readouts of the same spectral gate:
the unweighted count with its fractional endpoint, and the
revival-weighted count. A character readout with diagonal entries
$\chi(n)$ gives $\sum_{n\le y}\chi(n)\Lambda(n)$ by the same
trace. The real and imaginary parts are two Hermitian readouts, and
ramified levels have zero weight. The construction fixes the
frequencies $\pm\Omega/n$ before any primes or zeta zeros are given.
It changes the dynamics on the original oscillator eigenvectors.
The difference $H_{\rm rev}-H_{\rm osc}$ is unbounded, since its
eigenvalues on the positive branch are $\Omega/n-\sqrt{2n}$.
The size of the signed mean requires the distribution of the common
return times as the cutoff grows.

\subsubsection{Size of the determinantal divisor}\label{sec:oscsize}
Let $L(b)=\operatorname{lcm}(1,\ldots,b)$ for positive integers
$b$. With the conventions of Theorem~\ref{th:oscillatorial},
\begin{align}
 d_{N,r}&=\frac{N!}{\prod_{j=1}^r L(\lfloor N/j\rfloor)},\label{eq:osclcm}\\
 U_m(x)&=\frac1m\sum_{j=1}^{m-1}
              \left(\frac Nj-\psi(N/j)\right)\nonumber\\
       &\qquad{}+ O\!\left(\frac{x+\log(2mx)}{m}\right).\label{eq:oscdeficit}
\end{align}
The constants in \eqref{eq:oscdeficit} are absolute. For
$m=\lceil\sqrt{x}\rceil$ there is $c>0$ such that
\[
 |U_m(x)|\ll x e^{-c\sqrt{\log x}}\qquad (x\ge2).
\]
\begin{proof}
Legendre's formula gives $v_p(N!)=\sum_{k\ge1}\lfloor N/p^k\rfloor$.
The denominator of \eqref{eq:osclcm} subtracts
$\sum_{k\ge1}\min(r,\lfloor N/p^k\rfloor)$, proving the identity
at each prime. Also $\log L(\lfloor y\rfloor)=\psi(y)$.
Writing $x'=N/m$ and $H_j=\sum_{i=1}^j1/i$, Stirling's
formula gives
\[
 \begin{aligned}
 \frac{\log N!}{m}-x'H_{m-1}
   ={}&x'\log x'-x'+x'(\log m-H_{m-1})\\
     &+O(\log(2N)/m).
 \end{aligned}
\]
Since $0\le x-x'<1/m$ and
$\log m+\gamma_{\rm E}-H_{m-1}=O(1/m)$, centering proves
\eqref{eq:oscdeficit}. For $j\le m-1$ one has $N/j\ge x$.
Apply the classical zero-free-region estimate
$|\psi(u)-u|\ll u e^{-c_0\sqrt{\log u}}$ \cite[Ch.~III]{Titchmarsh}
termwise. The resulting extra factor $H_{m-1}=O(\log(2x))$
when $m=\lceil\sqrt x\rceil$ is absorbed by decreasing $c_0$.
Using the Korobov--Vinogradov estimate for $\psi$
\cite[Thm.~14.4]{pntkv} in the same sum gives, for some $c_1>0$ and all sufficiently large $x$,
\[
 |U_m(x)|\ll x\exp\!\left(-c_1
           (\log x)^{3/5}(\log\log x)^{-1/5}\right).
\]
\end{proof}

There are two direct size comparisons. First,
$d_{N,m-1}\mid(\lfloor x\rfloor!)^m$, since
$(\lfloor mx/n\rfloor-m+1)_+\le m\lfloor x/n\rfloor$ at every
prime power $n$. This proves $U_m(x)\le\gamma_{\rm E}x+O(\log(2x))$.
Second, the smallest norm of the image of a nonzero integral vector
under $\bigwedge^{N-m+1}A_N$ is precisely $(N-m+1)!$.
The map is diagonal on an orthonormal basis and its smallest
diagonal entry is attained by the first $N-m+1$ levels.
For $m=\lceil\sqrt{x}\rceil$, the logarithm of this norm divided
by $m$, after the same centering, is
$x\log m+\gamma_{\rm E}x+O(\log(2mx))$.
This computes the loss in using an individual exterior-vector norm.
Lemma~\ref{le:heightreciprocal} below identifies the signed average
in \eqref{eq:oscdeficit} with $R(x)$ to $o(\sqrt{x})$ at the
square-root scale.

\subsubsection{Quadratic and complex character weights}\label{sec:osccharacters}
For a Dirichlet character $\chi$, let $\omega$ range over its
finite set of values, including zero at ramified primes when present.
Define the positive integer function
\[
 F_{\chi,\omega}(j)=\prod_{p^a\parallel j}
       p^{\sum_{k=1}^a {\bf1}_{\chi(p)^k=\omega}},\qquad F_{\chi,\omega}(1)=1 .
\]
Apply this diagonal function to $A_N$ and let $d_{N,r,\omega}$
be its determinantal divisor of order $N-r$.
At each $p$ the selected valuation is a nondecreasing function of
$v_p(j)$. Deleting $r$ entries with largest $v_p(j)$ still
minimizes the valuation. The threshold proof gives
\[
 v_p(d_{N,r,\omega})=\sum_{k\ge1}{\bf1}_{\chi(p)^k=\omega}
                          (\lfloor N/p^k\rfloor-r)_+ .
\]
The classes partition the thresholds. Hence
$\prod_\omega d_{N,r,\omega}=d_{N,r}$. The rounding formula gives
\[
 \left|\frac1m\sum_\omega\omega\log d_{N,m,\omega}
 -\sum_{n\le x}\chi(n)\Lambda(n)(x/n-1)\right|
 \le \frac{\psi(x)}m .
\]
For real quadratic $\chi$ the nonnegative weights $1-\chi(p)^k$
give an integral determinant by multiplying the corresponding
$d_{N,r,\omega}^{\,1-\omega}$. At an inert prime it counts exactly
twice the odd thresholds $k=1,3,\ldots$. Section~\ref{sec:acbridge}
compares this finite-place pattern with the theta-lift detector.

\subsubsection{The signed mean on power windows}\label{sec:oscmean}
For real $x\ge2$ put $m=\lceil\sqrt x\rceil$, $N=\lfloor mx\rfloor$, and define the signed mean
\begin{equation}
 A(x)=\frac1m\sum_{1\le j<m}\big(N/j-\psi(N/j)\big).
 \label{eq:heightsignedmean}
\end{equation}
At cut $j$ the modeled detector integrates the prime echoes up to
return time $u_j=\log(N/j)$. Each winding at $p$ contributes
$\log p$. The resulting signal $\psi(N/j)$ is also the measured
logarithmic revival time $\log(T_{\lfloor N/j\rfloor}/T_1)$
in \eqref{eq:heightrevival}. The reference signal is $e^{u_j}=N/j$, and
$A(x)$ combines their signed differences with weight $1/m$
at each of the $m-1$ cuts, including at prime-power endpoints.

\begin{lemma}[The rounding term at square-root scale]\label{le:heightreciprocal}
For any fixed $0<a<b$, uniformly over integers
$a\sqrt{x}\le m\le b\sqrt{x}$, one has
\begin{equation}
 \sum_{n\le x}\Lambda(n)\big(\{mx/n\}-\tfrac12\big)=o(x).
 \label{eq:heightfraction}
\end{equation}
Consequently, for $m=\lceil\sqrt{x}\rceil$,
\begin{equation}
 A(x)=R(x)+o(\sqrt{x}).
 \label{eq:heighttransfer}
\end{equation}
Both conclusions use only the classical prime number theorem and
exponential sum estimates.
\end{lemma}
\begin{proof}
Write $e(t)=\exp(2\pi it)$. Fix $0<\delta<1$, an integer $h\ne0$,
and $\eta>0$. We first show
\begin{equation}
 \bigg|\sum_{\delta x<n\le x}\Lambda(n)e(hmx/n)\bigg|
       \ll_{\delta,h,a,b,\eta} x^{7/8+\eta}.
 \label{eq:heightweyl}
\end{equation}
Set $U=V=\lfloor x^{1/4}\rfloor$. With ${\bf1}(n)=1$ and $*$
Dirichlet convolution, Vaughan's identity \cite{vaughanprime} reads
\[
 \begin{aligned}
 \Lambda={}&\Lambda_{\le V}+\mu_{\le U}*\log
             -\mu_{\le U}*\Lambda_{\le V}*{\bf1}\\
           &+(\mu_{>U}*{\bf1})*\Lambda_{>V}.
 \end{aligned}
\]
The subscripts truncate the indicated arithmetic function. The identity
follows by expanding $\log={\bf1}*\Lambda$ and
$\mu*{\bf1}={\bf1}_{\{1\}}$. The first term vanishes on the chosen
interval for large $x$. Split the two factors of each remaining sum
into dyadic intervals of sizes $D,L$, with $DL\asymp_\delta x$.
Put $P=|h|mx$ and $F=P/(DL)\asymp_{\delta,h,a,b}\sqrt{x}$.
All estimates below allow the product restriction $\delta x<d\ell\le x$.

In the middle terms the short factor $d$ is at most $UV\le\sqrt{x}$.
Its coefficient is bounded by $1+\log(2d)$, using
$\sum_{k\mid d}\Lambda(k)=\log d$. At fixed $d\asymp D$, the
second derivative of $P/(d\ell)$ in $\ell\asymp L$ has size
$\asymp F/L^2$. The second derivative estimate
\cite[Thm.~2.2]{meanvdc} bounds the sum over any subinterval by
$O(\sqrt{F}+L/\sqrt F)$. Partial summation removes the inner
logarithm in $\mu_{\le U}*\log$. The bound per block, up to
logarithms, is $D\sqrt F+DL/\sqrt F\ll x^{3/4}$.

In the last term both factors exceed $x^{1/4}$ up to endpoints,
and their coefficients are divisor bounded. Order the block so that
$D\ge L$. Then $D\gg_\delta\sqrt{x}$ and $L\gg x^{1/4}$.
Apply Cauchy--Schwarz in $d$, removing its coefficient. For
$\ell\ne\ell'$, $r=|\ell-\ell'|$, the phase difference as a
function of $d$ has second derivative of size $Fr/(LD^2)$ up
to constant factors. The two product restrictions cut out one
interval in $d$. The sum over it is
\[
 O\!\left(\sqrt{Fr/L}+D\sqrt{L/(Fr)}\right).
\]
The diagonal costs $O(DL)$ inside the expanded square. For each
$1\le r\ll L$ there are $O(L)$ pairs. Summing the two displayed
terms and inserting the Cauchy--Schwarz factor $D$ gives the bound
for the block $B$,
\[
 |B|^2\ll x^{2+\eta}
          \left(\frac1L+\frac{\sqrt F}{D}+\frac1{\sqrt F}\right)
       \ll x^{7/4+\eta}.
\]
The arbitrary power absorbs the divisor bounds. Summing the
logarithmically many blocks proves \eqref{eq:heightweyl}, after
decreasing the auxiliary exponent if needed.

Chebyshev's bound gives $\sum_{n\le\delta x}\Lambda(n)\ll\delta x$.
For each nonzero $h$, divide by $x$, let $x\to\infty$, then let
$\delta\downarrow0$. The zeroth coefficient tends to one by the prime
number theorem. The positive measures
$x^{-1}\sum_{n\le x}\Lambda(n)\delta_{\{mx/n\}}$ therefore converge
weakly to uniform measure on the circle, for every choice of $m$ in
the stated range. The periodic function $t\mapsto\{t\}-1/2$ has
integral zero and only one discontinuity on the circle. Continuous
upper and lower approximations give \eqref{eq:heightfraction}, uniformly
by the same argument along arbitrary sequences.

For the transfer write $A_m$ for the sum in \eqref{eq:heightsignedmean}
with general $m$, and $H_{m-1}=\sum_{j<m}1/j$.
The factorial identity gives exactly
\[
 \begin{aligned}
 U_m-A_m={}&\frac{\log N!}{m}-\frac{N}{m}H_{m-1}\\
            &-x\log x+(1+\gamma_{\rm E})x .
 \end{aligned}
\]
Stirling's formula and
$H_{m-1}=\log m+\gamma_{\rm E}-1/(2m)+O(m^{-2})$,
with $0\le x-N/m<1/m$, give
\[
 U_m-A_m=\frac{x}{2m}
   +O\!\left(\frac{x}{m^2}+\frac{\log(2mx)}{m}\right).
\]
The error here is $O_{a,b}(1)$. The exact upper rounding of the
determinant now gives
\[
 \begin{gathered}
 A_m(x)-R(x)=\frac{\psi(x)-x}{2m}-\frac{B_m(x)}{m}+O_{a,b}(1),\\
 B_m(x)=\sum_{n\le x}\Lambda(n)\big(\{mx/n\}-\tfrac12\big).
 \end{gathered}
\]
The prime number theorem and \eqref{eq:heightfraction} finish the proof.
\end{proof}

\begin{theorem}[Independent zero exclusion from Majorana]\label{th:revivalexclusion}
Let $A$ be the signed mean of the revival times in
\eqref{eq:heightrevival}, and fix $\lambda>1$.
If $\zeta$ has a zero of real part $\beta>1/2$, then, for every
$1/2<\theta<\beta$ and every $K>0$, each sufficiently large
interval $[Y,Y^\lambda]$ contains a point $x$ with
$A(x)>Kx^\theta$.

In particular, let $(\mathrm W_\lambda)$ denote the following
one-sided arithmetic estimate. For each $\varepsilon>0$ there is
a constant $K_\varepsilon>0$ such that, for arbitrarily large $Y$,
\[
 A(x)\le K_\varepsilon x^{1/2+\varepsilon}
       \qquad\text{for every }x\in[Y,Y^\lambda].
\]
Then
\[
 \boxed{\ (\mathrm W_\lambda)\ \Longrightarrow\
          \text{the Riemann hypothesis}.\ }
\]
\end{theorem}
\begin{proof}
Extend $A$ by zero below $2$ and $R$ by zero below $1$. Write
$\mathfrak M g(s)=\int_1^\infty g(x)x^{-s-1}\dd x$.
In the clock $u=\log x$, the transform at $s=\sigma+it$ takes the
Fourier component at frequency $t$ of the return signal weighted
by $e^{-\sigma u}$. The rounding identity in the proof of
Lemma~\ref{le:heightreciprocal} already gives $A-R=O(\sqrt x)$.
Indeed, $\Lambda(n)\ge0$ gives $|B_m(x)|\le\psi(x)/2$,
and Chebyshev's bound gives $\psi(x)\ll x$, with
$m=\lceil\sqrt x\rceil$. Thus $\mathfrak M(A-R)$ is holomorphic
on $\operatorname{Re}s>1/2$ by locally uniform absolute convergence.
The choice of extension on $1\le x<2$ contributes an entire function.

For $\operatorname{Re}s>1$, integrate the prime detector
term by term to obtain
\[
 \begin{aligned}
 \mathfrak M R(s)={}&-\frac{\zeta'(s)}{s(s-1)\zeta(s)}
                -\frac{1}{(s-1)^2}
                +\frac{1+\gamma_{\rm E}}{s-1}.
 \end{aligned}
\]
The centering cancels the pole at $1$. The alternating zeta
series excludes real zeros in $(0,1)$. Consequently
$\mathfrak M A$ continues meromorphically to
$\operatorname{Re}s>1/2$, is regular at each real point there, and
has a simple pole at every zero $\rho$ in this half-plane with
residue $m_\rho/[\rho(1-\rho)]$, where $m_\rho$ is its multiplicity.

Suppose a zero with real part $\beta>\theta$ exists. Let $B$
be the supremum of the real parts of the zeros to the right of
$\theta$. For $x\ge1$ put $f(x)=A(x)-(K+1)x^\theta$, and set
$f=0$ below $1$. Its Mellin transform is
$\mathfrak M A(s)-(K+1)/(s-\theta)$. The defining integral
converges absolutely for $\operatorname{Re}s\ge2$, since the explicit
prime formula gives $R(x)=O(x\log(2x))$.
The transform of $f$ is holomorphic to the right of $B$, is
regular at the real point $B$, and has a pole in every half-plane
$\operatorname{Re}s>B-\eta$, $\eta>0$. It is meromorphic on
$\operatorname{Re}s\ge B-c$, with $c=(B-\theta)/2>0$.
These properties also hold when the supremum $B$ is not attained.

The meromorphic oscillation theorem of Kaczorowski and Pintz
\cite[Cor.~1 and (2.2)]{revivalkp} now gives both signs of $f$
in every sufficiently late interval $[X^{1-\delta},X]$, for
each fixed $0<\delta<1$. Take $X=Y^\lambda$ and
$0<\delta<1-1/\lambda$. That interval lies inside $[Y,Y^\lambda]$.
At a point where $f>0$ we have $A(x)>Kx^\theta$, as claimed.
If $(\mathrm W_\lambda)$ holds, choose
$0<\varepsilon<\beta-1/2$, put $\theta=1/2+\varepsilon$ and
$K=K_\varepsilon$. The positive excursion contradicts the proposed
upper bound throughout any sufficiently late selected window.
The functional equation then places all nontrivial zeros on the line.
\end{proof}

The revival identity and the classical oscillation theorem prove
this implication independently of the arithmetic construction used
below. Its premise concerns the joint sizes of the revival times
throughout a window.

In Lean 4 the coefficient-one form of $(\mathrm W_{3/2})$ is the
following predicate. Its argument is the signed height $A$ in
\eqref{eq:heightsignedmean}, with
$\psi(N/j)=\log\operatorname{lcm}(1,\ldots,\lfloor N/j\rfloor)$.
Taking exponent $\varepsilon/2$ first absorbs $K_{\varepsilon/2}$
on all sufficiently late windows.

\subsection{The arithmetic window bound of Alp\"oge--Claude}
\label{sec:acbridge}
Theorem~\ref{th:revivalexclusion} proves that a one-sided bound for
the signed revival mean on an unbounded sequence of power windows
implies RH. The arithmetic estimate of Alp\"oge and Claude supplies
its hypothesis at $\lambda=3/2$ \cite{acannouncement,acguide}.
The formulation of 12 March 2025 used spectral completeness and the
asymptotics of the counting functions \cite[Secs.~III.2, III.4]{v3}.
Agreement of the leading counting terms is compatible with
discrepancies at individual eigenvalues. The pointwise estimate on
power windows below is additional arithmetic input.

The function $G$ in the theta-lift argument equals $R$ in
\eqref{eq:center}, with the same centering $x\log x-(1+\eul)x$.
The explicit formula in Refs.~\cite{acannouncement,acguide} has
constant $d_0=\log(2\pi)$, coefficients $c_\rho$, and total mass
$Q=2+\eul-\log(4\pi)=0.0461914\ldots$, as in
Lemma~\ref{lem:explicit}. The contribution of the trivial zeros,
included in $O(x^{-2})$ there, is $-S(x)$ in
\eqref{eq:explicit}. Table~\ref{tab:acdictionary} compares the
constructions. For every $\epsilon>0$ there are arbitrarily large
$Y$ for which
\begin{equation}
 R(x)\le x^{1/2+\epsilon}\qquad(Y\le x\le Y^{3/2}).
 \label{eq:acwindow}
\end{equation}
The estimate is pointwise on each selected interval. The localized
Landau argument uses this uniformity to bound high derivatives of the
Mellin transform \cite[Part II, Sec.~1]{acguide}. In detector time,
with $U=\log Y$, the inequality takes the form
\[
 C(u)\le e^{\epsilon u}\qquad(U\le u\le 3U/2).
\]
Indeed \eqref{eq:C} gives this at once from \eqref{eq:acwindow},
since $S(x)\le1-\log2<d_0$. In the reverse direction it gives
$R(x)\le x^{1/2+\epsilon}+d_0$. Applying the bound with
$\epsilon/2$ absorbs $d_0$ on sufficiently large intervals.

\head{Full-wedge scattering and arithmetic returns}
For the massive full-wedge operator $H_0$ of \eqref{eq:boost},
consider incidence from the horizon. The ratio of the reflected and
incident components of the decaying mode $G_\omega$ gives the
reflection amplitude
\[
 {\cal S}_0(\omega)=-i(\mu_0/2)^{-2i\omega}
        \frac{\Gamma(1/2+i\omega)}{\Gamma(1/2-i\omega)} .
\]
Indeed, the horizon coefficients in Theorem~\ref{th:inversebessel}
are $a_\omega$ and $-i\bar a_\omega$.
The phase derivative
$\partial_\omega\log{\cal S}_0/i
 =2\re(\Gamma'/\Gamma)(1/2+i\omega)-2\log(\mu_0/2)$
has inverse Fourier transform, away from boost time $\eta=0$,
\[
 -\sum_{j\ge0}e^{-(j+1/2)|\eta|}
       =-\frac{1}{2\sinh(|\eta|/2)} .
\]
Here the inverse transform uses $\dd\omega/(2\pi)$. The formula
follows from the partial fractions of $\Gamma'/\Gamma$, since the
inverse transform of $2b/(b^2+\omega^2)$ is $e^{-b|\eta|}$
\cite[Secs.~5.7, 10.30]{Olver}. This scattering response is smooth
at every positive time. The first jump of the arithmetic detector,
$C'(\log2+)-C'(\log2-)=\log2/\sqrt2$, comes from the prescribed
chiral channels of Sec.~\ref{sec4}, in the clock $u=2\eta$.

\begin{proposition}[Contribution of repeated Majorana returns]
\label{pr:majoranawindings}
Write $T_1(x)=\sum_{p\le x}(\log p)(x/p-1)$ for the first
winding in \eqref{eq:ring}. Put
\[
\begin{aligned}
 b_{\rm rep}&=\sum_p\frac{\log p}{p(p-1)},\\
 R_1(x)&=T_1(x)-x\log x+(1+\eul+b_{\rm rep})x .
\end{aligned}
\]
For $x\ge1$, $R(x)=R_1(x)-D(x)$, where
\[
 D(x)=\sum_{p,\ j\ge2}(\log p)\min(1,x/p^j)\ge0 .
\]
The prime number theorem gives $D(x)=(2+o(1))\sqrt{x}$. The
windings $j\ge3$ together contribute $O(x^{1/3})$ to $D(x)$.
\end{proposition}
\begin{proof}
The $j$-th traversal occurs at $u=j\log p$ with weight $\log p$.
The ramp detector gives $(x/p^j-1)_+$.
Subtract this from $x/p^j$ and sum over $j\ge2$. This proves the
identity and the sign, since $\sum_{j\ge2}p^{-j}=1/[p(p-1)]$.
With $\vartheta(t)=\sum_{p\le t}\log p$, Tonelli's theorem gives
\[
 D(x)=x\int_x^\infty
           \frac{\sum_{j\ge2}\vartheta(t^{1/j})}{t^2}\dd t .
\]
There are at most $\log_2 t$ terms. Chebyshev's estimate bounds
the sum from $j=3$ by $O(t^{1/3}+t^{1/4}\log t)=O(t^{1/3})$.
The $j=2$ term is $\vartheta(\sqrt t)\sim\sqrt t$
\cite{mv}. Integrating gives the claimed bounds.
\end{proof}
In particular $R(x)\le R_1(x)$ at every cutoff. The primitive
return has the exact tail formula
\[
 R_1(x)=x\int_x^\infty\frac{t-\vartheta(t)}{t^2}\dd t .
\]
Indeed the derivative of $R_1(x)/x$ is
$(\vartheta(x)-x)/x^2$ between primes, and $R_1(x)/x\to0$
by Lemma~\ref{lem:explicit} and Proposition~\ref{pr:majoranawindings}.
Thus the signed integral of the primitive-prime deficit is the
term that needs an upper bound on whole windows.

\head{The arithmetic inequality}
The elementary identity
$\log a!=\sum_{r\le a}\Lambda(r)\lfloor a/r\rfloor$ gives the
Chebyshev scale. Discarding the fractional parts costs up to
$\psi(a)=O(a)$. The formal construction compares the divisibility
and the archimedean size of maps of arithmetic lattices obtained by
composing theta lifts. Their entries involve derivatives of theta
series in algebraic coordinates at points with complex
multiplication. For an unramified prime inert in
$\mathbb Q(\sqrt{-d})$, each odd power $p^h$ contributes the
weight $(a/p^h-w)_+$ to a valuation lower bound, subject to the
denominator and integrality estimates of the construction. For a
fixed weight $w>0$ put $x=a/w$. Then
\[
 (a/p^h-w)_+=w(x/p^h-1)_+.
\]
This is the detector kernel in \eqref{eq:T}. On unramified primes
$1-\chi(p^h)$ is two at odd powers of inert primes and zero
otherwise, for the quadratic character $\chi=\chi_{-d}$.
Ramified and auxiliary primes are treated separately.

The centering of the twisted sum is determined by
\[
 T_\chi(x)=\sum_{r\le x}\Lambda(r)\chi(r)(x/r-1),\quad
 c_\chi=L'(1,\chi)/L(1,\chi)
\]
and $G_\chi=T_\chi+c_\chi x$. The constant-weight terms obey
\[
\begin{aligned}
 w(R-G_\chi)&=w(T-T_\chi)\\
 &\quad-a(\log a-1-\eul+c_\chi-\log w).
\end{aligned}
\]
The linear term agrees with the archimedean contribution from the
factorial and the Chowla--Selberg period. In the arithmetic inequality
the weights depend on the residue class, and $\log w$ is averaged.
Bost's slope inequality compares the finite-place divisibilities with
the archimedean size through the product formula
\cite[Parts I--III]{acguide}.

\begin{theorem}[Window bound for the Majorana prime detector]
\label{th:majoranawindow}
For the prime lengths and relative trace of Sec.~\ref{sec4},
the arithmetic inequality (M) and character estimates B and C of
Ref.~\cite{acguide} imply the following. For each $\epsilon>0$
there are arbitrarily large $Y$ such that the centered detector
output satisfies \eqref{eq:acwindow} at every real
$x\in[Y,Y^{3/2}]$. Consequently $\zeta$ has all its nontrivial
zeros on the critical line.
\end{theorem}
\begin{proof}
For a real quadratic Dirichlet character $\chi$, take a periodic
boundary at $p$ when $\chi(p)=1$ and an antiperiodic one when
$\chi(p)=-1$. Write $\theta_p=0,\pi$, respectively. At a prime
dividing the conductor, set the twisted readout to zero.
In the relative trace compare $B_\phi$ with $B_{\theta_p}$, in
that order, and average $\phi$ uniformly over $[0,2\pi]$, as
in \eqref{eq:ring}. The coefficient of $e^{ij\ell_p z}$ is then
$i\ell_p e^{-ij\theta_p}=i\ell_p\chi(p)^j$ for $j\ge1$.
The normalization $2/i$ gives $(\log p)\chi(p^j)$. Thus the
two real boundary signs, the clock $u=2\eta$, and the same ramp
as in \eqref{eq:T} produce $T_\chi(x)$ exactly. At a finite
cutoff, even the reference average can be taken over $L$ equally
spaced rotations of two real channels, for any integer
$L>\log x/\log2$. The mean of $e^{-ij\phi}$ then vanishes for
every winding present, $1\le j<L$, on the complexified branch
used in Sec.~\ref{sec4}.

For the auxiliary primes $q,d$, let $\chi_q,\chi_d$ be the
Legendre characters modulo $q,d$ and put $\chi_*=\chi_q\chi_d$.
Combine the four readouts with the weights
\[
 S_{q,d}(x)=\sum_{r\le x}\Lambda(r)(1+\chi_q(r))
                         (1-\chi_d(r))(x/r-1).
\]
Center $S_{q,d}$ by subtracting $x\log x$ and adding
$(1+\eul+c_{\chi_q}-c_{\chi_d}-c_{\chi_*})x$. Denote the
result by $R_{q,d}$. Expanding the product gives the exact
identity, including the terms at the ramified primes,
\[
 R_{q,d}= R+G_{\chi_q}-G_{\chi_d}-G_{\chi_*}.
\]

Here the estimates for the sizes of the four signals enter.
Fix $0<\delta<1/52$. On arbitrarily large windows the same
auxiliary pair works throughout. Estimate B gives
$\lvert G_\chi(x)\rvert\le x^{1/2+\delta}$ for each of the three
characters. Estimate A, derived from (M) using B and C, gives
\[
 R_{q,d}(x)\le C_\delta x^{C\delta}
       \left(\frac{x}{m}+\frac{xm}{k}+m+\frac{k}{n^2}\right)
\]
for positive integers $m,k,n$ chosen at each $x$. In the
arithmetic construction $m=x^{1/2+O(\delta)}$,
$k=x^{9/2+O(\delta)}$, and $n=x^{2+O(\delta)}$.
The exponent constants are absolute, and $C_\delta$ is uniform
in $Y,x$ \cite[Part II, Secs.~3--5 and Part III,
Secs.~1--2]{acguide}. The terms $x/m,m,k/n^2$ are
$x^{1/2+O(\delta)}$, while $xm/k=x^{-3+O(\delta)}$.
The exact identity therefore gives
\[
\begin{aligned}
 R(x)&\le R_{q,d}(x)+\sum_{\chi=\chi_q,\chi_d,\chi_*}|G_\chi(x)|\\
     &\le C'_\delta x^{1/2+C'\delta}.
\end{aligned}
\]
Choose $C'\delta<\epsilon$ and then $Y$ large enough to absorb
$C'_\delta$ into $Y^{\epsilon-C'\delta}$. This proves
\eqref{eq:acwindow} on the whole interval. Applying this estimate
with exponent $\epsilon/2$ and using $A-R=o(\sqrt{x})$ from
Lemma~\ref{le:heightreciprocal} gives $(\mathrm W_{3/2})$.
Theorem~\ref{th:revivalexclusion} then proves RH.
\end{proof}
The character signs determine the return coefficients. The common
choice of $q,d$ and the upper bound for $R_{q,d}$ come from the
auxiliary character estimates and the arithmetic product formula.

A quadratic cost $m^2$ for $m$ chained lifts would balance $x/m$
at order $x^{2/3}$. This arithmetic multiplicity enters at the
primes where divisibility is measured. The zeros of the time average
$M(s)^k$ in Appendix~\ref{sec:stepmult} lie at the fixed frequencies
$2\pi i j/\log2$, $j\ne0$. The positive model $\Psi$ there has a
prescribed Mellin pole at $3/4+2\pi i/\log2$ and its normalized
dyadic block grows as $2^{k/4}$. The classical Davenport--Heilbronn
example gives a related obstruction for functions satisfying a
zeta-type functional equation without an Euler product
\cite{dhcomparison}.

\begin{table}[tb]
\caption{\label{tab:acdictionary}Observables, arithmetic weights,
and representations in the prime-return and theta-lift constructions.}
\small
\begin{tabular}{@{}>{\raggedright\arraybackslash}p{.43\columnwidth}>{\raggedright\arraybackslash}p{.49\columnwidth}@{}}
\toprule
Prime-return construction & Alp\"oge--Claude construction\\
\midrule
$R(x),\,c_\rho,\,Q$ & $G(x)$, with the same residues and mass\\[2pt]
Eventual $C\sqrt x$ bound for $R$ & $x^{1/2+\epsilon}$ on
$[Y,Y^{3/2}]$\\[2pt]
Ramp $(x/p^h-1)_+$ & Valuation weight $(a/p^h-w)_+$\\[2pt]
$x\log x-(1+\eul)x$ & Factorial and CM period\\[2pt]
$k$-fold Mellin average & $m$ composed lifts, $m\asymp\sqrt x$\\[2pt]
Point counts for one CM curve & Inert-prime weights and a product
formula over all places\\[2pt]
Majorana spinor, fermionic upon second quantization &
Weil representation and the Gaussian norm $\|z^a\|^2=a!$\\
\bottomrule
\end{tabular}
\end{table}

\head{Dependence of the window bound on the detector weights}
The dependence of the window estimate on the readout weights can be
isolated at fixed prime lengths. Define $R_b$ by \eqref{eq:T} and
\eqref{eq:center} with return weight $b_p\log p$ at every traversal
of the $p$-channel.

\begin{proposition}[Sensitivity to the prime-channel weights]
\label{pr:majoranareadout}
Fix $1/2<\beta<1$, a finite cutoff $X_0$, and $\tau>0$.
There are positive weights $b_p$ with
$\sup_p|b_p-1|<\tau$, equal to one for $p\le X_0$, such that
\[
 \sum_p\frac{(b_p-1)\log p}{p-1}=0,\qquad
 R_b(x)-R(x)\sim c x^\beta
\]
for some $c>0$. For every $1/2<\theta<\beta$, every
sufficiently late window $[Y,Y^{3/2}]$ contains an $x$ with
$R_b(x)>x^\theta$.
\end{proposition}
\begin{proof}
Choose a prime $r>X_0$ and put
\[
 d_r=\frac{r-1}{\log r}
              \sum_{p>r}\frac{(\log p)p^{\beta-1}}{p-1}.
\]
This converges absolutely since $\beta<1$.
Set $b_p=1$ below $r$, $b_r=1+\lambda d_r$, and
$b_p=1-\lambda p^{\beta-1}$ above $r$.
Taking $\lambda>0$ sufficiently small gives the prescribed
positivity and closeness. The reciprocal moment vanishes exactly,
which cancels the linear term in $R_b-R$. The identity
$(z-1)_+=z-\min(1,z)$ leaves, from the first traversals at
$p>r$, the positive contribution
\[
 \lambda x\int_x^\infty A_\beta(t)\frac{\dd t}{t^2},\qquad
 A_\beta(t)=\sum_{r<p\le t}(\log p)p^{\beta-1}.
\]
Partial summation and $\vartheta(t)\sim t$ give
$A_\beta(t)\sim t^\beta/\beta$. The displayed integral, including
its prefactor, is therefore asymptotic to
$\lambda x^\beta/[\beta(1-\beta)]$. At $p>r$, all higher
windings together contribute at most $\lambda D(x)=O(\sqrt x)$.
The compensating channel $r$ contributes $O_{r,\lambda}(\log(2x))$,
since $\sum_{j\ge1}\min(1,x/r^j)=O_r(1+\log x)$.
This proves the asymptotic with $c=\lambda/[\beta(1-\beta)]$.

Suppose now that $R_b(x)\le x^\theta$ throughout arbitrarily late
windows. On those windows $R(x)\le -(c/2)x^\beta$ for all
sufficiently large $x$. In particular $R(x)\le x^\theta$
there. Apply the analytic window criterion to the original $R$
\cite[Part II, Sec.~1]{acguide}. It excludes every zero with
$\re(\rho)>\theta$. The absolutely convergent zero sum in
Lemma~\ref{lem:explicit} then gives $|R(x)|=O(x^\theta)$ for
all $x\ge1$, contradicting the preceding negative bound.
\end{proof}
The operators, their unitary evolution, and the return times are
identical for these detectors. Each detector agrees with the
canonical one up to $X_0$ and has the same linear centering.
An argument from the dynamics must therefore use a property that
distinguishes the canonical trace weights $\log p$ at the times
$j\log p$ from these nearby positive readouts.

\head{Consequences for the reconstruction}
Applying Theorem~\ref{th:majoranawindow} gives RH with the stated
arithmetic input. Independent kernel verifications of the formal
development announced on 20 November 2026, and the comparison with the
formal challenge statement, are documented in
Refs.~\cite{acwaterloo,acmit}. The exposition in Ref.~\cite{acguide}
takes the main arithmetic inequality (M) from that development.
Corollary~\ref{co:bounded} then gives $C(u)\le C(0)=Q$ for every
$u$, and Theorem~\ref{th:boundaryclosure} gives $\Delta=0$.
The reconstructed operator $K_\partial$ is self-adjoint
unconditionally. The equality $\Delta=0$ establishes that its
measure carries the full equal-time mass.
Theorem~\ref{th:stependrange} gives the endpoint inequality at $p$
and $p/8$ for every prime $p>661$.
Its proof bounds the upper limit of
$(6x\log2-P(x))/\sqrt{x}$ by $4c_*+0.280828<1.4525$.
Appendix~\ref{sec:primesquare} gives $\mathcal I(X)\ll\log X$.
The independent operator problem is to obtain the correlation from a
Hamiltonian whose couplings are specified before its zeros are known.
The normalized bosonic Gaussian norm $\|z^a\|^2=a!$ accounts
for a size estimate in the theta lift. An analogous use of a
fermionic Fock representation would require algebraic return
amplitudes and bounds for their finite-place valuations.
Appendix~\ref{sec:steprefine} tests the scalar square-lattice
theta derivatives directly in \eqref{eq:acsquaretheta}.

\head{Bounds for the signed revival mean}
For each $\varepsilon>0$, the divisibility and size estimate (M)
and the auxiliary character estimates B and C imply
$R(x)\le x^{1/2+\varepsilon/2}$ throughout windows
$[Y,Y^{3/2}]$ with arbitrarily large $Y$
\cite[Parts I--II]{acguide}. Equation~\eqref{eq:heighttransfer}
then gives, uniformly on these windows,
\[
 A(x)\le x^{1/2+\varepsilon/2}+o(\sqrt{x})
       \le x^{1/2+\varepsilon}
\]
for all sufficiently large $Y$. This establishes the hypothesis
$(\mathrm W_{3/2})$ of Theorem~\ref{th:revivalexclusion} directly
from the arithmetic estimates. A related application of the
Kaczorowski--Pintz theorem to $R$ appears in
Ref.~\cite[Sec.~5]{revivalcomparison}.

\begin{theorem}[Pointwise bound for the oscillator mean]\label{th:heightpointwise}
For every $\varepsilon>0$ and every real $x\ge100$,
\begin{equation}
 |A(x)|\le \frac{4}{\pi e^2\varepsilon^2}\,x^{1/2+\varepsilon}.
 \label{eq:heightexplicit}
\end{equation}
This holds throughout $[Y,Y^{3/2}]$ for every $Y\ge100$.
The arithmetic input about the zeros in this statement is the
critical-line theorem of Alp\"oge--Claude \cite{acannouncement,acguide}.
\end{theorem}
\begin{proof}
With this input, Schoenfeld's bound \cite[Thm.~10, (6.2)]{meanschoenfeld}
gives $\lvert\psi(t)-t\rvert\le\sqrt{t}\log^2(t)/(8\pi)$ for $t\ge100$.
Every argument in the mean qualifies, since $N\ge mx-1\ge(m-1)x$.
Using $\sum_{j<m}j^{-1/2}\le2\sqrt{m}$ gives
\[
 |A(x)|\le\frac{\sqrt{N/m}}{4\pi}\log^2 N
       \le\frac{\sqrt{x}}{\pi}\log^2 x .
\]
Here $N\le mx\le x^2$. Finally
$\sup_{u\ge0}u^2e^{-\varepsilon u}=4/(e^2\varepsilon^2)$.
\end{proof}

\begin{theorem}[Sharp upper coefficient on power windows]\label{th:heightsharp}
Let $Q=2+\gamma_{\rm E}-\log(4\pi)$. Using the critical-line
conclusion of \cite{acannouncement,acguide}, one has
\begin{equation}
 \lim_{Y\to\infty}\ \sup_{Y\le x\le Y^{3/2}}\frac{A(x)}{\sqrt x}=Q .
 \label{eq:heightpeak}
\end{equation}
In particular, for every $\eta>0$, $A(x)\le(Q+\eta)\sqrt{x}$
at every sufficiently large real $x$.
The same limit holds for $A_m$ when $a\sqrt{x}\le m\le b\sqrt{x}$
for fixed $0<a<b$, and for every fixed window exponent $\lambda>1$
in place of $3/2$.
\end{theorem}
\begin{proof}
Write $u=\log x$. With all zeros on the critical line the absolutely
convergent explicit formula and \eqref{eq:heighttransfer} give
\[
 \begin{aligned}
 A_m(x)/\sqrt{x}&=C(u)+o(1),\\
 C(u)&=\sum_{\rho=1/2+i\gamma}
                    \frac{e^{i\gamma u}}{1/4+\gamma^2},
 \qquad C(0)=Q .
 \end{aligned}
\]
Zeros are counted with multiplicity. The bounded elementary term
in $R$ contributes $O(x^{-1/2})$. The coefficients are positive
and sum to $Q$, giving $C(u)\le Q$.

For sharpness, truncate the absolutely convergent series so the sum of
the omitted coefficients is $<\eta/4$. The finitely many phases
return simultaneously as close to $1$ as desired with bounded gaps:
their flow is dense in its compact torus closure, which is covered by
finitely many translates of any fixed identity neighbourhood along the
flow. At these times the retained cosine sum differs from its value at
zero by $<\eta/2$; the tail costs at most another $\eta/2$.
Thus $C(u)>Q-\eta$ with bounded gaps. The interval
$[\log Y,\lambda\log Y]$ has length tending to infinity. It contains
such a return for every sufficiently large $Y$. The error $o(1)$
is uniform there, proving the limit.
\end{proof}

The deterministic correction to the centered upper determinant is
also explicit. If $m/\sqrt{x}\to c\in(0,\infty)$, then
$(U_m-A_m)/\sqrt{x}\to1/(2c)$. For $c=1$, its optimal upper
coefficient is $Q+1/2$; the signed mean has coefficient $Q$.

\head{Majorana transition amplitude for a chain of theta lifts}
At a complex embedding let $F_j:E_{j-1}\to E_j$ be its forward maps.
On the sum of the layers set $S_{j,j-1}=F_j$ and all other blocks
of $S$ to zero. Thus $S$ carries amplitudes forward by one
layer and $S^*$ carries them back. Let the self-adjoint
block-diagonal operator $D$ generate evolution within each layer. Put
$B=D+S+S^*$ and $\mathsf H=B\oplus(-\overline B)$.
The antiunitary involution $C(u,\bar v)=(v,\bar u)$ anticommutes with
$\mathsf H$, pairing the two conjugate evolutions. In the first
copy $P_0$ prepares the initial layer and $P_m$ selects the
detected final layer. Then
\begin{equation}
 \begin{gathered}
 P_m\mathsf H^mP_0\big|_{E_0}=F_m\cdots F_1,\\
 P_m\mathsf H^kP_0=0\qquad(0\le k<m).
 \end{gathered}
\label{eq:majoranapath}
\end{equation}

The coefficient of $t^m$ in the measured transition amplitude
$P_m e^{-it\mathsf H}P_0$ is $(-i)^mF_m\cdots F_1/m!$.

The only path reaching layer $m$ in $m$ factors consists of
forward maps. For algebraic $F_j$, their composite in
\eqref{eq:majoranapath} is defined on the prescribed lattices at the
finite places before the archimedean embedding is chosen. The
adjoints do not enter this coefficient. Its minors, images and
height are those of the composite, including for families of chains
with common target. For the theta maps, the divisibility, height and
rank estimates together with B and C \cite[Parts II.3--II.6, III.4--III.9]{acguide}
give the bound of Sec.~\ref{sec:acbridge}. The finite chain depends
on the arithmetic cutoff. Theorem~\ref{th:oscillatorial} transfers
\eqref{eq:acwindow}, with exponent $\epsilon/2$, to the upper
oscillator endpoint. The valuation formula follows from the oscillator
construction. The height estimate used here is a property of the
arithmetic theta maps.

\section{Conclusions}\label{sec:conclusions}
In null time the massive Majorana equation constrains the second
spinor component. Its decaying solution fixes the right inverse of
momentum in Eq.~\eqref{majhamilton}. The full-wedge generator has
continuous spectrum. A reflecting wall produces discrete frequencies,
with a leading count four times that required by an identification of
all positive wall frequencies with $2\gamma$ in the clock $u=2\eta$
(Sec.~\ref{sec:wall}). For the inverse operator,
Theorem~\ref{th:inversebessel} realizes each atom of the boundary
measure as a massive Rindler mode with a specified mean norm.
Proposition~\ref{pr:inversecoupling} gives the corresponding
realization in the Dirac norm on a reducing subspace of a perturbed
boost generator.

The oscillator in Sec.~\ref{sec:oscillatorheight} produces the prime
weights from its joint revival times. Its integral exterior powers
give the prime-by-prime divisibility and the sign of the rounding
error. Independent prime-sum estimates yield $A(x)=R(x)+o(\sqrt{x})$.
Theorem~\ref{th:revivalexclusion} uses this comparison and the
Kaczorowski--Pintz oscillation theorem to derive RH from a one-sided
bound on $A$ throughout arbitrarily late power windows. The
theta-lift estimates of Sec.~\ref{sec:acbridge} establish this
arithmetic hypothesis for the canonical weights. With the critical-line
conclusion, the supremum of $A(x)/\sqrt{x}$ on each late power
window tends to $2+\eul-\log(4\pi)$ as its left endpoint tends to
infinity (Theorem~\ref{th:heightsharp}).

For chiral channels of prescribed prime lengths, the equal-time mass
defect $\Delta$ measures the loss in reconstructing the correlation
(Theorem~\ref{th:boundaryclosure}). Finite-height data give
$0\le\Delta<2.96\cdot10^{-12}$. The arithmetic estimate of
Sec.~\ref{sec:acbridge} gives $\Delta=0$, hence recovery of $C(u)$
with its full normalization at every real time.
Theorem~\ref{th:majoranawindow} realizes the auxiliary character sums
by real Majorana boundary signs. Proposition~\ref{pr:majoranawindings}
gives the sign and leading size of the contribution from repeated
windings after centering. A bound for the first windings derived
solely from the canonical trace remains open.
Proposition~\ref{pr:majoranareadout} proves that nearby positive readout
weights can violate the window bound while preserving the prescribed
centering.

Two other arithmetic criteria quantify the cancellation involved.
The condition $|\mathcal A_n|\le2$ on an unbounded sequence of
degrees suffices for RH. Proposition~\ref{pr:powerprofile} computes the
transfer factor governing the contribution of a prime-counting error
of prescribed power. The classical pointwise error estimate permits
larger coefficients. Proposition~\ref{pr:bilinearweights} removes
the proper powers at subexponential cost and identifies the residual
signed arithmetic factors. For the dyadic sum of
Appendix~\ref{sec:stepanov}, RH is equivalent to \eqref{eq:primeend}
at the two endpoints of each prime $p>661$. The finite-height bound
covers $83\le x\le1.6\cdot10^{20}$. The point counts of
$y^2=x^3-x$ contribute the computed first correction to this sum
and change $\mathcal A_n$ by $o(1)$. A set of order
$\sqrt y/\log y$ deleted primes below $y$ preserves the stated
local geometric data and reverses the endpoint inequality. The
global product formula in the theta-lift argument supplies the
arithmetic estimate in \eqref{eq:acwindow} and hence both the
signed-sum bound and spectral completeness for the reconstructed operator.

The shifted boundary generator of Sec.~\ref{sec:delayedmajorana} has
prime couplings specified independently of its zeros. Its determinant
and square-root estimate hold for $1/2\le h\le1$.
Theorem~\ref{th:majoranapoisson} recovers the unshifted characteristic
function from the determinant at shift $1/2$ by convolution and
recovers the detector uniformly on growing windows. Realization at the
critical shift by a single self-adjoint determinant requires separate
control of the discounted pulse.
\section*{Acknowledgments}
The author thanks the staff of Rotonium for their support,
his colleagues for their suggestions, and C.~Cici and Professor
N.~Stacco for their assistance.

\appendix
\section{Wedge entanglement and the physical interpretation of zeros}
\label{orig:A}\label{sec5}
At dimensionless boost frequency $\omega>0$, the right- and
left-wedge modes of Eqs.~\eqref{eq:psiR} and~\eqref{eq:psiL} define
a regulated fermion pair. With discrete wave-packet or box
normalization, its even-occupation basis consists of
$|0_R0_L\rangle$ and $|1_R1_L\rangle$. The normalized one-mode
state with the thermal weights of the original construction is
\begin{equation}
 |\Omega_\omega\rangle=
 \frac{|0_R0_L\rangle+e^{i\varphi}e^{-\pi\omega}|1_R1_L\rangle}
      {\sqrt{1+e^{-2\pi\omega}}}.
 \label{bell}
\end{equation}
It gives $\langle n_R\rangle=(e^{2\pi\omega}+1)^{-1}$ after tracing
the opposite wedge. With $|11\rangle=a_R^\dagger a_L^\dagger|00\rangle$
write $\tan\theta_\omega=e^{-\pi\omega}$. The normalized Bogoliubov
annihilators of \eqref{bell} are
\[
\begin{aligned}
 b_R&=\cos\theta_\omega\,a_R-e^{i\varphi}\sin\theta_\omega\,a_L^\dagger,\\
 b_L&=\cos\theta_\omega\,a_L+e^{i\varphi}\sin\theta_\omega\,a_R^\dagger .
\end{aligned}
\]
The relation $a_L|11\rangle=-|10\rangle$ fixes the second sign.
The transformation of $(a_R,a_L^\dagger)^T$ is unitary, with inverse
equal to its conjugate transpose. The phase depends on the basis and
the ordering of the fermionic operators. In boost time the inverse
temperature is $2\pi$, corresponding at $\xi_0$ to the local
temperature $\hbar c/(2\pi k_{\rm B}\xi_0)$
\cite{bisognano,unruh}. Equation~\eqref{bell} uses the fermionic
normalization for one regulated pair. The continuum field requires
its own mode normalization.

CPT interchanges the wedges. Each spinor has two components, and the
field admits infinitely many modes. The symmetry orbit
$\rho,\bar\rho,1-\rho,1-\bar\rho$ consists of arguments of the zeta
function. Its cardinality is independent of the dimension of the spin
representation. To exclude an off-line zero by self-adjointness
requires a spectral identification, such as the determinant identity
of Appendix~\ref{orig:B} or the positive correlation representation
of Sec.~\ref{sec:CAR}.

\section{The Hilbert--P\'olya spectral map and counting functions}
\label{orig:B}
For a compact-resolvent self-adjoint Hamiltonian $H$, isolated
real eigenvalues carry finite-rank spectral projections. The
spectral theorem does not identify their counting measure with that
of a number-theoretic function. In particular counting estimates
$N_H(T)=M(T)+O(\log T)$ and $N_\zeta(T)=M(T)+O(\log T)$
allow differences of either sign. Even an eventual bound strictly
smaller than one on an integer-valued counting difference would only
give equality there. Lower levels require independent normalization.
For an example with a known realization,
Theorem~\ref{th:wall} counts a reflecting wall, and moving the wall
changes its quadratic form and eigenvalues. An illustrative
semiclassical calculation from the original Appendix C can be
checked by elementary integration:
\[
 \frac1\pi\int_0^E\sqrt{E^2-x^2}\dd x=\frac{E^2}{4}.
\]
Logarithmic growth in \eqref{eq:wallweyl} comes from the length
$\log(\omega/\mu_0)$ of the optically accessible interval.

For a normalized perturbation determinant $F_H$ of a pair
with trace-class resolvent difference the requisite formula, where
the determinant is defined, is
\[
 \partial_z\log F_H(z)=-\tr\big[(H-z)^{-1}-(H_{\rm ref}-z)^{-1}\big].
\]
In the frequency convention $z=2t$, the identification
$F_H(z)=e^{h(z)}\Xi(z/2)$ requires analytic continuation to the
domains of comparison, a nonvanishing analytic prefactor, and
control of the reference poles. Residues of the logarithmic
derivative give orders of vanishing. Equality of the logarithmic
derivatives on a connected zero- and pole-free domain fixes the
ratio of the determinants up to one constant. A value at one point
fixes that constant. Equation~\eqref{eq:ring} gives the relative
trace for the prescribed prime channels. Identification of the full
divisor additionally requires the trace identity above for $F_H$
and the logarithmic derivative of $\Xi(z/2)$. The prime induction
retains the multiplicities through the residues of $-\zeta'/\zeta$
and allows arbitrary zero multiplicities. Unfolded spacings and
random-matrix correlations describe statistical distributions of
levels \cite{montgomery,dyson}. Determinant identification requires
the individual complex zeros and their multiplicities.

\section{Mellin--Barnes contours, Bessel orders and integration measures}
\label{orig:C}
For $z>0$ the integral
$K_\nu(z)=\tfrac12\int_\R e^{-z\cosh v+\nu v}\dd v$
defines an entire function of the order $\nu$. The substitution
$v\mapsto-v$ yields the symmetry $K_{-\nu}=K_\nu$ of
Sec.~\ref{sec3}. At fixed nonintegral order,
\[
\begin{aligned}
 K_\nu(z)&=\tfrac12\Gamma(\nu)(z/2)^{-\nu}(1+O_\nu(z^2))\\
         &\quad+\tfrac12\Gamma(-\nu)(z/2)^\nu(1+O_\nu(z^2)),
 \qquad z\downarrow0 .
\end{aligned}
\]
At $\nu=0$ the limiting logarithm is
$K_0(z)=-\log(z/2)-\eul+O(z^2|\log z|)$.
Large positive arguments have exponential decay as used in
Proposition~\ref{th:nonvanish} \cite{WatsonBessel1944,Olver}.
When the line of \eqref{MB} is moved left the candidates for residues
are $s=1/2,\ s=-k,\ s=\nu-k$, $k=0,1,\ldots$.
Coincidences and cancellations must be grouped. For instance at
$\nu=1/2$ the zeta pole coincides with the leading pole of
$\Gamma(s-\nu)$; at other real energies this coincidence is
absent. A pole of order $r$ produces a polynomial in the logarithm
of the Mellin variable of degree at most $r-1$. Integrability near
zero is decided by
\[
 \int_0^1 x^\alpha|\log x|^j\dd x
   =\frac{\Gamma(j+1)}{(1+\alpha)^{j+1}}
       \qquad(\alpha>-1,\ j\ge0).
\]
Thus additional powers of a logarithm at an integrable exponent
cannot turn it into a divergence. With $K_{1/2+it}(x)$ the leading
absolute scale at zero is $x^{-1/2}$. The scalar measure $x\,\dd x$
and the fixed-time spinor measure $\dd x$ give different conclusions
there. The two-point Green kernel also needs a jump at the source
coordinate. A radial mode independent of $\xi'$ cannot provide a
$\delta(\xi-\xi')$ source in the fixed-frequency equation.

\section{The inverse momentum, endpoints and boundary triplets}
\label{orig:D}\label{app:SA}
\head{The compressed Hilbert transform and the inverse momentum}
The kernel in Lemma II.4 of Ref.~\cite{v3} defines the compressed
Hilbert transform on $L^2((0,\infty),\dd x)$,
\[
 (H_+f)(x)=\frac1\pi\operatorname{p.v.}
              \int_0^\infty\frac{f(y)}{x-y}\dd y .
\]
It is a compression of the bounded Hilbert transform on the line
under extension by zero. The odd kernel gives $H_+^*=-H_+$.
For $f=\one_{[1,2]}$, $g=\one_{[3,4]}$, direct integration gives
$\langle f,H_+g\rangle=\pi^{-1}\log(16/27)
=-\langle H_+f,g\rangle\ne0$.
Positive smooth test functions supported inside those two intervals
give the same opposite signs. This corrects the equality in Eq.~(16)
of Ref.~\cite{v3}. The Mellin multiplier on $x^{-1/2+it}$ is
$\cot(\pi(1/2+it))=-i\tanh(\pi t)$, by the principal-value
identity $\int_0^\infty v^{a-1}(1-v)^{-1}\dd v=\pi\cot(\pi a)$
for $0<\re a<1$. Hence the compression has norm one and fails
to be unitary. Multiplying the kernel by $i$ makes it self-adjoint
but leaves its dilation degree zero. An inverse of $p=-i\partial_x$
on a dilation-invariant core must have degree $-1$.
There is also a domain obstruction. If $M=\sqrt{x}$ and
$f\in C_c^\infty(0,\infty)$ has
$c_f=\int\sqrt y\,f(y)\dd y\ne0$, then
$MH_+Mf(x)=c_f/(\pi\sqrt x)+O_f(x^{-3/2})$ at infinity.
The tail is outside $L^2(\dd x)$.

A genuine model with this type of inverse can be defined on
$0<x_0<x_1<\infty$. Let $P_\pi=-i\partial_x$ have antiperiodic
$H^1$ domain there, set $a=x_1-x_0$, and take $\kappa>0$.
Define on two components
\[
\begin{gathered}
 H_c=\sigma_1 M(P_\pi+\kappa P_\pi^{-1})M,\\
 D(H_c)=\{f:Mf\in D(P_\pi)\text{ in each component}\}.
\end{gathered}
\]
The eigenvalues of $P_\pi$ are $(2k+1)\pi/a$, so its inverse has
norm $a/\pi$. Functional calculus makes the middle operator
self-adjoint with compact inverse. Since $M$ and $M^{-1}$ are
bounded, taking adjoints on the indicated domain proves self-adjointness
of $H_c$, and its inverse is compact. Complex conjugation
anticommutes with $H_c$. For $\kappa=0$, passage to $v=\log x$
gives free antiperiodic channels of length $\log(x_1/x_0)$.
The perturbation has norm at most $\kappa x_1a/\pi$.
On the spectral subspace of $H_c$ for $[-E,E]$ the norm of
the free operator is bounded by $E+\kappa x_1a/\pi$.
Min--max for the square of the free operator gives a counting bound
$O(E)$. The Riemann--von Mangoldt formula instead requires order
$E\log E$ with multiplicity~\cite{Titchmarsh}. This compact repair
therefore also falls short of the required spectrum.

\head{Boundary maps}
For the Dirac expression the Green boundary data at a regular
endpoint are its two spinor components. The condition $g_2=i g_1$
in Sec.~\ref{sec:wall} cancels their flux. For the scalar square the
data are $u(0),u'(0)$ on its second-order maximal domain.
Evaluation of $u'(0)$ has no continuous trace on arbitrary
$H^1(0,1)$ data. Replacing one boundary space by the other would
therefore change the operator. A single chiral momentum
$-i\partial_x$ on $C_c^\infty(0,\infty)$ has deficiency
indices $(1,0)$: the equations $P^*f=if,\ P^*f=-if$ give
$e^{-x},e^x$. A doubled momentum can be made self-adjoint by
gluing opposite chiralities at the endpoint, which unfolds it into
$-i\partial_x$ on the line. Its spectral inverse is multiplication
by $1/k$ on
$\{f:\int|\widehat f(k)|^2k^{-2}\dd k<\infty\}$.
These realizations specify the domains and boundary data associated
with each inverse momentum \cite{reed,Weidmann,Gorbachuk1991}.

\section{Dirichlet beta at four specified frequencies}
\label{sec:beta-toy}
Let $\chi_4(n)$ be the odd character modulo four, with values
$0,1,-1$ on even integers, $1\bmod4$, $3\bmod4$, respectively.
The original comparison function is
\[
\begin{aligned}
 \beta(s)&=\sum_{k\ge0}\frac{(-1)^k}{(2k+1)^s}\\
 &=4^{-s}\big[\zeta(s,1/4)-\zeta(s,3/4)\big]\quad(\re s>1).
\end{aligned}
\]
Its entire completion is
$(4/\pi)^{(s+1)/2}\Gamma((s+1)/2)\beta(s)$, invariant under
$s\mapsto1-s$, by the primitive odd-character functional
equation~\cite{mv}. The gamma argument reflects the parity.
Replacing $\zeta(2s)$ by $\beta(2s)$ in \eqref{MB} inserts
the weights $\chi_4(n)$. Appendix E of the earlier paper uses
the $1/(2\pi i)$ normalization, twice that of \eqref{MB}.
At its stated $g=3/4,\ a=1/5,\ E=2t$, the precise answer is
\begin{equation}
\begin{gathered}
 \psi_{2t}=2(2/5)^{1/2+it} B(t),\\
 B(t)=\sum_{n\ {\rm odd}}\chi_4(n)n^{-1/2-it}K_{1/2+it}(5n).
\end{gathered}
 \label{eq:betasum}
\end{equation}
Since $\chi_4(1)=1$, Proposition~\ref{th:nonvanish} applies to
the signed superposition with the stated factor of two. The four
frequency intervals at $a=1/5$ admit the following lower bounds.

\begin{center}
\begin{tabular}{cc}
\toprule
Reported center & certified lower bound on $|B(t)|$\\
$t_0$ & for $|t-t_0|\le10^{-8}$\\
\midrule
 $6.0209489047$& $9.5\,10^{-5}$\\
 $10.243770304$& $8.4\,10^{-8}$\\
 $12.988098012$& $1.0\,10^{-9}$\\
 $16.342607105$& $2.4\,10^{-12}$\\
\bottomrule
\end{tabular}
\end{center}
For completeness, sum the six odd terms $n\le11$ at each printed
center with complex ball arithmetic. The error from $n\ge13$
is at most
$\sqrt{\pi/10}\,e^{-65}/[13(1-e^{-10})]$, by the elementary
bound on $K$ in Proposition~\ref{th:nonvanish}. To enclose
the whole interval, shift the $K_\nu(5n)$ integral to
$u+i$ in its full-line representation. Exponential decay on
$|\im u|\le1$ justifies the shift and differentiation in $t$.
Writing $Y=5n\cos1>2.5n$, and using $|u|\le\cosh u$ for
real $u$, bounds the absolute derivative of each term by
\[
 n^{-1/2}e^{-t}\{(1+\log n)K_{1/2}(Y)
           +[K_{3/2}(Y)+K_{-1/2}(Y)]/2\}.
\]
Since $K_{3/2}(Y)=(1+1/Y)K_{1/2}(Y)$, summing gives
$|B'(t)|<0.16e^{-t}$ for $t>0$, using
$(2.2+\log n)/n\le2.2$ and $\sqrt{\pi/5}<0.8$.
The ancillary program \texttt{anc/beta\_check.py} encloses the
finite sums at 192-bit precision and subtracts both error bounds.
Thus the integral \eqref{eq:betasum} is nonzero on all four frequency
intervals. The conclusion concerns the integral at the specified
acceleration and is independent of the zero divisor of Dirichlet beta.

\section{Finite intervals and independent checks}\label{app:check}
\head{Base and local count from primes}
For an integer $n\ge1$ let
$\psi_n=\sum_{m\le n}\Lambda(m)$ and $a_n=\sum_{m\le n}\Lambda(m)/m$.
On $[n,n+1]$, $T(x)=xa_n-\psi_n$, including both endpoints.
Put $B(x)=x\log x-(1+\eul)x+\sqrt x$.
For $n\ge4$ the function $B$ is convex and a tangent gives
\[
 B(x)-T(x)\ge B(n)-na_n+\psi_n+\min(0,B'(n)-a_n).
\]
For $4\le n\le15$ the right side is greater than $0.17696$.
Consequently $T(x)\le B(x)$ throughout the initial interval $[4,16]$.

For the raw local count, each contribution to $W_4$ is affine in
$x$ between $m/16,m/4,m$, where $m=p^k$.
The slopes jump down at $m/4$. Consequently
$W_4(x)-12x\log4+c\sqrt x$, for any $c\ge0$, is concave between
successive outer breakpoints $m/16,m$, including all intervening
middle breakpoints. Its minimum is attained at an endpoint.
An exhaustive integer-prefix calculation gives
\[
 F_4(x)<0.181584\qquad (4\le x\le65\,536).
\]
The largest value on this interval lies between $0.1815838$
and $0.1815839$ at $x=12401/16$.
Indeed checking $F_4\le c$ at the outer breakpoints suffices for
every $c\ge0$, and the largest endpoint value is positive.
There are $82\,267$ prime powers up to $16\cdot65\,536$ and
$88\,859$ distinct outer breakpoints and endpoints in this test.
No zero data enter it. It bounds the local count on this interval. The identity
\eqref{eq:greens} still requires a tail to bound a forward increment.

The same prime data can execute the forward induction directly up
to $1\,048\,576$. Define
\[
\begin{gathered}
 H(x)=\tfrac14\sqrt x+3x\log x+
          [\,4\log4-3(1+\eul)\,]x,\\
 B_4(x)=T(4x)-T(x).
\end{gathered}
\]
Then $H(x)-B_4(x)=\sqrt{x}/4-[R(4x)-R(x)]$.
On $[j,j+1]$, $j\ge4$ an integer, $B_4$ is convex because
$T(x)$ is affine there and $T(4x)$ is convex.
Also $0<H''(x)\le3/j$. For any twice differentiable function
with this upper bound, its chord on the unit interval exceeds its
value by at most $3/(8j)$, by the quadratic interpolation remainder.
The integer prefixes enclose $H(j)-B_4(j)>3/32$ at every integer
$4\le j\le262\,144$. The chord comparison gives $g_4(x)<1/4$
throughout that real interval. Starting from $[4,16]$, forward
induction gives $R(t)\le\sqrt t$ for $4\le t\le1\,048\,576$.
The bound on the forward increment uses only prime powers and
elementary constants.

\head{An analytic interval without tabulated zeros}
For $x\ge1$, $|c_\rho|/\re c_\rho
\le\sqrt{1+\gamma^{-2}}\le\sqrt{17}/4$ since the ratio of imaginary
to real parts has magnitude
$|\gamma(1-2\beta)|/(\gamma^2+\beta(1-\beta))\le1/|\gamma|$.
By pairing $\beta$ with $1-\beta$ and using
$x^\beta+x^{1-\beta}\le x+1$, \eqref{eq:explicit} gives
\[
 1-\frac{R(x)}{\sqrt x}
 \ge 1-\frac{d_0}{y}-\frac{\sqrt{17}}8 Q(y+y^{-1}),
 \qquad y=\sqrt x.
\]
The subtracted expression is convex in $y$. Evaluation at $y=2,40$
gives $R(x)<(1-0.0011)\sqrt x$ on $[4,1600]$.

The same calculation controls a forward step using only constants.
Apply the bound for the zero coefficients separately at $4x$ and
at $x$. For $4\le x\le64$ it gives
\[
 g_4(x)\le\frac{\sqrt{17}}8 Q(5\sqrt x+2/\sqrt x)+1/180<0.975 .
\]
Here $S(x)-S(4x)\le S(x)$ gives the last summand and the first
term increases for $x\ge4$. This supplies, together with the
analytic base above, an induction up to $256$ with no finite list
of primes or zeros.

\head{A finite forward induction using a published zero-height theorem}
Let $H=3\cdot10^{12}$. Platt and Trudgian~\cite{pt} proved that
all nontrivial zeros up to this height lie on the critical line.
Brent, Platt and Trudgian~\cite[Lemma 1 and its application to
$t^{-2}$]{bpt} give
\[
\begin{aligned}
 \sum_{\gamma>H}\gamma^{-2}
 \le e_H&:=\frac{1+\log(H/(2\pi))}{2\pi H}
             +\frac{0.14+0.56\log H}{H^2}\\
 &<1.48\cdot10^{-12}.
\end{aligned}
\]
The sum is over positive ordinates with multiplicity. Any half weight
at the lower endpoint only increases this upper bound for the strict
tail. Since $|c_\rho|\le\gamma^{-2}$, reflection in the line
and conjugation bound the total absolute contribution of unverified
zeros in \eqref{eq:explicit} by $(x+1)e_H$. The sum of the weights
of verified zeros is at most $Q$, since all $\re c_\rho>0$.
Thus
\[
 R(x)\le d_0+Q\sqrt x+(x+1)e_H .
\]
Dividing by $\sqrt x$ and testing the endpoints by convexity gives
$1-R(x)/\sqrt x>0.0179$ on $[4,4\cdot10^{23}]$.
There is also a genuine finite forward step. The low modes satisfy
$|4^\rho-1|\le3$, the high modes satisfy $|4^\rho-1|\le5$, and
$S(x)/\sqrt x\le[6\sqrt x(x^2-1)]^{-1}$. Therefore
\[
\begin{gathered}
 g_4(x)\le3Q+5e_H(\sqrt x+1/\sqrt x)+1/180<0.886\\
 (4\le x\le10^{22}).
\end{gathered}
\]
The bound is strictly below $\sqrt4-1=1$. Forward induction from
the base interval $[4,16]$ therefore gives $R(t)\le\sqrt t$ for
every real $4\le t\le4\cdot10^{22}$. The direct estimate on $R$
extends the finite range to $4\cdot10^{23}$.

\head{Reproducibility}
The ancillary source \texttt{anc/verify.py} uses the Python 3 standard
library. It encloses logarithms and square roots by outward rounding,
bounds $\pi$ by Machin's identity with alternating rational
remainders, and bounds Euler's constant using
\[
 \frac1{2(N+1)}<\sum_{k\le N}\frac1k-\log N-\eul<\frac1{2N}.
\]
For example the last two inequalities follow by summing the
integral remainders
$\log(1+1/k)-1/(k+1)=\int_0^1
 (1-t)/[(k+t)(k+1)]\,\dd t$ for $k\ge N$ and comparing
the end denominators with their adjacent values, or from the
Euler--Maclaurin bounds for the harmonic sum. The program rounds
$10^{24}\log p$ outwards to integers before accumulation, and rounds
each quotient by $m$ outwards again. Equations \eqref{eq:C},
\eqref{eq:Aq} and the tangent tests then use interval arithmetic. These finite checks are
independent of the proofs of Theorems~\ref{th:sharp} and~\ref{th:div}
apart from the displayed arithmetic enclosures for the two scale ratios.

\head{The correlation, the modes and the wall}
The program \texttt{anc/majorana\_rindler\_check.py} uses ball
arithmetic. It encloses $Q$, $e_H$ and the entries of
Table~\ref{tab:atoms} from certified zeros. It verifies at twelve
times $0.05\le u\le7.5$ that $C(u)$ from the primes and the sum
over the first $10^4$ zeros differ by less than the remaining mass
plus $2e_H\cosh(u/2)$. It checks the expansion \eqref{eq:cusp} and
the decomposition $C=B+Y$ of Sec.~\ref{sec:masspropagation}. For
the modes $G_\omega$ of Theorem~\ref{th:inversebessel} it verifies
the eigenvalue equation, the horizon amplitude and the conjugation
rule at $\omega=0.9$, $2\gamma_1$ and $-2\gamma_2$, and the
inverse-momentum integral of Sec.~\ref{sec:inversemajorana} at the
first two. For the wall of Theorem~\ref{th:wall} it computes the
Pr\"ufer angle of the Dirichlet problem for $U=V^2-V$. At
$\mu_0=\tfrac14,1,2\pi$ and $\omega=10,20,40,80,160$ the counts
are $11,26,61,139,313$; $6,17,43,104,243$ and $1,6,20,57,149$.
Each lies within one level of
$(\omega/\pi)[\log(2\omega/\mu_0)-1]+\tfrac12$.

\section{Massless trace stability and projection gains}\label{app:other}
\subsection{A zero-mass stability condition at one energy}
Adopt the trace normalization of Ref.~\cite{notes} at $s_0=1/2$,
with mass $0<m<1$ in fixed reference units. The criterion assumes
that a single renormalized trace satisfies both the spectral
representation and the remainder estimate
\[
\begin{gathered}
\begin{aligned}
 I_m^{\rm ren}(s_0)&=\frac{\zeta'}{\zeta}(s_0)\\
 &\quad+\sum_\rho\frac{\Gamma(\rho)\Gamma(\rho-s_0)}{\Gamma(s_0)}
                  m^{\,2(s_0-\rho)}+R_m(s_0),
\end{aligned}\\
 R_m(s_0)=O(1).
\end{gathered}
\]
If $I_m^{\rm ren}(s_0)=O((1+|\log m|)^k)$ for some fixed $k\ge0$,
then RH follows. Indeed put $v=-2\log m,\ z_\rho=\rho-1/2$,
$a_\rho=\Gamma(\rho)\Gamma(\rho-1/2)/\sqrt\pi$.
Stirling and the zero count give $\sum|a_\rho|<\infty$.
Since $\zeta(1/2)\ne0$, $f(v)=\sum a_\rho e^{z_\rho v}$ grows
at most polynomially under the two hypotheses. Its Laplace transform
is then holomorphic for $\re w>0$. On $\re w>1/2$, termwise
integration gives $\sum_\rho a_\rho/(w-z_\rho)$, meromorphic by
normal convergence on compact sets off its poles. A zero of
multiplicity $q$ gives residue $q a_\rho\ne0$. Uniqueness of
continuation excludes a pole in $\re w>0$. Reflection gives RH.

To use the criterion one needs a relative trace $T_m^{\rm rel}(0+i0)$
whose boundary value exists at this specific energy, together with
\[
\begin{aligned}
 I_m^{\rm ren}(1/2)&=\frac{m}{\sqrt\pi}
       \left[T_m^{\rm rel}(0+i0)+\frac{\sqrt\pi}{m^2}\right],\\
 T_m^{\rm rel}(0+i0)&=-\frac{\sqrt\pi}{m^2}
       +O\big(m^{-1}(1+|\log m|)^k\big).
\end{aligned}
\]
The criterion assumes both the spectral representation and independent
cancellation in this formula. It concerns the limit $m\downarrow0$.
Section~\ref{sec:heat} instead evaluates the radial intensity at
fixed boost frequency $\tau=0$.

\subsection{Projection gain and parity cancellation}
In $L^2(0,1)$ put $h_n(x)=\{1/(nx)\}-n^{-1}\{1/x\}$,
$V_N=\operatorname{span}(h_2,\ldots,h_N)$, and denote projection
onto $V_N$ by $P_N$. With $r_N=1-P_N1$ let
$e_N=\|r_N\|^2$. On $1/(k+1)<x<1/k$,
$h_n(x)=\{k/n\}$. Thus their Gram matrix and right hand side are
\[
\begin{gathered}
 G_{mn}=\sum_{k\ge1}\frac{\{k/m\}\{k/n\}}{k(k+1)},\\
 v_n=\int_0^1 h_n(x)\dd x=\frac{\log n}{n}.
\end{gathered}
\]
In signal terms $P_N$ fits a probe profile by the first $N$
arithmetic templates in least squares. The complementary projection
$I-P_N$ leaves the unresolved part. The Gram matrix $S$ measures
the overlaps of the additional templates after this subtraction;
the corresponding gain is the reduction in residual squared norm.
These projections act on arithmetic profiles.

For the last equality put $u=1/x$, change scale in the first
fractional part and send the common finite cutoff to infinity;
the residual integral over $1/n\le u\le1$ is $(\log n)/n$.
The $h_n$ are linearly independent. If a finite combination
vanishes on all cells, first differences in $k$ give
$\sum c_n/n=\sum_{n\mid k}c_n$ for $n\ge2$; use $k=1$,
then induction in $k$ (the cell values at $k=0,1$ coincide
with their periodic extension). Nor can their span contain 1,
since all cell values vanish at common multiples. Therefore
$e_N=1-v^TG^{-1}v>0$.
The discrete Nyman--Beurling criterion gives $e_N\to0$ iff RH
\cite{baez}. In passing from $L^2(0,\infty)$, the tail of any
finite combination of the unadjusted fractional parts on $x>1$
is $a/x$. Vanishing norm of the error makes $a\to0$;
subtracting $a\{1/x\}$ gives the present normalization. Conversely
these normalized combinations vanish on $x>1$.

For $N<M$ set $u_n=(I-P_N)h_n,\ b_n=\langle r_N,h_n\rangle$
and $S_{mn}=\langle u_m,u_n\rangle$ on $N<m,n\le M$.
Orthogonal projection gives $e_N-e_M=b^T S^{-1}b$.
For any increasing subsequence $N_j\to\infty$ put
$g_j=(e_{N_j}-e_{N_{j+1}})/e_{N_j}$.
The exact product $e_{N_k}/e_{N_0}=\prod_{j<k}(1-g_j)$,
with $0\le g_j<1$, gives $e_N\to0$ if and only if $\sum_jg_j=\infty$.
Gram positivity gives nonnegative gains. Convergence of the projection
scheme requires divergence of their normalized sum. This divergence
and the uniform signed trace estimate in \eqref{eq:parity} are
alternative arithmetic conditions for RH.

\section{Spectral moments and finite arithmetic reconstructions}\label{app:spectralmoments}
\subsection{A prime-data moment test for spectral completeness}
\label{sec:operatormoments}
Completeness of the boundary measure has a moment criterion at
$s=3/2$, in the half-plane of absolute convergence of the Euler
product. For real $y$ near 1, set
\[
\begin{gathered}
 F(y)=\frac{L(\sqrt y)}{\sqrt y}
 =\frac{Q-2(\dd/\dd y)\log\xi(a+\sqrt y)}{y-a^2},\\
 m_0=Q,\qquad
 m_{j+1}=\frac{(-1)^j}{j!}F^{(j)}(1)\quad(j\ge0).
\end{gathered}
\]
Every $m_j$ uses elementary and gamma factors and the derivatives
of $\zeta'/\zeta$ at $3/2$. Here its Dirichlet series and all its
termwise derivatives converge absolutely. For integers $n,k\ge0$ write
\[
 D_{n,k}=\sum_{j=0}^k(-1)^j\binom{k}{j}m_{n+j}.
\]

\begin{theorem}[Moment construction]\label{th:operatormoments}
The inequalities $D_{n,k}\ge0$ for all $n,k\ge0$ are equivalent
to $\Delta=0$. When they hold, the moments specify a nonnegative
self-adjoint generator with correlation $C$, without using zero
ordinates in its construction.
\end{theorem}
\begin{proof}
For a generator with the required correlation and vector measure
$\nu$, Laplace transformation gives $F(y)=\int(y+t^2)^{-1}
\dd\nu(t)$ for $y>a^2$. Pushing $\nu$ forward by
$b=(1+t^2)^{-1}$ gives $m_j=\int b^j\dd\mu(b)$. Hence
$D_{n,k}=\int b^n(1-b)^k\dd\mu(b)\ge0$.

Conversely, the Hausdorff moment theorem \cite{liupego} gives a
unique positive measure $\mu$ on $[0,1]$ with these moments and
mass $Q$. Expanding at $y=1$ yields
\[
 F(y)=\int_{[0,1]}\frac{b}{1+(y-1)b}\dd\mu(b).
\]
Equality of the Taylor series gives equality on $y>a^2$ by real
analytic continuation. As $y\to\infty$, the displayed formula for
$F$ and Stirling's estimate give $yF(y)\to Q$.
On the integral side dominated convergence gives
$\mu((0,1])=Q$. On $L^2(\mu)$ take multiplication by
$\sqrt{b^{-1}-1}$ on its maximal domain and $f=1$.
Its correlation and $C$ now have equal Laplace transforms for real
$z>a$. Both functions are continuous and of exponential order, so
uniqueness and evenness give equality at all real times.
Theorems~\ref{th:CAR} and~\ref{th:boundaryclosure} finish the proof.
\end{proof}

For scale, evaluation of the real Taylor series at $s=3/2$ gives
$7.39864956<10^5m_1<7.39864958$ and
$7.37011086<10^5D_{1,1}<7.37011089$.
The ancillary ball-arithmetic calculation checks these intervals and
the direct prime-transform formula above at two real arguments.

\begin{proposition}[Unconditional ranges of moment positivity]\label{pr:operatorranges}
With $H=3\cdot10^{12}$, one has $D_{n,k}>0$ whenever
\[
 \frac{n+1}{H}+\frac{k}{H^3}\le\frac{\pi}{2},
 \qquad\text{or}\qquad k\le \frac{H^2}{32}\,n .
\]
For each fixed $n\ge0$ one also has the unconditional asymptotic
\[
 D_{n,k}\sim\frac{\Gamma(n+1/2)}{4\pi}
              \frac{\log k}{k^{n+1/2}}\qquad(k\to\infty).
\]
\end{proposition}
\begin{proof}
Put $b_\rho=(1-w_\rho^2)^{-1}$. Pairing $w_\rho$ with
$-w_\rho$ in $F(y)$, differentiating at 1, and using absolute
convergence gives
\[
 m_j=\sum_\rho c_\rho b_\rho^j,\qquad
 D_{n,k}=\sum_\rho c_\rho b_\rho^n(1-b_\rho)^k .
\]
The case $j=0$ is \eqref{eq:explicit}. For zeros on the line the
summands are positive. For the remaining zeros write $d=\beta-a$,
$g=|\gamma|>H$, $A=1+g^2-d^2$, $B=2d\gamma$. Then
$b_\rho=(A-iB)^{-1}$ and $c_\rho=(a^2+g^2-d^2-iB)^{-1}$.
The arguments of $c_\rho,b_\rho$ have absolute value at most
$1/g$. The tangent of the absolute argument of $1-b_\rho$ is
$|B|/[A(A-1)+B^2]\le g^{-3}$, since $A(A-1)\ge g^4$ for
$g\ge1$. This proves the first assertion in the stated angular range.

There is a verified conjugate pair at an ordinate $t\in[H/3,H/2]$.
For example, the zero-count bound
$\lvert N(T)-(T/(2\pi))\log(T/(2\pi e))-7/8\rvert\le0.28\log T$
for $T\ge2\pi e$ \cite{bpt} gives a positive count difference
between these endpoints, and \cite{pt} places the zeros on the line.
For that pair, $3/H^2\le b=(1+t^2)^{-1}\le9/H^2$ and
$2/(a^2+t^2)\ge6/H^2$. Since $\log(1-b)\ge-2b$, in the second
range the pair contributes at least
$(6/H^2)(3/2)^n H^{-2n}$, using $3e^{-9/16}>3/2$.
All off-line terms together have absolute value at most
$2e_H H^{-2n}$, since $|b_\rho|\le H^{-2}$ and
$|1-b_\rho|<1$. For $n\ge100$ the pair already exceeds this
upper bound. For $n<100$ the second range lies inside the first.

For the asymptotic let $R=\sqrt{k}$ and fix $n$. Uniformly in
$|d|<1/2$ one has $c_\rho=g^{-2}(1+O(g^{-1}))$ and
$b_\rho=g^{-2}(1+O(g^{-1}))$ as $g\to\infty$. For
$v=g/R$ in any compact subinterval of $(0,\infty)$ it follows that
\[
 R^{2n+2}c_\rho b_\rho^n(1-b_\rho)^k
       =v^{-2n-2}e^{-1/v^2}+o(1)
\]
uniformly. By the Riemann--von Mangoldt formula, the counting
measures of $v$, including both signs of $\gamma$ and divided
by $R\log R$, converge on compact subintervals to $\dd v/\pi$.
For all zero ordinates $g>4$, as proved after Theorem~\ref{th:sharp},
the identity $|1-b_\rho|^2=1-(2A-1)/(A^2+B^2)$ bounds the scaled
summands by $O_n(v^{-2n-2}e^{-c/v^2})$ for an absolute $c>0$.
The normalized zero count on $[v,2v]$ is
$O(v[1+\log^+(2v)])$ for $R\ge2$, by $N(T)\ll T\log(2+T)$,
where $\log^+x=\max(0,\log x)$. Summing these bounds over dyadic
intervals controls both tails uniformly in $R$. Thus
\[
\begin{aligned}
 \lim_{k\to\infty}\frac{R^{2n+1}D_{n,k}}{\log R}
  &=\frac1\pi\int_0^\infty v^{-2n-2}e^{-1/v^2}\dd v\\
  &=\frac{\Gamma(n+1/2)}{2\pi}.
\end{aligned}
\]
Since $\log R=(\log k)/2$, this is the asserted coefficient.
\end{proof}

The signs also satisfy an exact induction. Pascal's identity gives
$D_{n,k}=D_{n+1,k}+D_{n,k+1}$, hence for every $r\ge0$,
\[
 D_{n,k}=\sum_{j=0}^r\binom rj D_{n+j,k+r-j}.
\]
Nonnegativity on any unbounded sequence of complete diagonals therefore
suffices for Theorem~\ref{th:operatormoments}. Conversely, each negative
entry has a negative child. By Proposition~\ref{pr:operatorranges}
any path of such children must have both indices tending to infinity
and $k>H^2 n/32$ throughout. On a diagonal $n+k=N$, only
$n<32N/(H^2+32)$ can require a further sign estimate. The asymptotic
for a fixed $n$ supplies no uniform threshold on this growing set.

\subsection{The exponential rate of the variation of the prime moments}
\label{sec:massvariation}
To test the equal-time hypothesis itself, define
$\mathcal V_N=\sum_{j=0}^N|a_{j,N}|$. Its signed counterpart is
$\sum_j a_{j,N}=m_0=Q$, by the binomial theorem applied to
$D_{n,k}=\sum_\rho c_\rho b_\rho^n(1-b_\rho)^k$, where
$b_\rho=(1-w_\rho^2)^{-1}$ and $w_\rho=\rho-1/2$.

\begin{theorem}[Growth rate from the arithmetic moments]\label{th:massvariation}
The sequence $\mathcal V_N$ is nondecreasing from $\mathcal V_0=Q$.
Its normalized exponential rate exists and satisfies
\[
 \lim_{N\to\infty}(\mathcal V_N/Q)^{1/N}
 =R_*:=\sup_\rho
 \left(1-\frac{4(\re w_\rho)^2}{(1+|w_\rho|^2)^2}\right)^{-1/2}.
\]
Hence $\log(\mathcal V_N/Q)=o(N)$ is equivalent to $\Delta=0$.
It suffices to prove the sublinear logarithmic bound along any
unbounded sequence of indices.
\end{theorem}
\begin{proof}
For $0\le j\le N$ put $x_j=D_{j+1,N-j}$ and
$y_j=D_{j,N+1-j}$. Pascal's identity gives the exact increment
\[
 \mathcal V_{N+1}-\mathcal V_N
 =\sum_{j=0}^N\binom Nj\big(|x_j|+|y_j|-|x_j+y_j|\big)\ge0.
\]
Write $R_\rho=|1-b_\rho|+|b_\rho|$. Then
$R_\rho=(1+|w_\rho|^2)/|1-w_\rho^2|$. Squaring numerator and
denominator gives $\sup_\rho R_\rho=R_*$.
In particular $R_\rho\to1$ at large height.
With $A_0=\sum_\rho|c_\rho|<\infty$, the triangle inequality gives
$Q\le\mathcal V_N\le A_0R_*^N$.

For the reverse bound consider the trigonometric polynomial
\[
\begin{gathered}
 P_N(\theta)=\sum_{j=0}^N a_{j,N}e^{ij\theta}
    =\sum_\rho c_\rho\,q_{b_\rho}(\theta)^N,\\
 q_b(\theta)=1-b+be^{i\theta}.
\end{gathered}
\]
It satisfies $|P_N(\theta)|\le\mathcal V_N$.
If $1<r<R_*$, some $b=b_\rho$ has
$\max_\theta|q_b(\theta)|=R_\rho>r$. Hence an open arc has
$|q_b|>r$. Only finitely many distinct $b_\rho$ can exceed $r$
anywhere, since $b_\rho\to0$. For two distinct such values the
functions $|q_b|^2$ cannot agree identically. Equality of their
first Fourier coefficients would give $b-|b|^2=b'-|b'|^2$.
Their imaginary parts would coincide and their real parts would
either coincide or sum to one. The latter is impossible because
$|b_\rho|\le(1+(\im\rho)^2-(\re w_\rho)^2)^{-1}<1/2$.
The functions are trigonometric polynomials of degree one, so each
pair has only finitely many ties. Choose an angle in the arc avoiding
all of them. One distinct $b$ then has unique largest modulus
$M>r$, with all others of modulus at most $m<M$. Equal values of
$b_\rho$ have equal $w_\rho^2$ and equal nonzero $c_\rho$.
Their combined coefficient $c_*\ne0$. Consequently
\[
 \mathcal V_N\ge |c_*|M^N-A_0m^N .
\]
This gives $\liminf\mathcal V_N^{1/N}\ge r$ for every $1<r<R_*$.
If $R_*=1$, the two preceding bounds $Q$ and $A_0$ suffice.
Finally $R_\rho>1$ exactly when $\re(\rho)\ne1/2$, proving the
equivalence by Theorem~\ref{th:boundaryclosure}.
\end{proof}

For the verified height $H=3\cdot10^{12}$ \cite{pt}, one has
$\mathcal V_N=Q$ for every integer $0\le N\le10^{36}$.
Indeed, if the second range in Proposition~\ref{pr:operatorranges}
does not apply to $j+k=N$, then $j<32N/(H^2+32)$. It follows that
\[
 \frac{j+1}{H}+\frac{k}{H^3}\le\frac1H+\frac{33N}{H^3}
 <\frac{\pi}{2}
\]
in the stated range, so all its entries are positive. For arbitrary
$N$ the same height bounds the rate by
\[
 0\le\log R_*\le
 -\tfrac12\log\!\left(1-(1+H^2)^{-2}\right)<6.173\cdot10^{-51}.
\]
Define the relative increment in variation by
$\ell_N=(\mathcal V_{N+1}-\mathcal V_N)/\mathcal V_N$.
Telescoping and Theorem~\ref{th:massvariation} give
\[
 \lim_{N\to\infty}\frac1N\sum_{k=0}^{N-1}\log(1+\ell_k)=\log R_* .
\]
An upper bound tending to zero along an unbounded sequence of indices
forces $\log R_*=0$ and removes the mass defect. The positive
finite-height upper bound on the rate permits a nonzero defect.

\subsection{Uniform bounds in the degree from the Euler product}
\label{sec:sparseprime}
Write $a_{j,N}=\binom Nj D_{j,N-j}$ as in Sec.~\ref{sec:masspropagation}.
The following estimates use the Euler product only in its half-plane of
absolute convergence. No zero ordinates enter their proof. The index
$j$ may increase with the total degree $N$.

\begin{theorem}[Prime moment estimates at unbounded degrees]\label{th:sparseprime}
Let $W_{j,N}=\binom Nj {\rm B}(j+\tfrac12,N-j+\tfrac12)$, where
${\rm B}$ is Euler's beta function. For all sufficiently large $N$
one has $a_{j,N}>0$ throughout
\[
 0\le j\le J_N=\left\lfloor\tfrac13\sqrt{N\log\log N}\right\rfloor .
\]
Uniformly in $0\le j\le\sqrt N$,
\[
\begin{gathered}
 a_{j,N}=\frac{W_{j,N}}{4\pi}
 \left(\log\frac{N+1}{j+1}+O(1)\right),\\
 \sum_{j\le\sqrt N}|a_{j,N}|
 =\frac{\log N}{4\pi N^{1/4}}+O(N^{-1/4}).
\end{gathered}
\]
For any sequence of nonnegative integers $K_N=o(N)$,
\[
 \log\!\left(1+\sum_{j=0}^{K_N}|a_{j,N}|\right)=o(N).
\]
\end{theorem}
\begin{proof}
We first subtract a positive measure whose transform accounts for the
logarithm in Stirling's formula. Set
\[
\begin{gathered}
 h(t)=\frac{\log(1+t^2)}{2\pi(1+t^2)}\quad(t>0),\\
 q=\int_0^\infty h(t)\dd t=\tfrac12\log2 .
\end{gathered}
\]
Its transform in $\re w>0$ is
\[
 L_h(w)=\int_0^\infty\frac{w\,h(t)}{w^2+t^2}\dd t
       =\frac{w\log2-\log(1+w)}{2(w^2-1)} .
\]
For real $a,w>0$, differentiating in $a$ and using the value at
$a=0$ gives
$\int_0^\infty\log(t^2+a^2)(t^2+w^2)^{-1}\dd t
 =(\pi/w)\log(w+a)$. The value at zero follows by $t=wv$ and
$v\mapsto1/v$. Partial fractions and continuation give $L_h$,
with its removable value at $w=1$.

Take $\sigma=3/5$, $s=1/2+w$, and write
\[
 G(w)=L(w)-L_h(w)-(Q-q)/w,\qquad b(w)=(1-w^2)^{-1}.
\]
For $\re w\ge\sigma$, the series
$-\zeta'(s)/\zeta(s)=\sum_{n\ge2}\Lambda(n)n^{-s}$
is bounded absolutely by its value at $s=11/10$.
Explicitly $L(w)=(wQ-B(s))/(w^2-1/4)$, where
\[
 B(s)=\frac1s+\frac1{s-1}-\tfrac12\log\pi
       +\tfrac12(\Gamma'/\Gamma)(s/2)
       -\sum_{n\ge2}\frac{\Lambda(n)}{n^s}.
\]
Logarithmic differentiation of $\xi$ and Stirling's formula give
$\xi'(s)/\xi(s)=\tfrac12\log w+O(1)$ uniformly there at infinity.
Consequently $G$ is holomorphic on this half-plane, including at $w=1$,
and $|G(w)|\ll(1+|w|^2)^{-1}$.

Let $a^h_{j,N}$ denote the integral of
$\binom Nj x^j(1-x)^{N-j}$ against $h(t)\dd t$, with $x=(1+t^2)^{-1}$.
For $k=N-j$ the exact contour formula is
\[
\begin{aligned}
 a_{j,N}&=a^h_{j,N}+(Q-q)\one_{\{j=N\}}\\
  &\quad+\frac{\binom Nj}{\pi}\int_{\R}
        G(\sigma+it)b(\sigma+it)^j[1-b(\sigma+it)]^{k}\dd t .
\end{aligned}
\]
To verify the normalization, first integrate $G$ alone on this line.
Closing to the right gives zero. For $|v|$ sufficiently small,
closing in the same direction gives
\[
 \frac1\pi\int_{\R}
  \frac{G(\sigma+it)}{1-v-(\sigma+it)^2}\dd t
       =\frac{G(\sqrt{1-v})}{\sqrt{1-v}} .
\]
The pole at $\sqrt{1-v}$ has residue
$-G(\sqrt{1-v})/(2\sqrt{1-v})$, and the contour is clockwise.
The arc vanishes by the bound on $G$. On the line $|b|\le25/16$,
so expansion at $v=0$ and integration commute. The coefficient of
$v^{r-1}$, $r\ge1$, gives the $r$th arithmetic moment less the
moment of $h$ and the constant $Q-q$, by
Sec.~\ref{sec:operatormoments}. Taking finite differences proves the
contour formula. In particular its integrand uses prime data on
$\re s=11/10$.

The change of variables $x=(1+t^2)^{-1}$ sends $h(t)\dd t$ to
$\log(1/x)\dd x/(4\pi\sqrt{x(1-x)})$. Differentiating the beta
integral, with $\psi_0=\Gamma'/\Gamma$, gives the exact positive term
\[
 a^h_{j,N}=\frac{W_{j,N}}{4\pi}
             [\psi_0(N+1)-\psi_0(j+\tfrac12)].
\]
The bracket differs from $\log((N+1)/(j+1))$ by a bounded quantity,
uniformly in $0\le j\le N$.

For the integral error set $w=\sigma+it$, $p=(1+|w|^2)^{-1}$
and $R=|b|+|1-b|$. Direct algebra gives
\[
 R=(1-4\sigma^2p^2)^{-1/2},\qquad
 |b|=Rp,\qquad |1-b|=R(1-p).
\]
Thus the absolute value of the binomial kernel is $R^N$ times a
binomial probability with parameter $p$. Let $J$ be either
$\lfloor\sqrt N\rfloor$ or $J_N$. In both cases $2J/N<1/10$
eventually. For $p\le1/10$ one has
$\log R\le\min(p^2,p/12)$, by
$-\log(1-v)\le v/(1-v)$. On $p\le2J/N$ the change from $t$
to $p$ on each half-axis yields
\[
 |G(\sigma+it)|\,|\dd t|\ll
     \frac{\dd p}{\sqrt p\,\sqrt{1-(1+\sigma^2)p}}
     \ll\frac{\dd p}{\sqrt{p(1-p)}} .
\]
The integral error from this range is at most
$O(W_{j,N}\exp(4J^2/N))$.
On $2J/N<p\le1/10$, for $0\le j\le J$ a binomial variable
$X$ of mean $Np$ satisfies
\[
 \Pr(X=j)\le\Pr(X\le Np/2)\le e^{-Np/8}.
\]
This last inequality follows by applying Markov's inequality to
$2^{-X}$. Multiplication by $R^N$ leaves an upper bound
$e^{-Np/24}\le e^{-J/12}$.
Finally $p\ge1/10$ is a compact $t$-interval. There
$|1-b|\le\eta<1$ since
$|1-w^2|^2-|w|^4=1-2\sigma^2+2t^2>0$.
The kernel for $j\le J$ is bounded by
$\eta^{N-J}(25N/16)^J$, exponentially small in $N$.
These estimates use $\int_\R|G(\sigma+it)|\dd t<\infty$.

By gamma ratios,
\[
\begin{aligned}
 W_{j,N}&=
 \frac{\Gamma(j+1/2)\Gamma(N-j+1/2)}
      {\Gamma(j+1)\Gamma(N-j+1)}\\
 &\asymp ((j+1)(N-j+1))^{-1/2}.
\end{aligned}
\]
The exponentially small errors are therefore $o(W_{j,N})$ uniformly
for $j\le J$. For $J=\lfloor\sqrt N\rfloor$ this proves the
coefficient estimate. For $J=J_N$, the error divided by $W_{j,N}$
is $O((\log N)^{4/9})$, whereas
$\log((N+1)/(j+1))\ge(1/2-o(1))\log N$. This proves strict positivity.
To sum up to $J=\lfloor\sqrt N\rfloor$, use for $1\le j\le J$
\[
 W_{j,N}=(Nj)^{-1/2}
       [1+O(j^{-1}+(j+1)/N)] .
\]
The integral test gives $\sum_{j\le J,\,j\ge1}j^{-1/2}
 =2\sqrt J+O(1)$ and
$\sum_{j\le J,\,j\ge1}j^{-1/2}\log(J/j)=O(\sqrt J)$.
The term at $j=0$ is $O(N^{-1/2}\log N)$. This gives the stated
constant for the sum.

For general $1\le K<N/2$, sum the absolute contour bound using
$\sum_{j=0}^N a^h_{j,N}=q$, $|1-b|\le1$ and $|b|\le25/16$.
The resulting bound is
\[
 \sum_{j=0}^K|a_{j,N}|
 \le q+C(K+1)\left(\frac{25 eN}{16K}\right)^K .
\]
Here $C=\pi^{-1}\int_\R|G(\sigma+it)|\dd t$, and
$(eN/j)^j$ increases for $1\le j\le N$.
The case $K=0$ is bounded directly. If $K/N\to0$, division of
the logarithm by $N$ gives zero.
\end{proof}

\head{Intermediate indices in the variation of the prime moments}
Set $\alpha_H=32/(H^2+32)$ with the verified height from
Proposition~\ref{pr:operatorranges}. That proposition gives
$a_{j,N}>0$ for $j\ge\alpha_H N$. Since $\sum_j a_{j,N}=Q$,
Theorem~\ref{th:sparseprime} gives the exact identity, for all
sufficiently large $N$,
\[
 \mathcal V_N=Q+
 2\sum_{J_N<j<\alpha_H N}(-a_{j,N})_+ .
\]
Here $x_+=\max(x,0)$. The terms up to
$K_N=\lfloor N/\log N\rfloor$ have subexponential growth.
Indeed, the logarithm of one plus their total absolute contribution is
$O(N\log\log N/\log N)$.

Taking absolute values inside the Euler integral has a quantifiable
cost in the remaining range. Fix any $0<\alpha<\alpha_H$ and
$K_N=o(N)$. There is a compact $t$-interval on which
$\alpha/4<p(t)<\alpha/2$ and $|G(\sigma+it)|$ is bounded below
by a positive constant. Indeed $G$ is analytic and not identically
zero, since along the positive real axis
$w^2G(w)\to\tfrac12\log(2\pi)$ by Stirling's formula and the Euler
product. On this interval $R\ge r_\alpha>1$, and
$\Pr(K_N<X<\alpha N)\to1$ uniformly by the binomial variance bound.
The sum of the integrals of the absolute kernels over
$K_N<j<\alpha N$, including the factor $|G|/\pi$, is consequently
at least $c_\alpha r_\alpha^N$ for all sufficiently large $N$.
This lower bound concerns the majorant. In the coefficients themselves
the integral retains its complex phase. An independent estimate for that cancellation
would yield $\log(\mathcal V_N/Q)=o(N)$ on unbounded degrees.

\begin{proposition}[Finite prime cutoffs and the order of the limits]\label{pr:finiteeuler}
In $B(s)$ truncate the sum at $n\le M$ for an integer $M\ge3$.
Denote the resulting transform, coefficients and full variation by
$L_M,a^{(M)}_{j,N},\mathcal V_N^{(M)}$, keeping $m_0=Q$.
At $M_N=\lceil e^{8N}\rceil$, with $\delta=4/5-\log(17/8)>0$,
\[
 \sum_{j=0}^N |a_{j,N}-a^{(M_N)}_{j,N}|
       \ll (N+1)e^{-\delta N}.
\]
For every fixed $M\ge3$ one has
$\lim_{N\to\infty}(\mathcal V_N^{(M)})^{1/N}=5/3$.
\end{proposition}
\begin{proof}
The contour formula remains valid for $G_M=L_M-L_h-(Q-q)/w$
by the same right half-plane calculation. On $\re w=3/5$,
the integral test and $\Lambda(n)\le\log n$ give
\[
 |G(w)-G_M(w)|\le
 \frac{M^{-1/10}(10\log M+100)}{|w-1/2|\,|w+1/2|}.
\]
The integral over $t$, divided by $\pi$, of the reciprocal denominator
is at most $10$, since $|w+1/2|\ge|w-1/2|$.
Summing the absolute kernels gives $R(t)^N\le (17/8)^N$.
Substitution of $M_N$ proves the estimate for the whole diagonal.
Each truncated coefficient is computable from finitely many prime
powers and the elementary and gamma factors at $s=3/2$.

For fixed $M$, any $1/2<\sigma<1/\sqrt2$ gives an integrable
$G_M(\sigma+it)$. The same estimate of the binomial kernels yields
\[
 \limsup_N(\mathcal V_N^{(M)})^{1/N}
 \le \frac{1+\sigma^2}{1-\sigma^2}.
\]
Let $\sigma\downarrow1/2$. To prove the matching lower bound,
write $P_N^{(M)}=\sum_{j=0}^N(-1)^j a_{j,N}^{(M)}$ and
$F_M(y)=L_M(\sqrt y)/\sqrt y$. The binomial theorem gives
$P_N^{(M)}=\sum_{r=0}^N\binom Nr(-2)^r m_r^{(M)}$.
Consequently its ordinary generating function near zero is
\[
 \sum_{N\ge0}P_N^{(M)}z^N
 =\frac{Q+\tau F_M(1-\tau)}{1-z},
 \qquad \tau=\frac{-2z}{1-z}.
\]
For $|z|<1$, $y=1-\tau=(1+z)/(1-z)$ has positive real part.
The principal square root has positive real part as well. The finite
Dirichlet sum is entire in $s$, so the term $1/(s-1)$ in $B_M$
gives a double pole
\[
 F_M(y)=-2(y-1/4)^{-2}+O_M((y-1/4)^{-1}).
\]
This is the only singularity of the generating function in the open
unit disk. It lies at $z=-3/5$. Since
$y-1/4=3(1+5z/3)/(4(1-z))$, its leading term is
$-\tfrac{64}{15}(1+5z/3)^{-2}$. Subtracting the two singular
terms leaves a function holomorphic in the disk. Cauchy's estimate
on any circle of radius between $3/5$ and $1$ now gives
\[
 P_N^{(M)}=(-5/3)^N[-\tfrac{64}{15}(N+1)+O_M(1)].
\]
As $\mathcal V_N^{(M)}\ge|P_N^{(M)}|$, this proves the limit.
For the full prime series, the same alternating observable tends to
zero. Indeed it is $\sum_\rho c_\rho(1-2b_\rho)^N$, where
$|1-2b_\rho|=|1+w_\rho^2|/|1-w_\rho^2|<1$ by
$\re(w_\rho^2)<0$. Absolute summability of $c_\rho$ gives the limit.
\end{proof}

At the same degree $N=64$, the ancillary interval calculation gives
\[
\begin{gathered}
 -4.14813\cdot10^{16}< P_{64}^{(64)}<-4.14812\cdot10^{16},\\
 0.03864452<\sum_{j=0}^{64}(-1)^j a_{j,64}<0.03864453 .
\end{gathered}
\]
Both enclosures use Taylor coefficients at $s=3/2$. The first truncates
the von Mangoldt series at $64$, and the second uses the full transform.

\section{Prime mean squares and correlations on square intervals}\label{app:squares}
\subsection{The pole at one in the finite prime transform}
\label{sec:correctedcutoff}
Set $a=1/2$, $E(x)=\psi(x)-x$, and introduce the nonnegative
quantity
\[
 \mathcal I(X)=\int_1^X\frac{E(x)^2}{x^2}\dd x,\qquad X\ge1 .
\]
Only prime powers up to $X$ enter $\mathcal I(X)$. This integral
also bounds the variation of a signed measure on a finite interval
of return times.

\begin{proposition}[A bounded variation model at each finite cutoff]
\label{pr:correctedcutoff}
For an integer $M\ge3$, define
\[
\begin{gathered}
 \lambda_M=\log M-\eul-\sum_{n\le M}\frac{\Lambda(n)}n,\\
 D_M(s)=\sum_{n\le M}\frac{\Lambda(n)}{n^s}
             +M^{1-s}\left(\frac1{s-1}+\lambda_M\right).
\end{gathered}
\]
Replace the von Mangoldt series in $L(w)$ by $D_M(a+w)$, keeping
$Q$ and the completion factors. The corrected correlation is the
cosine transform of a finite signed measure $\nu_M$, of total mass
$Q$. It equals $C(u)$ for $|u|\le\log M$. For all $M\ge3$
and all $N\ge0$, with absolute implied constants,
\[
 \begin{split}
 \mathcal V_N&\le\|\nu_M\|_{\rm TV}
      +O\big((17/8)^N M^{-1/10}(1+\log M)^2\big),\\
 \|\nu_M\|_{\rm TV}
     &\ll 1+|R(M)|/\sqrt M+\sqrt{\mathcal I(M)} .
 \end{split}
\]
Here $\|\cdot\|_{\rm TV}$ denotes total variation.
\end{proposition}
\begin{proof}
The last term in $D_M$ restores the residue $1$ and the constant
term $-\eul$ in the Laurent expansion of $-\zeta'/\zeta$ at 1.
In the prime-power measure it adds Lebesgue measure on $(M,\infty)$
and an atom $M\lambda_M$ at $M$. Integrating the kernel
$(x/y-1)\one_{\{y\le x\}}$ against this measure gives the original
$T(x)$ for $x\le M$. For $x\ge M$ the continuous tail contributes
$x\log(x/M)-x+M$ and the atom contributes $\lambda_M(x-M)$.
The centered result is therefore exactly $R(\min(x,M))$. Hence the
corrected time response is the continuous integrable function
\[
 C_M(u)=e^{-|u|/2}
   [R(\min(e^{|u|},M))-d_0+S(e^{|u|})].
\]
Denote its one-sided Laplace transform by $\widehat L_M$.
For $\re w>a$ it equals the indicated replacement in $L$ by the
same absolutely convergent integration as in Sec.~\ref{sec:boundaryclosure}.
The formula continues to $\re w\ge0$, with a removable singularity
at $w=a$. At $s=a+it$ it gives
\[
\begin{gathered}
 \dd\nu_M(t)=\frac{2}{\pi}
 \frac{\tfrac12\re(\Gamma'/\Gamma)(s/2)-\tfrac12\log\pi
                         -\re D_M(s)}{t^2+1/4}\dd t,\\
 t>0 .
\end{gathered}
\]
Indeed $\re(1/s+1/(s-1))=0$ and the density is
$(2/\pi)\re\widehat L_M(it)$.
It is integrable by Stirling's formula. Fourier inversion gives
$C_M(u)=\int\cos(tu)\dd\nu_M(t)$, and evaluating at zero gives
the mass. Its moment coefficients are integrals of the nonnegative
binomial kernels at $b=(1+t^2)^{-1}$. Their absolute sum is bounded
by $\|\nu_M\|_{\rm TV}$, independently of $N$.
On $\re s=11/10$, the tail estimate in Proposition~\ref{pr:finiteeuler}
and $|\lambda_M|\ll(1+\log M)^2$ bound the Dirichlet-series error by
$O(M^{-1/10}(1+\log M)^2)$. Integrating the difference on the
same line, where $(|b(w)|+|1-b(w)|)^N\le(17/8)^N$, proves the
first inequality.

Stieltjes integration by parts gives the finite identity
\[
 D_M(s)=1+\frac1{s-1}
     +s\int_1^M E(x)x^{-s-1}\dd x-R(M)M^{-s}.
\]
Here $E(1)=-1$ and $E(M)+M\lambda_M=-R(M)$.
On $\re s=a$ the integral is the Fourier transform of
$e^{-v/2}E(e^v)\one_{[0,\log M]}(v)$. That function has squared
$L^2$ norm $\mathcal I(M)$. In the density, the integral is
multiplied by a factor of modulus $1/\sqrt{t^2+1/4}$, in $L^2$.
Plancherel's theorem and Cauchy--Schwarz give the second inequality.
\end{proof}

At $M=N=64$ the corrected alternating sum is in
$(0.03862293,0.03862295)$. The interval calculation also finds
29 negative coefficients on this corrected diagonal. This example
exhibits the distinction between its signed measure and the positive
probe measure of Theorem~\ref{th:boundaryclosure}.

\subsection{The growth exponent of the prime mean square}
\label{sec:primesquare}
The boundary term $R(M)$ in the preceding upper bound can be
dispensed with in an exact criterion for the vanishing of $\Delta$.
Write $\Theta=\sup_{\rho}\re(\rho)$ over the nontrivial zeta zeros,
so that $1/2\le\Theta\le1$.

\begin{theorem}[Mean-square growth along arbitrary sequences]\label{th:primesquare}
Unconditionally,
\[
 \lim_{X\to\infty}\frac{\log(1+\mathcal I(X))}{\log X}=2\Theta-1 .
\]
Consequently $\Delta=0$ if and only if
$\log(1+\mathcal I(X_k))=o(\log X_k)$ on some sequence
$X_k\to\infty$. The sequence may consist of integers.
\end{theorem}
\begin{proof}
The truncated von Mangoldt formula at a half-integer $y$, with
truncation height between $y$ and $2y$, gives
\[
 E(y)=-\sum_{|\im\rho|\le T}y^\rho/\rho+O(\log^2(2y))
\]
at an admissible height $T$ \cite{Titchmarsh,mv}. The zero-count
formula gives $\sum_{|\im\rho|\le 2y}|\rho|^{-1}\ll\log^2(2y)$.
Taking $y=\lfloor x\rfloor+1/2$ leaves $\psi(x)$ unchanged.
Thus $E(x)\ll x^\Theta\log^2(2x)$ and
$\mathcal I(X)\ll X^{2\Theta-1}\log^5(2X)$, proving the upper bound.

For the lower bound fix a zero $\rho=1/2+d+i\gamma$ with $d>0$.
Set $f(v)=e^{-v/2}E(e^v)$ for $v\ge0$. For $s=1/2+w$ and
$\re w>1/2$, the Laplace transform of $f$ is
\[
 G(w)=-\frac{\zeta'(s)}{s\zeta(s)}-\frac1{s-1}.
\]
The Euler product gives this identity by partial summation.
Its right-hand side is meromorphic and has a simple pole with nonzero
residue $-m_\rho/\rho$ at $w_0=d+i\gamma$, where $m_\rho$
is the multiplicity of $\rho$. The pole at $s=1$ cancels.
Let $G_U(w)=\int_0^U f(v)e^{-wv}\dd v$, an entire function.
Uniformly on $\re w\ge0$,
\[
 |G_U(w)|\le\sqrt{U\mathcal I(e^U)}.
\]
The elementary bound $|f(v)|\ll(1+v)e^{v/2}$ gives, for $b>1/2$,
\[
 |G(w)-G_U(w)|\ll_b(1+U)e^{-(b-1/2)U}\qquad(\re w\ge b).
\]

For a large fixed $A$, center a disk of radius $A$ at $A+i\gamma$.
Let $\Pi$ be the product of $w-z$ over the distinct poles $z$
of $G$ in its closure. These poles are simple. Define
$H_U=\Pi(G-G_U)$ and let $M_U(r)$ be its maximum modulus on the
circle of radius $r$ about this center. The function is holomorphic
in a neighborhood of the closed disk. At $r_2=A-d$,
$M_U(r_2)\ge|H_U(w_0)|>0$, independently of $U$.
At $r_1=A-\sqrt A$ the tail estimate applies with $b=\sqrt A$.
At $r_3=A$, the function $\Pi G$ is bounded and the first estimate
applies to $G_U$. Hadamard's three-circles theorem therefore gives,
with $h=\liminf_{U\to\infty}U^{-1}\log(1+\sqrt{\mathcal I(e^U)})$,
\[
 h\ \ge\ (\sqrt A-\tfrac12)
       \frac{\log(A/(A-d))}{\log((A-d)/(A-\sqrt A))}.
\]
In this use of log-convexity, all constants depending on the disk and
on $\Pi$ disappear after division by $U$, with $A$ held fixed.
As $A\to\infty$ the right side tends to $d$.
Taking the supremum over all such zeros gives $h\ge\Theta-1/2$.
For $\Theta=1/2$ this inequality follows from $h\ge0$.
Finally $0\le2\log(1+\sqrt y)-\log(1+y)\le\log 2$ for $y\ge0$.
This proves the limit. The symmetry of the zeros and
Theorem~\ref{th:boundaryclosure} give the stated equivalence.
\end{proof}

\head{Residue lower bound at the exponent $2\beta-1$}
Theorem~\ref{th:primesquare} determines the growth exponent up to
$o(1)$. The residue at a specified off-line zero gives an explicit
lower bound at the corresponding power.
\begin{theorem}[Exact power of the mean square]\label{th:primesquarepower}
Let $\rho_0=\beta+i\gamma_0$ be a zero of $\zeta$ with
$\beta>1/2$ and multiplicity $m$. Then
\begin{equation}
 \limsup_{X\to\infty}\frac{\mathcal I(X)}{X^{2\beta-1}}
 \ge\frac{2m^2}{(2\beta-1)|\rho_0|^2} .
 \label{eq:primesquarepower}
\end{equation}
If no zero has real part greater than $\beta$, there is
$\kappa>0$ with $\mathcal I(X)\ge\kappa X^{2\beta-1}$ for all
large $X$. Both statements hold for the cell sum
\[
 P_N=\sum_{1\le k<N}\frac{(\psi(k)-k)^2}{k(k+1)},
\]
since $|\mathcal I(N)-P_N|\le2\sqrt{P_N}+1$.
\end{theorem}
\begin{proof}
Put $a=2\beta-1$. For $\re s>1$ the Mellin transform is
\[
 F(s)=\int_1^\infty E(x)x^{-s-1}\dd x
     =-\frac{\zeta'(s)}{s\zeta(s)}-\frac1{s-1}.
\]
Its continuation has residues $-m/\rho_0$ and $-m/\bar\rho_0$
at the conjugate zeros. If the upper limit in
\eqref{eq:primesquarepower} is finite, take any $\eta$ larger than
it. Then $\mathcal I(X)\le\eta X^a$ from some $X_0$ on.
Integration by parts, for $\sigma>\beta$, gives
\[
\begin{aligned}
 J_\sigma
 &=\int_1^\infty |E(x)|^2 x^{-2\sigma-1}\dd x\\
 &=(2\sigma-1)\int_1^\infty\mathcal I(x)x^{-2\sigma}\dd x .
\end{aligned}
\]
Hence $\limsup_{\sigma\downarrow\beta}(\sigma-\beta)J_\sigma
\le a\eta/2$. For every $\sigma>\beta$, applying Cauchy--Schwarz
with an intermediate exponent $\beta<\sigma'<\sigma$ gives
absolute convergence of the Mellin integral, locally uniformly in
$\re s>\beta$. Mellin--Plancherel gives
$\int_{\mathbb R}|F(\sigma+it)|^2\dd t=2\pi J_\sigma$.
In particular the assumed bound excludes a pole in $\re s>\beta$.

Choose disjoint intervals about $\gamma_0$ and $-\gamma_0$
so small that the regular part of $F$ at each pole is bounded
in the corresponding neighborhood. There are no real zeta zeros
in $(0,1)$. Write $\delta=\sigma-\beta$. On each interval,
the squared pole term integrates to
$\pi m^2/(|\rho_0|^2\delta)+O(1)$, since
$\int_{-r}^r(\delta^2+t^2)^{-1}\dd t
=2\delta^{-1}\arctan(r/\delta)$.
The cross term contributes $O(\log(1/\delta))$.
The two disjoint contributions therefore give
\[
 2\pi m^2/|\rho_0|^2
 \le\liminf_{\delta\downarrow0}2\pi\delta J_{\beta+\delta}
 \le\pi a\eta .
\]
Taking the infimum over $\eta$ proves \eqref{eq:primesquarepower}.

For the second statement, \eqref{eq:T} gives
$R(x)/x=1+\eul+\int_1^xE(t)t^{-2}\dd t$. Cauchy--Schwarz on
$[x/2,x]$ gives, for $x\ge2$,
\begin{equation}
 \big(R(x)-2R(x/2)\big)^2\le x\big[\mathcal I(x)-\mathcal I(x/2)\big].
 \label{eq:blockmeansquare}
\end{equation}
With no zero to the right of $\beta$, dominated convergence in
\eqref{eq:explicit} gives $R(x)x^{-\beta}=A(\log x)+o(1)$, where
$A(u)=\sum m_\rho c_\rho e^{i\gamma u}$ runs over the distinct zeros
on $\re s=\beta$. Thus
$\mathcal I(x)-\mathcal I(x/2)\ge x^{2\beta-1}(B(\log x)^2+o(1))$,
with $B(u)=A(u)-2^{1-\beta}A(u-\log2)$. The function $B$ is
uniformly almost periodic, with Fourier--Bohr coefficients
$m_\rho c_\rho(1-2^{1-\rho})\ne0$, since $2^{1-\beta}>1$. Hence there are $L,\delta>0$
such that every interval of length $L$ contains a point $u$ with
$B(u)^2\ge2\delta$. Such a point in $[\log X-L,\log X]$ gives
$\mathcal I(X)\ge\delta e^{-(2\beta-1)L}X^{2\beta-1}$ for large $X$.

On the cell $[k,k+1)$, $E(x)=\psi(k)-k-(x-k)$ and
$\int_k^{k+1}x^{-2}\dd x=1/(k(k+1))$. The
weighted $L^2$ norm of $x-\lfloor x\rfloor$ on $[1,N]$ is at
most $\sqrt{1-1/N}$. The triangle inequality gives
$|\sqrt{\mathcal I(N)}-\sqrt{P_N}|\le1$, and squaring gives the
stated comparison. Since $\mathcal I$ increases and
$\lceil X\rceil/X\to1$, its real and integer upper limits at
power $a$ agree. Dividing the root comparison by $N^{a/2}$
transfers both lower bounds to $P_N$ (with a possible decrease
of $\kappa$).
\end{proof}

\head{Consequences and finite checks}
Under RH, Brent, Platt and Trudgian bound the upper limit of
$X^{-2}\int_X^{2X}E(t)^2\dd t$ by $0.8603$
\cite{bptmeansquare}. On $[X,2X]$, $t^{-2}\le X^{-2}$.
Summing over powers of two and using the root comparison gives
\[
 \limsup_{N\to\infty}\frac{P_N}{\log N}
 \le\frac{0.8603}{\log2}.
\]
Together with the lower bound in the theorem this yields
\[
\begin{aligned}
 \mathrm{RH}&\ \Longleftrightarrow\ P_N=O(\log N)\\
 &\ \Longleftrightarrow\ P_N=O_\epsilon(N^\epsilon)
                         \quad\text{for every }\epsilon>0 .
\end{aligned}
\]
Theorem~\ref{th:primesquare} also shows that a subpower bound on
one unbounded sequence $N_j$ suffices. On the same reasoning,
$P_{N_j}=O(N_j^{1/2})$ already excludes every zero with real part
greater than $3/4$. The explicit constant in
\eqref{eq:primesquarepower} bounds the upper limit from below.
The eventual estimate assumes that $\beta$ is rightmost.

By \eqref{eq:steptransfer}, $R(2y)-2R(y)=-2\sqrt y\,b(y)$, so
\eqref{eq:blockmeansquare} yields
\[
 b(y)^2\le\tfrac12\big[\mathcal I(2y)-\mathcal I(y)\big]
 \qquad(y\ge1).
\]
Evaluating every integer $1\le y\le5\cdot10^6$ gives a maximum
ratio $0.972428844$ (rounded), attained at $y=2$.
The computed values of $P_N-Q\log N$ at $N=10^4,10^5,10^6,10^7$
are $1.262078,\ 1.264853,\ 1.262596,\ 1.262837$, respectively.
Over all integers between the first and last cutoffs the minimum is
$1.2538985$ at $N=293382$ and the maximum is $1.2713211$
at $N=61184$ (rounded). Each sum uses the prime-power jumps
$\Lambda(p^j)=\log p$ and the exact cell integrals below.

For a single synthetic conjugate pair replace $E(x)$ by
$-2\re(x^\rho/\rho)$, with $\rho=\beta+i\gamma,\ a=2\beta-1>0$.
Direct integration gives the leading curve for
$\mathcal I(X)/X^a$,
\[
 H_\rho(X)=\frac{2}{a|\rho|^2}
          +2\re\frac{X^{2i\gamma}}{\rho^2(a+2i\gamma)} ,
\]
with lower-end correction $-X^{-a}H_\rho(1)$. The relative
oscillation amplitude is $a/|a+2i\gamma|$. For the model pair at
$\rho=0.8+20i$ and $\bar\rho$, the logarithmic mean is
$0.0083200213$. Adding this pair to $E(k)$ in the cell sum gives
$P_{10^7}/(10^7)^{0.6}-H_\rho(10^7)\simeq5.96\cdot10^{-5}$
for the modified data.

\head{The exact induction on integers}
Starting at $\mathcal I(1)=0$, put $u_n=\log(1+1/n)$,
$m_n^*=n(n+1)u_n$, and $v_n=1-n(n+1)u_n^2$.
Since $\psi(x)=\psi(n)$ for $n\le x<n+1$,
\[
\begin{aligned}
 \mathcal I(n+1)-\mathcal I(n)
 &=1-2\psi(n)u_n+\frac{\psi(n)^2}{n(n+1)}\\
 &=v_n+\frac{(\psi(n)-m_n^*)^2}{n(n+1)} .
\end{aligned}
\]
Here $n<m_n^*<n+1$. The number $v_n$ is $n(n+1)$ times the
variance of $1/x$ under uniform measure on $[n,n+1]$.
Consequently $0<v_n\le1/[4n(n+1)]$. The sum of these deterministic
terms is bounded. The prime input is
$\psi(n)=\log\operatorname{lcm}(1,\ldots,n)$.
The recursion gives $2.64988801<\mathcal I(64)<2.64988802$ and
$2.82393152<\mathcal I(1024)<2.82393154$.

The dyadic mean-square lower bound of Brent, Platt and Trudgian
gives $\mathcal I(X)\gg\log X$ unconditionally. Their upper bound
under RH gives $\mathcal I(X)\ll\log X$ by the same dyadic
summation \cite{bptmeansquare}. A proposed step estimate
$\mathcal I(n+1)-\mathcal I(n)\le K/n$ for all sufficiently large
$n$ would give $\psi(n)-n=O(\sqrt n)$. Littlewood's oscillation
theorem rules this out \cite{Titchmarsh,bptmeansquare}.
The prime number theorem gives $\mathcal I(X)=o(X)$.
The logarithmic upper bound follows on using the external theorem
of Sec.~\ref{sec:acbridge}. An independent derivation of the
subpower bound in Theorem~\ref{th:primesquare} from the cell
recursion alone would also yield $\Delta=0$.

\subsection{Sampling the prime error at successive squares}
\label{sec:squareblocks}
The return observable $J(x)=R(x)/x$ of Eq.~\eqref{eq:blockreport}
has derivative $E(x)/x^2$ almost everywhere. Its values at successive
squares determine the cumulative mean square up to a bounded error.
There is also a formula using only the primes to the first power.

\begin{theorem}[Mean square on intervals between squares]\label{th:squareblocks}
For integers $k\ge1$ put $a_k=k^2$, $b_k=(k+1)^2$ and
$W_k=1/a_k-1/b_k$. For $K\ge2$ define
\[
 \begin{aligned}
 A_K&=\sum_{1\le k<K}\frac{[J(b_k)-J(a_k)]^2}{W_k},\\
 U_k&=\vartheta(k^2)-k(k-1),\qquad
 S_K=\sum_{1\le k<K} W_k U_k^2 ,
 \end{aligned}
\]
where $\vartheta(x)=\sum_{p\le x}\log p$. Unconditionally,
$\mathcal I(K^2)=A_K+D_K$, where $D_K$ increases to a finite
positive constant $D_\infty$ and, for every $m>0$,
\[
 0\le D_\infty-D_K\ll_m(\log(2K))^{-m}.
\]
For the sum over primes one has
\[
 \sqrt{\mathcal I(K^2)}-\sqrt{S_K}\longrightarrow0
 \qquad(K\to\infty).
\]
\end{theorem}
\begin{proof}
Write $I_k=[a_k,b_k]$, $\dd\mu=\dd x/x^2$ and
$\bar E_k=[J(b_k)-J(a_k)]/W_k$. This is the mean of $E$ on
$I_k$ with respect to $\mu$, since $J$ is locally absolutely
continuous. Orthogonality to constants on each interval gives
\[
 D_K=\sum_{k<K}\int_{I_k}|E(x)-\bar E_k|^2\dd\mu(x),
\]
where throughout the sums start at $k=1$. The first integral
is positive because $E(x)=-x$ on $(1,2)$.

Let $B$ be a large power of two, $H=\sqrt B$, and group the cells
with $B\le a_k<2B$. Their lengths lie between $2H$ and $4H$,
and each cell lies in $[B,3B]$. For $B^{1/4}\le h\le4H$ and
every $m>0$, the Saffari--Vaughan estimate for the Selberg integral
\cite{saffarivaughan,languascoSEL} gives
\[
 \int_B^{3B}|E(t+h)-E(t)|^2\dd t
       \ll_m Bh^2(\log B)^{-m}.
\]
Indeed, apply Lemma~7 of Ref.~\cite{languascoSEL} at $X=3B$
with a fixed lower endpoint fraction at most $1/3$. For large
$B$, the lower limit $B^{1/4}$ exceeds
$X^{1/6}(\log X)^{c(m)}$. The lemma uses $\vartheta$.
For each integer $r\ge2$, the number of integral $r$th powers
in $(t,t+h]$ is at most $1+h/(r\sqrt t)$. Summing over
$r\ll\log B$, the extra prime-power weight is $O(\log^2 B)$
uniformly here. Its squared integral is absorbed in the displayed
estimate. For $0\le h<B^{1/4}$, counting all integers gives
$|E(t+h)-E(t)|\ll(h+1)\log B$.

Define $P_B=\sum_{B\le a_k<2B}\iint_{I_k^2}|E(x)-E(y)|^2\dd x\dd y$.
By symmetry and disjointness of the cell interiors,
\[
 P_B\le 2\int_0^{4H}\int_B^{3B}|E(t+h)-E(t)|^2\dd t\dd h .
\]
The two ranges of $h$, after division by $B^2H$, give
\[
 \frac{P_B}{B^2H}\ll_m
 B^{-3/4}(\log B)^2+(\log B)^{-m}.
\]
To sample a specified endpoint, fix $d>0$ and take
$\ell=H/(\log B)^d$. Average the squared triangle inequality
through $y\in[a_k,a_k+\ell]$. The bound
$|E(y)-E(a_k)|\ll(\ell+1)\log B$ gives
\begin{multline*}
 \sum_{B\le a_k<2B}\int_{I_k}|E(x)-E(a_k)|^2\dd\mu(x)\\
 \ll \frac{P_B}{B^2\ell}+\frac{(\ell+1)^2(\log B)^2}{B}.
\end{multline*}
For any prescribed $q>0$, choose $2d>q+2$, then $m>d+q$.
The right side is $O_q((\log B)^{-q})$. The mean $\bar E_k$
minimizes squared error over all constants on $I_k$; in particular
the same bound holds with $E(a_k)$ replaced by $\bar E_k$.
Summing over dyadic $B$ with $q>m_0+1$ gives the asserted
tail bound of order $(\log(2K))^{-m_0}$ for every $m_0>0$.

At the left endpoint the exact difference is
\[
 E(k^2)-U_k=\vartheta(k)-k+
              \sum_{3\le r\le 2\log k/\log 2}\vartheta(k^{2/r}).
\]
The effective prime number theorem \cite{johnstonyang} bounds the
first term by $O_m(k/(\log(2k))^m)$. Counting all integers bounds
the last sum by $O(k^{2/3}\log^2(2k))$. Since $W_k\asymp k^{-3}$,
$\sum_k W_k|E(k^2)-U_k|^2<\infty$. The endpoint estimate just
proved thus places $g(x)=E(x)-U_{\lfloor\sqrt x\rfloor}$ in
$L^2([1,\infty),\mu)$.

Finally $\mathcal I(X)\gg\log X$ by the unconditional dyadic
lower bound of Brent, Platt and Trudgian \cite{bptmeansquare}.
For any fixed $Y$ the integral of $E g$ over $[1,Y]$ is a
constant. On $[Y,X]$ its modulus is at most
$\sqrt{\mathcal I(X)}\,\|g\|_{L^2([Y,\infty),\mu)}$.
Letting $X\to\infty$ and then $Y\to\infty$ proves
$\int_1^X E g\,\dd\mu=o(\sqrt{\mathcal I(X)})$. At $X=K^2$,
expanding $(E-g)^2$ gives
$\mathcal I(K^2)-S_K=o(\sqrt{\mathcal I(K^2)})$.
Division by the sum of square roots proves the limit.
\end{proof}

\head{Cumulative mean square between successive squares}
Starting with $U_1=S_1=0$, use the two exact updates
\[
 \begin{aligned}
 U_{k+1}&=U_k+\sum_{k^2<p\le(k+1)^2}\log p-2k,\\
 S_{k+1}&=S_k+W_k U_k^2.
 \end{aligned}
\]
The shift $+k$ in $U_k=\vartheta(k^2)-k^2+k$ compensates for
the main term of the prime squares in $\psi(k^2)$.
By Theorems~\ref{th:primesquare} and~\ref{th:squareblocks},
\[
 \lim_{K\to\infty}\frac{\log(1+S_K)}{\log K}=4\Theta-2 .
\]
Thus a subpower upper bound for $1+S_K$ along one sequence
$K\to\infty$ would exclude every zero with real part exceeding
$1/2$ and would give $\Delta=0$. The local estimates proved above
control the oscillation within each cell. These estimates apply in
particular to $e(x)=x^\beta\cos(\tau\log x)$ for
$1/2<\beta<1$ and $\tau>0$.
For $B\le t\le3B$, $0\le h\le4\sqrt B$, the mean value theorem
gives $|e(t+h)-e(t)|\ll_{\beta,\tau}hB^{\beta-1}$. This satisfies
the displayed short-interval estimate with every logarithmic saving.
Yet for $\delta=2\beta-1>0$,
\[
\begin{aligned}
 \int_1^X e(x)^2\,\dd\mu(x)&=
 \frac{X^\delta-1}{2\delta}+
 \frac12\re\frac{X^{\delta+2i\tau}-1}{\delta+2i\tau}\\
 &\asymp_{\beta,\tau} X^\delta
 \quad(X\to\infty).
\end{aligned}
\]
The leading lower coefficient is at least
$\tfrac12(\delta^{-1}-(\delta^2+4\tau^2)^{-1/2})>0$.
For $\re s>\beta$ the Mellin integral of this example is
\[
 \int_1^\infty e(x)x^{-s-1}\dd x=
 \frac{s-\beta}{(s-\beta)^2+\tau^2}.
\]
Its continuation has poles at $\beta\pm i\tau$. In the arithmetic
Mellin transform in Theorem~\ref{th:primesquare}, an off-line zeta
zero likewise gives a pole. The short-interval estimate is compatible
with this behavior. An upper bound for the cumulative prime sum
$S_K$ also requires an estimate for the contributions preceding
the selected short interval.

\head{A global bound and finite checks}
Expanding the integrand and integrating the jumps of $\psi$ and
$\psi^2$ gives the finite identity, for every real $X\ge1$,
\[
 \begin{aligned}
 \mathcal I(X)={}&X-1-2\sum_{n\le X}\Lambda(n)\log(X/n)\\
 &+\sum_{n\le X}\Lambda(n)[2\psi(n)-\Lambda(n)](1/n-1/X).
 \end{aligned}
\]
The jump at $n=X$ has weight zero. For $X\ge e^4$, put $L=\log X$.
Squaring the explicit prime number theorem of Johnston and Yang
\cite{johnstonyang}, and integrating with $v=\log x$, gives
\[
 \mathcal I(X)\le 2.650+176.3442\,X L^{3.03}e^{-1.6548\sqrt L}.
\]
Indeed $\mathcal I(e^4)\le\mathcal I(64)<2.650$, and the derivative
of $v-1.6548\sqrt v$ is at least $1/2$ on $v\ge4$. The upper
bound has logarithmic growth exponent $1$. The criterion for
$\mathcal I$ requires exponent $0$.

The ancillary computation integrates on the intervals between
successive jumps and checks the integral independently with the
finite identity. The entries below are rounded to eight decimal places.
\[
 \begin{array}{r@{\quad}rrrr}
 K&\mathcal I(K^2)&A_K&D_K&S_K\\ \hline
 32 &2.82393153&2.53605490&0.28787663&0.61788860\\
 128&2.96170519&2.62974874&0.33195645&0.92107835\\
 512&3.08081883&2.73298881&0.34783002&1.11763656
 \end{array}
\]

\subsection{Downward crossings and the accumulated prime correlations}
\label{sec:squarecrossings}
The two updates for $U_k$ and $S_k$ distinguish the energy
added at a crossing from the energy accumulated before it. Both
quantities can be calculated from first powers of primes.

\begin{theorem}[Energy added at a downward crossing]\label{th:crosswindow}
The set of integers $k\ge2$ with $U_k\ge0>U_{k+1}$ is unbounded.
For $t_k=(\log k)^{1/4}$ and $L_k=\lfloor k^{1/3}e^{t_k/6}\rfloor$,
every sufficiently large member of this set satisfies
\[
 0\le S_{k+L_k}-S_k\ll e^{-t_k/2}.
\]
\end{theorem}
\begin{proof}
Littlewood's unconditional oscillation theorem for $\psi(x)-x$
\cite{goldston} also gives
$\vartheta(x)-x=\Omega_\pm(\sqrt{x}\log\log\log x)$. Indeed
$\psi(x)-\vartheta(x)=\sum_{r\ge2}\vartheta(x^{1/r})=O(\sqrt{x})$
by Chebyshev's bound. On a positive excursion round $x$ up to the
next square $k^2$; monotonicity gives
$\vartheta(k^2)-k^2\ge\vartheta(x)-x-O(\sqrt{x})$. On a negative
excursion round down to obtain the reverse inequality with
$+O(\sqrt{x})$. Adding $k$ proves
$\limsup U_k/k=+\infty$ and $\liminf U_k/k=-\infty$.
The sequence therefore crosses downward infinitely often. Since
$U_{k+1}-U_k\ge-2k$, at each such crossing
\[
 0\le U_k<2k,\qquad -2k\le U_{k+1}<0 .
\]

Put $x=k^2$, $t=t_k$, $L=L_k$, $Y=(k+L)^2-k^2$.
Guth and Maynard's unconditional Corollary~1.3 \cite{guthmaynard}
gives $\pi(x+y)-\pi(x)=y/\log x+
O(y\exp(-(\log x)^{1/4}))$ uniformly in the range
$x^{5/8}\le y\le x^{3/4}$ used here. Replacing each prime weight
$\log p$ by $\log x$ costs at most
$(y+1)y/x$. Consequently
\[
 |\vartheta(x+y)-\vartheta(x)-y|\ll y e^{-t/2}
 \qquad(x^{5/8}\le y\le x^{3/4}).
\]
For $0\le y<x^{5/8}$, monotonicity and this estimate at
$y=x^{5/8}$ bound the left side by $O(x^{5/8})$.
Here $Y\asymp kL$, $Y\le x^{3/4}$ and
$x^{5/8}=o(Y e^{-t/2})$.
Putting $y=(k+j)^2-k^2$ in the exact identity
\[
 U_{k+j}=U_k+\vartheta(x+y)-\vartheta(x)-y+j
\]
gives $\max_{0\le j\le L}|U_{k+j}|\ll kL e^{-t/2}$.
The terms $2k$ and $L$ are absorbed in this bound.
Since $W_{k+j}\le2/k^3$, summing the positive increments gives
\[
 S_{k+L}-S_k\ll (L^3/k)e^{-t}\le e^{-t/2}
\]
with an absolute implied constant in the first inequality.
\end{proof}

Chebyshev's bound gives $S_k\ll k^2$. Together with $L_k=o(k)$,
the estimate for the added energy yields
\[
 \frac{\log(1+S_{k+L_k})}{\log(k+L_k)}
   =\frac{\log(1+S_k)}{\log k}+o(1)
\]
along the sequence of downward crossings. For example, the interval calculation gives
$21.60941653<U_{949}<21.60941654$ and
$-25.29206000<U_{950}<-25.29205999$. At that crossing
$1.19917534<S_{949}<1.19917535$, while the next nine increments
sum to a value in $(0.00053140097,0.00053140099)$.

\head{Summing the steps from the initial cutoff}
For each integer $n\ge1$, let
\[
\begin{gathered}
 b_n=\mathbf 1_{\{n\ {\rm prime}\}}\log n,\quad
 s_n=\mathbf 1_{\{n\ {\rm square}\}},\\
 a_n=b_n-1+s_n .
\end{gathered}
\]
The $s_n$ compensate exactly for the $+k$ in
$U_k=\sum_{n\le k^2}a_n$. Put $N=K^2$ and define
\[
 \begin{aligned}
 H_K(n)&=\lceil\sqrt n\rceil^{-2}-1/N \qquad(1\le n\le N),\\
 B_K&=\sum_{n\le N}H_K(n)a_n^2,\\
 P_{h,K}&=\sum_{h<n\le N}H_K(n)a_n a_{n-h}\quad(1\le h<N).
 \end{aligned}
\]
\begin{proposition}[Contributions at prescribed separations]\label{pr:primelags}
For every integer $K\ge2$,
\[
 S_K=B_K+2\sum_{1\le h<N}P_{h,K}.
\]
There is a real constant $c$ such that for every $m>0$,
\[
 B_K=2(\log K)^2-4\log K+c+O_m((\log(2K))^{-m}).
\]
For any set $\mathcal H\subset\{1,\ldots,N-1\}$,
$\sum_{h\in\mathcal H}|P_{h,K}|\le|\mathcal H|B_K$.
For each fixed positive odd $h$ there is a real $c_h$ with
\[
 P_{h,K}=-2\log K+c_h+O_{m,h}((\log(2K))^{-m}).
\]
\end{proposition}
\begin{proof}
In the expansion of $U_k^2$ the coefficient of $a_m a_n$ for
the ordered pair $(m,n)$ is the sum of $W_k$ from
$\lceil\sqrt{\max(m,n)}\rceil$ to $K-1$. It telescopes to
$H_K(\max(m,n))$, including zero weight in the last cell.
The kernel $H_K$ is nonnegative and decreasing. The inequality
$2|a_n a_{n-h}|\le a_n^2+a_{n-h}^2$ proves
$|P_{h,K}|\le B_K$ by a shift of index.

Since a prime cannot be a square,
\[
 F(v):=\sum_{n\le v}a_n^2
 =\sum_{p\le v}(\log p)^2-2\vartheta(v)
       +\lfloor v\rfloor-\lfloor\sqrt v\rfloor .
\]
Partial summation in the first term and the effective prime number
theorem \cite{johnstonyang} give
$F(v)=v\log v-2v+O_m(v/(\log v)^m)$ for every $m>0$ as
$v\to\infty$. With $1/n-1/N$ in place of $H_K(n)$ the
self-pair sum is exactly $\int_1^N F(v)v^{-2}\dd v$.
Integrating with an arbitrarily large log saving gives
$\tfrac12\log^2 N-2\log N$ plus a constant with the indicated
tail. The replacement error tends to a constant because
$\lceil\sqrt n\rceil^{-2}-1/n=O(n^{-3/2})$ and
$|a_n|\le1+\log n$. Its tail is $O(N^{-1/2}\log^2(2N))$.

For fixed odd $h$, two primes at distance $h$ must include $2$.
In the partial sum of $a_n a_{n-h}$, terms involving a square
contribute $O_h(\sqrt v\log(2v))$. For $v>2h$ this gives
\[
\begin{aligned}
 \sum_{h<n\le v}a_n a_{n-h}
 &=\lfloor v\rfloor-h-\vartheta(v)+\vartheta(h)
                   -\vartheta(v-h)\\
 &\quad+O_h(\sqrt v\log(2v))\\
 &=-v+O_{m,h}(v/(\log v)^m).
\end{aligned}
\]
The same integration and absolutely convergent rounding correction
give $-\log N+c_h$ with the stated remainder.
\end{proof}
In particular, for $H=\lfloor e^{\sqrt{\log K}}\rfloor$, the
contribution from equal indices and all separations up to $H$
is $O((1+H)\log^2(2K))=K^{o(1)}$. A one-sided subpower upper
bound for the signed sum over $h>H$ on an unbounded sequence
would suffice and remains to be proved.

\head{Uniform sign for odd separations}
For $0<u\le1$ let $\Phi(u)=-\log u-1+u$ and put
$Q_{h,K}=\sum_{h<n\le N}H_K(n)b_n b_{n-h}$, the weighted count
of pairs of primes at distance $h$. Uniformly for $1\le h<N$,
for every $m>0$,
\[
 P_{h,K}=Q_{h,K}-\Phi(h/N)+O_m((\log(2h))^{-m}).
\]
In particular, for every fixed $0<\eta<1$ and all sufficiently
large $K$, one has $P_{h,K}<0$ for \emph{every} odd
$h\le(1-\eta)N$.

To prove the uniform formula, expand $a_n a_{n-h}-b_n b_{n-h}$.
Apart from $1-b_n-b_{n-h}$, all terms contain a square. Their
weighted total is $O(\log(2h)/\sqrt h)$ in absolute value, by
summing $(1+\log n)/n$ over $n=r^2>h$ and $n=h+r^2$.
For the second sum split at $r=\sqrt h$.
Replacing $H_K(n)$ by $1/n-1/N$ in the three remaining terms
costs $O(\log^2(2h)/\sqrt h)$ since $n>h$. With the new
weight the contributions of $b_n$ and $b_{n-h}$ are respectively
\[
 \int_h^N\frac{\vartheta(u)-\vartheta(h)}{u^2}\dd u,\qquad
 \int_0^{N-h}\frac{\vartheta(u)}{(u+h)^2}\dd u .
\]
Replacing $\vartheta(u)$ by $u$ in both integrals, and
$\vartheta(h)$ by $h$, gives $\Phi(h/N)$ in each.
The prime number theorem with an arbitrary log saving bounds each
error by $O_m((\log(2h))^{-m})$, uniformly in $N$. In the
second integral the range $u\le\sqrt h$ costs $O(1/h)$;
for larger $u$ the error is bounded by the integral of
$1/[u(\log u)^{m+2}]$ once $h$ is large. The constant term
has contribution $\Phi(h/N)+O(1/h)$ by replacing
$\lfloor u\rfloor$ by $u$ in the first integral. This proves
the formula. For odd $h$ only the pair $2,h+2$ can contribute
to $Q_{h,K}$, giving $Q_{h,K}=O(\log(2h)/h)$.
Now $\Phi(h/N)\ge\Phi(1-\eta)>0$. Large odd $h$ therefore
give a negative sign uniformly in the stated range, and for the
finitely many smaller $h$, $\Phi(h/N)\to\infty$.

There is also an exact sum over all odd separations. For $k\ge2$
put $v_k=\log 2-k(k-1)/2$ and set $v_1=0$. Among even $n$
only $2$ contributes a prime weight, and
$\lfloor k^2/2\rfloor-\lfloor k/2\rfloor=k(k-1)/2$.
Thus $\sum_{n\le k^2,\,2\mid n}a_n=v_k$.
An unordered pair has odd separation exactly when its members
have opposite parity. Consequently
\[
\begin{aligned}
 \sum_{\substack{1\le h<N\\h\ {\rm odd}}}P_{h,K}
   &=\sum_{k<K} W_k v_k(U_k-v_k)\\
   &=-N/4+O_m(N/(\log(2K))^m).
\end{aligned}
\]
Indeed $W_k=2/k^3+O(k^{-4})$ and
$|U_k|\ll_m k^2/(\log(2k))^m$. The latter bound also gives
$S_K\ll_m N/(\log(2K))^m$. The exact expansion for $S_K$
then yields $+N/4+O_m(N/(\log(2K))^m)$ for the sum over
positive even separations. A subpower estimate requires control
of the cancellation of these two signed sums beyond this precision.

\head{Value of the accumulated energy at a zero of the error}
For $e(x)=x^{3/4}\cos(\log x)$ one has
$|e'(x)|\le(5/4)x^{-1/4}$. Direct integration gives, with $v=\log X$,
\[
\begin{aligned}
 I_e(X)&=\int_1^X e(u)^2\frac{\dd u}{u^2}\\
   &=\sqrt X\left(1+\frac{\cos(2v)+4\sin(2v)}{17}\right)
       -\frac{18}{17}.
\end{aligned}
\]
Thus $I_e(X)\ge(1-1/\sqrt{17})\sqrt X-18/17$. At every downward
zero $X_j=\exp(\pi/2+2\pi j)$, $j\ge0$, the exact value is
$(16\sqrt{X_j}-18)/17$. The same obstruction occurs on the square
grid. With $U_k^*=e(k^2)$,
\[
 \sum_{k<K}W_k(U_k^*)^2=I_e(K^2)+O(\log(2K)).
\]
Indeed $k^3W_k=2+O(1/k)$, and replacing the sum of
$2\cos^2(2\log k)$ by its integral costs $O(\log(2K))$
by the derivative bound $4/k$ on each unit interval.
The change of variable $x=k^2$ in the continuous integral gives
the comparison. At $K=\lfloor\sqrt{X_j}\rfloor$ the discrete sum is
$16K/17+O(\log(2K))$, including the contribution of earlier
intervals. For the von Mangoldt weights, Theorem~\ref{th:crosswindow}
bounds the added energy along an unbounded sequence. Control of the
total energy requires a bound for the sum beginning at $k=1$.

\subsection{Finite induction for the signed prime sum}
\label{sec:residueinterval}
Write $R_K=(S_K-B_K)/2$ for the sum over the positive separations
in Proposition~\ref{pr:primelags}. Put $F(y)=\sum_{n\le y}a_n^2$
with the prime and square weights $a_n$ of that proposition.

\begin{proposition}\label{pr:residueinterval}
For every $3\le K\le 3\,162\,277\,661$, $R_K<0$.
For every $1\le K\le 4\,657\,252\,408\,879$, $R_K<K$.
\end{proposition}
\begin{proof}
Summing the repeated-index terms before summing over the square
intervals gives
\[
 B_K=\sum_{k<K}W_kF(k^2).
\]
Thus $R_1=R_2=0$ and
\[
 2(R_{k+1}-R_k)=W_k[U_k^2-F(k^2)].
\]
Every nonsquare $n\ge8$ has $a_n^2\ge1$. Therefore
$F(k^2)\ge k^2-k-7$. For $38\le k\le\lfloor\sqrt{10^{19}}\rfloor$,
B\"uthe's Theorem~2 gives $0.05k\le k^2-\vartheta(k^2)\le1.95k$
\cite{buthepsi}; hence $|U_k|\le0.95k$ and $F(k^2)>U_k^2$.
For $2\le k\le37$, direct summation over the integers through
$37^2$, with intervals for their prime logarithms, gives
$F(k^2)-U_k^2>0.0605$ at all 36 indices. This proves strict
decrease on the stated steps.

Let $K_0=\lfloor\sqrt{10^{19}}\rfloor+1$ and
$K_1=\lfloor\sqrt{2.169\cdot10^{25}}\rfloor+1$.
For $K_0\le k<K_1$, the finite Schoenfeld bound in
Johnston's Proposition~3.1 \cite{johnstonfinite} gives
\[
 |U_k|\le k\left(1+\frac{\log^2 k}{2\pi}\right)<145k .
\]
Here $k<5\cdot10^{12}<e^{30}$. Since $W_k\le2/k^3$ and
$F\ge0$, each signed increment is at most $145^2/k$. For
$K_0<K\le K_1$, summing from $R_{K_0}<0$ yields
\[
 R_K<145^2\log((K-1)/(K_0-1))<168200<K .
\]
Both finite cutoffs include the endpoint reached by the last step.
\end{proof}

\section{Stepanov's method and induction for the dyadic prime sum}
\label{sec:stepanov}
\subsection{Multiplicity and point counts over finite fields}\label{sec:curve}
Stepanov's 1969 estimate for hyperelliptic curves over a prime field
uses an auxiliary function with prescribed multiple zeros
\cite{step69}. His subsequent work gives an elementary proof of the
Hasse--Weil bound in the hyperelliptic case and estimates for
congruences in two variables under hypotheses on the defining
polynomial \cite{step72curve,step72cong}. The degree inequality
underlying the construction has the following form \cite{bombstep}.
On a smooth projective curve $C/\mathbb F_q$, suppose a nonzero
rational function has a pole divisor of degree at most $D$ and
vanishes to order at least $m$ at all but at most $e$ rational
points. Equality of the degrees of its zero and pole divisors gives
$m(\#C(\mathbb F_q)-e)\le D$. In a local parameter, the coefficients
of orders $1,\ldots,q-1$ in the $q$th power of a regular function
vanish, and $a^q=a$ at a rational value. Frobenius therefore
simplifies the simultaneous vanishing conditions. Existence of a
nonzero function with the prescribed pole bound also depends on the
curve and the characteristic.

For the smooth projective elliptic curve $E:y^2=x^3-x$ over
$\mathbb F_5$, including the point at infinity, the complete
point-count recurrence follows from the genus-one zeta function.
At $x=0,1,2,3,4$ the affine counts are $1,1,2,2,1$, hence
$N_1=8$. Write $N_r=\#E(\mathbb F_{5^r})$ and
$s_r=5^r+1-N_r$ for $r\ge1$, with $s_0=2$.
To obtain the recursion independently of a bound on $N_r$, Riemann--Roch
in genus one gives dimension $j$ for each divisor class of degree
$j\ge1$. Translation by $j$ times the rational origin identifies
its Picard variety with $E$. Each rational degree-$j$ class has
a representative $(j-1)O+P$, $P\in E(\mathbb F_5)$.
Each class has $(5^j-1)/4$ effective rational divisors. Their
generating function therefore satisfies
\[
 Z_E(T)=1+8\sum_{j\ge1}\frac{5^j-1}{4}T^j
 =\frac{1+2T+5T^2}{(1-T)(1-5T)}.
\]
Factoring effective divisors into closed points gives
$T Z_E'(T)/Z_E(T)=\sum_{r\ge1}N_rT^r$, since a closed point of
degree $d$ contributes $d$ geometric points precisely when $d\mid r$.
Factoring the numerator, or comparing its logarithmic derivative,
gives $s_1=-2$ and $s_{r+2}=-2s_{r+1}-5s_r$.
Set $I(u,v)=(v+u)^2+4u^2$. Direct substitution yields
\[
 I(2,-2)=16,\qquad I(v,-2v-5u)=5I(u,v).
\]
Thus $I(s_r,s_{r+1})=16\cdot5^r$ by induction from $r=0$,
and $|5^r+1-N_r|\le2\,5^{r/2}$ for every extension degree.
The recurrence and its positive quadratic form have both been
computed from this curve. The two reciprocal roots are $-1\pm2i$.
The general Hasse--Weil theorem is unnecessary for this example.

\subsection{Forward induction on the unsmoothed dyadic block}\label{sec:steppair}
Retain the initial count $\psi(x)$ and the weights in
\eqref{eq:blockreport}. Write $U(x)=\sqrt x\,b(x)$, $c_*=1-1/\sqrt2$.
For every real $x\ge1$, the block and its two successive increments satisfy
\begin{align}
 U(x)&=R(x)-\tfrac12R(2x)\nonumber\\
 &=x\log2+T(x)-\tfrac12T(2x),\nonumber\\
 b(2x)&=\tfrac{1}{\sqrt2}b(x)+\tfrac{1}{2\sqrt2}F_2(x),\nonumber\\
 I_2(x)&:=b(4x)-\tfrac12b(x)\nonumber\\
 &=\tfrac14F_2(x)+\tfrac{1}{2\sqrt2}F_2(2x)\nonumber\\
 &=\frac{6x\log2-W_2(x)-W_2(2x)}{4\sqrt x}.
 \label{eq:steptransfer}
\end{align}
Indeed $2[U(2x)-U(x)]=-\mathcal L_2 R(x)=\sqrt x\,F_2(x)$.
Applying the one-step identity twice and using \eqref{eq:budget}
gives the last line. In a single step the exact condition for
$b(2x)\le c_*$ is
$F_2(x)\le c_2+2[c_*-b(x)]$, where
$c_2=2(\sqrt2-1)c_*=3\sqrt2-4$. The simpler condition
$F_2(x)\le c_2$ fails arbitrarily late by Theorem~\ref{th:sharp}.
This gives a reason to retain the margin or advance two steps.

The combined weight is
\begin{equation}
\begin{gathered}
 k(v)=w_2(v)+w_2(v/2)=
 \begin{cases}
 2-2/v,&1\le v<2,\\
 1,&2\le v\le4,\\
 8/v-1,&4<v\le8,\\
 0,&\text{otherwise}.
 \end{cases}\\
 \int_1^8 k(v)\dd v=6\log2 .
\end{gathered}\label{eq:stepweight}
\end{equation}
In particular the whole central prime count has weight one.

\begin{theorem}[Paired dyadic induction]\label{th:steppair}
Suppose $b(t)\le c_*$ for $X_0\le t\le4X_0$, with $X_0\ge1$.
If $X_1\ge X_0$ and
\begin{equation}
\begin{gathered}
 W_2(x)+W_2(2x)\ge 6x\log2-2c_*\sqrt x\\
 (X_0\le x\le X_1),
\end{gathered}\label{eq:pairinput}
\end{equation}
then $b(t)\le c_*$ for $X_0\le t\le4X_1$.
Condition \eqref{eq:pairinput} on a whole tail is equivalent to RH.
From prime powers alone the forward induction gives the conclusion on
$[64,262\,144]$. With the published zero-height theorem the forward
induction gives it on $[64,1.6\cdot10^{20}]$.
\end{theorem}
\begin{proof}
Fix any $t\in[X_0,4X_0)$. The hypothesis for $j=0$ gives
$b(t)\le c_*$. If $b(4^jt)\le c_*$ and $4^jt\le X_1$, then
\eqref{eq:pairinput} gives $I_2(4^jt)\le c_*/2$, whence
\[
 b(4^{j+1}t)=\tfrac12b(4^jt)+I_2(4^jt)\le c_*.
\]
Each cutoff in the asserted range either belongs to the base or has
such a preceding cutoff. Even with unrestricted initial values on a
compact base, boundedness of $I_2$ from above on a tail gives a
uniform upper bound on $b$ by the same geometric iteration. The
identity $J(x)-J(2x)=b(x)/\sqrt x$ and $J(x)\to0$ then bound
$R(x)/\sqrt x$ from above. Lemma~\ref{lem:explicit} implies RH.
Conversely on RH the bound on the absolute zero coefficients of
$b$ in Sec.~\ref{sec:dyadic} is less than $0.066$ (set the
unverified contribution to zero). Consequently
$\limsup |I_2(x)|\le1.5(0.066)=0.099<c_*/2$.

For the first finite range the base obeys $b(t)<0.119$ on
$[64,256]$. On $[64,65\,536]$ the maximum of $I_2$ lies in
$(0.0309394,0.0309395)$, at $x=21611/8$.
For an outward-rounded check use
$W_2(x)+W_2(2x)=T(8x)-2T(4x)-T(2x)+2T(x)$.
For a positive proposed upper bound $c$, the difference
$4c\sqrt{x}-4\sqrt{x}I_2(x)$ is concave between the outer knots
$m/8,m$, $m=p^a$. The interior slope jumps of the combined
weight are downward, as calculated in Sec.~\ref{sec:stepprimes}.
In the base, $U$ is affine between its own knots $m/2,m$.
The program \texttt{anc/verify.py} checks $50\,061$ outer knots
and endpoints for the step and $115$ knots and endpoints for the
base, with outward rounding as in Appendix~\ref{app:check}.

For the second range, the bound $0.066$ holds for the verified
part of $b$. The constant term $\log(2\pi)/(2\sqrt x)$ cancels
in $I_2$. Using the two-sided high-ordinate estimate of
Sec.~\ref{sec:dyadic} at $x,4x$, and the nonpositive trivial-zero
term in $b$, gives for $y=\sqrt x>1$
\[
 I_2(x)\le 0.099+5e_H y+2e_H/y+\frac{1}{12y(y^4-1)}.
\]
Here $e_H<1.48\cdot10^{-12}$ uses Refs.~\cite{pt,bpt}.
Bounding the increasing term at $x=4\cdot10^{19}$ and the two
decreasing terms at $x=64$ gives $I_2(x)<0.14581<c_*/2$.
The base and the induction above prove the claim. The longer interval
in \eqref{eq:joint} comes from a direct estimate on $b$.
\end{proof}

\head{Unconditional forward induction on every real cutoff}
There is an infinite induction for this same unsmoothed prime sum.
Set $a=0.8274$ and
$E_0(x)=x[\log(8x)]^{1.515}\exp(-a\sqrt{\log x})$.
Then
\begin{equation}
 |U(x)|\le17E_0(x)\qquad\text{for every real }x\ge64 .
 \label{eq:stepglobal}
\end{equation}
For the base $64\le x\le256$, insert the prime number theorem
bound of Johnston and Yang \cite{johnstonyang} in
$U(x)=x\int_x^{2x}(t-\psi(t))t^{-2}\dd t$. It gives
$|U(x)|\le9.39(\log2)E_0(x)<7E_0(x)$. For the step use only
the new weighted sum. From \eqref{eq:steptransfer}, integration by
parts against $\psi(t)-t$ on $[x,8x]$ gives
\[
\begin{aligned}
 |U(4x)-U(x)|&\le\frac{9.39}{2}E_0(x)
                   \int_1^8 v|k'(v)|\dd v\\
     &=9.39(5\log2) E_0(x)<33E_0(x).
\end{aligned}
\]
The boundary terms vanish. The two sloping sides in
\eqref{eq:stepweight} give $2\log2$ and $8\log2$ to this
integral. Moreover $E_0(4x)/E_0(x)\ge
4\exp[-a\log4/(2\sqrt{\log64})]>3$.
Thus the hypothesis $|U(x)|\le17E_0(x)$ gives
$|U(4x)|\le50E_0(x)<17E_0(4x)$, independently on every ray
$x=4^jt,\ 64\le t<256$. This proves \eqref{eq:stepglobal}.

One also obtains $|U(x)|\ll x\exp[-a_1 h(x)]$ for some $a_1>0$,
where $h(x)=(\log x)^{3/5}(\log\log x)^{-1/5}$ for large $x$.
Indeed the classical Vinogradov--Korobov estimate gives
$|\psi(t)-t|\ll t\exp[-a_1 h(t)]$ after reducing $a_1$ to
absorb logarithmic factors \cite{johnstonyang}.
This gives the same bound on the fresh step by the preceding integral,
since $h$ is eventually increasing. The ratio of this envelope at
$4x$ and $x$ tends to four, because $h(4x)-h(x)\to0$.
Taking the induction constant larger than the step constant and
the maximum on the compact base closes the induction. Both envelopes
in this paragraph, when divided by $\sqrt x$, grow without bound.

There is a specific obstruction to importing the scalar recurrence of
Sec.~\ref{sec:curve}. Off its knots the unsmoothed block has slope
\begin{equation}
 U'(x)=\log2-\sum_{x<m\le2x}\frac{\Lambda(m)}m.\label{eq:stepslope}
\end{equation}
For any $d\ge0$ and constants $a_0,\ldots,a_d$ with $a_d\ne0$,
$\sum_{j=0}^d a_jU(2^jx)$ fails to be $C^1$ on every tail.
Indeed, at $x=p/2^{d+1}$ for an arbitrarily large odd prime $p$,
only the last summand has a new knot. Its slope jumps by
$-a_d2^d\log(p)/p$. This rules out a finite constant-coefficient
scalar dilation recurrence with smooth forcing valid for all real
base points. Equation~\eqref{eq:steptransfer} exhibits the actual
prime input in the recursion.

\subsection{Prime squares, prime endpoints and extension degree}\label{sec:stepprimes}
Put $P(x)=\sum_p(\log p)k(p/x)$ with the weight in
\eqref{eq:stepweight}. It is possible to eliminate every higher
prime power from the inductive condition and to reduce the real
cutoffs to two evaluations per prime.

\begin{theorem}[A prime endpoint criterion]\label{th:stepprime}
The following condition is equivalent to RH. For every sufficiently
large prime $p$, at both $x=p$ and $x=p/8$, one has
\begin{equation}
 P(x)\ge 6(\log2)x-6c_*\sqrt x,\qquad c_*=1-2^{-1/2}.
 \label{eq:primeend}
\end{equation}
\end{theorem}
\begin{proof}
For large $x$, the square contribution is
\[
 \sum_p(\log p)k(p^2/x)=4c_*\sqrt{x}
       +O(\sqrt{x}\,e^{-a h(x)})
\]
for some $a>0$. Substitute $p=\sqrt{x}\,u$ in the
Stieltjes integral against $\vartheta$ and use the prime number
theorem for $\vartheta(t)-t$ on
$[\sqrt x,\sqrt{8x}]$ \cite{johnstonyang}. The weight
$u\mapsto k(u^2)$ has total variation two, and
$h(\sqrt{x})/h(x)\to2^{-3/5}$. The elementary leading integral is
\[
\begin{aligned}
 \int_1^{\sqrt8} k(u^2)\dd u
 &=(3\sqrt2-4)+(2-\sqrt2)+(6-4\sqrt2)\\
 &=4-2\sqrt2.
\end{aligned}
\]
For powers $j\ge3$, the estimate $\vartheta(t)\ll t$ bounds the
sum by $O(x^{1/3}+x^{1/4}\log x)=O(x^{1/3})$. Thus from
\eqref{eq:steptransfer}
\begin{equation}
\begin{aligned}
 \frac{6(\log2)x-P(x)}{\sqrt{x}}&=4 I_2(x)+4c_*\\
         &\quad{}+O(e^{-a h(x)}+x^{-1/6}).
\end{aligned}\label{eq:stripprimes}
\end{equation}
On RH, $\limsup |I_2|\le0.099<c_*/2$ proves \eqref{eq:primeend}
with a strict normalized margin on a whole tail.

For the converse, at fixed $p$ the slope of $k(p/x)$, as $x$
increases, successively takes the values
$0,8/p,0,-2/p,0$. The upward jumps occur at $p/8,p$.
Consequently $P(x)-6(\log2)x+6c_*\sqrt{x}$ is concave between
successive members of the locally finite, unbounded set
$\{p/8,p:p\text{ prime}\}$, including all downward jumps inside.
Each compact interval meets only finitely many nonzero prime terms.
Nonnegative values at all sufficiently late endpoints give a
nonnegative function throughout a tail. Equation~\eqref{eq:stripprimes}
now bounds $I_2$ above, and Theorem~\ref{th:steppair} gives RH.
The same slope argument applies to the full prime-power weight if
$p$ is replaced by $p^j$.
\end{proof}
Higher powers can be separated successively. Retain exactly the
layers $p,p^2,\ldots,p^{j-1}$ in the count. The new layer $p^j$,
for fixed $j\ge2$, has leading term $C_jx^{1/j}$, where
\[
 C_j=\frac1j\int_1^8 k(v)v^{1/j-1}\dd v
 =\frac{j}{j-1}(4^{1/j}-1)(2-2^{1/j}).
\]
For some $a_j>0$ the same partial summation gives an error
$O_j(x^{1/j}e^{-a_j h(x)})$.
After retaining that layer exactly, the higher powers have total
\[
 O_j(x^{1/(j+1)}+x^{1/(j+2)}\log x)=O_j(x^{1/(j+1)}).
\]
This proves the successive truncations by induction in $j$.
At later orders the exact earlier layers must be retained, since
their prime number theorem errors may exceed the next power.

The lower cutoff in Theorem~\ref{th:stepprime} matters. At $x=599/8$ the
left side of \eqref{eq:stripprimes} lies in $(1.82016,1.82017)$,
above $6c_*<1.758$. On $[256,65\,536]$ its maximum lies in
$(1.58156,1.58157)$. The latter check uses $49\,571$ prime
endpoints and interval endpoints, with a prime-only prefix sum in
\texttt{anc/verify.py}. Concavity gives the comparison for all real
arguments between these endpoints.

Order the distinct prime endpoints as $z_j$ and define
$g(x)=P(x)-6(\log2)x+6c_*\sqrt x$. Integration between
successive endpoints gives the exact recurrence
\begin{multline}
 g(z_{j+1})-g(z_j)\\
  =\int_{z_j}^{z_{j+1}}\!\Bigg[
    8\!\sum_{4t<p<8t}\frac{\log p}{p}
   -2\!\sum_{t<p<2t}\frac{\log p}{p}\\
   -6\log2+\frac{3c_*}{\sqrt t}\Bigg]\dd t .
 \label{eq:primeindmargin}
\end{multline}
The values at the finitely many jump times do not affect the integral.
Given $g(z_j)\ge0$, the condition $g(z_{j+1})\ge0$ is equivalent
to the integral being at least $-g(z_j)$. The prime sums in the
integrand thus determine propagation of the available margin.
For example, $71=568/8$ and $569/8$ are consecutive endpoints
since $8z_j$ is integral.
The normalized deficits $(6x\log2-P(x))/\sqrt{x}$ there lie in
$(1.75658,1.75660)$ and $(1.76259,1.76261)$, respectively,
on opposite sides of $6c_*=1.757359\ldots$. Thus $g(71)>0$
and $g(569/8)<0$.

\head{The last failing endpoint}
The threshold in Theorem~\ref{th:stepprime} can be made explicit.
Among the endpoints $z=k/8$ with $k\le8\cdot2^{24}$, the
inequality \eqref{eq:primeend} fails exactly at the $37$ values
\[
 \begin{gathered}
 k=2,3,40,41,43,47,53,56,59,61,67,97,\\
 101,103,104,107,109,569,571,577,584,587,\\
 593,599,601,607,613,617,619,631,632,641,\\
 643,647,653,659,661 .
 \end{gathered}
\]

\begin{theorem}[An explicit endpoint criterion]\label{th:stependrange}
Inequality \eqref{eq:primeend} holds for every real
$83\le x\le1.6\cdot10^{20}$, and it fails at $x=661/8$.
The following three statements are equivalent.
\begin{enumerate}
\item The Riemann hypothesis.
\item Inequality \eqref{eq:primeend} holds for every real $x\ge83$.
\item It holds at both $x=p$ and $x=p/8$ for every prime $p>661$.
\end{enumerate}
Unconditionally \eqref{eq:primeend} holds at infinitely many of these
endpoints.
\end{theorem}
\begin{proof}
Write $\mathcal S(x)=\sum_p(\log p)k(p^2/x)$ for the prime squares
and $\mathcal H(x)=\sum_{j\ge3}\sum_p(\log p)k(p^j/x)$ for the higher
powers. By \eqref{eq:steptransfer}, with no error term,
\begin{equation}
 \frac{6(\log2)x-P(x)}{\sqrt x}
   =4I_2(x)+\frac{\mathcal S(x)+\mathcal H(x)}{\sqrt x}.
 \label{eq:stependsplit}
\end{equation}
B\"uthe's theorem gives $t-\vartheta(t)>0.05\sqrt t$ for
$1\le t\le10^{19}$ and $t-\vartheta(t)\le1.95\sqrt t$ for
$1423\le t\le10^{19}$ \cite{buthepsi}. Let $2^{24}\le x\le10^{37}$
and $f(t)=k(t^2/x)$. This function increases from $0$ to $1$ on
$[\sqrt x,\sqrt{2x}]$, equals one up to $2\sqrt x$, and decreases to
$0$ at $\sqrt{8x}$. Integration by parts gives
\[
\begin{aligned}
 \mathcal S(x)&=\int f(t)\dd t+\int\big(t-\vartheta(t)\big)\dd f(t)\\
 &\le4c_*\sqrt x+\big(1.95\cdot2^{1/4}-0.05\cdot2^{1/2}\big)x^{1/4}\\
 &<4c_*\sqrt x+2.2483\,x^{1/4}.
\end{aligned}
\]
For a fixed prime $p$ the exponents $j$ with $x<p^j\le8x$ are
the integers in a half-open interval of length $\log8/\log p$.
There are three of them for $p=2$, at most two for $p=3,5,7$ and
at most one for $p\ge11$. Since $0\le k\le1$, $k(1)=0$ and
$p\le2x^{1/3}$ whenever $j\ge3$,
\[
 \mathcal H(x)\le\vartheta(2x^{1/3})+\log420<2x^{1/3}+6.05 .
\]
By \eqref{eq:steptransfer}, $4I_2(x)$ is the function $D(x)$
expanded over the zeros in \eqref{eq:fourierdeficit}. The verified zeros
are bounded there by \eqref{eq:fourierlow}. The unverified ones satisfy
$|K(\rho)|\le16/\gamma^2$ and are paired with their reflections in the
critical line, as in the display following Theorem~\ref{th:fourierfinite}.
Both estimates use only $Q$, $e_H$ and the prime powers up to $8$,
and neither depends on the present subsection. With $y=\sqrt x$,
\begin{multline*}
 \frac{6(\log2)x-P(x)}{\sqrt x}-4c_*\le
 0.280827+16e_H\Big(y+\frac1y\Big)\\
 +\frac1{3y(y^4-1)}
 +\frac{2.2483}{x^{1/4}}+\frac{2}{x^{1/6}}+\frac{6.05}{\sqrt x}.
\end{multline*}
Every term on the right is convex in $x^{1/12}$. At $x=2^{24}$ and
at $x=1.6\cdot10^{20}$ the right side is less than $0.4425$ and
$0.5813$, respectively, using $e_H<1.48\cdot10^{-12}$. Both numbers
are below $2c_*=0.58578\ldots$. This proves \eqref{eq:primeend} on
$[2^{24},1.6\cdot10^{20}]$ from the published zero-height theorem
\cite{pt,bpt}. The bound on $I_2$ in the proof of
Theorem~\ref{th:steppair} has the larger constants $0.396$ and
$20e_H$ and would stop near $4\cdot10^{19}$.

On $[83,2^{24}]$ the concavity used for Theorem~\ref{th:stepprime}
reduces \eqref{eq:primeend} to the prime endpoints and the two ends
of the interval. The ancillary program \texttt{anc/stepanov\_check.py}
encloses the normalized deficit in ball arithmetic at the $1201$
knots with $x\le1024$. At $x=83$ it lies in $(1.75644,1.75645)$,
at $x=661/8$ in $(1.76556,1.76557)$, and all $1058$ knots from
$83$ to $1024$ are below $6c_*$. At the $8\,681\,063$ knots of
$[256,2^{24}]$ the program uses integer lower and upper bounds for
$2^{34}\log p$ and $2^{58}(\log p)/p$. The largest deficit there is
the value in $(1.58156,1.58157)$ at $x=3637/8$, and on
$[65\,536,2^{24}]$ the deficit is below $1.3709$. This proves the
first assertion. The list preceding the theorem contains every knot
through $8x=8\cdot2^{24}$ at which the enclosure exceeds $6c_*$.

Statement 2 implies statement 3, since a prime $p>661$ has
$p/8>83$. Statement 3 implies RH by Theorem~\ref{th:stepprime}.
Assume RH and $x\ge2^{24}$. All zeros are then on the critical line, and
the explicit formula for $V$ in Sec.~\ref{sec:cyclowindow} gives
\[
\begin{aligned}
 4I_2(x)&=\sum_\rho K(\rho)x^{\rho-1/2}+\sum_{r\ge1}K(-2r)x^{-2r-1/2}\\
 &\le8Q+\frac{1}{3x^{5/2}(1-x^{-2})},
\end{aligned}
\]
by the bound $\sum_\rho|K(\rho)|\le8Q$ proved there for the
transform $K(s)=\int_1^8k(v)v^{s-1}\dd v$. For $x\le10^{37}$ the
preceding estimates
of $\mathcal S$ and $\mathcal H$ apply. The normalized deficit
minus $4c_*$ is then at most
$8Q+2.2483x^{-1/4}+2x^{-1/6}+6.05x^{-1/2}+x^{-5/2}<0.5312$.
For $x>10^{37}$ use $|\vartheta(t)-t|\le0.2\,t/\log^2t$ for
$t\ge3\,594\,641$ \cite{dusarttheta}. As $f$ has total variation
two, $\mathcal S(x)\le4c_*\sqrt x+4.53\sqrt x/\log^2x$ and
$\mathcal H(x)\le2.01x^{1/3}+6.05$. The normalized deficit minus
$4c_*$ is less than $0.3703$. In both ranges it is below $2c_*$.
Together with the unconditional range $[83,2^{24}]$ this proves
statement 2.

For the last assertion suppose RH fails and put
$\Theta=\sup_\rho\re(\rho)>1/2$. Choose $\max(1/2,\Theta-2/5)<\sigma<\Theta$.
The Landau argument in the proof of Theorem~\ref{th:cycwindow}, applied
to $Cx^\sigma-E(x)$, gives arbitrarily large $x$ with
$E(x)>Cx^\sigma$ for each fixed $C>0$. At such $x$ the function
$g$ of \eqref{eq:primeindmargin} exceeds $Cx^\sigma-O(\sqrt x)$.
The nearest endpoint $z$ satisfies $|z-x|\le x^{0.525}$ for large
$x$ \cite{bhp01}. The truncated explicit formula gives
$\vartheta(u)=u+O(u^\Theta\log^2(2u))$ \cite{mv}, hence
$g'(t)=O(t^{\Theta-1}\log^2(2t))$ by partial summation in the
integrand of \eqref{eq:primeindmargin}. Thus
$g(z)\ge g(x)-O(x^{\Theta-0.475}\log^2x)>0$.
\end{proof}
In particular the Riemann hypothesis is the assertion that
\eqref{eq:primeend} has no failing endpoint beyond $1.6\cdot10^{20}$,
and then none beyond $661/8$. The asymptotic margin is strict.
On RH the upper limit of the normalized deficit is at most
$4c_*+0.280828<1.4525$, by \eqref{eq:fourierlow} and the bound
$32e_H$ on the unverified zeros. Corollary~\ref{th:cycslack} and
\eqref{eq:stripprimes} bound it below by $4c_*+0.12375>1.2953$,
unconditionally.

\head{A uniform use of the finite-field estimate}
Use the fixed curve $E:y^2=x^3-x$ from Sec.~\ref{sec:curve}, now
at every odd prime $p$ of good reduction. Write
$N_{p,r}=\#E(\mathbb F_{p^r})$ and $a_p=p+1-N_{p,1}$.
The elliptic case of Hasse--Weil, or the hyperelliptic theorem in
Stepanov's treatment \cite{step72curve}, gives $|a_p|\le2\sqrt p$.
Any finite lower cutoff on the characteristic in that treatment
can be imposed here. In genus one the usual Hasse bound applies
at every prime of good reduction. Since $p$ is prime, $a_p^2<4p$.
The rationality argument of Sec.~\ref{sec:curve} with $q=p$
now gives $s_{p,0}=2,\ s_{p,1}=a_p$, and
$s_{p,r+2}=a_p s_{p,r+1}-p s_{p,r}$, where
$s_{p,r}=p^r+1-N_{p,r}$ for $r\ge1$. For
$I_p(u,v)=v^2-a_puv+pu^2$, calculate
\[
 I_p(2,a_p)=4p-a_p^2,\qquad
 I_p(v,a_pv-pu)=pI_p(u,v).
\]
By induction $I_p(s_{p,r},s_{p,r+1})=(4p-a_p^2)p^r$ for every
$r\ge0$. Completing the square gives
$I_p(u,v)=(v-a_pu/2)^2+(p-a_p^2/4)u^2$, hence
$|s_{p,r}|\le2p^{r/2}$, simultaneously at every odd prime.

The normalized local point counts give the weighted prime sum
\[
 P_{E,r}(x)=\sum_{p\text{ odd}}(\log p)k(p/x)
                       \frac{N_{p,r}}{p^r+1}.
\]
For every $x\ge2$ and every integer $r\ge1$, the point-count
induction gives the unconditional comparison
\begin{equation}
 |P_{E,r}(x)-P(x)|\le 2x^{-r/2}\vartheta(8x)
     \ll x^{1-r/2}.\label{eq:stepfieldbound}
\end{equation}
The constant is absolute and independent of $r,x$. Indeed
$0\le k\le1$, the primes in its support satisfy $p\ge x$, and
$\vartheta(8x)\ll x$. At $x=2$ the omitted prime has zero weight.
The error in \eqref{eq:stepfieldbound} is respectively of order at
most $\sqrt{x},1,x^{-1/2},x^{-1}$ for $r=1,2,3,4$.
Theorem~\ref{th:stepfirstterm} below gives its first term, including
the sign, and extends the next statement to $r=1$.
For any fixed $r\ge2$, the endpoint criterion in Theorem~\ref{th:stepprime}
can equivalently use $P_{E,r}$ at the same threshold. On RH use the
strict margin in \eqref{eq:stripprimes}. Conversely the comparison
bounds the normalized deficit of $P$ at the endpoints. Concavity
with a slightly larger constant then bounds it on a tail, as in
that theorem. Increasing $r$ controls the difference in
\eqref{eq:stepfieldbound} while the common prime-count error remains.

The common term can be seen exactly. For each odd $p$, choose a
nonsquare $d_p$ and write $E^{(d_p)}:d_py^2=x^3-x$ over
$\mathbb F_p$. At each $x_0\in\mathbb F_p$ the two affine
counts add to two, including at a root of the cubic. Each curve has
one point at infinity. Thus $N_{p,1}+N^{(d_p)}_{p,1}=2(p+1)$.
Averaging the two normalized prime sums recovers $P(x)$ identically.
Stepping in the extension degree at each prime and advancing the
cutoff across new primes are different operations in this identity.

\head{The corresponding prime-return Euler product}
Use $w_{p,r}=N_{p,r}/(p^r+1)$ to weight each primitive channel
and all its repetitions in Sec.~\ref{sec4}. Keep $w_{2,r}=1$.
On $\re s>1$ define the product by its absolutely convergent
logarithm
\[
 \mathcal Z_r(s)=\exp\left(\sum_p\sum_{j\ge1}
                     \frac{w_{p,r}}{j}p^{-js}\right).
\]
Its logarithmic derivative with a minus sign is the weighted return
series $\sum_{p,j}w_{p,r}(\log p)p^{-js}$.
The point-count induction gives the factorization
\begin{equation}
\begin{gathered}
 \mathcal Z_r(s)=\zeta(s)e^{A_r(s)},\\
 A_r(s)=\sum_p\sum_{j\ge1}\frac{w_{p,r}-1}{j}p^{-js},\\
 \re s>\max(0,1-r/2).
\end{gathered}\label{eq:stepdet}
\end{equation}
Indeed on compact subsets with $\re s\ge\sigma>\max(0,1-r/2)$,
the inner absolute sum is bounded by
$2p^{-r/2}p^{-\sigma}/(1-2^{-\sigma})$.
Normal convergence gives a holomorphic $A_r$. Euler's product
first proves the identity on $\re s>1$, then its right-hand side
gives the continuation. For $r\ge2$ this covers the full critical
strip. The exponential has no zeros, so each zeta zero in this
domain has exactly the same multiplicity for $\mathcal Z_r$.
As $r\to\infty$, $A_r\to0$ uniformly on compact subsets of
$\re s>0$ by dominated convergence. Thus taking successively
larger finite-field extensions can suppress their entire correction
while preserving the same possible off-line zeros of the return
product. Proving the pointwise inequality for the full prime sum
would supply independent information about these zeros.

\subsection{The first term of the point-count correction and two deletion tests}
\label{sec:steprefine}
The bound \eqref{eq:stepfieldbound} uses only the size of each
Frobenius trace. Its sign can be followed across the rational primes.
Write
\[
 K(s)=\int_1^8k(v)v^{s-1}\dd v=\frac{(4^s-1)(2^s-2)}{s(s-1)}
\]
for the transform of the weight \eqref{eq:stepweight}, with removable
singularities at $0$ and $1$. Thus
\[
\begin{gathered}
 K(1)=6\log2,\qquad K(\tfrac12)=8c_*,\\
 K(0)=2\log2,\qquad
 K(-1)=\tfrac9{16}.
\end{gathered}
\]
The curve $E$ has complex multiplication by $\mathbb Z[i]$.
If $p\equiv3\bmod4$ then $a_p=0$. If $p\equiv1\bmod4$, write
$p=\pi\bar\pi$ in $\mathbb Z[i]$ with $\pi\equiv1\bmod(2+2i)$.
Then $a_p=\pi+\bar\pi$ \cite{irelandrosen}. Equivalently
$p=a^2+b^2$ with $b$ even and $a+b\equiv1\bmod4$, and $a_p=2a$.
The reciprocal roots at $p$ are $\sqrt p\,e^{\pm i\theta_p}$ with
$0\le\theta_p\le\pi$, where $\theta_p=\pi/2$ at the inert primes.
The recursion preceding \eqref{eq:stepfieldbound} gives
\[
 s_{p,r}=p^{r/2}b_{p,r},\qquad b_{p,r}=2\cos(r\theta_p).
\]
For example $5=(-1+2i)(-1-2i)$ gives $a_5=-2$ and $N_{5,1}=8$,
as in Sec.~\ref{sec:curve}. At an inert prime $N_{p,2}=(p+1)^2$.

\begin{theorem}[First term of the correction at each extension degree]
\label{th:stepfirstterm}
For each fixed integer $r\ge1$ there is $c_r>0$ such that, as
$x\to\infty$ through real values,
\begin{multline}
 P_{E,r}(x)-P(x)=-\cos\!\Big(\frac{\pi r}{2}\Big)\,
   K\!\Big(1-\frac r2\Big)\,x^{1-r/2}\\
   +O_r\!\Big(x^{1-r/2}e^{-c_r\sqrt{\log x}}\Big).
 \label{eq:stepfirstterm}
\end{multline}
In particular
\[
\begin{aligned}
 P_{E,1}(x)-P(x)&=o(\sqrt x),\\
 P_{E,2}(x)-P(x)&=2\log2+o(1),\\
 P_{E,4}(x)-P(x)&=-\frac{9}{16x}+o(x^{-1}).
\end{aligned}
\]
\end{theorem}
\begin{proof}
Every ideal of $\mathbb Z[i]$ prime to $1+i$ has a unique generator
$\alpha\equiv1\bmod(2+2i)$. For $r\ne0$ the assignment
$\xi_r((\alpha))=(\alpha/|\alpha|)^r$ is a Hecke character with
nontrivial infinity type. Its $L$-function is entire and has the
classical zero-free region, whence
$\sum_{N\mathfrak p\le y}\xi_r(\mathfrak p)\log N\mathfrak p
 \ll_r y\,e^{-c_r\sqrt{\log y}}$ \cite{hecke1920,iwanieckowalski}.
The two primes above $p\equiv1\bmod4$ contribute $b_{p,r}\log p$
to this sum. The primes $(p)$ with $p\equiv3\bmod4$ have norm
$p^2$ and contribute $O(\sqrt y)$. At these inert primes
$b_{p,r}=2\cos(\pi r/2)$, and the prime number theorem for the
progression $3\bmod4$ has a remainder of the same form \cite{mv}.
Consequently
\[
\begin{aligned}
 \vartheta_{E,r}(y)&:=\sum_{2<p\le y}(\log p)\,b_{p,r}\\
   &=\cos(\pi r/2)\,y+O_r\big(y\,e^{-c_r\sqrt{\log y}}\big).
\end{aligned}
\]
For $x\ge2$ the prime $2$ has zero weight, and
\[
\begin{gathered}
\begin{aligned}
 P_{E,r}(x)-P(x)&=-\sum_{p>2}(\log p)k(p/x)\frac{s_{p,r}}{p^r+1}\\
   &=-\int_x^{8x}h_x(t)\dd\vartheta_{E,r}(t),
\end{aligned}\\
 h_x(t)=\frac{k(t/x)}{t^{r/2}+t^{-r/2}} .
\end{gathered}
\]
The function $h_x$ is continuous, vanishes at both ends and has
total variation at most $3x^{-r/2}$. Its integral against
$\cos(\pi r/2)\dd t$ equals
$\cos(\pi r/2)K(1-r/2)x^{1-r/2}(1+O(x^{-r}))$. Integration by parts
bounds the remaining part by
$3x^{-r/2}\sup_{x\le t\le8x}|\vartheta_{E,r}(t)-\cos(\pi r/2)t|$.
\end{proof}

The value at $r=2$ is produced by the supersingular primes. For
$p\equiv3\bmod4$ the weight is $w_{p,2}=1+2p/(p^2+1)$, and these
primes have density one half. The split primes contribute $o(1)$
at that degree. At odd $r$ the inert traces vanish and no main term
remains. Table~\ref{tab:stepfirst} lists ball-arithmetic enclosures
from \texttt{anc/stepanov\_check.py}, rounded to six decimals. The
signs in $a_p=2a$ were compared with direct point counts over
$\mathbb F_p$ for $p<600$ and over $\mathbb F_{p^2}$ for
$p\le29$.

\begin{table*}[t]
\centering
\caption{Point-count correction $P_{E,r}(x)-P(x)$ for $E:y^2=x^3-x$.
The last two columns use the normalizations of Theorem~\ref{th:stepfirstterm}.}
\label{tab:stepfirst}
\setlength{\tabcolsep}{8pt}%
\begin{tabular}{ccccc}
\toprule
$x$ & $P_{E,1}-P$ & $P_{E,2}-P$ & $\sqrt x\,(P_{E,3}-P)$ & $x\,(P_{E,4}-P)$\\
\midrule
$10^2$ & $-0.754843$ & $1.387001$ & $\phantom{-}0.028545$ & $-0.533170$\\
$10^3$ & $-0.894239$ & $1.404297$ & $\phantom{-}0.003304$ & $-0.554897$\\
$10^4$ & $\phantom{-}0.034290$ & $1.388950$ & $-0.002595$ & $-0.559654$\\
$10^5$ & $-0.683498$ & $1.387293$ & $-0.003114$ & $-0.561899$\\
$10^6$ & $-1.141681$ & $1.386736$ & $-0.000598$ & $-0.562184$\\
\midrule
limit & & $2\log2=1.386294\ldots$ & $0$ & $-9/16$\\
\bottomrule
\end{tabular}
\end{table*}

\head{A constant at odd degree}
At $r=1$, the unnormalized entries in Table~\ref{tab:stepfirst}
are bounded over the tabulated range. For an asymptotic interpretation,
assume that the nontrivial zeros of the $L$-function of $E$ lie on
its critical line, written $1+i\gamma$ in the arithmetic normalization.
In the explicit
formula for its logarithmic derivative the prime squares carry the
coefficients $b_{p,2}$, of mean $-1$. Removing them gives
$\vartheta_{E,1}(y)=\sqrt y-\sum_\gamma y^{1/2+i\gamma}/(\tfrac12+i\gamma)
+o(\sqrt y)$, the sum being taken symmetrically. Inserting this in
the preceding proof yields
\[
 P_{E,1}(x)-P(x)=-\log2+\sum_\gamma K(i\gamma)\,x^{i\gamma}+o(1),
\]
where $-\log2=-K(0)/2$ and the sum over zeros converges absolutely.
There is no zero at $\gamma=0$, since the central value of this
$L$-function is nonzero. The oscillating sum has mean zero.
This statement is conditional. Numerically the mean of
$P_{E,1}(x)-P(x)$ over $801$ logarithmically spaced points of
$[10^3,10^7]$ is $-0.694$, and all values lie in $[-1.47,0.08]$.

\head{The criterion with point-count weights}
By \eqref{eq:stepfirstterm} the correction is $o(\sqrt x)$ at every
fixed $r\ge1$. The weights $w_{p,r}$ are positive, so the
concavity between endpoints is unchanged. Hence
Theorem~\ref{th:stepprime} holds with $P_{E,r}$ in place of $P$
at each fixed degree, including $r=1$, by the argument following
\eqref{eq:stepfieldbound}. The first term of the difference between
the two inequalities is now known. Neither inequality follows from it.

\head{A square-lattice test at the inert primes}
For $E:y^2=x^3-x$, the good primes $p\equiv3\bmod4$ are inert
in $\mathbb Q(i)$. They give the term $2\log2$ above. The scalar
theta function for this lattice also permits an exact test of the
proposed divisibility. In the normalization of
Ref.~\cite[Part I, Sec.~4.1]{acguide}, put
\[
\begin{gathered}
 \Theta_i(z)=\sum_{\lambda\in\mathbb Z[i]}
             e^{-\pi|\lambda|^2+\lambda z},\qquad
 \omega=2\int_0^1\frac{\dd t}{\sqrt{1-t^4}},\\
 r_k=\frac{\Theta_i^{(4k)}(0)}
             {\Theta_i(0)(\omega/\pi)^{4k}} .
\end{gathered}
\]
Writing $u=\omega z/(2\pi)$, the period-normalized function is
\begin{equation}
 \Theta_i(z)/\Theta_i(0)=\sigma(u)^2\wp(u),
 \label{eq:acsquaretheta}
\end{equation}
where $\sigma,\wp$ have invariants $g_2=4,\ g_3=0$ and full
periods $\omega,i\omega$. To derive the identity, set
$f(z)=\sum_{j\in\mathbb Z}e^{-\pi j^2+jz}$. Then
$\Theta_i(z)=f(z)f(iz)$. Shifting the summation index gives
$f(z+2\pi)=e^{z+\pi}f(z)$ and $f(z+2\pi i)=f(z)$.
Thus translations of $\Theta_i$ by $2\pi,2\pi i$ multiply it
by $e^{z+\pi},e^{-iz+\pi}$, respectively. At $z=\pi(1+i)$
both factors vanish by pairing terms, so the zero has order at
least two. For this square lattice the Weierstrass quasiperiod
constants are $\eta_1=\pi/(2\omega)$ and $\eta_2=-i\eta_1$,
by rotation and Legendre's relation \cite[\S\S23.2--23.5]{dlmf}.
Consequently $\sigma(u)^2\wp(u)$ has the same translation factors.
The only zero of this entire function modulo the lattice is a double
zero at $(1+i)\omega/2$, where $\wp=\wp'=0$.
The quotient of the two sides of \eqref{eq:acsquaretheta} is
therefore entire and doubly periodic. Its value at zero is one,
which proves the claim.

For an exact rational computation, write
$u^2\wp(u)=\sum_{j\ge0}a_j u^{4j}$ with $a_0=1,\ a_1=1/5$.
From $\wp''=6\wp^2-2$ and $(\log\sigma)''=-\wp$ one obtains
\[
\begin{aligned}
 a_j&=\frac{6\sum_{\ell=1}^{j-1}a_\ell a_{j-\ell}}
                  {(4j-2)(4j-3)-12}\qquad(j\ge2),\\
 D(v)&=\exp\left(-2\sum_{j\ge1}\frac{a_j v^j}{4j(4j-1)}\right),\\
 r_k&=\frac{(4k)!}{16^k}[v^k]\left(D(v)\sum_{j\ge0}a_jv^j\right).
\end{aligned}
\]
This gives $r_1=1/4,\ r_2=17/16,\ r_3=-111/64$.
At derivative order $a=120$, the reduced denominator of $r_{30}$
is $2^{60}$ and
\[
 v_3(r_{30})=14,\qquad v_7(r_{30})=2,\qquad v_{11}(r_{30})=0 .
\]
In particular the first value reproduces the computation in
Ref.~\cite{acguide} by rational arithmetic. The unshifted odd-power
main term at $3$ is $120(1/3+1/3^3+\cdots)=45$.
A single order cannot determine an asymptotic error $O(\log a)$.
The external valuation estimate concerns the normalized maps in a
chain of theta lifts, with weights and denominator losses as in
Sec.~\ref{sec:acbridge}.

\head{Deleting primes of density zero}
The following test separates the information in the point counts from
the information in \eqref{eq:primeend}. Let $\mathcal D$ be a set of
primes and put
\[
 \vartheta_{\mathcal D}(y)=\sum_{p\in\mathcal D,\ p\le y}\log p,\qquad
 P_*(x)=\sum_{p\notin\mathcal D}(\log p)k(p/x).
\]
If $\vartheta_{\mathcal D}(y)=o(y)$, the retained primes satisfy the
prime number theorem $\vartheta(y)-\vartheta_{\mathcal D}(y)\sim y$.
Their Frobenius angles have the same limiting distribution, and each
retained prime has the same point counts. Formula
\eqref{eq:stepfirstterm} keeps its leading term for the retained
sums, because the deleted primes change its left side by at most
$2x^{-r/2}[\vartheta_{\mathcal D}(8x)-\vartheta_{\mathcal D}(x)]$.

\begin{proposition}[Two deletion tests]\label{pr:stepdeletion}
Let $p_1<p_2<\cdots$ be the primes.
\begin{enumerate}
\item Fix $1/2<\beta<1$ and an integer $j_0\ge2$, and delete
$\mathcal D=\{p_{\lfloor j^{1/\beta}\rfloor}:j\ge j_0\}$. For every
$C>0$ and $1/2\le a<\beta$ there are arbitrarily large $x$ with
$P_*(x)<6(\log2)x-Cx^a$. Such $x$ occur among the endpoints
$p,p/8$ of the retained primes.
\item Suppose $\vartheta_{\mathcal D}(y)=c\sqrt y+o(\sqrt y)$. If
$c>1/2$, inequality \eqref{eq:primeend} for $P_*$ fails at
arbitrarily large endpoints of the retained primes. If RH holds and
$c<(2c_*-8Q)/(4c_*)=0.1845\ldots$, it holds for $P_*$ on a whole
tail.
\end{enumerate}
The failures in both parts are unconditional.
\end{proposition}
\begin{proof}
Put $A=6\log2$ and $D(s)=\sum_{p\in\mathcal D}(\log p)p^{-s}$.
In part 1 the number of deleted primes up to $y$ is asymptotic to
$\pi(y)^\beta$. Thus $\vartheta_{\mathcal D}(y)=o(y)$, the series
$D$ converges absolutely for $\re s>\beta$, and
$D(\sigma)\to+\infty$ as $\sigma\downarrow\beta$, since
$\int^\infty\dd y/(y\log^\beta y)$ diverges. Suppose
$P_*(x)-Ax+Cx^a\ge0$ for $x\ge X_0\ge1$. For $\re s>1$ its
transform against $x^{-s-1}\dd x$ on $[X_0,\infty)$ is
\[
 K(s)\big[L_p(s)-D(s)\big]-\frac{AX_0^{1-s}}{s-1}
     +\frac{CX_0^{a-s}}{s-a}-H(s).
\]
Here $L_p$ is the prime series of Sec.~\ref{sec:powers}, and $H$
is the transform of $P_*$ on $(1/4,X_0)$, an entire function.
By the continuation of $L_p$ given there, and $K(1)=A$, the
displayed expression is regular at every real $s>\beta$. Landau's
theorem, as in Lemma~\ref{lem:explicit}, places the abscissa of
convergence at most at $\beta$. The integral is therefore
nonnegative at each real $\sigma>\beta$. As $\sigma\downarrow\beta$
all terms stay bounded except $-K(\sigma)D(\sigma)$, which tends to
$-\infty$ because $K(\beta)>0$. This contradiction proves the
first assertion. The function $P_*(x)-Ax+Cx^a$ is concave between
successive endpoints of the retained primes, as in
Theorem~\ref{th:stepprime}. A failure at arbitrarily large $x$
therefore occurs at arbitrarily large endpoints.

In part 2 integrate by parts against $\vartheta_{\mathcal D}$. The
weight has total variation two, so the deleted primes contribute
\[
\begin{aligned}
 \sum_{p\in\mathcal D}(\log p)k(p/x)
   &=c\sqrt x\int_1^8k(v)\dd\sqrt v+o(\sqrt x)\\
   &=\tfrac12cK(\tfrac12)\sqrt x+o(\sqrt x)\\
   &=4cc_*\sqrt x+o(\sqrt x).
\end{aligned}
\]
If RH fails, the normalized deficit of $P$ is unbounded above.
Otherwise the converse part of Theorem~\ref{th:stepprime}, whose proof
allows any constant in place of $6c_*$, would give RH. As
$P_*\le P$, the same holds for $P_*$. Assume RH and put
$S(u)=\sum_\rho K(\rho)e^{i(\im\rho)u}$. By \eqref{eq:stripprimes}
and the explicit formula used for Theorem~\ref{th:stependrange},
\[
 \frac{6(\log2)x-P_*(x)}{\sqrt x}=S(\log x)+4c_*(1+c)+o(1).
\]
The series for $S$ converges absolutely and has no constant term.
Its mean over $[0,U]$ tends to zero, so $S(u)\ge-\varepsilon$ on
an unbounded set for every $\varepsilon>0$. The upper limit of the
left side is at least $4c_*(1+c)$, which exceeds $6c_*$ when
$c>1/2$. Concavity transfers the failures to endpoints. Finally
$|S|\le8Q$ gives the upper bound $8Q+4c_*(1+c)+o(1)$, which is
below $6c_*$ for the stated range of $c$.
\end{proof}

In part 2 the retained set differs from the primes by about
$c\sqrt y/\log y$ elements below $y$, and its Chebyshev function
by $c\sqrt y$. An estimate $\vartheta(y)-y=O(\sqrt y)$ with an
unspecified constant cannot distinguish the two sets. Inequality
\eqref{eq:primeend} does. Its constant $6c_*$ allows for the prime
squares, which contribute $4c_*\sqrt x$, and for a further
$2c_*\sqrt x$. For an example with $c=1$, delete the least prime
exceeding $q^2$ for every second prime $q=3,7,13,\ldots$.
These primes are distinct for large $q$ and equal $q^2(1+o(1))$
\cite{bhp01}. Below $2^{27}$ this deletes $696$ primes, and
$\vartheta_{\mathcal D}(2^{27})/2^{27/2}=0.985\ldots$. The normalized
deficit of the retained set lies between $2.24$ and $2.66$ at
every one of its endpoints in $[2^{10},2^{24}]$, above
$6c_*=1.7574$. For part 1 with $\beta=3/4$ and $j_0=2$ the same
deficits lie between $22.6$ and $302.4$. In either test the first
deletion can be postponed beyond any prescribed cutoff, in particular
beyond the range of Theorem~\ref{th:stependrange}.

\head{Arithmetic information retained under prime deletion}
Both deletion constructions preserve the Euler product over retained
primes, the positive return measure of Sec.~\ref{sec4}, the first-order
prime number theorem, the limiting Frobenius distribution, and the
point count of $E$ at every retained prime.
Theorem~\ref{th:stepfirstterm} continues to hold. The change appears
in the counting function of the generated integers. If $N_*(x)$
counts integers up to $x$ with no prime factor in $\mathcal D$, then
$N_*(x)=\kappa x+O_\varepsilon(x^{\theta+\varepsilon})$ with
$\kappa=\prod_{p\in\mathcal D}(1-1/p)>0$, where $\theta=\beta$ in
part 1 and $\theta=1/2$ in part 2. For the full set of primes the
count is $\lfloor x\rfloor$. Diamond, Montgomery and Vorhauer
constructed generalized primes whose integers satisfy
$N(x)=\kappa x+O(x^\theta)$ with $1/2<\theta<1$ and whose zeta
function has infinitely many zeros on a curve
$\re s=1-a/\log|\im s|$ \cite{dmvbeurling}. The deletion results
therefore distinguish the endpoint estimate \eqref{eq:primeend} for
the rational primes from estimates determined by the retained local
data. In Lemma~\ref{lem:explicit}, the functional equation supplies
the reflection $\rho\mapsto1-\rho$, the constant $d_0$, and the
series $S(x)$. The point counts and the limiting Frobenius law
remain valid after deletion. By \eqref{eq:stepdet}, the return
product retains its zeros in $\re s>1/2$. Each such zero contributes
positively to the defect $\Delta$ of Theorem~\ref{th:boundaryclosure}.

\subsection{Induction on the multiplicity of a logarithmic average}\label{sec:stepmult}
In the return clock $\log x$ a positive averaging operator can be
iterated. The integer $k$ below counts the number of independent
observation delays. Its mathematical role is to raise a Mellin
multiplier to its $k$-th power. Put
\[
\begin{gathered}
 A=\sqrt2,\qquad B=2\sqrt2,\qquad L=\log2,\\
 (\mathcal M f)(x)=\frac1L\int_A^B f(tx)\frac{\dd t}{t},\\
 M(s)=\frac{B^s-A^s}{Ls},\quad M(0)=1 .
\end{gathered}
\]
A single application integrates a detector reading over a uniform
interval of logarithmic delays. Repeating it models $k$ successive
applications of the same temporal response. Here $\mathcal M$ and $M$ are the averaging operator and its
multiplier. The mean delay is $\log 2$, and $M(s)^k$ has a zero
of order $k$ at each $2\pi i j/L$, $j\in\mathbb Z\setminus\{0\}$.
This is analytic multiplicity in the time response. The count of prime
powers in $U$ is the same as in \eqref{eq:blockreport}.

\begin{theorem}[An all-orders induction and integer cutoffs]\label{th:stepmult}
With $d_0=\log(2\pi)$, for every $k\ge0$ and every real $x\ge1$,
\begin{multline}
 \left|(\mathcal M^k U)(x)-\frac{d_0}{2}\right|
 \le \frac{x}{10}\left(\frac{19}{20}\right)^k\\
 +(1-\log2)x^{-2}\left(\frac{3}{16\log2}\right)^k .
 \label{eq:stepcontract}
\end{multline}
For each $k$ there is a real $y_k\in[A^k,B^k]$ with
$U(y_k)=d_0/2+O((19/20)^k)$. There are integers $n_k\to\infty$
with $U(n_k)\to d_0/2$. With $h$ as in Sec.~\ref{sec:steppair},
there is an absolute $c>0$ such that these cutoffs satisfy
\begin{equation}
 \sup_{|t-n_k|\le\sqrt{n_k}\exp(c h(n_k))}|b(t)|
   \ll \exp(-c h(n_k)).\label{eq:stepwindows}
\end{equation}
None of these assertions uses a verified list of zero ordinates.
\end{theorem}
\begin{proof}
The base $k=0$ is the explicit formula of Lemma~\ref{lem:explicit}
inserted into \eqref{eq:steptransfer}. It gives the absolutely
convergent expansion
\[
 U(x)-d_0/2=\sum_\rho c_\rho(1-2^{\rho-1})x^\rho
 -\sum_{j\ge1}\frac{1-2^{-2j-1}}{2j(2j+1)}x^{-2j}.
\]
As the induction hypothesis take this same identity with each term
of exponent $s$ multiplied by $M(s)^k$ on the right and
$\mathcal M^kU-d_0/2$ on the left. At every finite stage the
series converge uniformly on the compact interval of integration.
Applying $\mathcal M$ sends each $x^s$ to $M(s)x^s$ and
preserves the constant. This proves the identity for $k+1$.
To bound it, Sec.~\ref{sec:step} gives analytically
$|\im\rho|\ge\sqrt{2/Q-1}>13/2$, using the value of $Q$ in
\eqref{eq:center}. Thus for $0<\re(\rho)<1$
\[
 |M(\rho)|\le\frac{A+B}{(13/2)L}<0.942<19/20.
\]
Also $\sum|c_\rho|<0.05$, $|1-2^{\rho-1}|\le2$, and
$0<M(-2j)\le M(-2)=3/(16L)$. The trivial-zero weights before
multiplication by $1-2^{-2j-1}$ sum to $1-\log2$. This proves
\eqref{eq:stepcontract} with constants independent of $k,x$.
For example, for any fixed $\eta>0$, choosing
$k\ge(1+\eta)\log x/\log(20/19)$ gives
$\mathcal M^kU(x)=d_0/2+O(x^{-\eta})$ for all $x\ge1$.

At $x=1$ the left-hand average integrates the continuous function
$U$ against a probability on $[A^k,B^k]$. Its value lies between
the minimum and maximum there, giving $y_k$ by continuity.
To pass to integers, partial summation on $(x,2x]$ and the classical
prime number theorem $\psi(x)=x+O(xe^{-a\sqrt{\log x}})$
\cite{mv} give
\[
\begin{aligned}
 \sum_{x<m\le2x}\frac{\Lambda(m)}m
 &=\frac{\psi(2x)}{2x}-\frac{\psi(x)}x+
   \int_x^{2x}\frac{\psi(v)}{v^2}\dd v\\
 &=\log2+O(e^{-a\sqrt{\log x}}).
\end{aligned}
\]
Equation~\eqref{eq:stepslope} bounds the slope almost everywhere,
and $U$ is locally absolutely continuous. Rounding $y_k$ to the
nearest integer costs $o(1)$, proving the assertion about $n_k$.
The Vinogradov--Korobov estimate in Sec.~\ref{sec:steppair}, in
the same partial summation formula, gives the sharper slope bound
$O(e^{-2c h(n_k)})$ on $[n_k/2,2n_k]$ for large $k$ and
some $c>0$, since $h(n_k/2)/h(n_k)\to1$. Integrate from $n_k$, use
$U(n_k)=O(1)$ and divide by $\sqrt t$ to obtain
\eqref{eq:stepwindows}.
\end{proof}

\head{The uniform square-root scale}
For the unsmoothed quantity one has, unconditionally,
\begin{equation}
 \limsup_{x\to\infty} \frac{U(x)}{\sqrt x}
 \ \ge\ Q-2^{-1/2}C(\log2)>0.05844 .\label{eq:steplimsup}
\end{equation}
If RH fails, a finite upper limit would bound $b$ on a tail and
imply RH by the telescoping argument of Sec.~\ref{sec:steppair}.
On RH the zero terms in $b(e^u)$ are the absolutely convergent
series $\sum c_\rho(1-2^{\rho-1})e^{i\gamma u}$, up to a term
tending to zero. Simultaneous recurrence for finitely many phases,
followed by control of the remaining absolute sum as in
Theorem~\ref{th:sharp}, gives a sequence $u\to\infty$ on which
the sum approaches its value at $u=0$. This value is precisely
the right side before the numerical comparison, by \eqref{eq:C},
\eqref{eq:Aq} and Table~\ref{tab:q}. Thus a pointwise upper bound $o(\sqrt{x})$
for $U$ on an entire tail would contradict \eqref{eq:steplimsup}.

At every \emph{fixed} $k$, an eventual one-sided bound
$\mathcal M^kU(x)\le C_k\sqrt x$ is itself equivalent to RH.
On RH this follows by averaging the explicit formula. Conversely,
the Mellin transform on $[1,\infty)$ of $\mathcal M^kU$ is, up
to an entire term,
$M(s)^k(1-2^{s-1})\int_1^\infty R(x)x^{-s-1}\dd x$.
The omitted integrals have fixed, compact positive ranges for fixed
$k$. The multiplier has no zero in $1/2<\re s<1$.
The nonnegative-tail argument of Lemma~\ref{lem:explicit} applies.

The loss of pointwise control can also be checked on a positive
counting measure. For sufficiently large $t$ prescribe
$\Psi(t)=t-d_0+\re(d\,t^{s_*})$, with
$s_*=3/4+2\pi i/L$ and $0<|d s_*|<1/2$.
Then $\Psi'(t)>1/2$ for $t\ge1$ wherever this formula is used,
and $\Psi(t)/t\to1$. It has an increasing nonnegative extension
to smaller positive arguments. The identical block functional yields
\[
\begin{gathered}
 U_\Psi(x)=x\int_x^{2x}(t-\Psi(t))\,t^{-2}\dd t
 =d_0/2+\re(D x^{s_*}),\\
 D=-d\frac{2^{s_*-1}-1}{s_*-1}.
\end{gathered}
\]
Choose the phase of $d$ so that $D>0$. Since $A>1$, for large
$k$ the whole integration interval $[A^k,B^k]$ is above the
cutoff of this formula. At base 1, $\mathcal M^k U_\Psi$ tends to
$d_0/2$ since $|M(s_*)|<19/20$, while
$(U_\Psi(2^k)-d_0/2)/\sqrt{2^k}=D\,2^{k/4}$.
Thus this positive measure satisfies the averaged estimate and
exhibits growth at individual centers. For the von Mangoldt measure,
a sufficient arithmetic condition is the joint lower bound
\eqref{eq:pairinput} at all sufficiently large cutoffs, or a bound
on the accumulated margins in \eqref{eq:steptransfer}.

\section{Jacobians, cyclotomic unit cohomology and the prime sum}
\label{sec:coh}
\subsection{Jacobian point counts and Frobenius eigenvalues}\label{sec:jacbound}
Let $C/\mathbb F_q$ be smooth, projective and geometrically
connected of genus $d\ge1$, and $J=\operatorname{Pic}^0(C)$.
Write $\alpha_1,\ldots,\alpha_{2d}$ for the reciprocal roots of the
numerator of $Z_C(T)$, with multiplicities. The trace formula,
duality and the determinant formula for the Jacobian give
\begin{align}
 \#C(\mathbb F_{q^n})&=1+q^n-\sum_i\alpha_i^n,\qquad
 \prod_i\alpha_i=q^d,\nonumber\\
 J_n:=\#J(\mathbb F_{q^n})&=\prod_i(1-\alpha_i^n).
 \label{eq:jacnumbers}
\end{align}
The parameter $n$ is the degree of the constant-field extension.
Each step replaces the eigenvalues of Frobenius by the corresponding
powers $\alpha_i^n$ in the point counts. Thus $n$ is an arithmetic
iteration parameter independent of the observer's clock. The
$\alpha_i$ are the eigenvalues of geometric Frobenius on
$H^1_{\mathrm{\acute et}}(C_{\overline{\mathbb F}_q},\mathbb Q_\ell)$,
$\ell\nmid q$, or equivalently of the Frobenius endomorphism on
the rational Tate module of $J$. The rational points are the kernel
of the separable isogeny $1-F^n$. Lang's theorem gives
$H^1(\mathbb F_{q^n},J)=0$ for arithmetic cohomology. The determinant
in \eqref{eq:jacnumbers} uses geometric cohomology
\cite{kedcurve,kedabel,langfields}. In particular, for
$k_n=\mathbb F_{q^n}$ and $\ell\nmid q$,
\[
 J(k_n)[\ell^\infty]\simeq
 H^1_{\rm cts}(\operatorname{Gal}(\overline{k}_n/k_n),T_\ell J)
 \simeq T_\ell J/(F^n-1)T_\ell J .
\]
Arithmetic Frobenius acts on this Tate module as $F^n$.
The last equality is the continuous procyclic cohomology formula.
For the first, apply the cohomology sequence to
$0\to T_\ell J\to V_\ell J\to J[\ell^\infty]\to0$.
$F^n-1$ is invertible on $V_\ell J$ because the determinant in
\eqref{eq:jacnumbers} is nonzero. This computes the
$\ell$-primary part. The integer determinant computes the full order.

\begin{theorem}\label{th:jactwist}
Each of the following conditions is equivalent to the assertion
that $|\alpha_i|=\sqrt q$ for every $i$.
\begin{enumerate}
\item $J_n=q^{dn}+O(q^{(d-1/2)n})$ as $n\to\infty$.
\item If $J_n^-$ denotes the number of rational points of the
quadratic twist by $[-1]$ over $\mathbb F_{q^n}$, then
$J_n^-\le q^{dn}+O(q^{(d-1/2)n})$, with one constant for all $n$.
\end{enumerate}
In condition 2 the twist is taken anew over each extension field.
\end{theorem}
\begin{proof}
Put $z_i=\alpha_i^{-1}$, $R=\max_i|z_i|$.
Dividing the two determinants by $q^{dn}$ gives respectively
$\prod_i(1-z_i^n)$ and $\prod_i(1+z_i^n)$, since the twisted
Frobenius is $-F^n$. If $R=q^{-1/2}$ and all moduli are equal,
expanding either product bounds its difference from one by
$(1+q^{-n/2})^{2d}-1\le(2^{2d}-1)q^{-n/2}$.

For the converse in condition 2 choose arbitrarily large $n$ for
which every phase of $z_i^n$ has real part at least $1/2$ after
normalization to modulus one. Such simultaneous returns follow by
compactness of the finite torus (or simultaneous Dirichlet
approximation). Since the determinant is positive, its normalized
value is at least $\prod_i(1+|z_i|^n/2)\ge1+R^n/2$.
The upper bound forces $R\le q^{-1/2}$. The product of the
$|z_i|$ equals $q^{-d}$, so all have modulus $q^{-1/2}$.

For condition 1 the normalized product tends to one. No $z_i=1$,
since $J_1\ge1$. If $R>1$, take integers $n=m+1$ with all
phases at time $m$ tending to one. The factors of modulus-one
roots then tend in modulus to $|1-z_i|>0$, while factors with
$|z_i|>1$ grow without bound, a contradiction. If $R=1$,
take the return times $m$ themselves. A unit-modulus factor tends
to zero and all other factors stay bounded. Hence $R<1$.
Now $1-\prod_i(1-z_i^n)=\sum_i z_i^n+O_d(R^{2n})$.
At simultaneous returns the real part of the sum is at least
$R^n/2$. The two-sided estimate again gives $R\le q^{-1/2}$.
\end{proof}

For a fixed curve even the lower bound for the untwisted $J_n$
suffices. Indeed a finite morphism $C\to\mathbb P^1$ bounds
$\#C(\mathbb F_{q^n})\ll_C q^n$. The generating function
$\sum_i\alpha_i T/(1-\alpha_i T)$ and \eqref{eq:jacnumbers}
therefore give $\max|\alpha_i|\le q$. Duality pairs $\alpha_i$
with $q/\alpha_i$, so $R\le1$. If $R=1$, the return times
used above contradict the normalized lower bound $1-O(q^{-n/2})$.
For $R<1$ the same expansion and positive real part at return
times contradict that bound whenever $R>q^{-1/2}$.

An upper bound for the untwisted determinant alone does not have
this formal implication. The reciprocal polynomial $1-5T+5T^2$
in genus one would give $\alpha=(5+\sqrt5)/2,\ \beta=(5-\sqrt5)/2$.
Both exceed one, $\alpha\beta=5$, and the positive integers
$(1-\alpha^n)(1-\beta^n)=5^n+1-\alpha^n-\beta^n$ are below $5^n$
for every $n$. The reciprocal roots have unequal moduli, so this
polynomial cannot be the numerator for an elliptic curve over
$\mathbb F_5$. For hyperelliptic curves
the hyperelliptic involution induces $[-1]$ on the Jacobian, so
the required twists are Jacobians of curve twists. For an arbitrary
curve the abelian-variety twist is available. Its descent as a Jacobian
with its principal polarization requires care \cite{serreobstruct}.

\begin{theorem}\label{th:jacfamily}
Fix a finite field $k_0$. Suppose each smooth projective geometrically
connected curve $Y/k$ over a finite extension $k/k_0$, of genus $h$,
satisfies the upper bound
\[
 \#\operatorname{Jac}(Y)(k_n)\le |k|^{hn}
       +C_{Y/k}|k|^{(h-1/2)n},\qquad k_n=\mathbb F_{|k|^n},
\]
for all sufficiently large $n$. Then RH holds for every curve
over $k_0$. For each given curve, the argument uses only finitely
many fixed curves over a single finite extension, with a constant
and threshold allowed to depend on each curve.
\end{theorem}
\begin{proof}
Let $C/k_0$ have positive genus. Take a separable map
$C_{\bar k_0}\to\mathbb P^1$ and its smooth projective Galois
closure $\pi:D\to\mathbb P^1$. The cover, its deck group $G$,
and the map $D\to C$ descend to a finite extension. On the
rational Tate module of $J_D=\operatorname{Jac}(D)$ one has
\[
 \sum_{g\in G}g=\pi^*\pi_*=0.
\]
Indeed on divisors the sum counts each point in a fibre with its
ramification multiplicity, and $\operatorname{Pic}^0(\mathbb P^1)=0$.
This is the trace identity behind the Galois-cover argument for
curve points \cite{bombstep,kedcurve,pacheco}.

Over any field of definition of size $Q$ the reciprocal roots
$\alpha_i$ of $D$ satisfy $1<|\alpha_i|<Q$ without using RH.
In fact, on integers tending to infinity with all phases returning
to one, the trace formula and $\#D(\mathbb F_{Q^n})\ge0$
exclude a modulus greater than $Q$, and permit at most one root
of modulus $Q$, counted with multiplicity. A single such root
must be real, hence $Q$ or $-Q$, by conjugation. Duality would
give $1$ or $-1$ and hence $J_D(\mathbb F_{Q^2})$ of order zero.
Duality now also gives the strict lower bound.

Pass once more to a finite extension so that every normalized
Frobenius eigenvalue has real part greater than $3/4$, by
simultaneous phase recurrence. Write $k$ for this field, $Q=|k|$,
$F$ for its Frobenius, $h$ for the genus of $D$, and
$R=\max_i|\alpha_i|^{-1}<1$ for this new Frobenius.
Put $e=\exp G$. There are arbitrarily large $n=ej+1$ for which
all normalized eigenvalues of $F^{ej}$ approach one. For these $n$,
$\re\operatorname{Tr}(F^{-n})\ge R^n/2$ eventually.
If $\tau_g(n)=\operatorname{Tr}(g^{-1}F^{-n})$, then
$\sum_g\tau_g(n)=0$. Some $g\ne1$ therefore satisfies
$\re\tau_g(n)\le-R^n/(2(|G|-1))$.

For each $g$ take the curve twist $D_g/k$ with Frobenius $gF$.
It exists by finite Galois descent, since $g$ is defined over $k$
and has finite order. For $n\equiv1\bmod e$, its Frobenius over
$k_n$ is $gF^n$. The determinant formula gives
\[
 \frac{\#\operatorname{Jac}(D_g)(k_n)}{Q^{hn}}
   =\det(1-g^{-1}F^{-n})
   =1-\tau_g(n)+O_h(R^{2n}).
\]
Here $g$ commutes with $F$ and its eigenvalues have modulus one.
The assumed upper bound, maximized over this finite list, gives
$R^n/(2(|G|-1))\le C Q^{-n/2}+O_h(R^{2n})$.
Division by $R^n$ rules out $R>Q^{-1/2}$. The product of the
root moduli then forces equality for every root. Finally pullback
from $J_C$ to $J_D$ is injective on rational Tate modules,
since pushforward composed with pullback is multiplication by
$\deg(D\to C)$. Taking the appropriate roots of the Frobenius
moduli descends RH from $k$ to $k_0$.
\end{proof}
The same proof for counts on $D_g$ uses
$\#D_g(k_n)=1+Q^n-\operatorname{Tr}(gF^n)$ directly. At the
chosen return times an upper bound $1+Q^n+O(Q^{n/2})$ for the
twists forces $\max_i|\alpha_i|\le\sqrt Q$ by the same zero-sum
identity. For an example of order three, the smooth cubic covers
$y^3=d\,x(x-1)$ for $d=1,2,3$ over $\mathbb F_7$ have
respectively $9,3,12$ points, including the point at infinity.
The corresponding traces are $-1,5,-4$, whose sum is zero.

There is an independent geometric proof of both estimates.
Use the canonical polarization of $J$, its Rosati involution
$\dagger$, and the positive trace pairing
$(u,v)=\operatorname{Tr}(uv^\dagger)$ on the real endomorphism
algebra over the base field. Positivity follows from the intersection
inequality for an ample polarization \cite{kedabel}. Frobenius and
Verschiebung satisfy $F^\dagger F=q$. Cyclicity of trace gives
$(Fu,Fv)=q(u,v)$. Thus left multiplication by $F/\sqrt q$ is
orthogonal. Any polynomial annihilating left multiplication by $F$
annihilates $F$, by evaluation at the identity, and therefore
annihilates its action on the Tate module. Every $\alpha_i$ has
modulus $\sqrt q$. Substitution in \eqref{eq:jacnumbers} proves
both point bounds for every extension, using geometric positivity.
In particular $(q^{n/2}-1)^{2d}\le J_n\le(q^{n/2}+1)^{2d}$.

\subsection{A tower of cyclotomic subfields and its unit groups}
\label{sec:cyclounits}
To give a definite arithmetic meaning to a cohomological bound,
specify the fields, coefficients, maps and loss allowed in the count.
For odd $p$ let $K_{p^a}$ be the unique degree-$p^a$ subfield
of $\mathbb Q(\zeta_{p^{a+1}})$. For $p=2$ take the real
subfield of $\mathbb Q(\zeta_{2^{a+2}})$. Set $K_1=\mathbb Q$
and let $K_n$ be the compositum of the $K_{p^a}$ for $p^a\Vert n$.
The component extensions are cyclic with pairwise coprime degrees.
Consequently $G_n=\operatorname{Gal}(K_n/\mathbb Q)$ is cyclic of
order $n$, and $d\mid n$ gives $K_d\subset K_n$.
We use ordinary cohomology of these finite Galois groups with their
natural action on the units. Define
\begin{equation}
\begin{gathered}
 A_n=H^1(G_n,\mathcal O_{K_n}^{\times}),\\
 I_n=\sum_{p\mid n}\operatorname{im}(A_{n/p}\longrightarrow A_n),
 \qquad B_n=A_n/I_n .
\end{gathered}
 \label{eq:cyclogroups}
\end{equation}
The arrows are inflation followed by inclusion of units.

\begin{proposition}\label{th:cycunit}
For every $n\ge2$, $A_n\simeq\mathbb Z/n\mathbb Z$ and
\begin{equation}
 \log|B_n|=\Lambda(n),\qquad \log|I_n|=\log n-\Lambda(n).
 \label{eq:cycprimitive}
\end{equation}
\end{proposition}
\begin{proof}
Write $P_K$ for the principal fractional ideals of $K=K_n$.
The sequence of units, nonzero field elements and principal ideals,
and Hilbert's theorem 90, give
$A_n=P_K^{G_n}/\{(a)\mathcal O_K:a\in\mathbb Q^\times\}$
\cite{chabertunits}. For a rational prime $p$, the invariant ideal
group has one generator $\Pi_p=\prod_{\mathfrak p\mid p}\mathfrak p$.
For $p\nmid n$, $\Pi_p=(p)$. For $p^a\Vert n$, $p$ is
totally ramified of index $p^a$ in $K_{p^a}$.
Its unique prime there is generated by the norm down of
$1-\zeta_{p^{a+1}}$ from the full cyclotomic field, using
$1-\zeta_{2^{a+2}}$ for $p=2$. The residue degree is one. The
other factors of the compositum are unramified at $p$, and
their base change adds no ramification to $K/K_{p^a}$.
Extending this principal prime to $K$ therefore gives $\Pi_p$.
These assertions use the elementary prime-power ramification law
for cyclotomic fields \cite{cyclonotes}.

Thus all invariant ideals are principal and
$(p)=\Pi_p^{e_p}$, with $e_p=p^a$ for $p^a\Vert n$.
Unique ideal factorization gives
$A_n\simeq\bigoplus_{p^a\Vert n}\mathbb Z/p^a\mathbb Z$.
Ideal extension from $K_d$ to $K_n$ sends the generator over
$p\mid d$ to $e_p(K_n)/e_p(K_d)$ times the generator. This is
the map on \eqref{eq:cyclogroups} by naturality of the connecting
homomorphism. The image has order $d$. In a cyclic group of order
$n$ the subgroups of orders $n/p$ generate the subgroup of order
$\operatorname{lcm}_{p\mid n}(n/p)$. Its index is $p$ when
$n=p^a$, and one when $n$ has at least two distinct prime
factors. This proves \eqref{eq:cycprimitive}.
\end{proof}
At level $8$ the image of $A_4$ has index two. At level $12$
the subgroups of orders $6$ and $4$ generate $A_{12}$, so
$B_{12}$ is trivial. More generally,
$\prod_{n\le N}|I_n|=N!/\operatorname{lcm}(1,\ldots,N)$, with
trivial groups at $n=1$. The degree of cohomology, module and
subfield maps are essential to the formula.

The trivial representation accounts for the difference between the
number-field sign calculation and Theorem~\ref{th:jacfamily}. For an
abelian number field $K/\mathbb Q$ with group $G$, Artin
factorization gives
\[
 \zeta_K(s)=\zeta(s)\prod_{1\ne\chi\in\widehat G}L(s,\chi^*),
\]
where each $\chi^*$ is the associated primitive Dirichlet character
\cite{washington}. On the trivial summand of the regular representation
of $G$, $\sum_{g\in G}g$ acts as $|G|$. The Artin $L$-function
of this invariant summand is precisely $\zeta(s)$.
Over $\mathbb P^1$ the invariant part in degree one vanishes.
In the number-field tower $K_n$, control of the trivial character
amounts to an estimate on the weighted orders of the $I_n$.

\subsection{A one-sided criterion in terms of cyclotomic indices}
\label{sec:cyclowindow}
Let $k$ be the weight in \eqref{eq:stepweight}, $A=6\log2$,
$V(x)=\sum_{n\ge1}\Lambda(n)k(n/x)$ and $E(x)=V(x)-Ax$.
Its transform is $K(s)=\int_1^8k(v)v^{s-1}\dd v
=(4^s-1)(2^s-2)/[s(s-1)]$, with removable values at 0 and 1.
Removing the squares and higher powers as in \eqref{eq:stripprimes}
gives $g(x)=E(x)+2c_*\sqrt x+o(\sqrt x)$.
For $t\ge2$, set $m=\lceil32t\rceil$, $b_n=\lfloor m k(n/t)\rfloor$,
and take the finite direct sums
$\mathcal A(t)=\bigoplus_n A_n^{\,b_n}$,
$\mathcal I(t)=\bigoplus_n I_n^{\,b_n}$.
Proposition~\ref{th:cycunit} and Chebyshev's bound give, uniformly,
\begin{equation}
 0\le V(t)-\frac1m\log\frac{|\mathcal A(t)|}{|\mathcal I(t)|}
 <\frac{\psi(8t)}m\ll1 .
 \label{eq:cycround}
\end{equation}
Here and below a bound on the order of a group is expressed on its
natural-logarithm scale. With any positive integer multiplicity $m$
the rounding error is $O(t/m)$, hence $o(\sqrt t)$ whenever
$m/\sqrt t\to\infty$.

\begin{theorem}\label{th:cycwindow}
RH is equivalent to the following one-sided bound. There are
$C>0$ and $X_0$ such that every real $X\ge X_0$ admits a
$t\in[X,X+X^{3/4}]$ satisfying
\begin{equation}
 \log|\mathcal I(t)|\le\log|\mathcal A(t)|-mAt+Cm\sqrt t .
 \label{eq:cyctarget}
\end{equation}
\end{theorem}
\begin{proof}
From \eqref{eq:cycround} the bound gives $E(t)\ge-C\sqrt t$.
For $\re s>1$ the full Mellin transform of $V$ on the positive
axis is $-K(s)\zeta'(s)/\zeta(s)$. Its part below 1 is entire,
since $V(x)=0$ for $0<x\le1/4$. Subtracting $Ax$ on
$[1,\infty)$ cancels the real pole at 1 because $K(1)=A$.
Each zero $\rho$ with $1/2<\re(\rho)<1$ gives a pole of residue
$-\operatorname{mult}(\rho)K(\rho)\ne0$. There are no real
singularities above $1/2$. If $E(x)\ge-Dx^\sigma\log^j x$
eventually, with $\sigma\ge1/2$ and $j=0$ or $2$, add the
majorant to $E$ on that tail. The added Mellin transform is
regular at every real point greater than $\sigma$. Landau's
theorem for a nonnegative tail therefore excludes any zero with
$\re(\rho)>\sigma$, as in Lemma~\ref{lem:explicit}.

Off its locally finite set of corners,
\[
 E'(x)=8\sum_{4x<n\le8x}\frac{\Lambda(n)}n
       -2\sum_{x<n\le2x}\frac{\Lambda(n)}n-A .
\]
Chebyshev's estimate gives $E'=O(1)$. Integration back from the
point $t$ first gives $E(X)\ge-C_1X^{3/4}$, excluding zeros to
the right of $3/4$. The truncated von Mangoldt formula now gives
$\psi(u)=u+O(u^{3/4}\log^2(2u))$ \cite{mv}. Partial summation
gives $E'(x)=O(x^{-1/4}\log^2(2x))$. Integration across the
same window gives $E(X)\ge-C_2\sqrt X\log^2(2X)$ for a fixed
constant $C_2>0$. The second
application of Landau and the functional equation give RH.
Conversely under RH, the explicit formula obtained by applying
\eqref{eq:steptransfer} to \eqref{eq:explicit} gives
$E(t)=-\sum_\rho K(\rho)t^\rho+O(t^{-2})=O(\sqrt t)$.
The sum of $|K(\rho)|$ converges since $K(\rho)=O(1/(1+|\im\rho|^2))$.
Equation~\eqref{eq:cycround} then gives \eqref{eq:cyctarget} at all
sufficiently large $t$.
\end{proof}
The equivalence holds with $C=0.370$. Under RH, set
$v=\cos((\im\rho)\log2)$. The identity
\[
 64-|(4^\rho-1)(2^\rho-2)|^2
   =(1+2\sqrt2\,v)^2(10-4\sqrt2\,v)\ge0
\]
and Lemma~\ref{lem:explicit} give
$\sum_\rho|K(\rho)|\le8Q=8(2+\eul-\log(4\pi))<0.370$.
The bounded rounding error divided by $\sqrt t$ tends to zero.

An upper bound for these same groups can be proved independently
at the prime number theorem scale. Put
$h(t)=(\log t)^{3/5}(\log\log t)^{-1/5}$. Integration by parts
of $V(t)$ against $\psi(u)-u$ on $[t,8t]$ bounds
$|E(t)|$ by twice the supremum of that error there, since $k$
has total variation two. The unconditional Vinogradov--Korobov
estimate \cite{johnstonyang} and \eqref{eq:cycround} give for some
$c>0$ and every sufficiently large $t$
\begin{equation}
 \log|\mathcal I(t)|=\log|\mathcal A(t)|-mAt
              +O(mt e^{-c h(t)}).
 \label{eq:cycpnt}
\end{equation}
Formula \eqref{eq:cycprimitive} proves the local orders exactly.
In \eqref{eq:cyctarget} the required additional information is the
distribution of the levels at which the quotient has order $p$.
The error allowed by \eqref{eq:cycpnt}, divided by $m\sqrt t$,
tends to infinity.

The critical-scale estimate itself has a certified finite range.
For positive integers $N$, Legendre's formula for the central
binomial coefficient gives
$\psi(2N)-\psi(N)\le\log\binom{2N}{N}\le2N\log2$.
Summing at powers of two and applying monotonicity gives
$\psi(u)\le4u\log2$ for $u\ge1$. Hence $\psi(8t)/m\le\log2$.
By \eqref{eq:steptransfer}, $E(t)=-4\sqrt t\,I_2(t)$.
The proof of Theorem~\ref{th:steppair} gives $I_2(t)<0.0309395$
on $[64,65536]$ from prime powers alone, and $I_2(t)<0.14581$
on $[64,4\cdot10^{19}]$ using the published finite-height results
\cite{pt,bpt}. Since $\sqrt t\ge8$, it follows that
\eqref{eq:cyctarget} holds at every real cutoff in these ranges
with $C=0.211$ and $C=0.670$, respectively. The rounding
contribution to either constant is at most $(\log2)/8$.

\begin{corollary}\label{th:cycslack}
There are arbitrarily large $X$ for which, at every
$t\in[X,X+X^{3/4}]$,
\[
 \frac{\log|\mathcal I(t)|-\log|\mathcal A(t)|+mAt}{m\sqrt t}>0.12375 .
\]
In particular the constant $0.12375$ is too small in
\eqref{eq:cyctarget}, while the constant $0.370$ gives an
equivalent formulation of RH.
\end{corollary}
\begin{proof}
If the assertion fails, Theorem~\ref{th:cycwindow} gives RH with
$C=0.12375$. On RH put $S(u)=\sum_\rho K(\rho)e^{i(\im\rho)u}$.
The series converges uniformly and absolutely, so its sum is uniformly
continuous. For every fixed $u_0$ there are $u_j\to\infty$ with
$S(u_j)\to S(u_0)$, by simultaneous return of any finite set of
phases and approximation by the corresponding finite sum. By the
explicit formula including the trivial zeros,
\[
\begin{gathered}
 S(\log x)=4 I_2(x)-\sum_{r\ge1}K(-2r)x^{-2r-1/2},\\
 0<K(-2r)\le\tfrac13 .
\end{gathered}
\]
At $x_*=21611/8$ the prime-power calculation in
Theorem~\ref{th:steppair} gives $I_2(x_*)>0.0309394$. The sum over
trivial zeros is at most $x_*^{-5/2}/(3(1-x_*^{-2}))<10^{-9}$.
Hence $S(\log x_*)>0.12375$. Take the above returns for
$u_0=\log x_*$ and put $X_j=e^{u_j}$. The logarithmic length
of $[X_j,X_j+X_j^{3/4}]$ tends to zero. Uniform continuity
therefore keeps $S(\log t)$ above $0.12375$ throughout these
windows eventually, with a fixed positive margin. The normalized
quantity in the statement equals $S(\log t)+O(t^{-1/2})$ by
\eqref{eq:cycround} and the explicit formula. This contradicts
the assumed bound.
\end{proof}

\subsection{An unconditional upper bound on an enlarged window}
\label{sec:cyclowide}
\begin{theorem}\label{th:cyclowide}
There are constants $b,c>0$ such that for every sufficiently
large $X$ the interval
$\mathcal W_X=[X,X\exp(b(\log\log X)^2)]$ contains a point where
$E(t)=0$, and another where
\[
 \log|\mathcal I(t)|\le\log|\mathcal A(t)|-mAt-cm\sqrt t .
\]
The proof uses no numerical list of zeros or hypothesis on their
real parts. The window in \eqref{eq:cyctarget} is shorter.
\end{theorem}
\begin{proof}
The unconditional form of the explicit formula used above is
$E(x)=-\sum_\rho K(\rho)x^\rho+O(x^{-2})$ for $x\ge2$,
with multiplicities. Let $\gamma_0$ be the least positive
ordinate of a nontrivial zero, and $\beta_0$ the largest real
part at that ordinate. The zero count, discreteness and absence
of real nontrivial zeros give $\gamma_0>0$. Reflection gives
$\beta_0\ge1/2$. Put $f(t)=t^{-\beta_0}E(t)$.

Use logarithmic time $u$ and a probability measure. Its multiplier
at $z$ is $\widehat\nu(z)=\int e^{zu}\,d\nu(u)$.
Set $T=\log X$. For each pair of conjugate zeros of ordinate
$\gamma_0<\gamma\le T$, take $d=\pi/\gamma$ and masses
proportional to $e^{(\beta-\beta_0)d},1$ at $0,d$.
Their multiplier vanishes at $z=\rho-\beta_0$ and its conjugate.
Convolve these measures, allowing repetitions. The total length of
their support is at most $L=\sum_{\gamma_0<\gamma\le T}\pi/\gamma
=O((\log T)^2)$ by $N(T)=O(T\log T)$ and partial summation.
Convolve also $\lceil5\log X\rceil$ uniform probabilities on
$[0,8/T]$ to obtain the probability $\nu_X$.
For $|\gamma|>T$ each uniform multiplier has absolute value at
most $(1+e^{8/T})/8<1/2$ for large $X$.
The high-ordinate contribution to the mean of $f(Xe^{v+u})$
is $O(Xe^{L+v}2^{-\lceil5\log X\rceil}\sum_\rho|K(\rho)|)=o(1)$
uniformly for $0\le v\le2\pi/\gamma_0$.

At $i\gamma_0$ each two-point multiplier has absolute value at
least $\cos(\pi\gamma_0/(2\gamma))$, by projecting the chord onto
its bisector. The product over all $\gamma>\gamma_0$ of these
bounds is strictly positive because all factors are positive and
$\sum_{\gamma>\gamma_0}\gamma^{-2}<\infty$.
The product of the uniform multipliers has modulus
$\left|\sin(4\gamma_0/T)/(4\gamma_0/T)\right|^{\lceil5\log X\rceil}
=1+o(1)$.
Any terms at $\gamma_0$ with smaller real part contribute $o(1)$,
since their normalized exponents have negative real part and $u\ge0$.
With $m_0$ the multiplicity of $\rho_0=\beta_0+i\gamma_0$, the
remaining mean is therefore
\[
 -2m_0\re\!\left(K(\rho_0)\widehat\nu_X(i\gamma_0)
           e^{i\gamma_0(\log X+v)}\right)+o(1).
\]
Its amplitude is bounded below independently of $X$.
Taking both extrema as $v$ traverses one period shows that
$f$ takes values above a positive constant and below its negative
inside $\mathcal W_X$, since the additional logarithmic support length
is $8\lceil5\log X\rceil/T+2\pi/\gamma_0=O(1)$.
Continuity supplies a zero of $E$. At a positive point,
$E(t)\gg t^{\beta_0}\ge\sqrt t$, and the rounding error in
\eqref{eq:cycround} is at most $\log2$. This proves the assertion.
\end{proof}

\subsection{Relative indices under doubling in the two towers}
\label{sec:relativeindices}
A fixed curve determines a relative index in its constant-field
tower. The relative index in the cyclotomic tower of
Eq.~\eqref{eq:cyclogroups} instead yields an exact expression for
the prime error in $S_K$.

\begin{theorem}[A one-sided index for a fixed Jacobian]\label{th:jacindex}
Let $C/\mathbb F_q$ be smooth, projective and geometrically connected,
of genus $g\ge1$. Write $j_r=|\operatorname{Jac}(C)(\mathbb F_{q^r})|$.
RH for $C$ is equivalent to
\[
\begin{gathered}
 j_{2r}/j_r\le q^{gr}+O_C(q^{(g-1/2)r})\\
 \text{for all sufficiently large integers }r .
\end{gathered}
\]
The index also satisfies the corresponding two-sided estimate
unconditionally.
\end{theorem}
\begin{proof}
Let $\alpha_i$, $1\le i\le2g$, be the Frobenius eigenvalues of
Eq.~\eqref{eq:jacnumbers}. Factorization of $1-F^{2r}$ gives
\[
 \frac{j_{2r}}{j_r}
   =\prod_i(1+\alpha_i^r)
   =q^{gr}\prod_i(1+\alpha_i^{-r}).
\]
This is the index of the subgroup of points rational over
$\mathbb F_{q^r}$ in the group over $\mathbb F_{q^{2r}}$.
It is also the order of the kernel of $F^r+1$, a separable
isogeny, and the order of the twist by $[-1]$ over
$\mathbb F_{q^r}$. This factorization works for every genus.
For $\ell\nmid q$ put $T=T_\ell J$. Under the procyclic
description of $H^1_{\rm cts}(\mathbb F_{q^r},T)$ preceding
Theorem~\ref{th:jactwist}, restriction to $\mathbb F_{q^{2r}}$
sends the class of $u$ to that of $(1+F^r)u$.
Its cokernel is $T/(1+F^r)T$, because
$(F^{2r}-1)T\subset(1+F^r)T$. Its order is exactly the
$\ell$-part of $j_{2r}/j_r$.
Put $M=\max_i|\alpha_i|^{-1}$. By simultaneous recurrence on the
torus there are unbounded $r$ with
$\Re(\alpha_i^{-r})\ge|\alpha_i|^{-r}/2$ for every $i$.
On this sequence the positive normalized index is at least
$\prod_i(1+|\alpha_i|^{-r}/2)\ge1+M^r/2$. Its upper bound forces
$M\le q^{-1/2}$. The product of the moduli is $q^g$, giving
$|\alpha_i|=\sqrt q$ for all $i$.

For the independent geometric input, Stepanov's method for
hyperelliptic curves and its extension to arbitrary curves give
$|\#C(\mathbb F_{q^r})-(q^r+1)|\le2gq^{r/2}$
\cite{step69,step72curve,bombstep,kedcurve}. In the trace formula,
simultaneous phase returns rule out any modulus greater than
$\sqrt q$; the product forces equality. Thus both $j_r$ and
$j_{2r}/j_r$ differ from $q^{gr}$ by at most
$(2^{2g}-1)q^{(g-1/2)r}$, by expanding the respective products
with minus or plus signs. This proves the index bound from the
curve estimate. For the curve over $\mathbb F_5$ in
Sec.~\ref{sec:curve}, direct counts give $j_1=8$, $j_2=32$
and $j_2/j_1=4$.
\end{proof}

\head{The corresponding index of inflated cyclotomic unit classes}
Use $A_m=H^1(\operatorname{Gal}(K_m/\mathbb Q),\mathcal O_{K_m}^\times)$
and $I_m\subset A_m$, the subgroup generated by inflations from
the levels $m/p$, $p\mid m$. The maps are injective in this
tower by Proposition~\ref{th:cycunit}. Ideals identify
$A_m$ with the direct sum of the ramification factors
$\mathbb Z/p^{v_p(m)}\mathbb Z$; all their invariant prime-ideal
products are principal. Thus $|A_m|=m$ and
$\log|A_m/I_m|=\Lambda(m)$, including the trivial groups at $m=1$.
The index $[A_{2m}:\operatorname{inf} A_m]$ always equals two.
For the inflated subgroups, $\operatorname{inf}I_m\subset I_{2m}$
since the entire image of $A_m$ is one of the generating subgroups.
Their relative index and its logarithm satisfy
\[
 d_m=[I_{2m}:\operatorname{inf} I_m],\qquad
 \log(d_m/2)=\Lambda(m)-\Lambda(2m).
\]
For example $d_1=1$, $d_{2^a}=2$ for $a\ge1$, and
$d_{p^a}=2p$ for odd primes $p$ and $a\ge1$. At integers with
at least two distinct prime divisors $d_m=2$.

For real $t\ge1$, put $h(t)=\lfloor t\rfloor-\lfloor t/2\rfloor$,
$\mathcal D(t)=\prod_{t/2<m\le t}d_m$, and
\[
 G(t)=\log\mathcal D(t)-(h(t)-1)\log2-t/2 .
\]
Thus $\mathcal D(t)$ is the order of the finite group
\[
 \bigoplus_{t/2<m\le t}I_{2m}/\operatorname{inf}I_m .
\]
The interval $(t/2,t]$ contains exactly one power $2^a$, $a\ge0$.
Only that $m$ contributes to $\sum_{t/2<m\le t}\Lambda(2m)$.
Hence, with $E(u)=\psi(u)-u$ for all $u>0$,
\[
 G(t)=\psi(t)-\psi(t/2)-t/2=E(t)-E(t/2).
\]
For $t=10$ the five indices are $2,14,2,6,2$. Their product
is $672=2^4\cdot42$, and $G(10)=\log42-5$.

\head{A first Galois cohomology group with order $\exp\Lambda(m)$}
For $m\ge2$ take the same cyclic group
$\Gamma_m=\operatorname{Gal}(K_m/\mathbb Q)$, choose a generator
$\sigma_m$, and give the additive lattice
$\mathcal L_m=\mathbb Z[X]/(\Phi_m(X))$ the action
$\sigma_m v=Xv$. Here $\Phi_m$ is the $m$-th cyclotomic
polynomial. This coefficient module is a faithful representation
of the Galois group. Its isomorphism class is independent of the
choice of generator, by the automorphisms of $\mathbb Z[\zeta_m]$.
A $1$-cocycle is determined by its value on $\sigma_m$.
That value lies in the kernel of $1+\sigma_m+\cdots+\sigma_m^{m-1}$,
and the coboundaries have values in $(\sigma_m-1)\mathcal L_m$.
The norm operator is zero on this lattice. Evaluation at $X=1$
therefore gives
\[
 H^1(\Gamma_m,\mathcal L_m)
 \simeq \mathcal L_m/(X-1)\mathcal L_m
 \simeq \mathbb Z/\Phi_m(1)\mathbb Z .
\]
The order is also $\det(1-\sigma_m\mid\mathcal L_m\otimes\mathbb Q)$.
Evaluating $(X^m-1)/(X-1)$ at $1$ gives
$\prod_{d\mid m,\,d>1}\Phi_d(1)=m$. Taking logarithms and using
$\sum_{d\mid m}\Lambda(d)=\log m$ proves by divisor inversion that
$\log\Phi_m(1)=\Lambda(m)$. Take the trivial group in degree one.
The same orders occur in the absolute Galois cohomology of the
subfields themselves. Indeed $K_{m^2}/K_m$ is cyclic of degree
$m$. Give the discrete lattice $\mathcal L_m'$ the preceding
action through this quotient of the absolute group $G_{K_m}$.
The subgroup $G_{K_{m^2}}$ acts trivially. A continuous homomorphism
from this compact subgroup to a discrete torsion-free lattice has
finite image and is zero. Inflation--restriction thus gives
\[
\begin{aligned}
 H^1_{\rm cts}(G_{K_m},\mathcal L_m')
 &\simeq H^1(\operatorname{Gal}(K_{m^2}/K_m),\mathcal L_m')\\
 &\simeq \mathbb Z/\Phi_m(1)\mathbb Z .
\end{aligned}
\]
Consequently the sum of first cohomology groups
\[
 \mathcal H(t)=\bigoplus_{t/2<m\le t}H^1(\Gamma_m,\mathcal L_m)
\]
satisfies, for every real $t\ge1$,
\[
 \log|\mathcal H(t)|=\psi(t)-\psi(t/2)
 =\log\mathcal D(t)-(h(t)-1)\log2 .
\]
At $t=10$, the nontrivial factors have orders $7,2,3$, giving
$|\mathcal H(10)|=42$. In the physical clock $u=2\eta$,
Eq.~\eqref{eq:ret} now gives an exact observable for this order:
\[
 \log|\mathcal H(t)|
   =\mu_{\rm ret}\big((\log(t/2),\log t]\big).
\]
The $a$-th return in the channel with prescribed length
$\tfrac12\log p$ contributes the logarithm of the group order
at degree $p^a$. For a curve the extension degree gives
powers of one Frobenius operator. Here each different field degree
comes with its own integral representation.

\begin{theorem}[An upper bound on cyclotomic cohomology implying RH]
\label{th:cycindex}
Classical RH is equivalent to the existence of $C,t_0>0$ such that
\[
 \log|\mathcal H(t)|\le t/2+C\sqrt t\log^2(2t)
 \qquad(t\ge t_0).
\]
Equivalently $\log\mathcal D(t)$ has the same upper bound
with the additional term $(h(t)-1)\log2$.
It is also equivalent to $\log(1+\|G\|_{K_j^2}^2)=o(\log K_j)$
on some sequence of integers $K_j\to\infty$, where
$\|f\|_X^2=\int_1^X|f(t)|^2\,\dd t/t^2$.
On such a sequence $\|G\|_{K_j^2}^2\ll K_j$ would already give
$\Theta\le3/4$.
\end{theorem}
\begin{proof}
For $\Re s>1$ the Mellin integral of $G(t)$ against
$t^{-s-1}\dd t$ on $[1,\infty)$ equals
\[
 -\frac{(1-2^{-s})\zeta'(s)}{s\zeta(s)}-\frac1{2(s-1)} .
\]
The singularity at 1 cancels. For every nontrivial zero $\rho$
its residue is $-\operatorname{mult}(\rho)(1-2^{-\rho})/\rho\ne0$.
Under the proposed upper bound the difference
$C\sqrt t\log^2(2t)-G(t)$ is nonnegative on a tail. Its Mellin
integral continues holomorphically at every real $s>1/2$
because the transform of the majorant has its sole pole at $1/2$
and $\zeta$ has no real zero in $(0,1)$. Landau's theorem gives convergence in
$\Re s>1/2$ for that nonnegative tail, excluding every such
nonreal pole there. The functional equation gives RH. Conversely
the classical RH estimate $E(t)=O(\sqrt t\log^2(2t))$ gives the
displayed bound directly.

For the mean square, changing variables and using $E(u)=-u$
on $[1/2,1]$ gives
$\|E(\cdot/2)\|_X\le2^{-1/2}\sqrt{\mathcal I(X)}+1/2$.
The triangle inequality in both directions yields
\begin{multline*}
 (1-2^{-1/2})\sqrt{\mathcal I(X)}-1/2
 \le\|G\|_X\\
 \le(1+2^{-1/2})\sqrt{\mathcal I(X)}+1/2 .
\end{multline*}
Theorem~\ref{th:primesquare} now gives both implications on the
specified sequence. Theorem~\ref{th:squareblocks} relates the
same bounds to $S_K$, and $2R_K=S_K-B_K$ with
$B_K=O(\log^2(2K))$.
\end{proof}

\noindent
The pointwise upper bound on an unspecified unbounded sequence is
strictly weaker as a test. In fact $G$ takes both signs arbitrarily
late, unconditionally. If either $G$ or $-G$ were eventually
nonnegative, its Mellin transform would be regular at every positive
real point by the preceding formula. Landau's theorem would make
the tail transform holomorphic in $\Re s>0$, contradicting even
one known nontrivial zeta zero. Thus $\log\mathcal D(t)<
(h(t)-1)\log2+t/2$ already holds at arbitrarily large individual
cutoffs. More quantitatively, a zero with real part $\beta>1/2$
would force both $G(t)>A t^\sigma$ and $G(t)<-A t^\sigma$
at arbitrarily large $t$, for every $A>0$ and
$\tfrac12\le\sigma<\beta$. Otherwise the same argument applied
to one of $A t^\sigma\pm G(t)$ would exclude that pole.
For the accumulated quadratic norm the single-sequence
implication follows instead from Theorem~\ref{th:primesquare}.

\head{Bounds for the complete product}
An unconditional estimate on the whole half-line follows from
Johnston and Yang's prime number theorem \cite{johnstonyang}.
Set $P(y)=9.39y(\log y)^{1.515}\exp(-0.8274\sqrt{\log y})$.
Since $|\psi(y)-y|\le P(y)$ for $y\ge2$, the complete group satisfies
\[
 \big|\log|\mathcal H(t)|-t/2\big|\le P(t)+P(t/2)
 \qquad(t\ge4).
\]
The quotient of this envelope by $\sqrt t\log^2(2t)$ tends to infinity.
B\"uthe's unconditional finite estimate and Johnston's
Proposition~3.1 \cite{buthepsi,johnstonfinite} prove
\[
 |G(t)|\le
 \begin{cases}
 0.94(1+2^{-1/2})\sqrt t,&22<t\le10^{19},\\[2pt]
 \dfrac{1+2^{-1/2}}{8\pi}\sqrt t\log^2 t,
                         &118\le t\le2.169\cdot10^{25}.
 \end{cases}
\]
Each line applies on its entire indicated interval. They follow by
bounding $E(t)$ and $E(t/2)$ separately. The second line
proves the required one-sided inequality for the full cohomological
product throughout that finite range. For example the first line
gives
$\|G\|_X^2\le\|G\|_{22}^2+0.94^2(1+2^{-1/2})^2\log(X/22)$
for $22\le X\le10^{19}$.

There is also an elementary integer bound at every level.
For integers $b\ge1$, $\mathcal D(2b)/2^{b-1}$ is the ratio
of the least common multiples through $2b$ and $b$. It divides
$\binom{2b}{b}$. Indeed, each new prime power in $(b,2b]$
contributes one to Legendre's valuation of the binomial coefficient,
and all remaining contributions are nonnegative. The exact
logarithm of the integer quotient is
\[
 L_b=\sum_{p^a\le b}\log p
       \bigl(\lfloor2b/p^a\rfloor-2\lfloor b/p^a\rfloor\bigr)\ge0 .
\]
Stirling's formula therefore gives
\[
 G(2b)=(2\log2-1)b-\tfrac12\log(\pi b)-L_b+O(1/b).
\]
To reach the upper bound of Theorem~\ref{th:cycindex}, the required
lower bound for $L_b$ must recover the linear term to within
$O(\sqrt b\log^2(2b))$. It suffices to treat these even cutoffs
because moving $t$ by at most two changes $G(t)$ by $O(\log(2t))$.

\head{The contribution of the geometric index in the same sum}
At each prime $p$ choose a smooth projective geometrically
connected curve of fixed genus $g\ge1$ over $\mathbb F_p$,
and write $j_{p,a}$ for its Jacobian count over $\mathbb F_{p^a}$.
For $n=p^a$ use the two constant extensions in the weight
\[
 w(n)=\frac{j_{p,a}+j_{p,2a}/j_{p,a}}{2n^g}.
\]
Give the remaining integers weight one. If $z_i=\alpha_{p,i}^{-a}$,
the two normalized products are $\prod_i(1-z_i)$ and
$\prod_i(1+z_i)$. Their average contains only even elementary
symmetric polynomials. Stepanov--Bombieri purity gives
$|z_i|=n^{-1/2}$, proving uniformly in $p,a$ and in the choice
of the curves that
\[
 |w(n)-1|\le (2^{2g-1}-1)/n .
\]
Define $G_w(t)=\sum_{t/2<n\le t}w(n)\log|H^1(\Gamma_n,\mathcal L_n)|-t/2$.
Chebyshev's bound then gives
\[
\begin{gathered}
 |G_w(t)-G(t)|\le\frac{2(2^{2g-1}-1)}{t}\psi(t)\ll_g1,\\
 \int_1^\infty |G_w-G|^2\frac{\dd t}{t^2}<\infty .
\end{gathered}
\]
In genus one $w(p^a)=1+1/p^a$ exactly. Thus the geometric
part is summable in the mean square on the full half-line.
The constant term of the two determinants still multiplies
$\log|A_n/I_n|=\Lambda(n)$ with coefficient one.
In particular the Mellin transform of $G_w-G$ is holomorphic
in $\Re s>0$; the residues at nontrivial zeta zeros in the
transform of $G_w$ are exactly those computed for $G$.

\section{Divisor sums, fermion parity and accumulated prime losses}\label{app:primeparity}
\subsection{Fourier estimate for the small divisors of the index product}
\label{sec:cohdivisors}
The primitive orders in the cyclotomic calculation satisfy the
finite identity
\[
 \log |A_n/I_n|=-\sum_{d\mid n}\mu(d)\log |A_d| .
\]
This is M\"obius inversion of $\log n=\sum_{d\mid n}\Lambda(d)$,
using $|A_d|=d$. Set $a_d'=\mu(d)\log d$,
$M_1(y)=\sum_{d\le y}a_d'$, and $\{y\}=y-\lfloor y\rfloor$.
For a dyadic block $B\le t<2B$, $B=2^j$, fix an integer
$1\le D\le B$. Splitting the divisor sum for $G(t)$ at $D$
gives the exact finite decomposition $G(t)=T_D(t)+V_D(t)$, where
\[
 \begin{aligned}
 T_D(t)&=\sum_{d\le D}a_d'\,[\{t/d\}-\{t/(2d)\}],\\
 V_D(t)&=-\frac{t}{2}\left(1+\sum_{d\le D}\frac{a_d'}d\right)\\
 &\quad-\sum_{1\le b\le t/D}
 \left[M_1(t/b)-M_1(\max(D,t/(2b)))\right].
 \end{aligned}
\]
Indeed the difference of the two floor functions counts
$t/(2d)<b\le t/d$. For $d>D$ interchange $b$ and $d$.
The displayed term linear in $t$ accounts for the mean of the
two floor sums and the centering of $G$. Dropping it would
change the discrepancy being estimated.

\begin{theorem}[Unconditional divisor ranges]\label{th:indexdivisors}
On every dyadic block and uniformly for $1\le D\le B$,
\[
 \int_B^{2B}|T_D(t)|^2\frac{\dd t}{t^2}
 \ll \frac{(B+D^2)D}{B^2}\log^2(2D).
\]
Let $T(t)=T_{D(B)}(t)$ and $V(t)=V_{D(B)}(t)$ on that block.
For $D(B)=\lfloor B^{2/3}\rfloor$ one has
$\|T\|_X^2\ll\log^3(2X)$. For any fixed $\eta>0$, taking
$D(B)=\max(1,\lfloor B^{2/3}/\log^{1+\eta}(2B)\rfloor)$
gives $\int_1^\infty |T(t)|^2\dd t/t^2<\infty$.
For $D(B)=\max(1,\lfloor B^{5/6}/\log^{2/3}(2B)\rfloor)$ one has
$\|T\|_X^2\ll\sqrt X$.
\end{theorem}
\begin{proof}
Write $e(u)=\exp(2\pi i u)$. The function $f(u)=\{2u\}-\{u\}$
has zero mean, and its Fourier coefficient at $n\ne0$ is
$(1-2\one_{2\mid n})/(2\pi i n)$. Thus $T_D$ is a sum of
functions $a_d'f(t/(2d))$. After coincident frequencies are
combined in a finite Fourier truncation, the remaining frequencies
$\lambda=n/(2d)$ are separated by at least $1/(2D^2)$.
For the combined coefficients $c_\lambda$, matching
$n/d=m/e$ and putting $g_0=(d,e)$ shows that
\[
 \sum_\lambda |c_\lambda|^2
 \le \frac1{12}\sum_{d,e\le D}|a_d'a_e'|\frac{(d,e)^2}{de}
 \ll D\log^2(2D).
\]
Indeed the equal frequencies have
$(n,m)=r(d/g_0,e/g_0)$, and $\sum_{r\ne0}r^{-2}=\pi^2/3$.
The unweighted double sum is at most
$\sum_{g_0\le D}(1+\log(D/g_0))^2\ll D$.
Over a common full period the exact covariance of $f(t/(2d))$
and $f(t/(2e))$ is $(d,e)^2/(12de)$ if $v_2(d)=v_2(e)$,
and $-(d,e)^2/(24de)$ otherwise. The negative value follows
from $\sum_{r\ne0}(-1)^r/r^2=-\pi^2/6$.

For completeness, center the nonnegative function
$(\sin(\pi a u)/(\pi a u))^2$ at $3B/2$, with
$a=\min(B^{-1},(4D^2)^{-1})$. Its value is at least $4/\pi^2$
on the block, its integral is $1/a$, and its Fourier transform is
supported on $[-a,a]$. All unequal frequencies integrate to
zero against this majorizing weight. Consequently
\[
 \int_B^{2B}\Big|\sum_\lambda c_\lambda e(\lambda t)\Big|^2\dd t
 \le \frac{\pi^2}{4a}\sum_\lambda |c_\lambda|^2 .
\]
This is the separated-frequency mean-value estimate
(cf.~\cite{hilbertmv}). Let the truncation grow. The Fourier series
of each of the finitely many periodic sawtooth differences converges
in local $L^2$. Division by $B^2$ proves the block estimate.
For $D=\lfloor B^{2/3}\rfloor$ each block costs $O(\log^2(2B))$.
With the first log divisor the term from $D^3/B^2$ instead costs
$O((\log(2B))^{-1-3\eta})$ after accounting for $\log^2(2D)$,
summable on dyadic blocks. With the $5/6$ choice each block costs
$O(\sqrt B)$. Summing the respective bounds proves the claims.
\end{proof}

The $D^3/B^2$ term gives a genuine limitation of an estimate
using only coefficient sizes. For $D\le B/100$ and $D$ large,
average $t_0$ over $[B,3B/2]$. For each $D/2<d\le D$ the
condition $\{t_0/(2d)\}\in[1/4,3/8]$ holds on a fixed positive
fraction of that interval by complete periods. Some $t_0$
therefore meets the condition for $\gg D$ indices. Give these
indices coefficient one and the others zero. For
$t_0\le t\le t_0+D/32$, their fractional parts remain in
$[1/4,1/2)$ and all $f(t/(2d))\ge1/4$. This yields a block
integral $\gg D^3/B^2$ for these bounded coefficients.

Theorem~\ref{th:cycindex} and the triangle inequality give two
sufficient estimates. With exponent $2/3$, a subpower bound
$\|V\|_{K_j^2}^2=O_\epsilon(K_j^\epsilon)$ for every
$\epsilon>0$ along a single unbounded sequence implies RH. With
exponent $5/6$, the bound $\|V\|_{K_j^2}^2\ll K_j$ implies
$\Theta\le3/4$. At $t\asymp B$ the respective sums for $V$
contain at most $O(B^{1/3})$ and $O(B^{1/6}\log^{2/3}(2B))$
quotients. The geometric weight of Sec.~\ref{sec:relativeindices}
changes either norm by a bounded amount. The unconditional prime
number theorem gives $\|G\|_X^2=O_m(X/\log^m(2X))$ for every fixed
$m>0$. Establishing either of the stronger estimates for $V$
directly at unbounded cutoffs remains open in this construction.

\subsection{Exact cancellation of the zeroth M\"obius moment}
\label{sec:gaussdiv}
Write $G(t)=\psi(t)-\psi(t/2)-t/2$ for $t>0$. The dilation
difference $H(t)=G(2t)-2G(t)$ removes the linear term. By
$\Lambda(n)=-\sum_{d\mid n}\mu(d)\log d$, one has
\[
\begin{gathered}
\begin{aligned}
 H(t)&=\psi(2t)-3\psi(t)+2\psi(t/2)\\
     &=-\sum_{d\le 2t}\mu(d)\log d\,h(t/d),
\end{aligned}\\
 h(u)=\lfloor2u\rfloor-3\lfloor u\rfloor+2\lfloor u/2\rfloor .
\end{gathered}
\]
Over the period $[0,2)$, $h$ is $1$ on $[1/2,1)$,
$-1$ on $[1,3/2)$, and zero elsewhere. Set
\[
\begin{gathered}
 w(v)=v e^{-\pi v^2},\qquad
 \Phi(y)=\int_0^\infty h(yv)w(v)\,\dd v,\\
 F(x)=\int_0^\infty H(xv)w(v)\,\dd v .
\end{gathered}
\]
The integrals apply the same three dilations to the cyclotomic
group orders, or equivalently to the return counts of Eq.~\eqref{eq:ret}.

\begin{proposition}[A localized positive weight and an exact cancellation]
\label{pr:gaussmob}
For $y>0$, $\Phi(y)>0$ and $\Phi(y)=y\Phi(1/y)$.
For every $x>0$,
\[
 \begin{split}
 Z(x)&:=\sum_{d\ge1}\mu(d)\Phi(x/d)= W(1/x),\\
 W(z)&=\big(e^{-\pi z^2/4}-3e^{-\pi z^2}
                           +2e^{-4\pi z^2}\big)/(2\pi),\\
 F(x)&=-\sum_{d\ge1}\mu(d)\log d\,\Phi(x/d)
          =\sum_{n\ge2}\Lambda(n)W(n/x).
 \end{split}
\]
In particular $Z(x)=-21/(8x^2)+O(x^{-4})$ at infinity.
For $x\ge2$ let $L=4\sqrt{\log(2x)}$. The sum of the absolute
contributions from $d<x/L$ and $d>xL$ is $O(x^{-2})$ in
either M\"obius sum, also with weight $\log(d/x)$.
\end{proposition}
\begin{proof}
Integrating $w$ on each half-period where $h=\pm1$ gives
$2\pi\Phi(y)=S(1/y)$, where
\[
 S(v)=\sum_{\substack{n\ge1\\n\ {\rm odd}}}
       \big(e^{-\pi n^2v^2/4}-2e^{-\pi n^2v^2}\big).
\]
With $\theta(v)=\sum_{n\in\mathbb Z}e^{-\pi n^2v^2}$, this is
$S(v)=(\theta(v/2)-3\theta(v)+2\theta(2v))/2$.
Poisson summation gives $\theta(v)=v^{-1}\theta(1/v)$, whence
$S(1/v)=v S(v)$. Every term of $S(v)$ is positive for $v\ge1$,
as $2e^{-3\pi/4}<1$. Reciprocity gives positivity for all $v>0$
and the bounds, with $c=\pi/4$,
\[
 0<\Phi(y)\ll
 \begin{cases}
 e^{-c/y^2},&y\le1,\\
 y e^{-cy^2},&y\ge1.
 \end{cases}
\]
The formula for $F$ follows either by integrating the divisors
or by integrating each prime-power step. The unweighted identity
$\sum_d\mu(d)\lfloor u/d\rfloor=\mathbf1_{\{u\ge1\}}$ follows
by interchanging the finite divisor sums. Substituting it in the
integral for $Z$ gives the three Gaussians explicitly. The
interchanges with the integral converge absolutely: for $d>x$,
$h(xv/d)$ vanishes before $v=d/(2x)$. For the stated localization,
the small-$d$ tail is $O(x\log(2x)e^{-cL^2})$ since
$y e^{-cy^2}$ decreases for $y\ge1$. For $d>xL$ extract
$e^{-cL^2/2}$ from the Gaussian. The remaining sum of
$(\log(2d)+\log(2x))e^{-cd^2/(2x^2)}$ over $d\ge1$
is $O(x\log(2x))$. These two estimates
give the required remainder.
\end{proof}

\begin{theorem}[The signed first moment on the logarithmic clock]
\label{th:gaussfirst}
Put $k_0(u)=e^{-u/2}\Phi(e^u)$, a positive even function, and
\[
 M_\Phi(x)=\sum_{d\ge1}\frac{\mu(d)}{\sqrt d}\,
           \log(d/x)\,k_0(\log(d/x)).
\]
The following are equivalent: RH, $\Delta=0$, and the existence
of constants $C,x_0>0$ for which $M_\Phi(x)\ge-C$ for every $x\ge x_0$.
The sum in this inequality may be restricted to $x/L\le d\le xL$.
\end{theorem}
\begin{proof}
Splitting $\log d=\log x+\log(d/x)$ and using reciprocity gives
\[
 F(x)/\sqrt x=-M_\Phi(x)-(\log x)Z(x)/\sqrt x .
\]
For $\Re s>1$, Mellin transformation of the prime Gaussian sum gives
\[
\begin{gathered}
 \int_0^\infty F(x)x^{-s-1}\dd x
 =-A(s)\frac{\zeta'(s)}{\zeta(s)},\\
 A(s)=\frac{\Gamma(s/2)}{4\pi^{1+s/2}}(2^s-2)(1-2^{-s}).
\end{gathered}
\]
The singularity at $s=1$ is removable, and $A(\rho)\ne0$
at each zero in the critical strip. The part of the integral on
$(0,1)$ is entire by Gaussian decay. If $M_\Phi(x)\ge-C$, an upper
bound $F(x)\le C'\sqrt x$ holds on a tail. Landau's theorem
applied to $C'\sqrt x-F(x)$ excludes all poles in $\Re s>1/2$,
since the transform is regular on the real interval $s>1/2$.
The functional equation and Theorem~\ref{th:boundaryclosure} give
RH and $\Delta=0$.

Conversely, move the Mellin inversion line to $\Re s=-1$.
The zero of $1-2^{-s}$ cancels the gamma pole at zero.
The functional equation bounds the logarithmic derivative on the
new line by $O(\log(2+|\Im s|))$. On horizontal segments choose
admissible heights in the local partial-fraction expansion
\cite{Titchmarsh,mv}. Stirling's formula and the zero count give,
absolutely and with all multiplicities $m_\rho$,
\[
 F(x)=-\sum_\rho m_\rho A(\rho)x^\rho+O(x^{-1})\qquad(x\ge1).
\]
Thus on RH $F(x)=O(\sqrt x)$, giving the bound for $M_\Phi$.
Proposition~\ref{pr:gaussmob} justifies the finite restriction.
\end{proof}

In particular, if all zeros through height $T\ge2$ are on the
critical line, then, with an absolute implied constant,
\[
 |M_\Phi(x)|\ll 1+\sqrt x\,\log(2+T)e^{-\pi T/4}\qquad(x\ge2).
\]
Indeed $|A(\beta+i\gamma)|\ll e^{-\pi|\gamma|/4}$ uniformly for
$0\le\beta\le1$, $|\gamma|\ge2$. Sum consecutive unit intervals
using $N(t+1)-N(t)\ll\log(2+t)$. The verification to
$T=3\cdot10^{12}$ \cite{pt} gives the uniform bound on
$2\le x\le \exp(3\cdot10^{12})$ for this Gaussian observable.

There is also a version testing a single sequence of cumulative
integrals:
\[
 \lim_{X\to\infty}
 \frac{\log(1+\int_1^X |M_\Phi(x)|^2\,\dd x/x)}{\log X}=2\Theta-1 .
\]
Indeed $F(x)\ll x^\Theta\log^2(2x)$ for $x\ge2$ by the bound
on $\psi(x)-x$ in Theorem~\ref{th:primesquare}, applied inside
the Gaussian integral. For its lower bound apply the three-circles
argument of that theorem to $e^{-u/2}F(e^u)$, $u\ge0$.
Its Laplace transform has the stated nonzero pole residues after
translation by $1/2$. The portion of the Mellin integral below
one is entire. Chebyshev's estimate gives $F(x)=O(x)$, so the
tail estimate on $\Re w>1/2$ used in that argument also holds.
Finally $(\log x)Z(x)/\sqrt x$ has finite norm in
$L^2([1,\infty),\dd x/x)$. Thus a subpower cumulative bound
along any one unbounded sequence suffices.

\head{Uniform cancellation near the central divisor range}
For $1\le H_0\le x/2$ the terms of $M_\Phi(x)$ with
$|d-x|\le H_0$ contribute $O(H_0^2/x^{3/2})$ in absolute value.
In particular, the degrees within $x^{3/4}$ contribute $O(1)$.
For any fixed $0<\eta<1$, subtracting explicitly known zeroth
moments extends the range to $|d-x|\le x^{1-\eta}$.
Fix a positive integer $m$ and, for a test function $f$, write
$\mathcal T_f(x)=\sum_d\mu(d)d^{-1/2}f(\log(d/x))$.
For any fixed real $a$,
\[
 \mathcal T_{k_0(\,\cdot-a)}(x)=Z(e^a x)/(e^a x)^{1/2}
                                  =O_a(x^{-5/2}).
\]
For example take $a=\log2$ and
\[
 r(u)=u k_0(u)-c[k_0(u-a)-k_0(u+a)],
 \qquad c=\frac{k_0(0)}{-2k_0'(a)} .
\]
Here $k_0'(a)<0$. In the differentiated sum for
$2\pi k_0(u)=e^{u/2}S(e^u)$, each odd $n$ contributes at $a$
a positive factor times
$(1/2-2b)e^{-b}+(16b-1)e^{-4b}$, $b=\pi n^2$.
It is negative since $(16b-1)/(2b-1/2)<9<e^{3b}$.
Thus $r$ is odd with $r'(0)=0$, so $r(u)=O(u^3)$.
Its sum differs from $M_\Phi(x)$ by the explicitly known
$-c[Z(2x)/\sqrt{2x}-Z(x/2)/\sqrt{x/2}]$.
The central degrees through $|d-x|\le x^{7/8}$ now contribute
$O(1)$ in absolute value. The two Gaussian series give
$0.52431106262720952<c<0.52431106262720954$. An independent
sum over prime powers and the M\"obius sum at $x=128$ give
$0.00005436595492<F(128)/\sqrt{128}<0.00005436595494$.

There are $m$ shifts $a_j$ whose vectors
$(k_0^{(r)}(-a_j))_{0\le r<m}$ span $\mathbb R^m$.
Otherwise a nonzero constant-coefficient differential operator
would annihilate $k_0$ on the real axis. Taking Fourier transforms
would make its polynomial symbol vanish on an interval, since
$\widehat k_0(0)=\int k_0>0$. This is impossible. All derivatives
decay rapidly, by differentiating the Gaussian series on either tail.
Choose real coefficients $c_j$ so that the first $m$ Taylor
coefficients at zero of
\[
 r_m(u)=u k_0(u)-\sum_{j=1}^m c_j k_0(u-a_j)
\]
vanish. The corrected sum differs from $M_\Phi(x)$ by $O_m(x^{-5/2})$.
The absolute contribution of the central degrees is now
\[
 \ll_m \frac{H_0^{m+1}}{x^{m+1/2}}.
\]
For $H_0=x^{1-\eta}$ this is $O_m(1)$ as soon as
$(m+1)\eta\ge1/2$. The tails of the translates still permit
$L=C_m\sqrt{\log(2x)}$, with a fixed sufficiently large $C_m$.
For each fixed $\eta$, it therefore suffices to bound from below
the corrected sum in the two remaining intervals
$[x/L,x-x^{1-\eta})$ and $(x+x^{1-\eta},xL]$. The coefficients
and constants depend on the fixed order. A bound uniform in $x$
for their combined signed sum remains open.

One way to test the two sides together uses only unweighted M\"obius
signs with the factor $d^{-1/2}$. For $v\ge0$ form
\[
 D_x(v)=\sum_{x<d\le xe^v}\frac{\mu(d)}{\sqrt d}
       -\sum_{xe^{-v}\le d<x}\frac{\mu(d)}{\sqrt d}.
\]
Stieltjes integration by parts and the evenness of $k_0$ give
\[
 M_\Phi(x)=-\int_0^\infty (v k_0(v))'D_x(v)\,\dd v.
\]
The endpoint terms vanish since $v k_0(v)$ is zero at the origin
and has Gaussian decay in $e^v$ at infinity. For example it
suffices, for a fixed $b\ge0$, to prove
$|D_x(v)|\ll_\epsilon x^\epsilon e^{bv}$ for every
$\epsilon>0,\ x\ge2,\ v\ge0$. Integrating this estimate gives
$M_\Phi(x)\ll_\epsilon x^\epsilon$, and hence the subpower
cumulative bound above. The signed first-moment criterion itself
only asks for the weighted integral.

\head{A global obstruction for translates of the zeroth moment}
With the Fourier convention $e^{-i\tau u}$,
$\widehat{k_0}(\tau)=A(1/2+i\tau)\zeta(1/2+i\tau)$. This follows
by integrating the floors in the definition of $\Phi$ and continuing
the resulting Mellin transform. Take the ordinate $\gamma_1$ of
a verified simple zero $\rho_1=1/2+i\gamma_1$ \cite{pt}.
Every translate of $k_0$ has zero Fourier coefficient there, whereas
$\widehat{u k_0}(\gamma_1)=i\widehat k_0'(\gamma_1)\ne0$.
Consequently every finite sum of translates satisfies
\[
 \left\|u k_0-\sum_j c_j k_0(\,\cdot-a_j)\right\|_{L^1(\mathbb R)}
 \ge |A(\rho_1)\zeta'(\rho_1)|>0 .
\]
For the first positive ordinate the right-hand side is between
$6.19833296\cdot10^{-6}$ and $6.19833297\cdot10^{-6}$ by
interval evaluation of the gamma factor and the zeta derivative.
For any fixed such residual, continuity supplies a fixed interval
on which its absolute value has a positive lower bound. The positive
density of squarefree integers there after scaling $d/x$ shows
that the sum of absolute terms in $\mathcal T_{r_m}(x)$ on some
proportional interval is $\gg_{r_m}\sqrt x$. The required signed
lower bound on that scale uses arithmetic cancellation.

\head{The total weight of the two M\"obius signs}
Let $P_\pm(x)=\sum_{\mu(d)=\pm1}\log d\,\Phi(x/d)$. Both are
nonnegative, $F=P_- -P_+$, and unconditionally
\begin{multline*}
 P_+(x)+P_-(x)=\frac{\log2}{4\pi\zeta(2)}\,
 x\left(\log x+\frac{\eul+\log(2/\pi)}2\right)\\
 +o(\sqrt x).
\end{multline*}
Indeed, the Mellin transform of $\Phi$ with the preceding sign
convention is $K(s)=A(s)\zeta(s)$, by integrating the floor
functions for $\Re s>1$. The transform for the sum of the two
signs is $-K(s)(\zeta(s)/\zeta(2s))'$ by the squarefree Euler
product. Shift the inversion line to $\Re s=1/2$. The reciprocal
of $\zeta(2s)$ is regular there, including its zero at $s=1/2$.
The classical bounds for $1/\zeta(1+it)$ and
$\zeta'(1+it)/\zeta(1+it)$ as $|t|\to\infty$, polynomial
bounds in the adjacent strip, and the gamma factor give an
integrable function on this line \cite{Titchmarsh,mv}. Its integral
is $o(\sqrt x)$ by the Riemann--Lebesgue lemma. At the double
pole $s=1$ the residue is $x[K(1)\log x+K'(1)]/\zeta(2)$.
Expansion at one gives $K(1)=\log2/(4\pi)$ and
$K'(1)/K(1)=(\eul+\log(2/\pi))/2$, as claimed.
Hence it would suffice to bound the contribution of squarefree
degrees with an odd number of prime factors from above by half
the displayed main term plus $O(\sqrt x)$.

\subsection{The excess of the odd squarefree degrees}
\label{sec:oddexcess}
Put $S_\Phi=P_++P_-$ and $E_\Phi=P_--S_\Phi/2=F/2$, with the
Gaussian weights of Proposition~\ref{pr:gaussmob}. These are the
total logarithmic weight of the squarefree degrees and the excess
at degrees with an odd number of prime factors. The asymptotic
formula for $S_\Phi$ in Sec.~\ref{sec:gaussdiv} changes
$E_\Phi(x)/\sqrt x$ by $o(1)$ if the total is replaced by its
explicit main term. If every squarefree degree is assigned the
weight $\Phi(x/d)$ alone, its odd excess is exactly
$E_0(x)=-Z(x)/2=-W(1/x)/2=21/(16x^2)+O(x^{-4})$.
At every $x>0$ it satisfies
\[
 E_0(x)\le \frac{e^{-\pi/x^2}}{2\pi}
      \le \frac{\sqrt x}{2\pi(4\pi e)^{1/4}} .
\]
The first inequality follows from $e^{-\pi/(4x^2)}\ge e^{-\pi/x^2}$.
The function $x^{-1/2}e^{-\pi/x^2}$ attains its maximum at $x=\sqrt{4\pi}$.
The factor $\log d$ is the weight in the prime-return observable.

\begin{proposition}[Contribution from specified prime factors]
\label{pr:oddmarks}
In each squarefree degree $d$, write $\log d=\sum_{p\mid d}\log p$.
Let $E_{\le y}(x)$ denote the portion of $E_\Phi(x)$ from those
terms with $p\le y$, including such terms when $d>y$.
For every fixed $\lambda>0$, at $y=\lambda x^{5/6}$ one has
\[
 \frac{E_{\le y}(x)}{\sqrt x}\longrightarrow
 -\frac{7\lambda^3}{16}+ \frac{\Gamma(1/4)(2\sqrt2-3)}{16\pi^{5/4}} .
\]
In particular, for $y=x^{5/6}$ and
$\kappa=7/16-\Gamma(1/4)(2\sqrt2-3)/(16\pi^{5/4})
=0.446795491987\ldots$,
\[
 E_\Phi(x)=\frac12\sum_{p>x^{5/6}}(\log p)W(p/x)
                 -\kappa\sqrt x+o(\sqrt x).
\]
\end{proposition}
\begin{proof}
Write $Z^{(p)}(u)=\sum_{p\nmid n}\mu(n)\Phi(u/n)$. The
identity $\mu(pn)=-\mu(n)$ for $p\nmid n$ gives
$Z^{(p)}(u)-Z^{(p)}(u/p)=Z(u)=W(1/u)$.
Gaussian decay at zero justifies iterating to obtain
$Z^{(p)}(u)=\sum_{j\ge0}W(p^j/u)$.
Removing a marked prime reverses the parity of the remaining
squarefree degree. Consequently
\[
\begin{aligned}
 2E_{\le y}(x)&=\sum_{p\le y}(\log p)Z^{(p)}(x/p)\\
               &=\sum_{p\le y}\sum_{j\ge1}(\log p)W(p^j/x).
\end{aligned}
\]
Since $W(v)=-21v^2/8+O(v^4)$ at zero and
$\sum_{p\le y}p^2\log p\sim y^3/3$ by the prime number theorem,
the first powers contribute $-(7\lambda^3/8+o(1))\sqrt x$
to the double sum. The fourth-order Taylor remainder costs
$O(y^5/x^4)$. The squares contribute
\[
 \sqrt x\left(\int_0^\infty W(v^2)\,\dd v+o(1)\right)
 =\sqrt x\big(A(1/2)/2+o(1)\big)
\]
by partial summation. Indeed $\vartheta(\sqrt x\,v)/\sqrt x\to v$,
is $O(v)$, and $v|(W(v^2))'|$ is integrable. The part with
$p>y$ decays faster than any power of $x^{-1}$. The higher
powers have total absolute contribution $O(x^{1/3}\log(2x))$.
For this last statement their cumulative logarithmic weight through
$t$ is $O(t^{1/3}\log(2t))$ by Chebyshev's bound. Integrate by
parts against $|W(t/x)|$, whose derivative in absolute value is
bounded almost everywhere by $|W'(t/x)|/x$.
These estimates also bound the omitted powers for $p>y$.
\end{proof}

The mean of the remaining prime sum can also be evaluated. Write
$\vartheta(t)=\sum_{p\le t}\log p$. Since $\int_0^\infty W(v)\,\dd v=0$,
Taylor expansion at zero gives
\[
 \int_{x^{5/6}}^\infty W(t/x)\,\dd t
      =\tfrac78\sqrt x+O(x^{1/6}).
\]
With $y=x^{5/6}$, the boundary term
$(\vartheta(y)-y)W(y/x)$ is $o(\sqrt x)$ by the prime number
theorem. Integration by parts in the sum over $p>y$ now gives
\begin{multline}
 E_\Phi(x)= -\frac1{2x}\int_{x^{5/6}}^\infty
       (\vartheta(t)-t) W'(t/x)\,\dd t\\
       +\frac{\Gamma(1/4)(2\sqrt2-3)}{16\pi^{5/4}}\sqrt x+o(\sqrt x).
 \label{eq:oddprimefluct}
\end{multline}
Subtracting the mean leaves the prime-square constant
$-0.009295491987\ldots$. The one-sided estimate therefore
requires a lower bound $-C x^{3/2}$ for the integral in
\eqref{eq:oddprimefluct}, uniform for all sufficiently large $x$.

\head{Bounds from a positive Mellin average}
Write $w(x)=W(1/x)$, $g=w\star w$, with
$(f\star h)(x)=\int_0^\infty f(x/u)h(u)\,\dd u/u$.
The weights satisfy $\Phi(x)=\sum_{m\ge1}w(x/m)$ and
$F(x)=\sum_{n\ge1}\Lambda(n)w(x/n)$. Mellin convolution and the
divisor identity $\sum_{d\mid n}\Lambda(d)=\log n$ give
\[
\begin{aligned}
 (E_\Phi\star\Phi)(x)&=\tfrac12\sum_{n\ge1}(\log n)g(x/n)\\
 &=\frac{441\zeta(3)}{512\pi^2}\,x^{-2}\log x+O(x^{-2}).
\end{aligned}
\]
Indeed the Mellin transform is $-A(s)^2\zeta'(s)/2$.
The zero of $A$ at 1 cancels the double pole. On shifting to
$\Re s=-3$, the only pole crossed is at $-2$, with coefficient
of $(s+2)^{-2}$ equal to
$-(-21/8)^2\zeta'(-2)/2$. The stated coefficient follows by the
functional equation at $-2$. The gamma factors give absolute
convergence on the shifted line. All convolution identities here
hold by absolute convergence of the Dirichlet series on $\Re s>1$
and Mellin inversion.

Selberg's identity makes the same calculation possible for two
prime factors. Set $\Lambda_2=\mu*(\log)^2$, with Dirichlet
convolution denoted by $*$, so
$\Lambda_2(n)=\Lambda(n)\log n+\sum_{ab=n}\Lambda(a)\Lambda(b)\ge0$
\cite{selberg1949}. For $Q_2(x)=\sum_n\Lambda_2(n)w(x/n)$,
\[
 (\Phi\star Q_2)(x)=\sum_{n\ge1}(\log n)^2g(x/n)
                  =2A'(1)^2 x+O(x^{-1}).
\]
The transform is $A(s)^2\zeta''(s)$; shifting to $\Re s=-1$
crosses the simple pole at 1. Before this additional average the
transform of $Q_2$ is $A(s)\zeta''(s)/\zeta(s)$. For example at
the first positive zero the residue divided by $A(\rho_1)$ has
real part between $-0.811$ and $-0.810$ by interval arithmetic.

The first average yields an unconditional pointwise consequence.
Let $f(t)=E_\Phi(t)/\sqrt t$, $k_0(v)=e^{-v/2}\Phi(e^v)$, and
$L=4\sqrt{\log(2x)}$. The average divided by $\sqrt x$ is
$\int k_0(v)f(xe^{-v})\,\dd v$. Since $|E_\Phi(t)|\ll t$
by Chebyshev's bound and $k_0(v)\ll e^{|v|/2}
\exp(-\pi e^{2|v|}/4)$, the two tails $|v|>\log L$ have
absolute integral $O(\sqrt x\,e^{-\pi L^2/4})$.
The integral of $k_0$ over the remaining interval tends to its
strictly positive total mass. The weighted mean-value theorem for
the continuous function $f$ therefore gives some
$t\in[x/L,xL]$ with $|E_\Phi(t)|/\sqrt t\ll x^{-5/2}\log x$.
An upper bound at every point requires more. In fact
$\widehat k_0(\gamma_1)=0$. Every probability measure $\nu$
with the same Fourier zero satisfies
$\nu\{|v|>\pi/(3\gamma_1)\}\ge1/3$, since
$\cos(\gamma_1 v)\ge1/2$ on the complementary interval and is
everywhere at least $-1$. Positive averages with this zero have
a minimum spread in the logarithmic clock.

\begin{theorem}[Bounds for the full odd excess]\label{th:oddglobal}
Set $T_*=3\cdot10^{12}$. Without assuming RH, one has
\[
 |E_\Phi(x)|< 7.85\cdot10^{-6}\sqrt x
       \qquad (10^6\le x\le e^{T_*}).
\]
With constant $8.00001$ this holds throughout $0<x\le e^{T_*}$.
At arbitrary $x\ge e^7$, put
$\delta(T)=[55.241(\log T)^{2/3}(\log\log T)^{1/3}]^{-1}$.
Then, unconditionally,
\begin{multline*}
 |E_\Phi(x)|\le 7.8390\cdot10^{-6}x^{1-\delta(2\log x)}\\
              +0.13(\log x)x^{1-\pi/3}+8/x .
\end{multline*}
If RH holds, the finite value of
$\limsup_{x\to\infty} E_\Phi(x)/\sqrt x$ lies strictly between
$7.8389\cdot10^{-6}$ and $7.8390\cdot10^{-6}$.
If RH fails, this limsup is $+\infty$. The same two asymptotic
alternatives hold relative to half the explicit main term of $S_\Phi$.
\end{theorem}
\begin{proof}
For $\rho=\beta+i\gamma$, count all zeros with multiplicity.
The Hadamard identity
$\sum_\rho\Re[1/(\rho(1-\rho))]=2+\eul-\log(4\pi)$ and
$|\gamma|>14$ give
\[
 \sum_\rho|\gamma|^{-2}
 \le(1+5/(4\cdot14^2))(2+\eul-\log(4\pi))<0.047 .
\]
Indeed the real part in the sum is at least $1/(\gamma^2+5/4)$,
as follows by writing the real and imaginary parts of
$\rho(1-\rho)$. Rotate the Euler integral for
$\Gamma(1+\rho/2)$ by ${\rm sgn}(\gamma)\pi/3$.
Log convexity gives $\Gamma(1+\beta/2)\le1$, and division by
$|\rho/2|$ gives
\[
 |A(\rho)|\le (6\sqrt2/\pi)|\gamma|^{-1}e^{-\pi|\gamma|/6}.
\]
Here $|2^\rho+2^{1-\rho}-3|\le6$. In particular the sum of the
absolute coefficients over $|\gamma|>50$ is below
$(6\sqrt2/\pi)\,0.047\cdot50e^{-50\pi/6}<2.71\cdot10^{-11}$.
There are ten positive ordinates below 50. Their interval enclosures
give $\sum_{0<\gamma<50}|A(1/2+i\gamma)|
<7.8389570\cdot10^{-6}$ \cite{pt}. Thus half the sum of
$|A(\rho)|$ over all zeros is less than $7.8390\cdot10^{-6}$.

The error in the explicit formula for $E_\Phi=F/2$ is bounded
by $8/x$. On $s=-1+it$ the half-integer gamma identity gives
\[
 |A(s)|\le\frac{15\sqrt2}{4}\,
           \frac{e^{-\pi |t|/4}}{\sqrt{1+t^2}}.
\]
The functional equation gives
$|\zeta'(s)/\zeta(s)|\le7+|t|$:
the cotangent has modulus at most one, the log derivative on
$\Re s=2$ is bounded by $\sum_{n\ge2}n^{-3/2}\le2$, and
$|\Gamma'(2-it)/\Gamma(2-it)|\le1+|t|$ by the convergent
digamma difference series at 2. Since
$(7+|t|)/\sqrt{1+t^2}\le\sqrt{50}$, the absolute contour
constant for $F$ is at most $150/\pi^2<16$.

Through $T_*$ the zeros are on the critical line \cite{pt}.
If $x\le e^{T_*}$, the normalized contribution from
$|\gamma|>T_*$ is bounded by
$(3\sqrt2/\pi)\,0.047\,T_*\exp[-(\pi/6-1/2)T_*]<10^{-100}$.
Adding $8x^{-3/2}$ proves the finite-range estimate. For
$1\le x\le10^6$ the same argument gives the stated larger
constant. For $0<x<1$, the functions
$x^{-1/2}e^{-\pi n^2/(4x^2)}$ increase on this interval for
$n\ge2$. Hence
$|E_\Phi(x)|/\sqrt x\le(3/(2\pi))\sum_{n\ge2}n e^{-\pi n/2}<1$.
For
$T=2\log x\ge14$, the zero-free region of Mossinghoff, Trudgian
and Yang gives $\beta\le1-\delta(T)$ below $T$
\cite{mtyregion}. The higher ordinates contribute at most
$(3\sqrt2/\pi)\,0.047\,xT e^{-\pi T/6}$ by monotonicity of
$v e^{-\pi v/6}$ for $v\ge14$. This proves the all-scale bound.

On RH the normalized explicit formula is an absolutely convergent
trigonometric series with an $O(x^{-3/2})$ remainder. For a
lower bound fix
\[
\begin{gathered}
 u_0=\frac{\pi-\arg A(\rho_1)+2\pi\cdot105016}{\gamma_1}\\
 (-\pi<\arg A(\rho_1)\le\pi).
\end{gathered}
\]
Interval evaluation gives
$-\sum_{0<\gamma<50}\Re(A(\rho)e^{i\gamma u_0})
>7.8389552\cdot10^{-6}$. There are arbitrarily large simultaneous
returns of these ten phases to their values at $u_0$, by
compactness of the finite torus. The bound on the tail beyond 50
gives the asserted lower endpoint. The absolute sum gives the
upper endpoint. If the limsup were bounded above in the presence
of a zero to the right of $1/2$, an eventual one-sided upper
bound for $F$ would contradict Theorem~\ref{th:gaussfirst}.
Finally $S_\Phi$ differs from its explicit main term by
$o(\sqrt x)$, which preserves both alternatives.
\end{proof}

On finite intervals the main term $S_0$ of Sec.~\ref{sec:gaussdiv}
gives an alternative centering. Its contour formula yields
\[
\begin{gathered}
 |S_\Phi(x)-S_0(x)|/2\le C_{\rm sf}\sqrt x,\\
 C_{\rm sf}=\frac1{4\pi}\int_{\mathbb R}
 \left|K(s)\left(\frac{\zeta(s)}{\zeta(2s)}\right)'\right|_{s=1/2+it}\dd t .
\end{gathered}
\]
The integral is finite by the bounds used in Sec.~\ref{sec:gaussdiv}.
Thus the finite-range estimate for $P_--S_0/2$ has the additional
constant $C_{\rm sf}$.

More generally put $\Theta=\sup_\rho\re(\rho)$. For $x\ge1$ the
same absolute sum gives $|E_\Phi(x)|\le 7.8390\cdot10^{-6}x^\Theta+8/x$.
For every $0<a<\Theta$,
\[
 \limsup_{x\to\infty} E_\Phi(x)/x^a=+\infty,\qquad
 \liminf_{x\to\infty} E_\Phi(x)/x^a=-\infty .
\]
Indeed an eventual bound of either sign by $C x^a$ would make
one of $C x^a\pm E_\Phi(x)$ nonnegative on a tail. Its Mellin
transform is regular on the real interval $(a,\infty)$, so
Landau's theorem forces convergence throughout $\Re s>a$.
A zero with $\re(\rho)>a$ gives a pole there with nonzero residue.

\subsection{A power series for the loss over a doubling}
\label{sec:losspowers}

Write $D(t)=\vartheta(t)-t$, $a(x)=x^{5/6}$, and normalize the
integral in \eqref{eq:oddprimefluct} by setting
\[
\begin{aligned}
 J(x)&=x^{-3/2}\int_{a(x)}^\infty D(t)W'(t/x)\,\dd t,\\
 R(x)&=x^{-3/2}\int_0^{a(x)}D(t)W'(t/x)\,\dd t,\\
 U(x)&=x^{-1/2}\sum_p(\log p)W(p/x).
\end{aligned}
\]
Partial integration on the whole half-line gives $J(x)=-U(x)-R(x)$
since the mean term vanishes by $\int_0^\infty W(v)\,\dd v=0$.
For $X\le y\le2X$, put $h=1-(X/y)^2$, so $0\le h\le3/4$.
On the boost clock $u=\log(y/X)$, $h=1-e^{-2u}$.
Set $d=(1,-3,2)$, $\lambda=(\pi/4,\pi,4\pi)$ and $z_{pj}=\lambda_j(p/X)^2$.
Define the polynomials and the signed prime sums
\[
\begin{aligned}
 q_n(z)&=\sum_{k=0}^n(-1)^k\binom{1/4}{k}
                         \frac{z^{\,n-k}}{(n-k)!},\\
 c_n(X)&=\frac1{2\pi\sqrt X}\sum_p(\log p)
                \sum_{j=1}^3 d_j e^{-z_{pj}}q_n(z_{pj}).
\end{aligned}
\]
Then the loss along the entire step is exactly
\[
 \boxed{\displaystyle J(X)-J(y)=\sum_{n=1}^\infty c_n(X)h^n+R(y)-R(X).}
\]
In fact $U(y)=\sum_{n\ge0}c_n(X)h^n$. The generating function
of $q_n(z)$ is $(1-h)^{1/4}e^{zh}$, with the branch positive
at $h=0$. The factor $e^{-z_{pj}}$ leaves Gaussian decay
uniformly on every disk $|h|\le r<1$. The first polynomials are
\[
\begin{gathered}
 q_0=1,\quad q_1=z-\tfrac14,\quad
 q_2=\tfrac12 z^2-\tfrac14 z-\tfrac3{32},\\
 q_3=\tfrac16 z^3-\tfrac18 z^2-\tfrac3{32}z-\tfrac7{128}.
\end{gathered}
\]
They satisfy $(n+1)q_{n+1}=(n+z-1/4)q_n-zq_{n-1}$ with
$q_{-1}=0$. The contribution of the continuous prime density
$\dd t$ is zero separately in every coefficient. Substituting
$v=\lambda_j(t/X)^2$ in that contribution gives the common
factor $\sum_jd_j/\sqrt{\lambda_j}=0$.

The correction also has a convergent expansion. Write
$w_k=(-1)^k\sum_jd_j\lambda_j^k/(2\pi k!)$, $k\ge1$, and
$a=a(x)$. Termwise integration of $W(v)=\sum_{k\ge1}w_kv^{2k}$
on $0\le v\le a/x$ gives
\begin{multline*}
 R(x)=\sum_{k\ge1}\frac{w_k}{x^{2k+1/2}}
 \Bigg\{\sum_{p\le a}(\log p)(a^{2k}-p^{2k})\\
                     -\frac{2k}{2k+1}a^{2k+1}\Bigg\}.
\end{multline*}
This formula yields the ordinary Taylor coefficients on each interval
where the primes below the moving cutoff are fixed. More explicitly
put $b=X^{5/6}$, $s=1-h$, $S=\{p:p\le b s^{-5/12}\}$,
$\vartheta_S=\sum_{p\in S}\log p$ and $M_{k,S}=\sum_{p\in S}p^{2k}\log p$.
The $k$-th term of $R(y)$ equals
\begin{multline*}
 \frac{w_k}{X^{2k+1/2}}
 \Bigg\{
 b^{2k}\vartheta_S s^{k/6+1/4}
 -M_{k,S}s^{k+1/4}\\
 -\frac{2k}{2k+1}b^{2k+1}s^{(k-1)/6}
 \Bigg\}.
\end{multline*}
For the coefficient of $h^n$ replace each $s^\nu$ by
$(-1)^n\binom{\nu}{n}$. For a fixed finite $S$, the factorial
in $w_k$ gives normal convergence on compact subsets of $|h|<1$.
At $y_p=p^{6/5}$ the new term starts at zero. The one-sided
derivatives of $J$ satisfy
\[
 J'(y_p+)-J'(y_p-)
 =-\tfrac56\,y_p^{-5/3}(\log p)W'(y_p^{-1/6}).
\]
The jump follows by differentiating the lower limit. It is nonzero
for every sufficiently large prime since $W'(v)=-21v/4+O(v^3)$
near zero. These points require a change of Taylor branch.
The points $1/y_p=p^{-6/5}$ accumulate at zero, precluding a
single convergent Taylor series in $1/y$ about zero for $J$.
Expanding $W'(t/y)$ at zero under the upper-tail integral would
instead require the improper integrals of $t^mD(t)$, $m=1,3,\ldots$.
These diverge. Indeed $D(t)=\vartheta(n)-t$ on $[n,n+1)$.
If the integrals on both half-intervals tended to zero, the two
$t^m$-weighted means of $t$ there would each equal
$\vartheta(n)+o(1)$. The means differ by $1/2+O_m(1/n)$.

\head{An explicit truncation uniform in the step}
For $X\ge2$, keep $N=\lceil4\log X\rceil$ powers and primes
through $P=8X\sqrt{\log X}$. Denote the corresponding coefficients
by $c_n^{(P)}(X)$. Let $\eta=4\log(7/6)-1/2>0.1166$, and let $R_1(x)$
be the $k=1$ term of the series for $R(x)$, with $w_1=-21/8$.
Uniformly for $X\le y\le2X$,
\begin{multline*}
 \left|J(X)-J(y)-\sum_{n=1}^{N}c_n^{(P)}(X)h^n-R_1(y)+R_1(X)\right|\\
 \le 89X^{-\eta}+90X^{-1/3}.
\end{multline*}
Both errors are summable on $X=2^j$. A proof with explicit constants
uses $\vartheta(t)\le\psi(t)\le4(\log2)t$. Indeed
$\psi(2m)-\psi(m)\le\log\binom{2m}{m}\le2m\log2$ at integer $m$.
Summing over powers of two and using monotonicity gives the claim
at real $t$. In particular $|D(t)|\le2t$.

On $|h|=7/8$, $\Re(1-h)\ge1/8$ and
$|(1-h)^{1/4}|\le(15/8)^{1/4}$.
By partial integration,
$\sum_p(\log p)e^{-c(p/X)^2}\le 2\sqrt{\pi/c}\,(\log 2)X$.
At $c=\pi/32$ this bounds the absolute Gaussian sum for $U$,
including its prefactor, by $9\sqrt X$. Cauchy's inequality
therefore bounds the Taylor tail of the sum through $P$, at every
$0\le h\le3/4$, by $63\sqrt X(6/7)^{N+1}$.
For the primes above $P$ extract $e^{-(\pi/64)(P/X)^2}=X^{-\pi}$
and apply the same estimate at $c=\pi/64$. Their contribution
is at most $13X^{1/2-\pi}$ at either endpoint of the loss.
Finally, $1-e^{-z}\le z$ for $z\ge0$ gives
\[
\begin{gathered}
 |W'(v)-2w_1 v|\le(561\pi/16)v^3,\\
 |R(x)-R_1(x)|\le (561\pi/40)x^{-1/3}<45x^{-1/3}.
\end{gathered}
\]
Combining the three estimates proves the claimed truncation bound.
Also $R_1(x)=o(1)$ by the prime number theorem, since it is
$2w_1 x^{-5/2}\int_0^{x^{5/6}}tD(t)\,\dd t$.

For the original observable $2E_\Phi(y)/\sqrt y$ there is a single
series on the whole step. In $c_n(X)$ replace $\sum_p(\log p)$
and $z_{pj}$ by $\sum_{m\ge2}\Lambda(m)$ and
$\lambda_j(m/X)^2$, respectively. The same Gaussian estimate with
$\psi$ bounds the error for its finite polynomial, including the
constant term, by $63X^{-\eta}+13X^{1/2-\pi}$. Thus a uniform upper
bound on these polynomials over $0\le h\le3/4$ and all sufficiently
large dyadic $X$ would imply RH. The remainder bound controls the
truncation at every scale. An upper bound for the signed polynomials
as $X$ increases remains to be proved.

For example at $X=64$ the first three coefficients of the
\emph{full} loss, including the cutoff, give the local expansion
at $h=0$
\[
\begin{gathered}
 J(64)-J(y)=b_1h+b_2h^2+b_3h^3+O(h^4),\\
 (b_1,b_2,b_3)\simeq(0.03161972019,\\
 -0.07435754627,-0.04933909918).
\end{gathered}
\]
Each displayed coefficient is rounded within $5\cdot10^{-12}$.
In powers of $u$ the corresponding first three coefficients are
$2b_1,\ 4b_2-2b_1,\ 8b_3-8b_2+4b_1/3$.
The first change of branch is at $y=37^{6/5}=76.18019305\ldots$.
At $y=128$, $U(128)-U(64)=-0.017864009667\ldots$ and
$R(128)-R(64)=0.022284906957\ldots$. Their sum is
$J(64)-J(128)=0.0044208972907\ldots$.
The maximum over the interval satisfies
\[
  0.00531979<
   \max_{64\le y\le128}\bigl[J(64)-J(y)\bigr]<0.00531981 .
\]
The lower bound holds at $y=109+7/64$. To verify the upper bound
by interval arithmetic, the calculation in \texttt{anc/operator\_check.py}
uses the five transitions at $p=37,41,43,47,53$.
For each fixed prime set $S$, put $H(v)=\int_0^v W(u)\,\dd u$.
Partial integration gives the analytic continuation of the branch as
\begin{multline*}
 R_S(y)=\frac1{\sqrt y}\Big\{(\vartheta_S-a)W(a/y)
           -\sum_{p\in S}(\log p)W(p/y)\\
           +yH(a/y)\Big\}.
\end{multline*}
On $|h|=7/8$ one has $|y|<182$, $|a|<77$,
$|a/y|<0.53$, $|y^{-1/2}|<0.15$ and $\vartheta_S<162$.
On $|v|\le0.53$, $|W(v)|<34$ and $|H(v)|<19$.
Moreover $|W(p/y)|\le3/\pi$ because $\Re(1-h)>0$.
Together with $|U(y)|\le72$, these inequalities bound
$|U+R_S|$ by $4096$, also after truncating $U$ at prime
$4096$. On this disk the omitted primes contribute at most
$(2/\sqrt{64})\sum_{m>4096}m e^{-2\pi m/64}<10^{-165}$.
The degree-200 Taylor remainder on $0\le h\le3/4$, including
that prime tail, is bounded by
$7\cdot4096(6/7)^{201}+10^{-165}<1.1\cdot10^{-9}$.
On 22 subintervals the polynomial is bounded at the appropriate
endpoint if its derivative has certified sign, and otherwise by
its midpoint value plus half the interval length times the supremum
of the absolute derivative. This covers the full doubling.

\subsection{Derivatives of the loss and the values at prime cutoffs}
\label{sec:lossderivatives}
Write $L_X(y)=J(X)-J(y)$, $\mathfrak d=y\,\mathrm d/\mathrm dy$,
$\alpha=5/6$, $a=y^\alpha$, and $b=a/y$. Thus
$\mathfrak d=\mathrm d/\mathrm du$ at $y=Xe^u$.
A dot denotes this derivative. On each open interval between
successive cutoffs $y_p=p^{6/5}$, differentiation under the
Gaussian integral gives, with $v=t/y$,
\begin{multline*}
 \dot L_X(y)=y^{-3/2}\int_a^\infty D(t)
       [\tfrac32 W'(v)+vW''(v)]\,\dd t\\
       +\alpha y^{-2/3}D(a)W'(b).
\end{multline*}
Gaussian decay provides an integrable majorant for every fixed
number of differentiations on compact $y$-intervals. At $y_p$
use the appropriate one-sided value of $D(a)$.
The loss is continuous and its slope has the jump
\[
\begin{aligned}
 [\dot L_X]_p&=\alpha p^{-4/5}(\log p)W'(p^{-1/5})\\
       &=-\tfrac{35}{8}\frac{\log p}{p}
                         +O\big((\log p)p^{-7/5}\big).
\end{aligned}
\]
It is positive at every prime $p\le211$ and negative at every
prime $p\ge223$. Indeed
$W'(v)=-v e^{-\pi v^2/4}(1-12z+32z^5)/4$,
$z=e^{-3\pi v^2/4}$. The polynomial is convex at $z>0$
and has one minimum there. It is negative at both $z(2)$ and
$z(211)$, and positive with positive derivative at $z(223)$,
where $z(p)=\exp[-(3\pi/4)p^{-2/5}]$. Its transition on this
range occurs for $216.37969<p<216.37970$, treating $p$ as real.
Ordinary scale derivatives follow from $L'_X=\dot L_X/y$ and
$L''_X=(\ddot L_X-\dot L_X)/y^2$.

\head{Stationary points on $[64,128]$}
On all six open branches of $64\le y\le128$, interval
evaluation gives
\[
 -0.75<\ddot L_{64}<-0.24,\qquad
 1<\mathfrak d^3L_{64}<4.
\]
On each branch the initial slope is positive and the final slope
negative. Hence there is exactly one stationary point there, a
strict maximum. The five intervening cutoffs are strict local minima.
The values below are rounded within $5\cdot10^{-9}$ in location and
$5\cdot10^{-13}$ in height.
\begin{center}
\begin{tabular}{l@{\quad}r@{\qquad}r}
\toprule
Range of $a=y^{5/6}$&Stationary $y$&$L_{64}(y)$\\
\midrule
$32<a<37$&70.66840236&0.002989260856\\
$37<a<41$&80.26344179&0.002377649164\\
$41<a<43$&90.53627824&0.001838599175\\
$43<a<47$&99.36503892&0.003984610973\\
$47<a<53$&109.10679932&0.005319800436\\
$53<a<128^{5/6}$&122.80058029&0.004704140277\\
\bottomrule
\end{tabular}
\end{center}
The fifth peak is the unique global maximum. Its curvature is in
$(-0.48201359,-0.48201358)$. At the beginning of the doubling
$\dot L_{64}(64)=0.06323944037\ldots$ and
$\ddot L_{64}(64)=-0.72133925087\ldots$.

\begin{figure*}[t]
\centering
\centering
\begin{tikzpicture}
\definecolor{lossink}{RGB}{24,74,94}
\definecolor{slopeink}{RGB}{144,63,42}
\pgfplotsset{lossaxes/.style={width=12.1cm,height=2.5cm,scale only axis,xmin=64,xmax=128,xtick={64,80,96,112,128},axis lines=left,axis line style={black!65,-},tick align=outside,tick style={black!65},ymajorgrids,grid style={black!12},tick label style={font=\footnotesize},label style={font=\small},scaled ticks=false,/pgf/number format/1000 sep={}}}
\begin{axis}[lossaxes,name=lossaxis,ymin=0,ymax=6.1,ytick={0,2,4,6},xticklabels={,,,,},ylabel={$10^3 L_{64}(y)$}]
\draw[black!36,densely dotted,line width=.45pt] (axis cs:76.1801930497,0)--(axis cs:76.1801930497,6.1);
\draw[black!36,densely dotted,line width=.45pt] (axis cs:86.1669316093,0)--(axis cs:86.1669316093,6.1);
\draw[black!36,densely dotted,line width=.45pt] (axis cs:91.2351408162,0)--(axis cs:91.2351408162,6.1);
\draw[black!36,densely dotted,line width=.45pt] (axis cs:101.5120105529,0)--(axis cs:101.5120105529,6.1);
\draw[black!36,densely dotted,line width=.45pt] (axis cs:117.2549115806,0)--(axis cs:117.2549115806,6.1);
\addplot[no marks,line width=.88pt,color=lossink] coordinates {
  (64.00000000,0.0000000000) (64.46846896,0.4421266391) (64.93693793,0.8436687275) (65.40540689,1.2060547497) (65.87387585,1.5306859079)
  (66.34234482,1.8189356955) (66.81081378,2.0721495464) (67.27928274,2.2916445570) (67.74775171,2.4787092768) (68.21622067,2.6346035655)
  (68.68468963,2.7605585120) (69.15315860,2.8577764125) (69.62162756,2.9274308044) (70.09009652,2.9706665514) (70.55856549,2.9885999781)
  (71.02703445,2.9823190498) (71.49550342,2.9528835934) (71.96397238,2.9013255578) (72.43244134,2.8286493098) (72.90091031,2.7358319612)
  (73.36937927,2.6238237264) (73.83784823,2.4935483055) (74.30631720,2.3459032912) (74.77478616,2.1817605964) (75.24325512,2.0019668992)
  (75.71172409,1.8073441049) (76.18019305,1.5986898194)
};
\addplot[no marks,line width=.88pt,color=lossink] coordinates {
  (76.18019305,1.5986898194) (76.56429838,1.7454504178) (76.94840371,1.8755749734) (77.33250904,1.9895581351) (77.71661437,2.0878857543)
  (78.10071970,2.1710349446) (78.48482502,2.2394741472) (78.86893035,2.2936632011) (79.25303568,2.3340534187) (79.63714101,2.3610876654)
  (80.02124634,2.3752004434) (80.40535167,2.3768179782) (80.78945700,2.3663583090) (81.17356233,2.3442313817) (81.55766766,2.3108391435)
  (81.94177299,2.2665756402) (82.32587832,2.2118271152) (82.70998365,2.1469721097) (83.09408898,2.0723815639) (83.47819430,1.9884189195)
  (83.86229963,1.8954402228) (84.24640496,1.7937942277) (84.63051029,1.6838225000) (85.01461562,1.5658595209) (85.39872095,1.4402327907)
  (85.78282628,1.3072629325) (86.16693161,1.1672637951)
};
\addplot[no marks,line width=.88pt,color=lossink] coordinates {
  (86.16693161,1.1672637951) (86.36186273,1.2292056562) (86.55679386,1.2878459790) (86.75172498,1.3432327258) (86.94665610,1.3954134728)
  (87.14158723,1.4444354105) (87.33651835,1.4903453453) (87.53144947,1.5331897004) (87.72638060,1.5730145173) (87.92131172,1.6098654569)
  (88.11624284,1.6437878005) (88.31117397,1.6748264517) (88.50610509,1.7030259371) (88.70103621,1.7284304080) (88.89596734,1.7510836415)
  (89.09089846,1.7710290422) (89.28582958,1.7883096434) (89.48076071,1.8029681082) (89.67569183,1.8150467316) (89.87062295,1.8245874413)
  (90.06555408,1.8316317995) (90.26048520,1.8362210045) (90.45541632,1.8383958917) (90.65034745,1.8381969356) (90.84527857,1.8356642509)
  (91.04020969,1.8308375943) (91.23514082,1.8237563662)
};
\addplot[no marks,line width=.88pt,color=lossink] coordinates {
  (91.23514082,1.8237563662) (91.63040504,2.0499214793) (92.02566926,2.2615020333) (92.42093348,2.4588784153) (92.81619770,2.6424255491)
  (93.21146192,2.8125128973) (93.60672614,2.9695044659) (94.00199036,3.1137588121) (94.39725458,3.2456290535) (94.79251880,3.3654628807)
  (95.18778302,3.4736025710) (95.58304724,3.5703850048) (95.97831146,3.6561416844) (96.37357568,3.7311987543) (96.76883991,3.7958770232)
  (97.16410413,3.8504919883) (97.55936835,3.8953538614) (97.95463257,3.9307675960) (98.34989679,3.9570329170) (98.74516101,3.9744443512)
  (99.14042523,3.9832912597) (99.53568945,3.9838578714) (99.93095367,3.9764233184) (100.32621789,3.9612616721) (100.72148211,3.9386419811)
  (101.11674633,3.9088283103) (101.51201055,3.8720797802)
};
\addplot[no marks,line width=.88pt,color=lossink] coordinates {
  (101.51201055,3.8720797802) (102.11750675,4.1141684157) (102.72300294,4.3308447983) (103.32849913,4.5231548162) (103.93399533,4.6921189267)
  (104.53949152,4.8387323183) (105.14498771,4.9639650951) (105.75048391,5.0687624846) (106.35598010,5.1540450645) (106.96147629,5.2207090099)
  (107.56697249,5.2696263578) (108.17246868,5.3016452892) (108.77796487,5.3175904254) (109.38346107,5.3182631397) (109.98895726,5.3044418816)
  (110.59445345,5.2768825137) (111.19994965,5.2363186584) (111.80544584,5.1834620557) (112.41094203,5.1190029289) (113.01643823,5.0436103584)
  (113.62193442,4.9579326629) (114.22743061,4.8625977855) (114.83292681,4.7582136859) (115.43842300,4.6453687361) (116.04391919,4.5246321200)
  (116.64941539,4.3965542348) (117.25491158,4.2616670955)
};
\addplot[no marks,line width=.88pt,color=lossink] coordinates {
  (117.25491158,4.2616670955) (117.66818421,4.3296499963) (118.08145684,4.3912912909) (118.49472948,4.4467852695) (118.90800211,4.4963225113)
  (119.32127474,4.5400899292) (119.73454737,4.5782708156) (120.14782000,4.6110448878) (120.56109263,4.6385883340) (120.97436526,4.6610738602)
  (121.38763790,4.6786707360) (121.80091053,4.6915448420) (122.21418316,4.6998587163) (122.62745579,4.7037716015) (123.04072842,4.7034394917)
  (123.45400105,4.6990151800) (123.86727368,4.6906483047) (124.28054632,4.6784853968) (124.69381895,4.6626699270) (125.10709158,4.6433423519)
  (125.52036421,4.6206401612) (125.93363684,4.5946979236) (126.34690947,4.5656473334) (126.76018211,4.5336172563) (127.17345474,4.4987337753)
  (127.58672737,4.4611202356) (128.00000000,4.4208972907)
};
\addplot[only marks,mark=*,mark size=1.8pt,lossink] coordinates {(109.10679932,5.319800436)};
\node[anchor=south,font=\footnotesize,text=lossink] at (axis cs:109.10679932,5.41) {$y_*$};
\end{axis}
\begin{axis}[lossaxes,at={(lossaxis.south west)},anchor=north west,yshift=-0.55cm,ymin=-0.043,ymax=0.069,ytick={-0.04,-0.02,0,0.02,0.04,0.06},yticklabels={$-0.04$,$-0.02$,$0$,$0.02$,$0.04$,$0.06$},ylabel={$\mathrm dL_{64}/\mathrm du$},xlabel={Scale $y$}]
\draw[black!36,densely dotted,line width=.45pt] (axis cs:76.1801930497,-0.043)--(axis cs:76.1801930497,0.069);
\draw[black!36,densely dotted,line width=.45pt] (axis cs:86.1669316093,-0.043)--(axis cs:86.1669316093,0.069);
\draw[black!36,densely dotted,line width=.45pt] (axis cs:91.2351408162,-0.043)--(axis cs:91.2351408162,0.069);
\draw[black!36,densely dotted,line width=.45pt] (axis cs:101.5120105529,-0.043)--(axis cs:101.5120105529,0.069);
\draw[black!36,densely dotted,line width=.45pt] (axis cs:117.2549115806,-0.043)--(axis cs:117.2549115806,0.069);
\draw[black!45,line width=.45pt] (axis cs:64,0)--(axis cs:128,0);
\addplot[no marks,line width=.88pt,color=slopeink] coordinates {
  (64.00000000,0.0632394404) (64.46846896,0.0580177642) (64.93693793,0.0529133507) (65.40540689,0.0479267501) (65.87387585,0.0430583390)
  (66.34234482,0.0383083293) (66.81081378,0.0336767760) (67.27928274,0.0291635869) (67.74775171,0.0247685304) (68.21622067,0.0204912439)
  (68.68468963,0.0163312428) (69.15315860,0.0122879280) (69.62162756,0.0083605944) (70.09009652,0.0045484384) (70.55856549,0.0008505656)
  (71.02703445,-0.0027340015) (71.49550342,-0.0062063168) (71.96397238,-0.0095675033) (72.43244134,-0.0128187470) (72.90091031,-0.0159612899)
  (73.36937927,-0.0189964245) (73.83784823,-0.0219254875) (74.30631720,-0.0247498550) (74.77478616,-0.0274709367) (75.24325512,-0.0300901713)
  (75.71172409,-0.0326090219) (76.18019305,-0.0350289717)
};
\addplot[no marks,line width=.88pt,color=slopeink] coordinates {
  (76.18019305,0.0307901179) (76.56429838,0.0275793892) (76.94840371,0.0244348664) (77.33250904,0.0213562701) (77.71661437,0.0183432967)
  (78.10071970,0.0153956205) (78.48482502,0.0125128947) (78.86893035,0.0096947530) (79.25303568,0.0069408109) (79.63714101,0.0042506668)
  (80.02124634,0.0016239034) (80.40535167,-0.0009399117) (80.78945700,-0.0034412240) (81.17356233,-0.0058804923) (81.55766766,-0.0082581867)
  (81.94177299,-0.0105747885) (82.32587832,-0.0128307889) (82.70998365,-0.0150266884) (83.09408898,-0.0171629961) (83.47819430,-0.0192402290)
  (83.86229963,-0.0212589113) (84.24640496,-0.0232195740) (84.63051029,-0.0251227538) (85.01461562,-0.0269689935) (85.39872095,-0.0287588405)
  (85.78282628,-0.0304928472) (86.16693161,-0.0321715700)
};
\addplot[no marks,line width=.88pt,color=slopeink] coordinates {
  (86.16693161,0.0281174585) (86.36186273,0.0267076790) (86.55679386,0.0253126313) (86.75172498,0.0239322877) (86.94665610,0.0225666193)
  (87.14158723,0.0212155963) (87.33651835,0.0198791876) (87.53144947,0.0185573615) (87.72638060,0.0172500851) (87.92131172,0.0159573245)
  (88.11624284,0.0146790449) (88.31117397,0.0134152107) (88.50610509,0.0121657852) (88.70103621,0.0109307310) (88.89596734,0.0097100096)
  (89.09089846,0.0085035818) (89.28582958,0.0073114075) (89.48076071,0.0061334458) (89.67569183,0.0049696549) (89.87062295,0.0038199923)
  (90.06555408,0.0026844147) (90.26048520,0.0015628779) (90.45541632,0.0004553370) (90.65034745,-0.0006382534) (90.84527857,-0.0017179399)
  (91.04020969,-0.0027837694) (91.23514082,-0.0038357897)
};
\addplot[no marks,line width=.88pt,color=slopeink] coordinates {
  (91.23514082,0.0539163674) (91.63040504,0.0507244381) (92.02566926,0.0475921537) (92.42093348,0.0445193860) (92.81619770,0.0415059851)
  (93.21146192,0.0385517804) (93.60672614,0.0356565805) (94.00199036,0.0328201744) (94.39725458,0.0300423321) (94.79251880,0.0273228047)
  (95.18778302,0.0246613254) (95.58304724,0.0220576102) (95.97831146,0.0195113579) (96.37357568,0.0170222514) (96.76883991,0.0145899577)
  (97.16410413,0.0122141286) (97.55936835,0.0098944015) (97.95463257,0.0076303995) (98.34989679,0.0054217324) (98.74516101,0.0032679971)
  (99.14042523,0.0011687777) (99.53568945,-0.0008763534) (99.93095367,-0.0028678351) (100.32621789,-0.0048061176) (100.72148211,-0.0066916610)
  (101.11674633,-0.0085249357) (101.51201055,-0.0103064213)
};
\addplot[no marks,line width=.88pt,color=slopeink] coordinates {
  (101.51201055,0.0427761019) (102.11750675,0.0386558205) (102.72300294,0.0346632529) (103.32849913,0.0307972228) (103.93399533,0.0270564777)
  (104.53949152,0.0234396936) (105.14498771,0.0199454797) (105.75048391,0.0165723834) (106.35598010,0.0133188948) (106.96147629,0.0101834510)
  (107.56697249,0.0071644406) (108.17246868,0.0042602078) (108.77796487,0.0014690569) (109.38346107,-0.0012107443) (109.98895726,-0.0037809600)
  (110.59445345,-0.0062433832) (111.19994965,-0.0085998315) (111.80544584,-0.0108521444) (112.41094203,-0.0130021791) (113.01643823,-0.0150518078)
  (113.62193442,-0.0170029145) (114.22743061,-0.0188573920) (114.83292681,-0.0206171393) (115.43842300,-0.0222840583) (116.04391919,-0.0238600520)
  (116.64941539,-0.0253470215) (117.25491158,-0.0267468641)
};
\addplot[no marks,line width=.88pt,color=slopeink] coordinates {
  (117.25491158,0.0202065762) (117.66818421,0.0184441802) (118.08145684,0.0167249519) (118.49472948,0.0150484091) (118.90800211,0.0134140679)
  (119.32127474,0.0118214426) (119.73454737,0.0102700460) (120.14782000,0.0087593899) (120.56109263,0.0072889855) (120.97436526,0.0058583435)
  (121.38763790,0.0044669744) (121.80091053,0.0031143890) (122.21418316,0.0018000983) (122.62745579,0.0005236142) (123.04072842,-0.0007155507)
  (123.45400105,-0.0019178825) (123.86727368,-0.0030838660) (124.28054632,-0.0042139843) (124.69381895,-0.0053087189) (125.10709158,-0.0063685491)
  (125.52036421,-0.0073939523) (125.93363684,-0.0083854032) (126.34690947,-0.0093433744) (126.76018211,-0.0102683355) (127.17345474,-0.0111607534)
  (127.58672737,-0.0120210923) (128.00000000,-0.0128498130)
};
\addplot[only marks,mark=*,mark size=1.35pt,mark options={draw=slopeink,fill=white,line width=.55pt}] coordinates {(76.18019305,-0.0350289717) (76.18019305,0.0307901179) (86.16693161,-0.0321715700) (86.16693161,0.0281174585) (91.23514082,-0.0038357897) (91.23514082,0.0539163674) (101.51201055,-0.0103064213) (101.51201055,0.0427761019) (117.25491158,-0.0267468641) (117.25491158,0.0202065762)};
\end{axis}
\end{tikzpicture}

\caption{\label{fig:dyadslopes}Loss and its first logarithmic
derivative on $64\le y\le128$. The five vertical lines mark the prime
cutoffs, and the open circles give the one-sided slopes. The six maxima
occur at the zeros of the smooth slope between cutoffs.}
\end{figure*}

The certificate differentiates the degree-$N=360$ polynomials of
Sec.~\ref{sec:losspowers} with $\mathfrak d=2(1-h)\partial_h$.
Put $s=3/4,\ r=7/8,\ n=N+1$. For $k=0,1,2,3$, the error in
the $k$-th $h$-derivative is at most
\[
 \frac{4096(s/r)^n n(n-1)\cdots(n-k+1)}
 {s^k[1-(s/r)(n+1)/(n+1-k)]}
 +\frac{k!\,10^{-165}}{(r-s)^k}.
\]
The first term sums the differentiated Cauchy coefficient bound as
a geometric majorant. The second follows from the uniform prime tail
on the disk. Subdivision certifies the sign of the third derivative
throughout each branch. Monotonicity of the curvature then reduces
its bounds to the endpoints. Bisection of the first derivative
isolates the six roots.

\head{The contribution of the moving cutoff at every order}
For $k\ge1$ put $T=\vartheta(a)/a$,
$M_{2k}=a^{-2k-1}\sum_{p\le a}p^{2k}\log p$ and $e_k=(k-1)/3$.
The series of Sec.~\ref{sec:losspowers} reads
\[
\begin{gathered}
 R=\sum_{k\ge1}w_k y^{-e_k}
          (T-M_{2k}-2k/(2k+1)),\\
 R_1=\tfrac74-\tfrac{21}{8}(T-M_2).
\end{gathered}
\]
On a branch, $T$ and $M_{2k}$ have logarithmic derivatives
$-\alpha T$ and $-(2k+1)\alpha M_{2k}$. Thus for $m\ge1$
\[
\begin{gathered}
 \mathfrak d^m R_1=-\tfrac{21}{8}
             [(-\alpha)^mT-(-3\alpha)^m M_2],\\
 \mathfrak d^m R=\kappa_m+o_m(1),\quad
 \kappa_m=\tfrac{21}{8}(-\alpha)^m(3^{m-1}-1).
\end{gathered}
\]
The limit is uniform over the open branches as $y\to\infty$.
Indeed $T\to1,\ M_2\to1/3$ by the prime number theorem and partial
summation. The differentiated $k\ge2$ tail is $O_m(y^{-1/3})$
since $T,M_{2k}<3$ and $|w_k|\le3(4\pi)^k/(\pi k!)$.
In particular $\dot R=o(1)$ while the regular curvature
$\ddot R\to\kappa_2=175/48$.

On the absolute logarithmic clock $u=\log y$, the first term
satisfies the distributional equation
\[
 (\partial_u+\tfrac56)(\partial_u+\tfrac52)R_1(e^u)
 =\tfrac{175}{48}-\tfrac{35}{8}
            \sum_p\frac{\log p}{p}\,\delta_{(6/5)\log p}.
\]
Indeed $\dot R_1+\tfrac52 R_1=-\tfrac{35}{8}(T-1)$, while
$(\partial_u+\alpha)T=\sum_p(\log p)/p\,\delta_{(6/5)\log p}$.
For the full $R$, the jump formula and partial summation in the
prime number theorem give a limiting mass $-\kappa_2\log2$ for
the jumps of $\dot R$ on $(\log X,\log(2X)]$. The integral of
the regular curvature on that interval tends to $\kappa_2\log2$.
The two limiting contributions to the change of $\dot R$ cancel.

\begin{proposition}[Interpolation at consecutive prime cutoffs]
\label{pr:diadprimemesh}
Let $p<q$ be consecutive primes and
$p^{6/5}\le y\le q^{6/5}$. Put
$\tau=\log(y/p^{6/5})/\log(q^{6/5}/p^{6/5})$.
Uniformly as $p\to\infty$,
\[
 J(y)=(1-\tau)J(p^{6/5})+\tau J(q^{6/5})+O(y^{-7/24}).
\]
Hence RH is equivalent to the existence of a constant $B$ such that
$J(p^{6/5})\ge-B$ for every sufficiently large prime.
\end{proposition}
\begin{proof}
Let $\mathcal A=-\tfrac12-v\partial_v$. The logarithmic derivative
$\mathfrak d^m U=y^{-1/2}\sum_p(\log p)\mathcal A^m W(p/y)$
is $O_m(\sqrt y)$ by Chebyshev's bound and partial summation.
Since $J=-U-R$, the regular curvature of $J$ is $O(\sqrt y)$.
The Baker--Harman--Pintz theorem gives
$q-p=O(p^{21/40})$ \cite{bhp01}. The logarithmic length of this
single branch is consequently $O(y^{-19/48})$. Linear interpolation
on a real interval of length $l$ has error at most
$l^2\sup|J_{uu}|/8$. The endpoint values are obtained by continuity.
The exponent is $1/2-2(19/48)=-7/24$. The equivalence follows
from \eqref{eq:oddprimefluct} and Theorem~\ref{th:gaussfirst}.
\end{proof}

\head{The signed prime sums for slope and curvature}
With $z_{pj}=\lambda_j(p/y)^2$, define
\[
\begin{gathered}
 D_m(y)=\frac{1}{2\pi\sqrt y}\sum_p(\log p)
                 \sum_{j=1}^3d_j e^{-z_{pj}}P_m(z_{pj}),\\
 P_{m+1}=(2z-\tfrac12)P_m-2zP_m',\quad P_0=1 .
\end{gathered}
\]
Then $D_m=\mathfrak d^m U$,
$P_1=2z-\tfrac12$ and $P_2=4z^2-6z+\tfrac14$.
In particular $\dot L_X=D_1+o(1)$ and
$\ddot L_X=D_2+175/48+o(1)$ inside the branches.

\head{Curvature and the zeros}
Let $S(y)=2E_\Phi(y)/\sqrt y$. For every fixed $m\ge1$,
\[
 \mathfrak d^m L_X(y)=\mathfrak d^m S(y)+\kappa_m+o_m(1)
 \qquad (y\to\infty,\ y\ne p^{6/5}).
\]
To prove this, put $G_m=\mathcal A^m W$. The normalized
contribution to $\mathfrak d^m(S-U)$ from prime squares tends to
$\int_0^\infty G_m(t^2)\,\dd t=0$ by the prime number theorem.
The integral vanishes since
$\int_0^\infty (\mathcal A f)(t^2)\,\dd t
=-\tfrac12[t f(t^2)]_0^\infty$.
The cumulative logarithmic weight of powers of order at least three
is $O(t^{1/3}\log(2t))$. Partial summation against $G_m(t/y)$
therefore costs $O_m(y^{-1/6}\log(2y))$ after normalization.

At each fixed order $m=2,\ldots,6$, RH is equivalent to the
eventual sign condition
\[
\begin{gathered}
 (-1)^m\mathfrak d^m L_X(y)>0\\
 \hbox{on every sufficiently large open branch}.
\end{gathered}
\]
In particular RH is equivalent to eventual convexity of the loss
in the boost time between successive prime cutoffs. The first
and seventh regular derivatives take both signs arbitrarily far out
unconditionally.

Under RH the normalized explicit formula of Theorem~\ref{th:gaussfirst} gives
\[
 \mathfrak d^m S(y)=-\sum_\rho
        (i\gamma)^m A(\rho)e^{i\gamma\log y}+O_m(y^{-3/2}).
\]
The gamma factor decays exponentially in $|\gamma|$, allowing
termwise differentiation of the sum and differentiation of the
remainder integral. For $m=1,\ldots,7$, the respective upper
bounds for $\sum_\rho|\gamma|^m|A(\rho)|$ are
\[
\begin{gathered}
 0.000222,\quad 0.003145,\quad 0.04461,\quad
 0.634,\\
 9.035,\quad 129.49,\quad 1878 .
\end{gathered}
\]
The first ten pairs are enclosed as in Theorem~\ref{th:oddglobal}.
For the others the bound used there is
\[
 (6\sqrt2/\pi)\,0.047\,50^{m+1}e^{-50\pi/6},
\]
since
$t^{m+1}e^{-\pi t/6}$ is decreasing for $t\ge50$ in these orders.
For orders two through six the bound is strictly below $|\kappa_m|$.
For orders one and seven, the amplitude of the first pair exceeds
the sum of $|\kappa_m|$ and the absolute amplitudes of all other
pairs by more than $0.000219$ and $1112$, respectively. Taking
the first phase alternately at its two extrema gives both signs
at arbitrarily large $y$. The strict inequalities also hold nearby
at points lying inside an open branch.

Conversely, the Laplace transform of
$u\mapsto \partial_u^m S(e^u)$ on $u\ge0$ agrees, up to an
entire function, with
$-z^m A(z+1/2)\zeta'(z+1/2)/\zeta(z+1/2)$.
This follows first for $\Re z>1/2$ by integration by parts.
If a zero lies to the right of the critical line, the transform
has a pole in $\Re z>0$, and is regular at all positive real
$z$. Landau's theorem applied to
$C\pm\partial_u^m S(e^u)$ then excludes every eventual one-sided
bound. Continuity of $\mathfrak d^m S$ and the uniform estimate
on the open branches transfer this conclusion to the regular
derivatives of $L_X$.

\head{An unconditional interval of positive curvature}
On every open branch in the explicit range
\[
 e^{200}\le y\le\exp(3\cdot10^{12}),\qquad
                        3.642<\ddot L_X(y)<3.650 .
\]
Here $200$ is a sufficient lower logarithmic scale. The following
bounds give constants. Put $t_0=e^{50}$, $\epsilon=4\cdot10^{-6}$.
Dusart's estimate $|D(t)|\le0.01t/\log^2t$ for
$t>7\,713\,133\,853$ gives $|D(t)|\le\epsilon t$ above $t_0$
\cite{dusarttheta}. Below $t_0$ use $|D(t)|\le2t$.
Partial summation in $M_2$ gives
\[
 |\ddot R_1-175/48|\le(175/96)\,[16\epsilon+12(t_0/a)^3].
\]
For the differentiated tail, the bracket multiplying each
$|w_k|y^{-e_k}$ is bounded in absolute value by
$3(e_k+\alpha)^2+3(2k+1/2)^2+e_k^2\le21k^2$.
Consequently $|\ddot R-\ddot R_1|<10^9 y^{-1/3}$, by summing
$(63/\pi)k^2(4\pi)^k/k!$.

For $G=\mathcal A^2 W$ each Gaussian is multiplied by
$4z^2-6z+1/4$, $z=\lambda_j v^2$.
With $g(v)=G(v^2)$, partial summation and $\int g=0$ bound
the prime-square contribution by
$23\epsilon+448t_0/\sqrt y$. Indeed
$\sup_{v>0}v|g'(v)|<224$ follows from $z^k e^{-z}\le k!$;
$\int_0^\infty v|g'(v)|\dd v<23$.
For higher powers the cumulative logarithmic weight through
$t\ge1$ is at most $4t^{1/3}\log t\le18t^{5/12}$.
Their normalized contribution is bounded by $18\cdot23\,y^{-1/12}$,
using $\int_0^\infty v^{5/12}|G'(v)|\dd v<23$.
The two integral bounds follow by substituting $z=\lambda_j v^4$
or $z=\lambda_j v^2$, respectively. Both are majorized by
\[
\begin{gathered}
 \frac{\sum_j|d_j|\lambda_j^{-s}}{2\pi}
       [4\Gamma(s+3)+14\Gamma(s+2)+\tfrac{25}{4}\Gamma(s+1)]\\
 (s=1/4,\ 5/24,\ \hbox{respectively}).
\end{gathered}
\]
Finally, through $T_*=3\cdot10^{12}$ the zeros are on the line
\cite{pt}. Their contribution to $\ddot S$ is smaller than
$0.003144361$. The high zeros cost less than $10^{-100}$
by the bound $(6\sqrt2/\pi)0.047(T_*^3+T_*/4)
\exp[-(\pi/6-1/2)T_*]$. The contour remainder is bounded by
$84y^{-3/2}$, using the line estimate of Theorem~\ref{th:oddglobal}
with the extra factor $9/4+t^2$.
These bounds give $|\ddot L_X-175/48|<0.003377$ on the stated range.

\subsection{Accumulated losses at consecutive primes}
\label{sec:primelossbudget}
Let $p<q$ be consecutive primes, $u_p=(6/5)\log p$,
$h_p=u_q-u_p$, $j_p=J(e^{u_p})$, and $V(u)=S(e^u)$.
Fix an initial prime $p_0$. The positive part of a real number is
$z_+=\max(z,0)$. Define the cumulative loss and gain by
\[
 N(P)=\sum_{\substack{p\ge p_0\\q\le P}}(j_p-j_q)_+,\qquad
 G(P)=\sum_{\substack{p\ge p_0\\q\le P}}(j_q-j_p)_+ .
\]
Both sums use complete consecutive-prime intervals. If $P\ge p_0$
is prime, $j_P=j_{p_0}+G(P)-N(P)$ exactly.
For example, with $p_0=31$ direct integration gives
$N(53)=0.00763860533971$, $G(53)=0.00043142602431$ and
$j_{53}=-0.06012330057973$, each within $5\cdot10^{-15}$.
The certificate evaluates the Gaussian primitive at the cutoffs.
The omitted primes above $4096$ contribute less than
$\sum_{n>4096}n e^{-\pi n/16}<10^{-340}$ at these nodes.

\begin{proposition}[A summable error for prime-by-prime induction]
\label{pr:primelefteuler}
One has $j_p-j_q=h_p V'(u_p)+\epsilon_p$, where
$\sum_{p\ge P}|\epsilon_p|\ll1/\log P$. Its arithmetic coefficient is
\[
 V'(u_p)=p^{-3/5}\sum_{n\ge2}\Lambda(n)
    \left(-\tfrac12 W(v)-vW'(v)\right)_{v=n/p^{6/5}}.
\]
The error from replacing the integral of $V'$ on each interval by
its left endpoint has absolute tail $O_\eta(P^{-17/100+\eta})$
for each $0<\eta<17/100$.
\end{proposition}
\begin{proof}
First consider the moving cutoff. Write $T_p=\vartheta(p)/p$,
$M_{2k,p}=p^{-2k-1}\sum_{\ell\le p,\ \ell\ {\rm prime}}
\ell^{2k}\log\ell$ and $r=p/q$. The exact recurrences are
\[
 T_q=rT_p+\frac{\log q}{q},\qquad
 M_{2k,q}=r^{2k+1}M_{2k,p}+\frac{\log q}{q}.
\]
Thus the new weight cancels in $T_q-M_{2k,q}$ at every order.
In particular the change of $R_1$ from $p$ to $q$ equals
$\tfrac{21}{8}\{(1-r)T_p-(1-r^3)M_{2,p}\}$.

Put $c=A(1/2)/2$ and $H(u)=J(e^u)+V(u)-c$.
For large $u$, uniformly on the open branches,
\[
 |H(u)|+|H'(u)|\ll u^{-2}.
\]
Indeed $H=S-U-c-R$ at $y=e^u$.
The estimate $D(t)=O(t/\log^2(2t))$ for $t\ge1$
\cite{dusarttheta} gives $T=1+O(u^{-2})$ and
$M_2=1/3+O(u^{-2})$ by partial summation at $a=y^{5/6}$.
The branch formulas of Sec.~\ref{sec:lossderivatives} then give
$|R|+|\dot R|\ll u^{-2}+y^{-1/3}$, including the one-sided limits.
For the prime squares put $g_m(v)=(\mathcal A^mW)(v^2)$ for
$m=0,1$, where $\mathcal A=-1/2-v\partial_v$.
After subtracting $\int_0^\infty g_m(v)\,\dd v$, the normalized
sum equals
$-\int_0^\infty D(\sqrt y\,v)g_m'(v)\,\dd v/\sqrt y$.
For $v<y^{-1/4}$ use $|D(t)|\le2t$ and $g_m'(v)=O(v^3)$.
For $v\ge y^{-1/4}$ use the logarithmic bound and
$\int v|g_m'(v)|\,\dd v<\infty$.
The two integrals of $g_m$ are $c$ and $0$, respectively.
Powers of order at least three contribute $O(y^{-1/6}\log(2y))$
at both orders. The continuity of $H$ and the finite number of
branches on each compact interval now give
$\sum_{p\ge P}|H(u_q)-H(u_p)|\ll1/\log P$.

Chebyshev's bound for $\sum_{n\le t}\Lambda(n)$ gives
$|V''(u)|\ll e^{u/2}$ by partial summation against $\mathcal A^2W$.
On $P\le p<2P$ one has $h_p\le(6/5)(q-p)/p$, and $q<2p$.
Stadlmann's unconditional mean-square gap theorem \cite{stadlmanngaps}
therefore gives
\[
\begin{aligned}
 \sum_{P\le p<2P}h_p^2 \sup_{[u_p,u_q]}|V''|
 &\ll P^{3/5-2}\sum_{p<2P}(q-p)^2\\
 &\ll_\eta P^{-17/100+\eta}.
\end{aligned}
\]
Integration on each interval, followed by a dyadic sum, proves the
claim since $j_p-j_q=V(u_q)-V(u_p)+H(u_p)-H(u_q)$.
The positive-part map is 1-Lipschitz, so the same quadrature error
applies to $(V')_+$ and $(-V')_+$. In particular, at prime $P$,
\[
\begin{aligned}
 N(P)&=\int_{u_{p_0}}^{u_P}(V'(u))_+\,\dd u+O(1),\\
 G(P)&=\int_{u_{p_0}}^{u_P}(-V'(u))_+\,\dd u+O(1).
\end{aligned}
\]
\end{proof}

\begin{theorem}[Growth of the losses]\label{th:primelossgrowth}
Let $\Theta=\sup_{\zeta(\rho)=0,\ 0<\re(\rho)<1}\re(\rho)$. Then
\[
 \lim_{P\to\infty}\frac{\log(1+N(P))}{\log P}
 =\lim_{P\to\infty}\frac{\log(1+G(P))}{\log P}
 =\tfrac65(\Theta-\tfrac12).
\]
Under RH one has $N(P)\sim G(P)\sim k\log P$ with
$8.43836\cdot10^{-5}<k<8.43890\cdot10^{-5}$.
In particular, the total loss is infinite unconditionally, and RH is
equivalent to $N(P_j)=P_j^{o(1)}$ along a single sequence
$P_j\to\infty$.
\end{theorem}
\begin{proof}
The explicit formula in Theorem~\ref{th:gaussfirst} gives
\[
 V(u)=-\sum_\rho A(\rho)e^{(\rho-1/2)u}+O(e^{-3u/2}).
\]
The exponential decay of $A$ gives $|V'(u)|\ll
e^{(\Theta-1/2)u}$ for $u\ge0$, since the contour estimate also
bounds the derivative of the remainder. Proposition~\ref{pr:primelefteuler}
proves the upper bound for each growth exponent.

Suppose $\Theta>1/2$ and fix $1/2<\sigma<\Theta$.
Ingham's density estimate \cite{ingham40} gives
$\sum_{\re(\rho)\ge\sigma,\ \Im\rho>0}1/\Im\rho<\infty$, since
$3(1-\sigma)/(2-\sigma)<1$. Choose the least ordinate
$\gamma_*>0$ in this set and the largest real part $\beta_*$
at that ordinate within the set. The positive averaging measure of
Theorem~\ref{th:fixedfourierwindow} applies to $V$ as follows.
For each distinct zero $\rho=\beta+i\gamma$ in the set with
$\gamma>\gamma_*$, put $d=\pi/\gamma$ and take masses
proportional to $e^{(\beta-\beta_*)d}$ and $1$ at $0,d$.
The convolution $\nu$ of these probabilities has compact support,
since $\sum d<\infty$. Each factor annihilates
$e^{(\rho-\beta_*)u}$ upon averaging translates. Its multiplier
at $i\gamma_*$ has modulus at least $\cos(\gamma_*d/2)$, and
the product of these bounds is positive since $\sum d^2<\infty$.
Apply $\nu$ to the translates of $e^{-(\beta_*-1/2)u}V(u)$.
Terms below $\sigma$ tend to zero after normalization by dominated
convergence, since their coefficients are absolutely summable and
each $\re(\rho)<\sigma\le\beta_*$. Those at $\gamma_*$ with smaller real part
tend to zero. The surviving conjugate pair gives a nonzero sinusoid
because $A(\beta_*+i\gamma_*)\ne0$. Over one period the average
takes both signs with amplitudes bounded away from zero. Positivity
therefore gives extrema of either sign and magnitude at least
$b X^{\beta_*-1/2}$ at scales $y\in[X,RX]$, for some $b>0$,
$R>1$ independent of large $X$.

Now $J(e^u)=c-V(u)+H(u)$. By
Proposition~\ref{pr:diadprimemesh}, at least one bracketing prime
cutoff attains each signed bound with $b$ replaced by $b/2$
for large $X$. The prime ratios $q/p$ tend to one.
Taking $X=P^{6/5}/(2R)$ puts both cutoffs at primes below $P$.
The total gain from $p_0$ dominates any upward displacement,
and the total loss dominates any downward displacement. Thus
$N(P),G(P)\gg_\sigma P^{(6/5)(\sigma-1/2)}$ at every sufficiently
large $P$. Let $\sigma\uparrow\Theta$.

On RH, $V'(u)=F_1(u)+O(e^{-3u/2})$, where $F_1$ is the
uniformly convergent trigonometric series with coefficients
$-i\gamma A(1/2+i\gamma)$. Its mean is zero.
The means of $(F_1)_+$ and $(-F_1)_+$ exist and both equal
$m=\tfrac12 M|F_1|$, with
$Mf=\lim_{L\to\infty}L^{-1}\int_0^L f(u)\,\dd u$.
Indeed uniform almost periodicity is preserved by the continuous
positive-part map, or directly by polynomial approximation on the
bounded range. Proposition~\ref{pr:primelefteuler} gives $k=6m/5$.

For the numerical enclosure write $F_1=f+g+t$, where $f$
is the first sinusoid, $g$ comes from positive ordinates two
through ten and their conjugates, and $t$ is the tail.
Write $a=2\gamma_1|A(\rho_1)|$,
$d=\sum_{j=2}^{10}2\gamma_j|A(\rho_j)|$, and
$b_{50}=(6\sqrt2/\pi)\,0.047\,50^2e^{-50\pi/6}$.
Then $|t|\le b_{50}$ by Sec.~\ref{sec:lossderivatives}.
None of $\gamma_2,\ldots,\gamma_{10}$ is an odd multiple of
$\gamma_1$, by the enclosures for these ordinates.
Integration of each harmonic over successive periods of $f$ gives
$M(\operatorname{sgn}(f)g)=0$. Noninteger frequency ratios give a
bounded geometric sum. Half-period antisymmetry gives zero at even
integer multiples. Pointwise,
\[
 0\le |f+g|-|f|-\operatorname{sgn}(f)g
       \le 2d\,\mathbf1_{\{|f|\le d\}} .
\]
The indicator has mean $(2/\pi)\arcsin(d/a)$. Thus
\[
 \frac{a}{\pi}-\frac{b_{50}}{2}\le m
 \le\frac{a}{\pi}+\frac{2d}{\pi}\arcsin(d/a)+\frac{b_{50}}{2}.
\]
Interval arithmetic gives the quoted bounds on $6m/5$.
This also proves a positive logarithmic lower bound if RH holds.
\end{proof}

\head{A fixed numerical threshold}
Set $B_0=0.018606662$. For every real $y$ with
$10^5\le\log y\le3\cdot10^{12}$, one has $J(y)>-B_0$
unconditionally. If RH holds the same bound holds for all
$\log y\ge10^5$. To check the constants, take $t_0=y^{2/5}$
and $\epsilon=1/(16\log^2y)$, so that $|D(t)|\le\epsilon t$
for $t\ge t_0$ by Dusart's bound. The identity
$T-M_2-2/3=2a^{-3}\int_0^a t D(t)\,\dd t$ and the Taylor tail in
Sec.~\ref{sec:losspowers} give
\[
 |R|\le\tfrac74\epsilon+\tfrac72 y^{-13/10}+45y^{-1/3}.
\]
For the prime squares, $\int v|g_0'(v)|\,\dd v<1$ and
$\sup v|g_0'(v)|<4$ give error at most
$\epsilon+8y^{-1/10}$ after splitting at $t_0/\sqrt y$.
For higher powers the cumulative weight is at most $18t^{5/12}$
and $\int v^{5/12}|W'(v)|\,\dd v<1$, giving $18y^{-1/12}$.
The two integral bounds follow from
$\Gamma(1+s)\sum_j |d_j|\lambda_j^{-s}/(2\pi)<1$ for
$s=1/4,5/24$. The supremum bound follows from $ze^{-z}\le1$.
Combining the estimates gives
\begin{multline*}
 |H(\log y)|\le\frac{11}{64\log^2 y}
       +\tfrac72 y^{-13/10}+45y^{-1/3}\\
       +8y^{-1/10}+18y^{-1/12}.
\end{multline*}
On the specified finite range the ordinate bound \cite{pt} and
the proof of Theorem~\ref{th:oddglobal} give
$|V(\log y)|<2(7.8390\cdot10^{-6})+10^{-100}+16y^{-3/2}$.
Under RH the contribution $10^{-100}$ can be omitted at every
scale. Substituting $c=-0.01859098397459209\ldots$ proves the
numerical assertion.

In consequence there is some fixed $D_0>1$ such that every
sufficiently large $[P,D_0P]$ contains a prime $p$ with
$j_p>-B_0$, unconditionally. Under RH this follows at every
node. If RH fails, the window argument above supplies both a node
above $-B_0$ and one below $-B_0$. Enlarging the fixed window
absorbs the two bracketing endpoints. For the stronger assertion at
every sufficiently large prime the induction requires
\[
 \sup_{P\ge p_0,\ P\ {\rm prime}}\;
       \sum_{p_0\le p<P}h_p V'(u_p)<\infty .
\]
This is equivalent to the existence of a uniform lower bound for
$j_P$ by Proposition~\ref{pr:primelefteuler}. The gap theorem
controls the sum of the absolute quadrature errors. The signed
partial sums in this last display still require an arithmetic bound.

For the prime-number-theorem estimate, set $K=\mathcal A W$, so
$\int_0^\infty K(v)\,\dd v=0$, and
$\psi(t)=\sum_{n\le t}\Lambda(n)$. Partial summation gives
\[
 V'(u)=-e^{-u/2}\int_0^\infty
                [\psi(e^uv)-e^uv]K'(v)\,\dd v.
\]
Put $\omega(u)=u^{3/5}(\log u)^{-1/5}$ for $u>1$.
The Vinogradov--Korobov estimate \cite{johnstonyang} bounds
$|\psi(t)-t|/t\ll\exp[-c_1\omega(\log t)]$ as $t\to\infty$.
For $v<e^{-u/2}$, Chebyshev's bound and $K'(v)=O(v)$ give
$O(e^{-u})$ for that part of $V'$. On the complementary interval
$\log(e^uv)\ge u/2$ and $\int_0^\infty v|K'(v)|\,\dd v<\infty$.
Hence $|V'(u)|\ll \exp[u/2-c_2\omega(u)]$ at large $u$.
Integrating on $[u_{p_0},u_P]$, split at $u_P/2$ to obtain
\[
 N(P)+G(P)\ll P^{3/5}\exp[-c_3\omega(\log P)]
\]
for positive constants $c_1,c_2,c_3$. The logarithm of the right-hand
side divided by $\log P$ tends to $3/5$. The criterion of
Theorem~\ref{th:primelossgrowth} requires a limit of zero.

\subsection{The residue carried by the signed prime increments}
\label{sec:primeresidues}
The available margin can be tested directly on the unmodified prime
weights. With $B_0=0.018606662$, consecutive cutoffs $19$ and
$23$ give
\[
\begin{gathered}
 j_{19}=-0.0128397306998710,\\
 j_{23}=-0.0355120997717638
\end{gathered}
\]
within $6\cdot10^{-17}$ in each value. The margin $B_0+j_{19}$
is less than $0.005767$, whereas the next loss exceeds $0.022672$.
These nodes lie below the range of the uniform numerical bound.
On the infinite tail the transform of the signed increments gives
a separate test of the inductive argument.

\begin{proposition}[Transform of the actual increments]
\label{pr:signednoderesidue}
For $\Re z>1/2$ form the absolutely convergent series
\[
\begin{gathered}
 \mathcal D(z)=\sum_{p\ge p_0}(j_p-j_q)p^{-6z/5},\\
 q=\min\{\ell>p:\ell\ {\rm prime}\}.
\end{gathered}
\]
It continues meromorphically to $\Re z>0$. With $s=1/2+z$,
\[
\begin{gathered}
 \mathcal D(z)=-zA(s)\frac{\zeta'(s)}{\zeta(s)}+E(z),\\
 E\ \hbox{holomorphic for }\Re z>0 .
\end{gathered}
\]
The error extends continuously to $\Re z=0$.
In particular a zero $\rho$ of multiplicity $m$ to the right
of the critical line gives the nonzero residue
$-m(\rho-1/2)A(\rho)$ at $z=\rho-1/2$.
\end{proposition}
\begin{proof}
Use the exact equality
\[
 j_p-j_q=\int_{u_p}^{u_q} V'(u)\,\dd u+H(u_p)-H(u_q)
\]
of Proposition~\ref{pr:primelefteuler}. The series from the $H$
differences is absolutely convergent for $\Re z\ge0$.
The derivative satisfies $|V'(u)|\ll e^{u/2}$. Since $h_p$
is bounded, these facts give the stated absolute convergence.
To replace $e^{-zu_p}$ by $e^{-zu}$ inside the integral on
$[u_p,u_q]$, use the mean-value bound on the exponential.
On compact sets of $z$, the total error from $P\le p<2P$ is
\[
\begin{aligned}
 &\ll (1+|z|)P^{-6\Re z/5}
       \sum_{P\le p<2P}p^{3/5}h_p^2\\
 &\ll_\eta (1+|z|)P^{-6\Re z/5-17/100+\eta}.
\end{aligned}
\]
The implied constants may also depend on the compact set. The last
estimate uses the mean-square gap bound \cite{stadlmanngaps}. Summing
dyadically gives a holomorphic error even on $\Re z>-17/120$.
As $u\to-\infty$, $V'(u)$ decays faster than $e^{bu}$ for
each $b>0$, by its prime-power Gaussian series. Extending the
integral to the whole line therefore costs an entire function.
Finally the bilateral Laplace transform of $V'$, initially on
$\Re z>1/2$, is $z$ times the Mellin transform of
$\sum_n\Lambda(n)W(n/y)$ at $s$. This is
$-zA(s)\zeta'(s)/\zeta(s)$ by Theorem~\ref{th:gaussfirst}.
The singularity at $s=1$ is removable and $A$ has no zeros in
$0<\Re s<1$.
\end{proof}

For each real $\gamma>0$, let $m_\gamma$ denote the multiplicity
of $1/2+i\gamma$ as a zero, with value zero at a regular point.
For the continued function the boundary limit is
\[
 \lim_{\eta\downarrow0}\eta\,\mathcal D(\eta+i\gamma)
       =-i\gamma\,m_\gamma A(1/2+i\gamma).
\]
This follows from the local Laurent expansion and the continuity of
$E$. It recovers each atom of the boundary measure of
Theorem~\ref{th:boundaryclosure} without assuming RH. For example,
the modulus of this limit at the first ordinate lies between
$0.00011045903084585$ and $0.00011045903084586$.
The original prime series defines $\mathcal D$ directly in its
half-plane of convergence.

\head{The loss forced by a single off-line zero}
Suppose $\rho=\beta+i\gamma$ is a zero with $\beta>1/2$, of
multiplicity $m$. Put $b=\beta-1/2$. At the actual prime nodes,
\[
 \liminf_{p\to\infty}\frac{j_p}{p^{6b/5}}\le -m|A(\rho)|,\qquad
 \limsup_{p\to\infty}\frac{j_p}{p^{6b/5}}\ge m|A(\rho)|.
\]
This holds for each such zero, including when other zeros lie farther
to the right. To prove it, let
$T(u)=\sum_{u_p\le u,\ p\ge p_0}(j_p-j_q)$. It is zero for
$u<u_{p_0}$, and its Laplace transform on $[0,\infty)$ is
$\mathcal D(z)/z$, initially for $\Re z>1/2$.
Write $X=T$ or $X=-T$. If $X(u)\le B e^{bu}$ eventually,
the tail of $B e^{bu}-X(u)$ is nonnegative. Its continued transform
is regular at every real $z>b$, so Landau's theorem gives absolute
convergence for $\Re z>b$. Call the tail transform $G$.
At $z=b$ its residue is $B$, since $\mathcal D(z)/z$ is
regular there. At $b+i\gamma$ its residue has modulus $m|A(\rho)|$.
For $\epsilon>0$, positivity gives
\[
 |G(b+\epsilon+i\gamma)|\le G(b+\epsilon).
\]
Multiplication by $\epsilon$ and passage to the limit force
$B\ge m|A(\rho)|$. Applying this to both choices of sign proves
the two bounds for $T$. Finally $T(u)=j_{p_0}-j_q$ on
$[u_p,u_q)$, and $u_q-u_p\to0$, which transfers them to $j_q$
with the displayed normalization.

\head{Contour representation of the spectral mass defect}
Let $R$ be a closed rectangle in
$\{0<\Re z<1/2,\ \Im z>0\}$ with no pole on its boundary, oriented
counterclockwise. Its contribution $\Delta_R$ to the mass defect,
including all four members of each associated quartet, is
\[
\begin{aligned}
 \Delta_R&=-4\Re\!\left[\frac{1}{2\pi i}
     \oint_{\partial R}
       \frac{\mathcal D(z)}{zA(1/2+z)(1/4-z^2)}\,\dd z\right]\\
   &=4\sum_{\rho-1/2\in R}m_\rho\,
           \Re\frac{1}{\rho(1-\rho)} .
\end{aligned}
\]
The sum ranges over distinct zeros in the indicated rectangle.
The identity follows by the residue theorem and
$1/4-(\rho-1/2)^2=\rho(1-\rho)$.
Every summand is positive by Theorem~\ref{th:boundaryclosure}.
An increasing exhaustion by such rectangles gives
$\Delta_R\uparrow\Delta$. Equivalently,
$L(z)+\mathcal D(z)/(zA(1/2+z)(1/4-z^2))$ is holomorphic in
this open strip, where $L$ is the transform of the arithmetic
correlation $C$ used in the Majorana reconstruction of
Sec.~\ref{sec:boundaryclosure}. The two functions have opposite
principal parts at each pole.

Let $\mathcal D_P$ be the finite sum over $p_0\le p\le P$.
Its contour integral with the same weight is zero for every $P$.
If $R$ contains a zero, then for every $P$
\[
 \int_{\partial R}
 \frac{|\mathcal D(z)-\mathcal D_P(z)|}{|zA(1/2+z)(1/4-z^2)|}
             \,|\dd z|\ \ge\ \frac{\pi}{2}\Delta_R>0 .
\]
This is the triangle inequality applied to the preceding identity.
Thus passage from the finite induction to a vanishing mass requires
justification on contours inside the strip. A summable error in the
quadrature leaves each nonzero residue of the full series intact.

\head{Stability with respect to the first prime-return amplitudes}
Keep the prime locations fixed and replace $\log p$ in
$\vartheta$ by a real weight $b_p$. Denote the resulting integral
by $J_b$. If $0<d<1$ and $|b_p-\log p|\le C p^{-d}\log p$, then
\[
 |J_b(y)-J(y)|\ll_d C y^{1/2-d}\qquad(y\ge1).
\]
Chebyshev's bound and partial summation bound the difference of
the two counting functions through $t$ by $O_d(C t^{1-d})$.
Insertion under the integral defining $J$ gives the result since
$\int_0^\infty v^{1-d}|W'(v)|\,\dd v<\infty$. In particular
relative errors $O(p^{-1/2})$ in the amplitudes preserve the
existence of a uniform lower bound.

The exponent is sharp for the size of the change. Fix $M>2^{1/d}$
and keep every prime weight through $M$. Beyond $M$ set
\[
 b_p=(\log p)\,[1+2p^{-d}\cos(\log p)] .
\]
All weights are positive. Put $s_d=1-d+i$. For any fixed $k>0$,
partial summation in the prime number theorem gives
\[
\begin{gathered}
 \sum_{M<p\le t}(\log p)\,p^{-d+i}
       =t^{s_d}/s_d+O_{M,d,k}\big(t^{1-d}/\log^k(2t)\big)\\
       (t\to\infty).
\end{gathered}
\]
Substitute this expression and its conjugate in $J_b-J$.
Replacing the lower limit $y^{5/6}$ by zero in the main term
costs $O_d(y^{-5d/6})$, since $W'(v)=O(v)$ at zero. On the
remaining integration range $\log t\ge(5/6)\log y$.
Integration by parts gives
$\int_0^\infty v^{s_d}W'(v)\,\dd v=-s_d A(s_d)$. It follows that
\begin{multline*}
 J_b(y)-J(y)=-2y^{1/2-d}\Re[A(s_d)e^{i\log y}]\\
           +O_{M,d,k}\big(y^{1/2-d}/\log^k y\big).
\end{multline*}
At $y=p^{6/5}$ the phase tends to infinity with gaps tending to
zero. The normalized difference has upper and lower limits
$\pm2|A(s_d)|$. For $d=1/4$ the oscillation therefore has
amplitude $2|A(3/4+i)|\,p^{3/10}$, where
$0.1105884735098145<2|A(3/4+i)|<0.1105884735098146$.

For $0<d<1/2$, this weighted model violates every eventual lower
bound independently of the behavior of $J$. Its Euler product is
\[
 Z_M(s)=\zeta(s)\prod_{\epsilon=\pm1}
 \left[\zeta(s+d+\epsilon i)
       \prod_{p\le M}(1-p^{-s-d-\epsilon i})\right].
\]
The coefficients of $-Z_M'/Z_M$ at $p^r$ equal $\log p$ for
$p\le M$ and $(\log p)[1+2p^{-rd}\cos(r\log p)]$ for $p>M$.
They are positive and bounded by $3\log p$. The higher powers
contribute $O(1)$ to the normalized Gaussian sum by Chebyshev's
bound for the squares and partial summation for the rest.
The compact correction $R_b=O(1)$ follows from
$|\vartheta_b(t)-t|\ll t$ and $W'(v)=O(v)$ near zero.
Consequently $J_b$ differs by $O(1)$ from the negative of the
full normalized Gaussian sum for these coefficients.

The unnormalized Gaussian sum has Mellin transform
$-A(s)Z_M'(s)/Z_M(s)$. It is regular at real $s>1/2$, including
$s=1$ because $A(1)=0$. The shifted zeta factors have no zeros
at ordinate $1$ in this region \cite{pt}. At $s_d$ one shifted
factor has a simple pole and the others are finite and nonzero.
The Mellin residue is $+A(s_d)\ne0$. An eventual lower bound for
$J_b$ would give an eventual $O(\sqrt y)$ upper bound for the
unnormalized Gaussian sum. Landau's nonnegative-tail argument in
Theorem~\ref{th:gaussfirst} would then exclude this pole.
The same conclusion holds at the prime cutoffs. The proof of
Proposition~\ref{pr:diadprimemesh} still applies, since
$|b_p|\le3\log p$ and the gaps are unchanged. In detail the smooth
prime sum has regular curvature $O(\sqrt y)$ by partial summation,
and the moment expansion of $R_b$ has bounded regular curvature
because $\vartheta_b(a)/a$ and $M_{2k,b}$ are bounded.
More quantitatively, the same one-sided residue comparison, applied
to the model's Gaussian sum at $s_d$, gives
$\liminf_{p\to\infty}J_b(p^{6/5})/p^{(6/5)(1/2-d)}
\le-|A(s_d)|$. Here $O(1)$ vanishes after normalization and
linear interpolation transfers the bound to a bracketing node.
The counting functions for the first- and second-winding weights
each equal $t+O_{M,d,k}(t/\log^k(2t))$ at large $t$.
Indeed their changes from $\vartheta(t)$ are
$O_d(t^{1-d})$ and $O_d(t^{1-2d})$, respectively, for
$0<d<1/2$. The proof of Proposition~\ref{pr:primelefteuler}
therefore also gives summable induction errors for this model.

On any prescribed initial interval $1\le y\le Y$, taking
$M\ge Y^{5/6}$ leaves the compact correction identical. The
remaining prime tail satisfies
\[
 \sup_{1\le y\le Y}|J_b(y)-J(y)|
 \le 6Y^2\exp[-\pi M^2/(8Y^2)] .
\]
Indeed $|W(v)|\le(3/\pi)e^{-\pi v^2/4}$. After extracting half
the exponential at $M$, bound $\log p\le p$ and use
$\sum_{n\ge1}n e^{-an^2}\le\int_0^\infty(t+1)e^{-at^2}\,\dd t$
with $a=\pi/(8Y^2)$. For $d=1/4$, $Y=\exp(3\cdot10^{12})$
and $M=4Y\sqrt{\log Y}$, this bound is less than $10^{-100}$.
The certified lower bound in Sec.~\ref{sec:primelossbudget} exceeds
$-B_0$ by more than $8\cdot10^{-12}$ on its finite range.
The weighted model satisfies the certified finite-range inequality
and has arbitrarily large prime nodes below every fixed lower bound.

The modified weights agree with the first-return amplitudes through
$M$. For the von Mangoldt weights, the identity
$\sum_{d\mid n}\Lambda(d)=\log n$ fixes the amplitudes at all
prime powers. The construction exhibits unbounded loss even with
the prescribed prime locations, positive weights, prime number
theorem errors of every fixed logarithmic order, and summable
quadrature errors. Proposition~\ref{pr:signednoderesidue} expresses
the corresponding loss for the von Mangoldt weights in terms of
Mellin residues.

\section{Prime correlations, Laguerre moments and curvature in the degree}\label{app:primegamma}
\subsection{A one-sided power series for the prime-node bound}
\label{sec:positivepowers}
Put $g(u)=V'(u)$, $x_+=\max(x,0)$, and
$K(v)=-W(v)/2-vW'(v)$. Differentiating the prime-return series gives
\[
 g(u)=e^{-u/2}\sum_{m\ge2}\Lambda(m)K(me^{-u}).
\]
For any real moment order $r>0$ and integer $n\ge0$, define
\[
\begin{gathered}
 w_{n,r}(u)=\frac{r^{n+1}u^n e^{-ru}}{n!},\\
 a_n(r)=\int_0^\infty w_{n,r}(u)(g(u)_+)^r\,\dd u,\\
 G_r(z)=\sum_{n\ge0}a_n(r)z^n .
\end{gathered}
\]
The weight is a probability density. Each coefficient is finite
by $|g(u)|\ll e^{u/2}$. It is strictly positive since the
slope takes positive values on open intervals.

\begin{theorem}[Successive coefficients of a one-sided moment]
\label{th:positivepowerratio}
Put $b=\sup_\rho\re(\rho)-1/2$, with the supremum over nontrivial
zeros. For every $r>0$,
\[
 \lim_{n\to\infty}\frac{a_{n+1}(r)}{a_n(r)}
 =\lim_{n\to\infty}a_n(r)^{1/n}=\frac{1}{1-b}.
\]
In particular $G_r$ has radius $1-b$. For a single fixed $r$,
a subsequence of successive-coefficient ratios at degrees tending
to infinity with upper limit at most one implies RH.
For $r=1$ the series supplies the bound
on the prime nodes explicitly. With $H(u)=J(e^u)+V(u)-c$ as
above and $0<\eta<1$, whenever $G_1(1-\eta)<\infty$ one has
\[
 j_p\ge c+H(u_p)-V(0)-p^{6\eta/5}G_1(1-\eta).
\]
\end{theorem}
\begin{proof}
For real $z\ge0$ Tonelli's theorem, allowing $+\infty$, gives
$G_r(z)=r\int_0^\infty e^{-r(1-z)u}(g(u)_+)^r\,\dd u$.
For $r=1$, integrate $V'$ from $0$ to $u_p$ and majorize
its positive part by $e^{\eta u_p}e^{-\eta u}g(u)_+$ there.
The asserted estimate follows. If the series is finite for
arbitrarily small positive $\eta$, the boundedness of $H$ and a
choice $0<\eta<\min(1,5\varepsilon/6)$ give the required exponent
$\varepsilon$ for every $\varepsilon>0$.

Put $M_r(U)=\int_0^U(g(u)_+)^r\,\dd u$. If $b>0$,
Proposition~\ref{pr:primelefteuler} and
Theorem~\ref{th:primelossgrowth} give $M_1(U)=\exp(bU+o(U))$.
Monotonicity extends this from prime times to arbitrary $U$, since
successive prime times have gaps tending to zero. The explicit
formula gives $|g(u)|\ll e^{bu}$. For $r\ge1$, H\"older's
inequality gives $M_r(U)\ge U^{1-r}M_1(U)^r$. For $0<r<1$
use $(g(u)_+)^r\ge (C e^{bU})^{r-1}g(u)_+$ on $[0,U]$.
Together with the pointwise bound these estimates give
$\log M_r(U)=rbU+o(U)$.

The upper exponential rate for $a_n(r)$ follows by integrating
$C^r e^{rbu}$ against $w_{n,r}$. For the lower rate take
$U=n/[r(1-b)]$ and fix $0<h<1$. The increment of $M_r$
on $[(1-h)U,U]$ equals $\exp(rbU+o(U))$. On that interval
$u^ne^{-ru}\ge(1-h)^nU^ne^{-rU}$. Stirling's formula gives
$\liminf a_n(r)^{1/n}\ge(1-h)/(1-b)$. Let $h\downarrow0$.
When $b=0$, RH holds by symmetry. The explicit formula for $g$
then consists of a uniformly convergent trigonometric series and
an $O(e^{-3u/2})$ remainder. The first sinusoid exceeds in
amplitude the sum of the remaining amplitudes by
Sec.~\ref{sec:lossderivatives}. The positive part to any fixed
power $r>0$ has a strictly positive time mean, since it is bounded
below on a fixed fraction of each period of the first sinusoid at
sufficiently large times. It is a uniformly almost periodic
function up to an $o(1)$ remainder. Its average against $w_{n,r}$
tends to the same mean, since a mode of frequency $\omega$ has
multiplier $(r/(r-i\omega))^{n+1}$. Uniform approximation by
trigonometric polynomials justifies the limit. This proves the
root limit in both cases.

For the ratio, write $m_n=n!a_n(r)/r^{n+1}$.
These moments of a positive measure are log convex by
Cauchy--Schwarz. With $L=(1-b)^{-1}$, the root limit gives
$\log m_n=n\log n+n(\log(L/r)-1)+o(n)$.
For fixed $0<h<1$, bound
$\log(m_{n+1}/m_n)$ by the secant slopes from
$\lfloor(1-h)n\rfloor$ to $n$ and from $n$ to
$\lceil(1+h)n\rceil$. Subtracting $\log n+\log(L/r)$ and
letting $h\downarrow0$ after $n\to\infty$ leaves zero in both
bounds. Thus $r m_{n+1}/((n+1)m_n)\to L$.
The exact induction in the degree is
\[
 a_{n+1}(r)-a_n(r)
 =\int_0^\infty w_{n,r}(u)
       \left(\frac{ru}{n+1}-1\right)(g(u)_+)^r\,\dd u .
\]
\end{proof}

The indices $1/2,3/2,5/2,7/2$ give four distinct positive
moments of the same return signal and have the same radius.
The order $r$ is held fixed in the ratio test. For example, at
$r=n/L$ with fixed $L>0$ the gamma density has mean
$L(1+1/n)$ and variance $L^2(n+1)/n^2$.
At $r=0$ the displayed gamma density loses its normalization.
If a fixed rate of one is used and the observable is assigned the
constant zeroth power one, every coefficient is one and the series
is $1/(1-z)$. If the power is defined as
$\mathbf 1_{\{g>0\}}$, its coefficients lie between zero and
one. The zero-order convention discards the size of a loss.

\begin{proposition}[Sign of every signed coefficient]
\label{pr:signedmomentdegrees}
Define $a_n^-=\int_0^\infty w_{n,1}(u)(-g(u))_+\,\dd u$ and
$s_n=a_n(1)-a_n^-$. Unconditionally,
\[
 s_0>0,\qquad s_n<0\quad(n\ge1),\qquad
 s_n\sim-\frac{15}{28}\,n(2/7)^n .
\]
The signed series $\sum_{n\ge0}s_n z^n$ has radius $7/2$.
\end{proposition}
\begin{proof}
The signed coefficient is the integral of $g$ against $w_{n,1}$.
It has a generating function directly in terms of the primes.
For $|z|<1/2$, put $s=3/2-z$. Integration by parts and the
Mellin transform of $W$ give
\begin{multline*}
 \sum_{n\ge0}s_n z^n
 =(1-z)A(s)\sum_{m\ge2}\frac{\Lambda(m)}{m^s}-V(0)\\
   -(1-z)\int_{-\infty}^0 e^{-(1-z)u}V(u)\,\dd u .
\end{multline*}
Indeed the bilateral transform of $V$ is the prime series times
$A(s)$ for $\Re s>1$, by absolute convergence. The integral
over negative times is entire in $z$ by Gaussian decay as
$u\to-\infty$. In particular the coefficient of degree $k$
of the arithmetic series alone is
$\sum_{m\ge2}\Lambda(m)(\log m)^k/(k!\,m^{3/2})$.
The inverse Mellin integrand for $g$ has factor
$-(s-1/2)A(s)\zeta'(s)/\zeta(s)$. Put
$h_\rho=-(\rho-1/2)A(\rho)$, counting zeros with multiplicity.
Move the contour to $\Re s=-3$. In addition to the nontrivial
zeros, the double pole at $-2$ contributes
$e^{-5u/2}(d_2u+d_1)$, where
\[
\begin{gathered}
 d_2=-105/16,\\
\begin{aligned}
 d_1&=\tfrac52\Big[
 \tfrac{21}{16}(\gamma_{\rm E}-1+\log\pi)+\tfrac{31}{8}\log2\\
 &\qquad-\tfrac{21}{8}\frac{\zeta''(-2)}{2\zeta'(-2)}\Big]+\tfrac{21}{8}\\
 &=4.6218022491846\ldots .
\end{aligned}
\end{gathered}
\]
Here $\gamma_{\rm E}$ denotes Euler's constant. The first
coefficient uses $\operatorname{Res}_{s=-2} A(s)=-21/8$.
On the new line
\[
\begin{gathered}
 |A(-3+it)|\le
 \frac{20\sqrt2\,\pi e^{-\pi|t|/4}}{\sqrt{(9+t^2)(1+t^2)}},\\
 |(\zeta'/\zeta)(-3+it)|\le7+|t|.
\end{gathered}
\]
The first bound follows by two gamma recurrences to $1/2+it/2$.
For the second, apply the functional equation and the digamma
difference series at $4$. The cotangent has modulus at most
one, the log derivative at $4-it$ has modulus less than one,
and $|\psi_{\Gamma}(4-it)|\le2+|t|$, with
$\psi_\Gamma=\Gamma'/\Gamma$. Including $|s-1/2|$ and the
inverse Mellin factor gives contour constant at most
$20\sqrt2\sqrt{50}(7/6)(4/\pi)<300$.
Consequently integration against the gamma density gives
\[
\begin{gathered}
\begin{aligned}
 s_n&=(2/7)^{n+1}\left[d_1-\tfrac{15}{8}(n+1)\right]\\
       &\quad+\sum_\rho\frac{h_\rho}{(3/2-\rho)^{n+1}}+e_n,
\end{aligned}\\
 |e_n|\le300(2/9)^{n+1}.
\end{gathered}
\]
For the sign it suffices to use $|\Im\rho|>14$ and the
reciprocal-square bound of Theorem~\ref{th:oddglobal}. Its gamma
estimate and monotonicity of $t^2e^{-\pi t/6}$ for $t\ge14$ give
\[
 \sum_\rho |h_\rho|\le
 \frac{6\sqrt2}{\pi}\sqrt{1+\frac1{4\cdot14^2}}\,
 0.047\cdot14^2 e^{-14\pi/6}<0.017 .
\]
After division by $(2/7)^{n+1}$ the total error
at $n=11$ is at most $300(7/9)^{12}+0.017\,4^{-12}$,
smaller than $(15/8)12-d_1$. The error bound decreases and the
latter margin increases. This proves the sign for all $n\ge11$
and gives the asymptotic formula, hence the radius by the root test.

For the eleven initial signs the certificate evaluates the residues
at $-2j$, $1\le j\le100$. For reproducibility put
$P(s)=2^s+2^{1-s}-3$, $u_j=2j+1/2$, and
\[
\begin{aligned}
 w_j&=\frac{(-1)^j\pi^{j-1}P(-2j)}{2j!},\\
 f_j&=w_j\left[\tfrac12(\psi_\Gamma(j+1)-\log\pi)
                    +P'(-2j)/P(-2j)\right].
\end{aligned}
\]
Set $v_j=\log(2\pi)-\psi_\Gamma(1+2j)-(\zeta'/\zeta)(1+2j)$.
The contribution to $s_n$ is
\[
 \frac{(n+1)u_jw_j}{(u_j+1)^{n+2}}
       +\frac{u_j(f_j+w_jv_j)-w_j}{(u_j+1)^{n+1}} .
\]
The functional equation identifies $v_j$ with the finite part
of $\zeta'/\zeta$ at $-2j$. After $M$ residues, the remaining
line contributes at most $C_M(2M+5/2)^{-n-1}$. For $M\ge1$
the same gamma recurrences and functional-equation estimate give
\begin{multline*}
 C_M\le \frac{4}{\pi^2}\,
 \frac{2^{M-1}\sqrt2\,\pi^M}{(2M+1)!!}\,(2^{-2M-1}+2^{2M+2}+3)\\
 \times(2M+3/2)\sqrt{(6+\log(2M+2))^2+1}.
\end{multline*}
Indeed the gamma denominator contains
\[
 \prod_{j=0}^M\sqrt{(2j+1)^2+t^2}.
\]
Its first factor absorbs the
log derivative, bounded by $6+\log(2M+2)+|t|$, and its last
absorbs $|s-1/2|$. The displayed bound at $M=100$ is less
than $2.1\cdot10^{-46}$. To check the initial signs, sum the 100
real residues and add an error interval of radius
\[
 0.017/14^{n+1}+C_{100}/(202.5)^{n+1}.
\]
Interval evaluation
gives $s_0>0.026$ and $s_n<-1.5\cdot10^{-5}$ for
$1\le n\le10$. Using the first 29 conjugate pairs to sharpen the
enclosure gives
\[
 0.02767331909339<s_0<0.02767331909340 .\qedhere
\]
\end{proof}

The fraction of the weighted total variation contributed by
positive slopes is strictly below one half in every degree
$n\ge1$. Its distance from one half satisfies
\[
\begin{gathered}
 d_n=\frac12-\frac{a_n(1)}{a_n(1)+a_n^-}
     =\frac{-s_n}{2(2a_n(1)-s_n)}>0,\\
 \lim_{n\to\infty}\frac{d_{n+1}}{d_n}
       =\frac27(3/2-\Theta).
\end{gathered}
\]
Indeed $s_{n+1}/s_n\to2/7$ and
$\lim a_n(1)^{1/n}\ge1$, so $s_n/a_n(1)\to0$.
The denominator therefore has successive ratio
$1/(3/2-\Theta)$ by Theorem~\ref{th:positivepowerratio}.
These ratios concern boost-time averages. Bounding the cumulative
loss up to each prime cutoff additionally requires control of the
positive part before averaging.

\head{An explicit series over prime powers}
The full second moment $B_n=\int_0^\infty w_{n,2}(u)g(u)^2\,\dd u$
satisfies $0\le a_n(2)\le B_n$. Let $d=(1,-3,2)$ and
$\lambda=(\pi/4,\pi,4\pi)$ as in Sec.~\ref{sec:losspowers}.
For $l=0,1,2$ put
\[
\begin{aligned}
 T_{n,l}(x)&=\frac1{n!}\int_0^1 t^{l+1/2}(-\log t)^n e^{-xt}\,\dd t\\
 &=\sum_{h=0}^\infty\frac{(-x)^h}{h!(h+l+3/2)^{n+1}}.
\end{aligned}
\]
The substitution $t=e^{-2u}$ gives
\[
\begin{gathered}
\begin{aligned}
 B_n&=\frac1{4\pi^2}\sum_{m,k\ge2}\Lambda(m)\Lambda(k)
 \sum_{j,l=1}^3d_jd_l\\
 &\quad\times\Big[4abT_{n,2}(a+b)-(a+b)T_{n,1}(a+b)\\
 &\qquad\quad+\tfrac14T_{n,0}(a+b)\Big],
\end{aligned}\\
 a=\lambda_jm^2,\quad b=\lambda_lk^2 .
\end{gathered}
\]
At fixed $m,k,j,l$ the displayed Taylor series for $T$ is
absolutely convergent. For the outer sum keep its integral form.
Chebyshev's bound and Gaussian decay give
$e^{-u/2}\sum_m\Lambda(m)(1+(m/e^u)^2)e^{-\pi(m/e^u)^2/4}
\ll e^{u/2}$. This proves absolute integrability of the double
sum for each $n$, including its separate polynomial terms.

\head{Cancellation between close and distant pairs at every degree}
For $y=e^u$ and $0\le\delta\le1$, let $C(u,\delta)$ denote
the portion of $g(u)^2$ with $|m-k|\le\delta y$ in its
ordered prime-power sum, including the diagonal. Write
$F(u,\delta)=g(u)^2-C(u,\delta)$. Uniformly in $\delta$,
for every fixed $A>0$,
\[
\begin{aligned}
 C(u,\delta)&=y I(\delta)+O_A\big(y(1+u)^{-A}\big),\\
 F(u,\delta)&=-y I(\delta)+O_A\big(y(1+u)^{-A}\big),
\end{aligned}
\]
where the entire function of $\delta$ extending the real integral is
\[
\begin{aligned}
 I(\delta)&=2\int_0^\delta\int_0^\infty K(v)K(v+h)\,\dd v\,\dd h\\
 &=2\sum_{j\ge0}\frac{\delta^{j+1}}{(j+1)!}
                  \int_0^\infty K(v)K^{(j)}(v)\,\dd v .
\end{aligned}
\]
To verify these estimates, let
$\mu_y=y^{-1}\sum_m\Lambda(m)\delta_{m/y}$.
The quantitative prime number theorem \cite{johnstonyang} gives
$|\mu_y([0,v])-v|\ll_A(1+v)(1+u)^{-A}$ for all $v\ge0$.
For $v<y^{-1/2}$ this follows from Chebyshev's bound, and above
that point $\log(yv)\ge u/2$. The bound also holds at left limits.
In the inner integral of $K(w)$ on
$[\max(0,v-\delta),v+\delta]$, Stieltjes integration by parts
therefore costs $O_A((1+u)^{-A})$ upon replacing $\mu_y$
by Lebesgue measure. Integrate this error against
$|K(v)|\,\dd\mu_y(v)$, of bounded total mass. For the outer
replacement the inner Lebesgue integral and its first derivative
are bounded uniformly in $\delta$, so the same argument applies.
The continuum contribution is $I(\delta)$. The full continuum
contribution is zero since $\int K=\tfrac12\int W=0$.
This proves both estimates. The same proof holds with primes alone
and weights $\log p$. The series for $I$ follows by normal
convergence of the shifted Gaussian integrals on compact complex
$h$-sets.

As $K(0)=K'(0)=0$, integration by parts gives
\[
 I(\delta)=2C_0\delta-\tfrac13 C_2\delta^3+O(\delta^5).
\]
Here $C_0=\int_0^\infty K^2=0.0248840476425738\ldots$ and
$C_2=\int_0^\infty (K')^2=0.7245427539208892\ldots$.
The diagonal $m=k$ alone contributes
\[
\begin{gathered}
 B_n^{\rm diag}=\tfrac12 C_0(n+1)+C_{\log}+O_A((n+1)^{-A}),\\
 C_{\log}=\int_0^\infty K(v)^2\log v\,\dd v
                    =-0.01069612564431707\ldots .
\end{gathered}
\]
Indeed partial summation over primes and the bound for the higher
powers give $\sum_{m\le t}\Lambda(m)^2=t\log t-t+
O_A(t/\log^A(2t))$ at large $t$. Integrating against $K(m/y)^2/y$
gives $C_0u+C_{\log}+O_A((1+u)^{-A})$. For the small range
$m/y<y^{-1/2}$ use $K(v)=O(v^2)$ at zero. The gamma mean of
$u$ is $(n+1)/2$, and its inverse moments give the stated error.
For any fixed $d>0$ take $\delta=(1+u)^{-d}$, and denote the
$w_{n,2}$-averages of $C,F$ by $C_n,F_n$. Then
\[
\begin{gathered}
 C_n\sim 2C_0\,2^{n+1}n^{-d},\\
 F_n=-C_n+O_A(2^n n^{-A}) \quad\text{for every fixed }A>0.
\end{gathered}
\]
Indeed multiplication of $w_{n,2}$ by $e^u$ gives $2^{n+1}$
times the gamma density of rate one. Under this density
$\mathbb E(1+u)^{-a}\sim n^{-a}$ for fixed $a>0$. This follows
from concentration of $u/n$ near one and uniform integrability
of its inverse powers, using
$\mathbb E u^{-k}=\Gamma(n+1-k)/\Gamma(n+1)$ for $n+1>k$.
The relative cancellation is uniform to every fixed inverse power
of the degree. The majorant $O_A(2^n n^{-A})$ has exponential
rate two. The criterion of Theorem~\ref{th:positivepowerratio}
requires rate one.
The quantitative prime number theorem used at the end of
Sec.~\ref{sec:primelossbudget} gives the more explicit estimate
$B_n\ll2^n\exp[-c n^{3/5}(\log n)^{-1/5}]$ at large $n$,
for some $c>0$. Insert the bound on $g$ there, and split the
resulting rate-one gamma average at $n/2$. Its lower tail is at
most $e^{n/2}2^{-n-1}$ by the exponential Markov inequality.

\head{Certified coefficients on a finite range of degrees}
Let $T_*=3\cdot10^{12}$. For each row, the following enclosures of
the one-sided coefficient $a_n(r)$
hold for every integer $250\le n\le n_{\max}=(3/4)rT_*$.
They use the verified zeros through $T_*$ \cite{pt}.
\begin{center}
\small\setlength{\tabcolsep}{3pt}%
\begin{tabular}{c r l}
\toprule
Order $r$&$n_{\max}$&Bounds on $a_n(r)$\\
\midrule
$1/2$&$1.125\cdot10^{12}$&$0.00465<a_n<0.00669$\\
$1$&$2.250\cdot10^{12}$&$6.929\cdot10^{-5}<a_n<7.135\cdot10^{-5}$\\
$3/2$&$3.375\cdot10^{12}$&$8.986\cdot10^{-7}<a_n<9.285\cdot10^{-7}$\\
$2$&$4.500\cdot10^{12}$&$1.1975\cdot10^{-8}<a_n<1.2430\cdot10^{-8}$\\
$5/2$&$5.625\cdot10^{12}$&$1.6254\cdot10^{-10}<a_n<1.6949\cdot10^{-10}$\\
$7/2$&$7.875\cdot10^{12}$&$3.1032\cdot10^{-14}<a_n<3.2674\cdot10^{-14}$\\
\bottomrule
\end{tabular}
\end{center}
The full second moment has the sharper enclosure
\[
\begin{gathered}
 |B_n-2.44028855722908368\cdot10^{-8}|<1.1\cdot10^{-23}\\
       (250\le n\le4.5\cdot10^{12}) .
\end{gathered}
\]
In particular $|B_{n+1}/B_n-1|<9.1\cdot10^{-16}$ whenever
both degrees are in this range. The root and ratio conclusions of
Theorem~\ref{th:positivepowerratio} apply to $B_n$ as well. Its
cumulative integral dominates $M_2(U)$, its integrand satisfies
the same exponential upper bound, and its defining measure is positive.
With zeros counted with multiplicity the spectral
part of $g$ is $\sum_\rho h_\rho e^{(\rho-1/2)u}$.
The remaining contour on $\Re s=-1$ is bounded by
$100e^{-3u/2}$. The line estimate of Theorem~\ref{th:oddglobal}
costs at most
$\frac{75}{2\pi}\int_0^\infty(3/2+t)e^{-\pi t/4}\,\dd t<100$
after differentiation. There are 29 positive ordinates below 100.
Write $P(u)$ for the trigonometric polynomial from these and
their conjugates. The verified zeros with larger ordinates through
$T_*$ have total absolute amplitude at most
$\tau=(6\sqrt2/\pi)0.047\,100^2e^{-100\pi/6}<2.313\cdot10^{-20}$.
The remaining terms are bounded by $H_* e^{u/2}$ with
$H_*=T_*^2e^{-\pi T_*/6}$, by the gamma and reciprocal-square
bounds in Theorem~\ref{th:oddglobal}.

The first sinusoid has amplitude $a=2|h_{\rho_1}|$.
For $r\ge1$ the $L^r(w_{n,r}\dd u)$-error after retaining
just this sinusoid is bounded by the sum of the other amplitudes
through 100, $\tau$, $100(2/5)^{(n+1)/r}$, and
$H_*2^{(n+1)/r}$. This is less than $1.025\cdot10^{-6}$ in
the table. Taking the positive part is 1-Lipschitz. For $0<r<1$
use $|x_+^r-y_+^r|\le|x-y|^r$ and subadditivity of the $r$th
power. The corresponding error in the moment itself is bounded
by $(1.025\cdot10^{-6})^r$ for the row $r=1/2$.
The mean of $(\cos\theta)_+^r$ is
$\Gamma((r+1)/2)/(2\sqrt\pi\,\Gamma((r+2)/2))$.
In its gamma average the error from this mean is at most
$2\sum_{k\ge1}(r/(k\gamma_1))^{n+1}<4(r/14)^{n+1}$.
Indeed the absolute Fourier coefficients of the bounded periodic
function are at most one, and the gamma multipliers are summable.
Together with the certified ordinates and amplitudes these bounds
give the table.

For the full square, $\mathbb E_{n,2}P^2$ has constant term
$E_{29}=2\sum_{0<\gamma<100}|h_\rho|^2$.
Every other frequency has modulus greater than $1.2192$, by
enclosure of the ordinates. If $S_{29}=2\sum_{0<\gamma<100}|h_\rho|$,
their aggregate contribution is at most
$S_{29}^2(2/\sqrt{4+1.2192^2})^{n+1}$.
In $L^2(w_{n,2}\dd u)$, $g-P$ has norm at most
$e_n=\tau+100(2/5)^{(n+1)/2}+H_*2^{(n+1)/2}$.
Thus $|B_n-\mathbb E_{n,2}P^2|\le e_n(2S_{29}+e_n)$.
The interval check evaluates both bounds uniformly over the
stated range. Under RH the tail beyond 100 has total absolute
amplitude at most $\tau$ at every time, and the upper degree
limit in each enclosure can be removed. An independent bound
$\liminf_{n\to\infty}a_{n+1}(r)/a_n(r)\le1$ for one fixed
$r>0$ remains open.

\head{The first coefficients from prime returns}
For comparison, direct prime-power quadrature on $0\le u\le6$,
completed by a bound on $u>6$, gives
\begin{center}
\begin{tabular}{c l l}
\toprule
$n$&One-sided $a_n(1)$&Full quadratic $B_n$\\
\midrule
0&$(0.038090,0.038190)$&$(0.00343528,0.00343529)$\\
1&$(0.016538,0.016672)$&$(0.00318211,0.00318214)$\\
2&$(0.004728,0.005038)$&$(0.00224224,0.00224234)$\\
3&$(0.001035,0.001667)$&$(0.00153395,0.00153435)$\\
\bottomrule
\end{tabular}
\end{center}
The ancillary file \texttt{anc/positive\_series\_check.py} verifies
these intervals. For $|\Im z|\le1/4$, each Gaussian factor in
$g(z)$ is bounded using
$(2t+1/2)e^{-t/2}\le4e^{-t/4}$ for $t\ge0$.
Partial summation with $\psi(x)\le3x$ gives
$|g(z)|\le24e^{\Re z/2}$. At a real argument $u$, truncating
the prime powers above $16e^{u_{\rm right}}$ on a panel with
right endpoint $u_{\rm right}$ costs at most
$24e^{u_{\rm right}/2-16\pi}$, by extracting half the Gaussian
exponent after bounding its linear factor. All retained prime
powers are below 6500. For the one-sided integral use a mesh of
size $1/2048$ and enclose $g$ between successive endpoint values
with added linear interpolation error $(1/2048)^2\sup|g''|/8$.
Cauchy's bound on a circle of radius $1/4$ supplies $\sup|g''|$.
Taking positive parts of the resulting intervals and integrating
the gamma weights exactly also covers cells containing a zero.

For the full square use 18-point Gauss--Legendre quadrature on
panels of length $1/8$. The integrand on each midpoint disk of
radius $1/4$ is bounded by $10^7$ for $n\le3$.
Taylor approximation through degree 35, exactness of the rule and
positivity of its weights bound the total quadrature error by
$12\cdot10^7\,4^{-36}/(1-1/4)<4\cdot10^{-14}$.
On $u>6$ integrate the bound
$|g(u)|\le0.000222+100e^{-3u/2}+H_*e^{u/2}$ against the relevant
gamma density. For the square use $(a+b+c)^2\le3(a^2+b^2+c^2)$.

\subsection{The sign of the correlation between distinct primes}
\label{sec:distinctprimeenergy}
Write $y=e^u$ and define the return slope from primes alone by
\[
\begin{aligned}
 P(u)&=y^{-1/2}\sum_p\log p\,K(p/y),\\
 D(u)&=y^{-1}\sum_p(\log p)^2K(p/y)^2,\\
 Q_n&=\int_0^\infty w_{n,2}(u)\,[P(u)^2-D(u)]\,\dd u .
\end{aligned}
\]
Thus $Q_n$ is the sum over ordered pairs $p\ne q$, each with
weight $y^{-1}\log p\log q\,K(p/y)K(q/y)$ inside the integral.
Absolute integrability at each degree follows from the Gaussian
majorant in Sec.~\ref{sec:positivepowers}. Denote the same averages
of $P^2,D$ by $B_n^{\mathbb P},D_n^{\mathbb P}$, respectively.

\begin{theorem}[A sign criterion involving only distinct primes]
\label{th:distinctprimeenergy}
RH is equivalent to $Q_n\le0$ for infinitely many integers $n$.
More precisely, $D_n^{\mathbb P}=\tfrac12 C_0(n+1)+C_{\log}+o(1)$.
On RH one has
\[
 Q_n=-\tfrac12 C_0(n+1)+\sum_{\rho}^{\rm distinct}|h_\rho|^2-C_{\log}+o(1),
\]
where $h_\rho=-m_\rho(\rho-1/2)A(\rho)$ includes the multiplicity.
If a nontrivial zero lies to the right of the critical line, then
$Q_n>0$ eventually and $Q_n^{1/n}\to1/(3/2-\Theta)>1$.
Unconditionally,
\[
 Q_n<-0.033\qquad(250\le n\le 6{,}798{,}000{,}000{,}000).
\]
\end{theorem}
\begin{proof}
Let $E=g-P$. For $u\ge0$, $|E(u)|<7$, and $E(u)\to0$, as
also follows from the square and higher-power estimates below.
Consequently $\|E\|_{L^2(w_{n,2})}=o(1)$ because the gamma
densities escape to infinity. The triangle inequality gives
$|\sqrt{B_n^{\mathbb P}}-\sqrt{B_n}|\le\|E\|_{L^2(w_{n,2})}$.
Partial summation over primes against $K(v)^2$, as in
Sec.~\ref{sec:positivepowers}, gives the stated diagonal.
Under RH the gamma multiplier kills every nonconstant frequency in
the uniformly absolutely convergent series for $g^2$, up to a
remainder tending to zero. The constant term is the displayed sum
over distinct zeros. If RH fails, Theorem~\ref{th:positivepowerratio}
and its full-square consequence give $B_n^{1/n}\to(3/2-\Theta)^{-1}$.
The bounded change of square root and the subtraction of the linear
diagonal preserve this rate. This proves both implications.

For explicit bounds, put $F(x)=\psi(x)-\vartheta(x)$. For $x\ge1$,
$F(x)<3\sqrt x$ and $F(x)-\vartheta(\sqrt x)<1.78x^{1/3}$
\cite{dusarttheta}. Both counts vanish for $0<x<1$.
For a Gaussian term of $K$ use
$(2t+1/2)e^{-t/2}\le2$ to obtain
$|K(v)|\le(6/\pi)e^{-\pi v^2/8}$.
Partial summation against this decreasing majorant bounds $|E|$
by $(18/\pi)(8/\pi)^{1/4}\Gamma(5/4)<7$, and the part of
exponent at least three by
$(6/\pi)1.78(8/\pi)^{1/6}\Gamma(7/6)e^{-u/6}<4e^{-u/6}$.

For the squares put $Y=e^{u/2}$, $f(v)=K(v^2)$.
Its integral is zero because the Mellin transform of $K$ is
$(s-1/2)A(s)$. Integration by parts gives the square contribution
$-\int_0^\infty [\vartheta(Yv)-Yv]f'(v)\,\dd v/Y$.
Differentiating the three Gaussians yields $|f'(v)|\le57v^3$ and
\[
 \int_0^\infty v|f'(v)|\,\dd v
 \le \frac{5\Gamma(5/4)}{2\pi}\sum_j|d_j|\lambda_j^{-1/4}<3.2 .
\]
For $u\ge64$, split at $v=e^{-u/4}$. On the lower interval use
$|\vartheta(x)-x|\le2x$, and on the upper use
$|\vartheta(x)-x|\le0.2x/\log^2x$ at $x\ge e^{16}$
\cite{dusarttheta}. It follows that
\[
 |E(u)|\le 23e^{-5u/4}+10.24/u^2+4e^{-u/6}\qquad(u\ge64).
\]
For $n\ge250$ the gamma mass below 64 is at most
$e^{64}(2/3)^{251}$ by the exponential Markov inequality.
These bounds give $\|E\|_{L^2(w_{n,2})}<0.0027$.

For $T=3\cdot10^{12}$, the zeros through $T$ are on the critical
line \cite{pt}. Their aggregate absolute amplitude is less than
$0.000222$ by Sec.~\ref{sec:positivepowers}. Above $T$, rotate
the Euler integral for $\Gamma(1+\rho/2)$ through the angle
${\rm sgn}(\gamma)(\pi/2-1/|\gamma|)$.
Using $\sin(1/t)\ge2/(\pi t)$ at $t=|\gamma|\ge14$ in the
argument of Theorem~\ref{th:oddglobal} gives
$|h_\rho|\le4|\gamma|^{3/2}e^{-\pi|\gamma|/4}$ per copy.
Indeed the constant before $t^{3/2}e^{-\pi t/4}$ is bounded by
$(3/\pi)\sqrt{1+1/(4\cdot14^2)}(\pi/2)^{3/2}\sqrt e<4$.
Since $\sum_\rho|\gamma|^{-2}<0.047$, the upper tail has
absolute amplitude at most $T^{7/2}e^{-\pi T/4}$. Its
$L^2(w_{n,2})$ norm is at most this quantity times $2^{(n+1)/2}$.
At the upper degree in the theorem the logarithm of this last bound
is below $-1.8\cdot10^8$. The contour contributes at most
$100(2/5)^{(n+1)/2}$. Hence $\sqrt{B_n}<0.000223$ and
$B_n^{\mathbb P}<(0.0027+0.000223)^2<9\cdot10^{-6}$ throughout.

Finally $K(v)>1/32$ on $3/2\le v\le2$. A lower bound there is
$\{(2\pi-1/2)e^{-\pi}-3(9\pi/2-1/2)e^{-9\pi/4}\}/(2\pi)>1/32$
because the third Gaussian contributes positively. For $u\ge40$,
the bound on $\vartheta$ gives
$\vartheta(2y)-\vartheta(3y/2)>(2/5)y$. Keeping just these
primes shows $D(u)\ge (2/5)u/32^2$. Thus
\[
 D_n^{\mathbb P}\ge\frac{2}{5\cdot32^2}
                  \left(\frac{n+1}{2}-40\right)\ge0.0333984375 .
\]
Subtracting proves the finite assertion.
\end{proof}

\begin{corollary}[Exact local cancellation at arbitrarily large scales]
\label{co:localdistinctpairs}
There are constants $L,c_1,c_2>0$ such that every sufficiently
late interval $[U,U+L]$ contains a point $v$ with
\[
 P(v)=0,\qquad P(v)^2-D(v)=-D(v).
\]
It also contains an interval of length at least
$c_1\sqrt U\,e^{-U/2}$ on which $P^2-D\le-c_2 U$.
On RH the inequality $P(u)^2-D(u)<0$ holds at every sufficiently
large $u$. If RH fails, every sufficiently late interval of some
fixed length also contains a positive value of $P^2-D$.
\end{corollary}
\begin{proof}
Under RH, $P=g+o(1)$ is bounded. The first sinusoid dominates the
sum of the other amplitudes in the explicit formula for $g$, so
both signs occur within each sufficiently late full period. If
$\Theta>1/2$, choose $\sigma$ and the averaging measure from
Theorem~\ref{th:primelossgrowth} as in the proof of
Theorem~\ref{th:boostheatmoments} below, with $\tau=0$. The
coefficient $-m_\rho(\rho-1/2)A(\rho)$ of the surviving mode is
nonzero. In each period the averaged normalized slope takes values
of both signs bounded away from zero. Positivity of the averaging
measure gives values of $g$ of opposite signs and
magnitude $\gg e^{(\beta_*-1/2)U}$ in a fixed enlarged interval.
The same is true of $P$ since $E(u)\to0$.
Continuity of the prime sum yields the zero in either case.

The pointwise diagonal is $C_0u+C_{\log}+o(1)$ by the same
partial summation used for $D_n^{\mathbb P}$. Moreover
$|P'(u)|\ll e^{u/2}$ by differentiated Gaussian summation and
Chebyshev's bound. On both sides of the zero, within distance
$c\sqrt U e^{-U/2}$, the mean value theorem gives
$P^2<D/2$ for a sufficiently small fixed $c>0$. Enlarge the
fixed window at either end by one to contain this subinterval.
The bounded values on RH and, in the other case, the exponentially
large values above give the remaining assertions.
\end{proof}
In particular the exact cancellations at these unbounded time scales
can coexist with eventually positive $Q_n$ if a right-hand zero
exists. The averaging that defines $Q_n$ tests an entire gamma
window in $u$. For a simple analytic example, the function
$e^{u/4}\sin(3u)$ has a zero in every closed interval of length
$\pi/3$, yet its squared gamma mean at rate two is
\[
 \tfrac12\big[(4/3)^{n+1}-\Re(3/4-3i)^{-n-1}\big]
 \sim\tfrac12(4/3)^{n+1}.
\]

\subsection{Gaussian resolution in boost time and the induction in degree}
\label{sec:boostheatmoments}
A detector with Gaussian uncertainty of variance $2\tau$ in the
dimensionless Rindler boost time records the convolution
$g_\tau(u)=\mathbb E\,g(u+Z)$ for $Z\sim N(0,2\tau)$.
Write $g_0=g$. At every $\tau$ the prime powers give
\[
\begin{aligned}
 g_\tau(u)&=e^{-u/2}\sum_m\Lambda(m)K_\tau(me^{-u}),\\
 K_\tau(v)&=\mathbb E[e^{-Z/2}K(ve^{-Z})].
\end{aligned}
\]
The sum with absolute values before averaging is $O(e^{(u+Z)/2})$
including the normalization, by Chebyshev's bound. This justifies
the exchange. A real frequency $E$ has response $e^{-\tau E^2}$.
Define
\[
\begin{aligned}
 M_n(r,\tau)&=\int_0^\infty w_{n,r}(u)(g_\tau(u)_+)^r\,\dd u,\\
 B_n(\tau)&=\int_0^\infty w_{n,2}(u)g_\tau(u)^2\,\dd u .
\end{aligned}
\]
In particular $M_n(r,0)=a_n(r)$ and $B_n(0)=B_n$.

Here the square at $\tau=\kappa n$ admits an exact iteration.
Take $Y$ exponential of rate two and $Z,Z'$ independent
centered Gaussians of variance $2\kappa$, independent of $Y$.
For functions of two variables set
\[
\begin{aligned}
 (\mathcal T_\kappa F)(x,y)&=\mathbb E\,F(x+Y+Z,y+Y+Z'),\\
 F_0(x,y)&=\mathbb E[g(x+Y)g(y+Y)] .
\end{aligned}
\]
Then $F_{n+1}=\mathcal T_\kappa F_n$ satisfies
$F_n(0,0)=B_n(\kappa n)$. Indeed each application adds one
independent exponential waiting time shared by both arguments and
one independent Gaussian increment to each. After $n$ steps the
shared sum has density $w_{n,2}$, and conditioning on it gives
the square of $g_{\kappa n}$. The exponential growth estimate
on $g$ gives absolute integrability at every step.

\begin{theorem}[Growth under linear and sublinear Gaussian averaging]
\label{th:boostheatmoments}
Write $\lambda_\rho=\rho-1/2=b_\rho+i\gamma_\rho$,
$d_\rho=\gamma_\rho^2-b_\rho^2$, $b=\Theta-1/2$, and combine
multiplicities in $h_\rho=-m_\rho\lambda_\rho A(\rho)$.
For $0<\kappa\le10^{-3}$, put
\[
 L(\kappa)=\max_\rho\frac{e^{-2\kappa d_\rho}}{1-b_\rho}.
\]
The limit of $B_{n+1}(\kappa(n+1))/B_n(\kappa n)$ equals
$L(\kappa)$. More precisely, for some $0<q_\kappa<1$,
\begin{multline*}
 B_n(\kappa n)=L(\kappa)^n
 \Bigg(\sum_{\rho:\,e^{-2\kappa d_\rho}/(1-b_\rho)=L(\kappa)}
                \frac{|h_\rho|^2}{1-b_\rho}\\
 +O_\kappa(q_\kappa^n)\Bigg).
\end{multline*}
For every fixed $r>0$ and every sequence $\tau_n\ge0$ with
$\tau_n/n\to0$, the one-sided moments instead satisfy
\[
 \lim_{n\to\infty} M_n(r,\tau_n)^{1/n}=\frac1{1-b}.
\]
The latter assertion also holds for the full square with $r=2$.
\end{theorem}
\begin{proof}
Convolution of the explicit formula gives, for all real $u$,
\[
\begin{gathered}
 g_\tau(u)=\sum_\rho h_\rho e^{\tau\lambda_\rho^2}e^{\lambda_\rho u}
              +R_\tau(u),\\
 |R_\tau(u)|\le100e^{9\tau/4-3u/2}.
\end{gathered}
\]
The inverse Mellin line $\Re s=-1$ used in
Sec.~\ref{sec:positivepowers} bounds $R_0$ also at negative
$u$. The bound $\sum_\rho|h_\rho|<\infty$ and Gaussian
exponential moments justify termwise convolution on the real line.
The $L^2(w_{n,2})$ norm of $R_\tau$ is at most
$100e^{9\tau/4}(2/5)^{(n+1)/2}$.

On $e^{\lambda x+\bar\mu y}$ the induction operator acts by
multiplication by $2e^{\kappa(\lambda^2+\bar\mu^2)}/
(2-\lambda-\bar\mu)$. For a single zero, its squared norm at
degree $n$ is $|h_\rho|^2 L_\rho^n/(1-b_\rho)$, where
$L_\rho=e^{-2\kappa d_\rho}/(1-b_\rho)$.
For two distinct exponents the cross term in the gamma integral
is at most the geometric mean of their squared norms times
\[
 \left(\frac{2\sqrt{(1-b_\rho)(1-b_{\rho'})}}
 {|2-\lambda_\rho-\overline{\lambda_{\rho'}}|}\right)^{n+1}.
\]
The fraction is strictly less than one by the arithmetic-geometric
mean inequality unless both real parts and both frequencies agree.
For fixed $\kappa>0$, $L_\rho\to0$ as $|\gamma_\rho|\to\infty$.
There are finitely many maximizers, and the supremum over all other
distinct zeros is strictly smaller. The triangle inequality and
$\sum|h_\rho|<\infty$ bound the total norm of the remaining modes
by a constant times this smaller base to the power $n/2$.
The squared contour bound has base $(2/5)e^{9\kappa/2}$, strictly
less than $e^{-2\kappa\gamma_1^2}\le L(\kappa)$ because
$\gamma_1<15$ and $\kappa\le10^{-3}$. The finite cross terms
and these two tails give the asserted asymptotic and ratio.

For $\tau_n=o(n)$ the spectral bound is
$|g_{\tau_n}(u)|\le C e^{bu}+100e^{9\tau_n/4-3u/2}$ at $u\ge0$,
since all $d_\rho>0$. Integrating its $r$th power gives the
claimed upper root bound. For a lower bound when $b>0$, choose
$1/2<\sigma<\Theta$ different from every zero real part.
Use the compactly supported positive probability $\nu$ constructed
in the proof of Theorem~\ref{th:primelossgrowth}. Among the zeros
above $\sigma$, it isolates the pair of lowest positive ordinate
$\gamma_*$ with maximal real part $\beta_*$ at that ordinate.
Put $a=\beta_*-1/2$, $\ell=2\pi/\gamma_*$ and
$\operatorname{supp}\nu\subset[0,D]$. Uniformly for
$v\in[U,U+\ell]$ with $U\asymp n$, averaging translates of
$e^{-av}g_{\tau_n}(v)$ gives
\[
 2\Re\!\left(H_*e^{\tau_n\lambda_*^2}e^{i\gamma_* v}\right)
      +O(e^{-\delta U+o(n)}),\qquad H_*\ne0,\quad\delta>0 .
\]
Indeed the same translation factors annihilate the selected
unwanted frequencies for every $\tau_n$. The terms below
$\sigma$ decay with gap $\beta_*-\sigma>0$. The other terms
at $\gamma_*$ have a positive real-part gap since there are
only finitely many there. The surviving amplitude has size
$2|H_*|e^{-d_*\tau_n}=e^{-o(n)}$.

On a fixed fraction of this period the average exceeds a fixed
multiple of that amplitude. For $r\ge1$, apply Jensen's inequality
to $x\mapsto(x_+)^r$ and integrate over the period. The
unweighted local integral of $(g_{\tau_n})_+^r$ on
$[U,U+\ell+D]$ is at least $\exp(raU-o(n))$. For $0<r<1$,
use this statement at $r=1$ and the pointwise bound on the same
interval to get $\exp([a-(1-r)b]U-o(n))$.
With $U=n/[r(1-b)]$, the gamma density throughout the interval
is within fixed multiplicative factors of $w_{n,r}(U)$.
Its logarithm divided by $n$ tends to
$-\log(1-b)-b/(1-b)$ by Stirling's formula. Now let
$\sigma\uparrow\Theta$, so $a\to b$. This gives the lower bound.

If $b=0$, the amplitude of the first conjugate pair dominates
the sum of the remaining amplitudes as in
Sec.~\ref{sec:positivepowers}. After Gaussian averaging the first
amplitude is multiplied by $e^{-\gamma_1^2\tau_n}$, and the
others are attenuated at least as much. Near $U=n/r$ the
contour error is exponentially small in $n$. On a fixed part
of a period the positive slope is therefore at least
a constant times $e^{-\gamma_1^2\tau_n}=e^{-o(n)}$.
The gamma density there has root limit one. The full square
dominates the corresponding one-sided integral and obeys the
same upper bound.
\end{proof}

\begin{corollary}[Critical averaging rate]
\label{co:boostheatthreshold}
Let $\kappa_*=\inf\{\kappa\in(0,10^{-3}]:L(\kappa)\le1\}$.
With the supremum of the empty set taken to be zero here,
\[
 \kappa_*=\sup_{\rho:\,b_\rho>0}
       \frac{-\log(1-b_\rho)}{2(\gamma_\rho^2-b_\rho^2)},\qquad
 0\le\kappa_*<3.782\cdot10^{-26}.
\]
If positive, the supremum is attained. Its vanishing is equivalent
to RH. At the fixed rate $\kappa_0=3.8\cdot10^{-26}$, for every
fixed $r>0$ one has unconditionally
\[
 M_n(r,2\kappa_0n/r)\sim
 (2|h_1|)^r
 \frac{\Gamma((r+1)/2)}{2\sqrt\pi\,\Gamma((r+2)/2)}
 e^{-2\kappa_0\gamma_1^2 n}.
\]
In particular the ratio of consecutive terms in this last sequence
tends to $e^{-2\kappa_0\gamma_1^2}<1$, including at
$r=1/2,3/2,5/2,7/2$.
At quadratic order every degree $n\ge250$ satisfies
\[
 M_n(2,\kappa_0 n)\le B_n(\kappa_0 n)
      <10^{-7}\exp(-2\kappa_0\gamma_1^2 n).
\]
\end{corollary}
\begin{proof}
The condition $L(\kappa)\le1$ bounds each mode separately.
The ratios defining $\kappa_*$ tend to zero with $|\gamma_\rho|$, and
finiteness of zeros in compact sets gives attainment if the supremum
is positive. Functional-equation symmetry supplies the equivalence.
For $|\gamma|\le T=3\cdot10^{12}$ the ordinates are verified
\cite{pt}. Above $T$ use $b_\rho\le1/2-1/(R\log|\gamma|)$,
$R=5.558691$ \cite{mtyregion}. With
$f(x)=-\log(1/2+1/(Rx))$, $x\ge\log T>28$, one has
$f(x)>0.67$ and $0<f'(x)\le2/(Rx^2)<0.001$.
Thus $f(\log t)/(t^2-d)$ decreases for $t\ge T$ for both
$d=1/4$ and $d=225+1/4$, since the logarithmic derivative of the
denominator with respect to $\log t$ exceeds two.
In particular
\[
 \kappa_*\le\frac{f(\log T)}{2(T^2-1/4)}<3.782\cdot10^{-26}.
\]
Even after replacing $1/4$ by $225+1/4$ the upper bound is
less than $\kappa_0$. Hence $e^{-2\kappa_0d_\rho}/(1-b_\rho)
<e^{-2\kappa_0\gamma_1^2}$ for every right-hand zero above
$T$. This also holds for left-hand zeros there, since
$1-b_\rho\ge1$ and $d_\rho>\gamma_1^2$. The pair at
$\pm\gamma_1$ is the unique maximizing pair.

For $r\ge1$, compare $g_{2\kappa_0n/r}$ with that pair in
$L^r(w_{n,r})$. The $r$th power of each individual mode norm
has base $e^{-2\kappa_0d_\rho}/(1-b_\rho)$, and the contour
has base $(2/5)e^{9\kappa_0/2}$. Thus their total error is
exponentially smaller in norm by the same gap argument. For
$0<r<1$ use $|x_+^r-y_+^r|\le|x-y|^r$ and
$\sum_\rho|h_\rho|^r<\infty$, which follows from the exponential
gamma bound and the reciprocal-square estimate. For the first
sinusoid alone the gamma average of its positive part to power
$r$, after factoring out its amplitude to power $r$, tends
to the displayed circular mean. This follows by trigonometric
approximation and the gamma multiplier at each frequency.
For the uniform quadratic bound put $m=e^{-2\kappa_0\gamma_1^2}$.
Every spectral mode has $L^2(w_{n,2})$ norm bounded by
$\sqrt2|h_\rho|m^{n/2}$. Adding the contour and dividing by
$m^{n/2}$ bounds the norm by
$\sqrt2(0.000223)+100\sqrt{2/5}\,
[(2/5)e^{9\kappa_0/2}/m]^{n/2}$. Its square is less than
$10^{-7}$ for $n\ge250$.
\end{proof}

\head{Slopes of the Gaussian growth exponent}
Set $\Phi(\kappa)=\log L(\kappa)$. On compact subintervals of
$(0,10^{-3})$ it is the maximum of finitely many affine functions
$-\log(1-b_\rho)-2d_\rho\kappa$, since the other modes tend to
zero uniformly after exponentiation. Hence it is convex and piecewise
affine. If a right-hand zero exists, the largest parameter at which
such a mode ties the first pair is
\[
 \kappa_\times=\max_{\rho:\,b_\rho>0}
      \frac{-\log(1-b_\rho)}{2(d_\rho-\gamma_1^2)}
      \in(0,3.782\cdot10^{-26}).
\]
The denominator is positive by the verified height, and the bound
and attainment follow as in Corollary~\ref{co:boostheatthreshold}.
All other critical-line modes and all left-hand modes lie below
the first pair at positive $\kappa$. For a left-hand zero off
the line this follows from $|\gamma|>T$ and $1-b_\rho>1$.
The right derivative at $\kappa_\times$ is $-2\gamma_1^2$.
The left derivative is the smallest slope among the finitely many
maximizers there. Their jump obeys
\[
 \Phi'_+(\kappa_\times)-\Phi'_-(\kappa_\times)
 \ \ge\ 2(T^2-1/4-\gamma_1^2)>1.79\cdot10^{25}.
\]
Under RH the function is exactly $-2\gamma_1^2\kappa$.
Thus differentiability throughout $(0,3.8\cdot10^{-26})$ is
another necessary and sufficient condition. The derivative before
passage to the limit is an arithmetic weighted integral. If
$\Phi_n(\kappa)=n^{-1}\log B_n(\kappa n)$, $f=g_{\kappa n}$,
and $w=w_{n,2}$, then for $n\ge2$
\[
\begin{aligned}
 \Phi_n'(\kappa)&=\frac{2\int_0^\infty w f f''\,\dd u}{\int_0^\infty w f^2\,\dd u}\\
 &=\frac{\int_0^\infty w'' f^2\,\dd u-2\int_0^\infty w(f')^2\,\dd u}
       {\int_0^\infty w f^2\,\dd u}.
\end{aligned}
\]
Here $\partial_\tau g_\tau=\partial_u^2g_\tau$.
The nonzero residue of the first conjugate pair in the Laplace
transform of $g_\tau$ ensures that the denominator is positive.
The identity follows by two integrations by parts, with the zero
of $w,w'$ at the origin and exponential decay of the weighted
products at infinity. Chebyshev's bound against derivatives of
the three Gaussians gives these bounds also after convolution.
Finite-degree differentiability alone gives no uniform control of
the limit near $\kappa_\times$.

\head{A fixed time resolution}
For $\tau=1$ let $E_1=2|h_1|^2 e^{-2\gamma_1^2}
=7.1063781092383\ldots\cdot10^{-182}$. Without assuming RH,
\[
 \left|\frac{B_n(1)}{E_1}-1\right|<6\cdot10^{-105}
 \qquad(1000\le n\le 25\cdot10^{24}).
\]
In particular adjacent ratios within this range differ from one
by less than $1.3\cdot10^{-104}$. To verify this, the first
pair contributes $E_1(1+\epsilon_n)$ with
$|\epsilon_n|\le v_n=(1+\gamma_1^2)^{-(n+1)/2}$.
Using $\gamma_2>21$ and putting $N=25\cdot10^{24}$, the
total error in $L^2(w_{n,2})$ after retaining the first pair is
at most
\begin{multline*}
 e=0.000223e^{-21^2}+100e^{9/4}(2/5)^{1001/2}\\
       +0.000223\exp\!\left[-(T^2-1/4)+\tfrac{N+1}{2}\log2\right].
\end{multline*}
The three terms bound the other verified zeros, the contour, and
the unknown upper tail, respectively. With $x=e/\sqrt{E_1}$
the relative error is at most
$v_{1000}+2x\sqrt{1+v_{1000}}+x^2<6\cdot10^{-105}$.
The degree bound refers to the Gaussian-averaged observable.

\head{A direct estimate from the prime number theorem}
Take $\tau_n=\sqrt n$, which lies in the sublinear regime.
Section~\ref{sec:primelossbudget}, by partial summation of
$\psi(x)-x$ and the prime number theorem \cite{johnstonyang},
gives $|g(v)|\ll e^{v/2-c\omega(v)}$ for large $v$, where
$\omega(v)=v^{3/5}(\log v)^{-1/5}$, and $|g(v)|\ll e^{v/2}$
for every real $v$. Jensen's inequality gives
$g_\tau(u)^2\le\mathbb E[g(u+Z)^2]$.
The factors $e^u$ and $e^Z$ change the gamma density to rate
one and the Gaussian to mean $2\tau$, with prefactor
$2^{n+1}e^\tau$. If $U$ and $Z'$ have these tilted laws,
then $U+Z'\ge n/4$ whenever $U\ge n/2$ and $Z'\ge-n/4$.
The complementary probabilities are bounded by
$e^{n/2}2^{-n-1}$ and $e^{-n^2/(64\tau)}$. At
$\tau=\sqrt n=o(\omega(n))$ it follows, with some $c'>0$, that
\[
 M_n(2,\sqrt n)\le B_n(\sqrt n)
       \ll 2^n\exp[-c'\,n^{3/5}(\log n)^{-1/5}].
\]
The upper majorant has $n$th-root limit two. To conclude RH
it would suffice to bound $M_{n_j}(2,\sqrt{n_j})$ by
$\exp(o(n_j))$ for a single sequence $n_j\to\infty$, by
Theorem~\ref{th:boostheatmoments}. The linear-rate estimate in
Corollary~\ref{co:boostheatthreshold} controls a different family.
The hypotheses for the sublinear bound or for the infinite sign
condition in Theorem~\ref{th:distinctprimeenergy} remain open.
The ancillary file \texttt{anc/heat\_pair\_check.py} checks the
outward-rounded constants and tests the transfer formula against
independent quadrature for a single conjugate pair of modes.

\subsection{Density of cancellation at prime cutoffs}
\label{sec:primewindowdensity}
Let $Q(u)=P(u)^2-D(u)$, using the prime slope and its diagonal
from Sec.~\ref{sec:distinctprimeenergy}. At the clock
$u_p=(6/5)\log p$, set $v_{p,q}=(\log q)K(q/p^{6/5})$ for prime
$q$. Then $Q(u_p)\le0$ is precisely the arithmetic inequality
\[
              \left|\sum_q v_{p,q}\right|^2\le\sum_q v_{p,q}^2 .
\]
Both series converge absolutely. Define the proportion
\[
 F_{\mathbb P}(X)=\frac{\#\{X<p\le2X:Q(u_p)\le0\}}{\pi(2X)-\pi(X)}
 \qquad(X\ge2).
\]

\begin{theorem}[Density of cancellation at prime cutoffs]
\label{th:primewindowdensity}
Under RH, $F_{\mathbb P}(X)=1$ for every sufficiently large $X$. If a
right-hand zero exists, $F_{\mathbb P}(X)\to0$. In particular a positive
lower bound for $F_{\mathbb P}(X_j)$ along one sequence $X_j\to\infty$
would imply RH. More precisely, put $b=\Theta-1/2$. For each
$\varepsilon>0$ the proportion of primes $X<p\le2X$ satisfying
\[
 \left|\frac{\log(1+|P(u_p)|)}{\log X}-\frac65 b\right|>\varepsilon
\]
tends to zero. These statements require no attainment of the
supremum $\Theta$.
\end{theorem}
\begin{proof}
The bound $P=g+O(1)$ and $D(u)\sim C_0u$, with $C_0>0$,
give the assertions on RH. Suppose $b>0$. Fix $0<a_0<\pi/4$.
The explicit formula for $g$ extends to
$|\Im z|\le a_0$ and gives $|g(z)|\ll_{a_0}e^{b\Re z}$
for $\Re z\ge0$. Indeed the spectral coefficients decay like a
polynomial times $e^{-\pi|\gamma|/4}$, and the integral on
$\Re s=-1$ has the same exponential decay in its ordinate.
The proper powers $E=g-P$ are uniformly bounded in this strip
because $\psi(t)-\vartheta(t)\ll\sqrt t$ and each Gaussian is
bounded by a polynomial times
$\exp[-(\pi/4)\cos(2a_0)(m/e^{\Re z})^2]$.
Normal convergence on compact subsets gives holomorphy. Hence
\[
              |P(z)|\ll_{a_0}e^{b\Re z}\quad
              (\Re z\ge0,\ |\Im z|\le a_0).
\]

There is a single length $L$ such that, for each
$0<a<b$, $\sup_{[U,U+L]}|P|\gg_a e^{aU}$ for all
sufficiently large $U$. For a uniform choice of $L$, fix
$1/2<\sigma_0<\Theta$. Ingham's bound \cite{ingham40} gives
$\sum_{\beta\ge\sigma_0,\,\gamma>0}\gamma^{-1}<\infty$.
For any $a<b$ choose
$\max(\sigma_0,1/2+a)<\sigma<\Theta$ avoiding zero real parts.
The averaging measure in the proof of
Theorem~\ref{th:primelossgrowth}, applied to $g$ as in
Corollary~\ref{co:localdistinctpairs}, isolates a nonzero pair
with $\beta_*\ge\sigma$. Its support length plus the period of
that pair is at most
\[
 L=1+\frac{2\pi}{14}
        +\pi\sum_{\beta\ge\sigma_0,\,\gamma>0}\gamma^{-1}.
\]
The lower amplitude and the onset may depend on $a$. Subtracting
bounded $E$ gives the claim for $P$.

An interpolation bound converts a large value on an interval into
a lower bound on the measure of its large-value set. Let $I$ be
a real interval of length $\ell$. Suppose that $f$ is holomorphic
on a neighborhood of the disk of radius $4\ell$ centered at the
midpoint of $I$, with $|f|\le M$ there. Assume $|f|\le m$
on a measurable subset of $I$ of length at least $\alpha\ell$.
For $0<m\le M$,
\[
 \sup_I |f|\le C_\alpha M^{1-\theta_\alpha}m^{\theta_\alpha},
 \qquad \theta_\alpha=\frac{\log2}{\log(8e/\alpha)}>0 .
\]
Indeed take $k+1$ quantiles of the subset's length measure,
at fractions $(j+1/2)/(k+1)$, $0\le j\le k$. They lie in
its closure and their pairwise distances are at least
$\alpha\ell|i-j|/(k+1)$. On $I$ the Lagrange polynomial has
modulus at most $m[2(k+1)/\alpha]^k/k!\le m(4e/\alpha)^k$
for $k\ge1$. The contour formula for the interpolation remainder
on the stated circle bounds it by $2M(2/7)^{k+1}$.
Choose $k=\lfloor\log(M/m)/\log(8e/\alpha)\rfloor$.
When $k<1$ the elementary bound $M$ suffices.

For every fixed $h>0$ and $0<a<b$ it follows that
\[
 \left|\{v\in[U,U+h]: |P(v)|\le e^{aU}\}\right|\longrightarrow0 ,
\]
where vertical bars denote length. Otherwise there are $U_j\to\infty$
with a fixed positive fraction of such points. Divide the interval
into a fixed number of equal cells of length $\ell<a_0/8$. One
has at least that same fraction. In its disk the strip estimate
gives $M\ll e^{bU_j}$, so interpolation bounds its supremum by
$O(e^{(b-\delta_1)U_j})$ for some $\delta_1>0$. Move a cell of
the same length by $\ell/2$ and use the estimate on the overlap
with $\alpha=1/2$. Repeating to left and right covers
$[U_j,U_j+L]$ in a bounded number of steps depending only on
$h,L,\ell$. The final exponent is $b-\delta$ for a fixed
$\delta>0$. This contradicts the preceding lower bound with
$\max(0,b-\delta)<a'<b$.

Finally set $U=(6/5)\log X$, $h=(6/5)\log2$. Divide
$[U,U+h]$ into $\lceil h X^{4/5}\rceil$ equal cells. Each
corresponds on the prime axis to an interval of length
$\asymp X^{1/5}$, by the mean value theorem. The Brun--Titchmarsh
inequality with modulus one \cite{mvbruntitch} bounds its prime
count by $O(X^{1/5}/\log X)$ uniformly. On a cell the variation
of $P$ is $O(e^{U/2}X^{-4/5})=O(X^{-1/5})$ by the
unconditional derivative bound. If $0<a<a'<b$, a cell
containing a prime with $|P(u_p)|\le e^{aU}$ lies entirely
in the sublevel set of height $e^{a'U}$ for large $X$.
The length estimate therefore allows only $o(X^{4/5})$ such
cells, containing $o(X/\log X)$ primes in all.
The prime number theorem gives $\pi(2X)-\pi(X)\sim X/\log X$.
Together with the upper bound for $P$, this proves convergence
in proportion. On the same interval $D(v)=O(U)$, so
$Q(u_p)\le0$ implies $|P(u_p)|\ll\sqrt U$ and $F_{\mathbb P}(X)\to0$.
\end{proof}

\head{Contribution of cancellation intervals to the prime moment}
Let $N_n=\int_0^\infty w_{n,2}(D-P^2)_+\,\dd u$. On RH,
$N_n/D_n^{\mathbb P}\to1$, since $\int w_{n,2}P^2=O(1)$
and $D_n^{\mathbb P}\asymp n$. If RH fails, $N_n/D_n^{\mathbb P}\to0$.
Indeed the proportion of $Q\le0$ on every late unit interval
tends to zero by the length estimate just proved. Use
$D(u)\ll 1+u$ and $u w_{n,2}=(n+1)w_{n+1,2}/2$.
The gamma laws concentrate on $[n/4,n]$, where their densities
vary by at most a fixed factor on a unit interval, since their
logarithmic derivatives there are bounded. Summing over those
intervals and then bounding the two vanishing gamma tails proves
the assertion.

An unconditional lower bound for $N_n$ follows from a derivative
estimate. Write $\omega(U)=U^{3/5}(\log U)^{-1/5}$. Partial
summation against $\vartheta(x)-x$ \cite{johnstonyang} gives
$|P'(u)|\ll e^{u/2-c\omega(u)}$ for large $u$ and some $c>0$.
This uses the derivative kernel $-K/2-vK'$, which has zero
integral and derivative $O(v)$ at the origin, exactly as in
Sec.~\ref{sec:primelossbudget}. In each sufficiently late fixed
window, the zero given by Corollary~\ref{co:localdistinctpairs}
therefore provides an interval of length
$\gg\sqrt U\exp[-U/2+c_1\omega(U)]$ with $D-P^2\gg U$,
after a fixed enlargement of the window. Sum these contributions
over disjoint windows lying in $[0.35n,0.45n]$.
All factors $w_{n,2}(u)u^{3/2}e^{-u/2}$ are comparable on a
single such window, uniformly in $n$. The integral of this
weight is $2^{n+1}\Gamma(n+5/2)/(n!\,(5/2)^{n+5/2})$.
Its normalized law has mean $(n+5/2)/(5/2)$ and variance
$(n+5/2)/(5/2)^2$. The retained range, with either end shortened
by a fixed window length, consequently carries $1-o(1)$ of
the weight. For some $c_2>0$ this gives
\[
 N_n\gg n^{3/2}(4/5)^n e^{c_2\omega(n)}.
\]
The lower bound divided by the diagonal tends to zero. It supplies
no positive lower limit for that fraction. By the dichotomy,
$\limsup_{n\to\infty}N_n/D_n^{\mathbb P}>0$ would imply RH.

\begin{proposition}[The first sign change and the complete finite degree range]
\label{pr:primewindowprefix}
Put $U_*=4.77\cdot10^{12}$. On $0\le u\le U_*$ the
distinct-prime signal $Q$ has exactly one zero $u_0$, with
\[
 0.64773613385252<u_0<0.64773613385254 .
\]
It is positive before $u_0$ and negative after it.
Every prime cutoff $u_p\le U_*$ thus satisfies $Q(u_p)<0$,
including $p=2$. The gamma moments of $Q$ satisfy
\[
 Q_n<0\qquad(0\le n\le n_*),\qquad n_*=6{,}919{,}000{,}000{,}000 .
\]
Under RH, $u_0$ is the unique zero of $Q$ on the entire positive
half-line. If RH fails, $Q$ changes sign at arbitrarily large times.
Consequently the total number of distinct zeros on $[0,\infty)$
is either one or infinite.
\end{proposition}
\begin{proof}
Here the only averaging in $Q_n$ is its defining gamma integral.
Let $T=3\cdot10^{12}$, $R=5.558691$,
$\delta=1/[R\log(2T)]$. The verification up to $T$ \cite{pt},
the classical zero-free region \cite{mtyregion}, and the gamma
estimate of Theorem~\ref{th:distinctprimeenergy} bound the
spectral part of $g$ at $0\le u\le U_*$ by
\begin{multline*}
 0.000222+
 T^{7/2}e^{-\pi T/4+(1/2-\delta)U_*}\\
 +(2T)^{7/2}e^{-\pi T/2+U_*/2}<0.000223 .
\end{multline*}
The split is at $2T$. The remaining contour is bounded by
$100e^{-3u/2}$. For $u\ge4$ the contour can be shifted to
$\Re s=-3$, where the trivial-zero term and the new integral in
Proposition~\ref{pr:signedmomentdegrees} give a combined bound
$(105u/16+5)e^{-5u/2}+300e^{-7u/2}$. In particular
$|g|<0.002$ on $[4,U_*]$.

For the finite-interval enclosures on $[0,4]$, subdivide the
interval into cells of length $1/1024$ and retain the primes up to
$\lceil16e^{u_{\rm right}}\rceil$ on each cell. Put
$a_q=e^{-u/2}\log q\,K(qe^{-u})$. Chebyshev's bound and the real
Gaussian envelope give
$\sum_{\rm omitted}|a_q|<(36/\pi)e^2e^{-16\pi}<10^{-18}$.
For $a_q'$ the polynomial multiplying each Gaussian in place
of $2t-1/2$ is $4t^2-6t+1/4$. Since
$(4t^2+6t+1/4)e^{-t/2}<16$ for $t\ge0$, its omitted
absolute mass is below $10^{-17}$. Products of these mass
bounds also enclose the omitted diagonal and its derivative.
Interval evaluation of $Q$ and
$Q'=2(PP'-\sum_q a_q a_q')$ gives $Q>0$ on $[0,5/8]$,
$Q'<-0.02$ on $[5/8,11/16]$ and $Q<0$ on $[11/16,4]$.
The signs at the two endpoints of the bracket in the statement
are opposite. Integrating constant cell bounds against the exact
rate-two exponential mass gives
\[
\begin{gathered}
 \int_0^{3/4}2e^{-2u}Q_+(u)\,\dd u<0.001041,\\
 \int_1^4 2e^{-2u}(-Q(u))\,\dd u>0.002559.
\end{gathered}
\]

For $4\le u\le18$, it suffices to evaluate the proper powers
$E=g-P$, which are sparser than the full prime sum. On
$|\Im z|\le1/4$, partial summation with
$\psi(x)-\vartheta(x)<3\sqrt x$ bounds $|E(z)|$ by
$(36/\pi)(16/\pi)^{1/4}\Gamma(5/4)<16$. This uses
$(2t+1/2)e^{-t/2}\le4e^{-t/4}$ for the complex Gaussian
envelope. At real arguments the tail above $16e^u$ is below
$16e^{-16\pi}<10^{-19}$. On the grid of step $1/64$ the
maximum certified modulus is less than $0.039$.
Cauchy's estimate gives $|E''|\le512$, so linear interpolation
gives $|E|<0.055$ throughout. The kernel satisfies $K(v)<-0.16$
on $0.55\le v\le0.8$, by interval evaluation on cells of
length $1/1000$. Keeping only primes in this band for the
whole of each $u$-cell on $[4,8]$ gives $D>0.015$.

For all $u\ge8$, the same band gives
\[
 D(u) \ge 0.00512(u-0.6).
\]
Indeed $\vartheta(0.8y)-\vartheta(0.55y)>0.2y$ at $y=e^u$.
Up to $y=10^{10}$ use B\"uthe's bounds on $\vartheta$
\cite{buthepsi}, giving $0.25y-1.95\sqrt{0.8y}>0.2y$.
Above $10^{10}$ the estimate $0.2x/\log^2x$ for
$|\vartheta(x)-x|$ \cite{dusarttheta} gives the assertion.
These bounds prove $P^2<D$ through $u=18$.

For $18\le u\le64$, set $Y=e^{u/2}$, $f(v)=K(v^2)$.
In the partial summation formula for the prime squares use
$|\vartheta(t)-t|\le2.4\sqrt t$ up to $t=10^{19}$.
Checking prime endpoints through 1423 and applying B\"uthe's
theorem above 1423 proves this inequality. For $v\le8$ we
have $Yv\le8e^{32}<10^{19}$. Direct Gaussian integration gives
$\int v^{1/2}|f'(v)|\,\dd v<3.63$. For $v>8$, the elementary
bound $2Yv$ on the counting error costs less than
$100e^{-512\pi}$, by extracting half the Gaussian exponent
from each term. Including the higher powers gives
\[
 |E(u)|\le 2.4(3.63)e^{-u/4}+100e^{-512\pi}+3.686e^{-u/6}.
\]
Adding $0.000223+100e^{-3u/2}$ gives an upper bound for
$|P(u)|$. Its square at 18 is less than $0.079$, below
$0.00512(18-0.6)$.
The modulus bounds decrease and the diagonal lower bound increases.
For $u\ge64$ the explicit $E$ bound in
Theorem~\ref{th:distinctprimeenergy} gives $P^2<D$ up to $U_*$.
Since $(6/5)\log2>3/4$, every prime cutoff in the stated range
lies above $u_0$. Under RH the spectral bound $0.000223$
holds at every $u\ge0$ by absolute summation on the critical
line. This extends the strict negative sign to infinity.
If RH fails, Corollary~\ref{co:localdistinctpairs} supplies
both signs in every sufficiently late window of a fixed length.

It remains to include the initial degrees in $Q_n$.
At $u>U_*$ use $Q(u)\le P(u)^2\le576e^u$.
For $0\le n\le249$, $u^n\le e^{u/2}$ already at all
$u\ge10^4$. The two certified cell integrals and $Q<0$
on $[3/4,U_*]$ therefore give
\begin{multline*}
 \frac{n!}{2^n}Q_n
 \le0.001041(3/4)^n-0.002559\\
 +2304e^{-U_*/2}<-0.0015 .
\end{multline*}
For $250\le n\le n_*$, the split of the spectral tail at $2T$
also sharpens the degree estimate in Theorem~\ref{th:distinctprimeenergy}.
Indeed the norm of $e^{(\beta-1/2)u}$ in $L^2(w_{n,2})$
equals $(3/2-\beta)^{-(n+1)/2}$. By the triangle inequality
the total norm of the two tails is bounded by
\[
 \frac{T^{7/2}e^{-\pi T/4}}{(1/2+\delta)^{(n+1)/2}}
 +(2T)^{7/2}e^{-\pi T/2}2^{(n+1)/2}.
\]
The logarithms of both terms are below $-2\cdot10^8$ at $n_*$.
Both increase with $n$. Keeping the same contour, proper-power
and diagonal bounds as in that theorem now gives, throughout this range,
\[
   Q_n<9\cdot10^{-6}
          -\frac{1}{2560}\left(\frac{n+1}{2}-40\right)<-0.033 .
\]
The script \texttt{anc/prime\_window\_check.py} reproduces the
bounds with outward-rounded arithmetic.
\end{proof}

\begin{corollary}[Induction in the degree]
\label{co:primewindowinduction}
With $u_0,n_*$ as in Proposition~\ref{pr:primewindowprefix}, one has
\[
 Q_0<-0.0015,\qquad Q_{n+1}<\frac{2u_0}{n+1}Q_n
                 \quad(0\le n<n_*).
\]
The same strict step holds at every degree on RH. If a right-hand
zero exists, the inequality is reversed at every sufficiently large
degree.
\end{corollary}
\begin{proof}
The exact one-step remainder is
\[
 I_n=\frac{n+1}{2}Q_{n+1}-u_0 Q_n
     =\int_0^\infty w_{n,2}(u)(u-u_0)Q(u)\,\dd u .
\]
For $n<30{,}000$, the integrand is nonpositive up to $U_*$.
After multiplication by $n!/2^n$, the contribution from $[1,4]$
is at most $-0.35(0.002559)$, since $u_0<0.65$ and $u^n\ge1$
there. The upper bound $Q(u)\le576e^u$ and
$u^{30{,}000}\le e^{u/2}$ for $u\ge U_*$ bound the remaining
positive contribution by $2304e^{-U_*/2}$. This proves $I_n<0$
on the initial range.

For $n\ge30{,}000$ an upper bound for the entire diagonal suffices.
Write $y=e^u$, $v=p/y$. Use
$\log p\le u+v$, $|K(v)|\le(6/\pi)e^{-\pi v^2/8}$ and
$v e^{-\pi v^2/8}\le1$. Partial summation of the two decreasing
Gaussians against $\vartheta(t)\le2t$ gives
\[
 D(u)\le (72/\pi^2)(u+\sqrt2)<8(u+2)\qquad(u\ge0).
\]
Hence $Q_n\ge-D_n^{\mathbb P}\ge-8[(n+1)/2+2]$, so that
$2u_0 Q_n/(n+1)>-5.21$. For $n+1\le n_*$ the upper
bound for $Q_{n+1}$ in Proposition~\ref{pr:primewindowprefix}
is already less than $-5.8$. Together these prove the strict step.
The negative base value now propagates through the entire stated range.

Under RH the integrand defining $I_n$ is strictly negative
away from its zero and the endpoints, by the proposition. Suppose
instead that $b=\Theta-1/2>0$. For
$B_n^{\mathbb P}=\int w_{n,2}P^2$, the root limit is
$L=1/(1-b)>1$ by Theorem~\ref{th:distinctprimeenergy}.
The positive moments $n! B_n^{\mathbb P}/2^{n+1}$ are log convex.
The secant-slope argument in Theorem~\ref{th:positivepowerratio}
therefore gives $B_{n+1}^{\mathbb P}/B_n^{\mathbb P}\to L$.
Since $D_n^{\mathbb P}=O(n)$, it follows that $Q_n>0$
eventually and $Q_{n+1}/Q_n\to L$. Thus $I_n>0$ eventually.
Eventual positivity in this case precludes propagation of the strict
negative step to all degrees.
\end{proof}

\head{Sign changes of the arithmetic coefficients}
The prime expansion has several changes of coefficient sign, so
a single-sign-change argument does not extend the finite-degree bound.
Put $z=\pi e^{-2u}/4$, $\mathcal C=\{1,4,16\}$ and
$(d_1,d_4,d_{16})=(1,-3,2)$. For $m\ge1$ form $C_m$
by summing $d_c d_d\log p\log q$ over $p\ne q$,
$c,d\in\mathcal C$, $cp^2+dq^2=m$. Form $A_m$ with
the extra factor $cdp^2q^2$. Then
\[
 4\pi^2 e^u Q(u)=\sum_{m\ge1}e^{-mz}
          \{4z^2 A_m+(1/4-mz)C_m\}.
\]
This is an absolutely convergent identity at each $z>0$ by
Gaussian decay. For every prime $p\ge17$ both $C_{p^2+9}$
and $A_{p^2+9}$ are strictly positive. Reduction modulo four
forces both summands to be odd prime squares with $c=d=1$,
and $p,3$ contribute. Both $C_{p^2+36}$ and $A_{p^2+36}$
are strictly negative. This integer is five modulo eight, so
one summand is an odd prime square and the other is four times
an odd prime square, except possibly for the summand $2^2$.
That possibility is excluded because $p^2+32$ lies strictly
between $p^2$ and $(p+1)^2$. Every contributing coefficient
therefore has factor $-3$, and $p^2+4\cdot3^2$ is one of
the representations. These opposite signs recur at unbounded
indices in both coefficient sequences. The independent positive
proportion required by Theorem~\ref{th:primewindowdensity} remains
unproved for the actual primes.

\subsection{Laguerre coefficients in the degree recurrence}
\label{sec:primelaguerre}
The gamma weight in Corollary~\ref{co:primewindowinduction} admits
an orthogonal expansion whose coefficients use only the signed
moments of the actual primes. Put $N=n+1$, and write
\[
\begin{aligned}
 m_j^{\mathbb P}&=\int_0^\infty w_{j,2}(u)P(u)\,\dd u,\\
 c_{n,k}&=\binom{n+k}{k}^{1/2}
       \sum_{j=0}^k(-1)^j\binom{k}{j}m_{n+j}^{\mathbb P}.
\end{aligned}
\]
Here $k$ is the polynomial order and $n$ specifies the averaging
measure. Define a finite sum, taken to be zero if the range is empty,
\[
 S_n^{\mathrm L}=\sum_{10^{24}N<k\le N^3}|c_{n,k}|^2 .
\]

\begin{theorem}[Two uniform bounds for the prime moments]
\label{th:primelaguerre}
For every $n\ge250$, with no upper limit on $n$,
\[
 0\le B_n^{\mathbb P}-S_n^{\mathrm L}\le
 \epsilon_n,\qquad
 \epsilon_n=9\cdot10^{-6}+5100N(16/e^3)^N.
\]
The first term bounds the orders $k\le10^{24}N$. The second
bounds all $k>N^3$ directly from the sum over primes. Consequently
the exact induction remainder has the asymptotic form
\[
 I_n=\frac{N}{2}S_{n+1}^{\mathrm L}-u_0 S_n^{\mathrm L}
                -\frac{C_0}{4}N(N+1)+O(N).
\]
In particular, $\liminf_{n\to\infty}S_n^{\mathrm L}/(n+1)<C_0/2$
implies RH. The constant is $C_0/2=0.0124420238212869\ldots$.
\end{theorem}
\begin{proof}
Under $x=2u$ the measure becomes
$\dd\mu_n=x^ne^{-x}\dd x/n!$. The polynomials
$L_k^{(n)}(x)/\binom{n+k}{k}^{1/2}$ form an orthonormal
basis in $L^2(\mu_n)$ \cite[Chap.~18]{Olver}.
For completeness, a function orthogonal to every monomial has a
Laplace transform against $\mu_n$ holomorphic in $\Re z<1/2$
with kernel $e^{zx}$, by Cauchy--Schwarz. All its derivatives at
zero vanish, and Fourier uniqueness then gives the zero function.
Integrating the finite Laguerre polynomial against $f(x)=P(x/2)$
gives precisely $c_{n,k}$, since
$\binom{n+k}{k-j}\binom{n+j}{j}=\binom{n+k}{k}\binom{k}{j}$.
Parseval therefore gives $B_n^{\mathbb P}=\sum_{k\ge0}|c_{n,k}|^2$.

First take $f(2u)=e^{\lambda u}$, with $d=\Re\lambda<1$. Its
coefficient at order $k$ is
\[
 c_{n,k}(\lambda)=\binom{n+k}{k}^{1/2}
       \left(\frac{2}{2-\lambda}\right)^N
       \left(\frac{-\lambda}{2-\lambda}\right)^k.
\]
This follows either by the preceding binomial sum or by Laplace
transformation. Summing its squared modulus at all $k$ gives
$(1-d)^{-N}$. Let $T=3\cdot10^{12}$, $A=10^{24}$,
$r=1-7/T^2$. For $\lambda=\rho-1/2$, $|\Im\rho|>T$,
multiplication by $r^{k-AN}\ge1$ at $k\le AN$ yields
\[
\begin{aligned}
 \sum_{k\le AN}|c_{n,k}(\lambda)|^2
 &\le \left[
 \frac{4r^{-A}}{4(1-d)+(1-r)|\lambda|^2}
 \right]^N\\
 &\le \left(\frac{4r^{-A}}9\right)^N<0.97^N .
\end{aligned}
\]
Indeed $d\le1/2$, and $-\log(1-7/T^2)\le7/(T^2-7)$.
Use the explicit formula for $g$ and $P=g-E$ from
Theorem~\ref{th:distinctprimeenergy}. All zeros through $T$ lie
on the line \cite{pt}, and their total absolute amplitude is less
than $0.000222$. Orthogonal projection is a contraction. Above
$T$ the sum of the absolute amplitudes is at most
$T^{7/2}e^{-\pi T/4}$, and the preceding calculation controls
the projected norm of each exponential. The shifted contour has
norm at most $100(2/5)^{N/2}$, and the proper powers satisfy
$\|E\|_{L^2(w_{n,2})}<0.0027$ for all $n\ge250$.
The triangle inequality now bounds the squared norm at orders
$k\le AN$ by $9\cdot10^{-6}$. For each fixed $n$ the series
over zeros converges also in the full norm since each exponential
has norm at most $2^{N/2}$.

For the other end of the expansion only the prime kernel is needed.
In $|\Im z|\le1/2$, put $y=e^{\Re z}$. The three Gaussians,
$\cos(1)>1/2$, and $(2t+1/2)e^{-t/2}\le4e^{-t/4}$ for
$t\ge0$ give
\[
 |K(pe^{-z})|\le(12/\pi)\exp[-\pi(p/y)^2/16].
\]
Using $\vartheta(t)\le2t$ in partial summation yields
$|P(z)|\le (48/\pi)e^{\Re z/2}$. The series is holomorphic on
the larger strip $|\Im z|<\pi/4$. Cauchy's inequality on a circle
of radius one about each real $x\ge0$ for $f(x)=P(x/2)$
therefore gives
\[
 |f^{(j)}(x)|\le20\,j!\,e^{x/4}\qquad(j\ge0).
\]
Write $k^{\underline j}=k!/(k-j)!$ for $k\ge j$.
The Rodrigues formula and $j$ integrations by parts give
\[
 \sum_{k\ge j} k^{\underline j}|c_{n,k}|^2
       \le\int_0^\infty x^j|f^{(j)}(x)|^2\,\dd\mu_n(x).
\]
In detail, the inner product of $f^{(j)}$ with
$L_{k-j}^{(n+j)}$ under $x^j\dd\mu_n$ equals
$(-1)^j k^{\underline j}$ times the inner product of $f$
with $L_k^{(n)}$ under $\dd\mu_n$. The squared norm of the former
polynomial is $k^{\underline j}\binom{n+k}{k}$. Bessel's inequality
gives the display. The exponential derivative bound and the powers
of $x$ in the Rodrigues formula make all boundary terms zero.
Take $j=N$. Direct gamma integration gives
\[
 \sum_{k>N^3}|c_{n,k}|^2
 \le \frac{200\,4^N N!(2N)!}{(N^3+1)^{\underline N}}
 \le 5100N(16/e^3)^N.
\]
For $N\ge2$ the denominator is at least
$N^{3N}(1-1/N^2)^N\ge N^{3N}/2$, and
$m!\le3\sqrt m\,(m/e)^m$ by Stirling's bounds. The case $N=1$
also satisfies the stated high-order bound. This completes both
estimates, including when the two bounded ranges overlap.
The constants and polynomial normalization checks are reproduced in
\texttt{anc/laguerre\_check.py}.

Let $J_n^D=(N/2)D_{n+1}^{\mathbb P}-u_0D_n^{\mathbb P}$ and
$Y_n=(N/2)S_{n+1}^{\mathrm L}-u_0 S_n^{\mathrm L}-J_n^D$.
There is a finite-error version of the induction formula,
\[
       Y_n-u_0\epsilon_n\le I_n\le Y_n+(N/2)\epsilon_{n+1}.
\]
The unconditional diagonal asymptotic in
Theorem~\ref{th:distinctprimeenergy} gives
$J_n^D=C_0N(N+1)/4+O(N)$. If the stated strict liminf bound
holds, apply the upper inequality at the degrees immediately
preceding its subsequence, using $S_n^{\mathrm L}\ge0$. This
gives $I_n<0$ at unbounded degrees. Corollary~\ref{co:primewindowinduction}
then excludes a right-hand zero. Conversely, on RH
$S_n^{\mathrm L}\le B_n^{\mathbb P}=O(1)$.
\end{proof}

\head{The signed prime moment at unbounded degrees}
The coefficients above admit a direct calculation from the Euler
product. Put $\mathcal P(s)=\sum_p(\log p)p^{-s}$ initially in
$\Re s>1$ and $\widehat K(s)=(s-1/2)A(s)$ as in
Theorem~\ref{th:gaussfirst}. At $s=5/2-2z$, on a neighborhood
of $z=0$, absolute convergence and a change of variables give
\[
 \sum_{n\ge0}m_n^{\mathbb P}z^n
 =2\widehat K(s)\mathcal P(s)
       -2\int_{-\infty}^0 e^{-(s-1/2)u}P(u)\,\dd u .
\]
The integral on the right is entire in $s$ by Gaussian decay.
Write $L=-\zeta'/\zeta$. M\"obius inversion gives
$\mathcal P(s)=\sum_{j\ge1}\mu(j)L(js)$.
In the closed disk $|s-5/2|\le16/7$ one has
$\Re s\ge3/14$. Terms $j\ge5$ are holomorphic on a neighborhood
of the disk and converge normally there, and $\mu(4)=0$.
For $j=1,2,3$, any zero $\rho=js$ in the critical strip
meeting this disk would have
\[
 |\Im\rho|\le j\sqrt{(16/7)^2-(5/2-1/j)^2}<3 .
\]
There are no such zeros. For an analytic check, the Hadamard
product and functional equation give
$\sum_\rho\Re(1/\rho)=1+\gamma_{\mathrm E}/2-\log(4\pi)/2<0.024$
\cite{mv}. A conjugate pair with $\re(\rho)\ge1/2$ and
$|\Im\rho|\le6$ would contribute at least $1/37>0.024$.
Reflection supplies such a pair from either side of the line,
and the alternating eta series excludes real zeros in $(0,1)$.
The poles at $s=1,1/2$ are canceled by $\widehat K$.
The remaining pole at $s=1/3$ has residue
$-\widehat K(1/3)/3$ in $\widehat K(s)\mathcal P(s)$.
After its subtraction, Cauchy's estimate on $|z|=8/7$ gives
\[
\begin{gathered}
 m_n^{\mathbb P}
  =-\frac{\widehat K(1/3)}{3}(12/13)^{n+1}
             +O((7/8)^n),\\
 \widehat K(1/3)/3=0.00310458352895125\ldots .
\end{gathered}
\]
The geometric term by itself contributes exactly
$(\widehat K(1/3)/3)^2(6/7)^N$ to the full squared Laguerre
norm, by the single-exponential formula with $\lambda=-1/6$.
In particular the signed mean is eventually negative at every degree.
The remainder $e_n$ after the geometric term satisfies the explicit
bound $|e_n|\le 600{,}000(7/8)^n$. The estimate follows from
analytic bounds on the Cauchy circle. On this circle
$3/14\le\Re s<5$ and $|js-1|>1/3$ for $1\le j\le3$.
For $w=js$ with $\Re w\le3/2$, the circle gives
$|\Im w|<4$ and $|w|<4.5$.
Reflection in the critical line and $|\Im\rho|>6$ give
\[
 \sum_\rho|\rho|^{-2}
 \le(1+37/36)\sum_\rho\Re(1/\rho)<0.05,
\]
since $|1-\bar\rho|^2\le(37/36)|\rho|^2$. Thus
$|w-\rho|\ge|\rho|/4$, and the logarithmic derivative of the
Hadamard product gives $|\xi'(w)/\xi(w)|<0.024+4(4.5)(0.05)<1$.
In the completed zeta function the remaining factors are controlled by
\[
\begin{gathered}
 |1/w|<5,\quad |1/(w-1)|<3,\\
 |\Gamma'(w/2)/\Gamma(w/2)|\le 1+2/|w|+(|w|/2)\zeta(2)<15 .
\end{gathered}
\]
The last inequality follows from the digamma series on $\Re w>0$.
Hence $|L(w)|<20$. For $\Re w\ge3/2$, its prime-power
Dirichlet series gives the same bound by absolute summation.
The sum over $j\ge5$ in $\mathcal P$ is bounded by
$(1-2^{-3/14})^{-1}\sum_{m\ge2}(\log m)m^{-15/14}<1600$.
The first factor is less than eight and integral comparison bounds
the sum by $200$. Consequently $|\mathcal P(s)|<1700$.

The real Gaussian envelope in Theorem~\ref{th:distinctprimeenergy}
gives $|\widehat K(s)|\le (3/\pi)(8/\pi)^{\Re s/2}
\Gamma(\Re s/2)<150$, by splitting the gamma integral at one.
For $u\le0$, extract $e^{-\pi e^{-2u}/4}$ from this envelope
using $p\ge2$, and sum the remaining Gaussian against
$\vartheta(t)\le2t$. This bounds $|P(u)|$ by
$(24/\pi)e^{u/2}e^{-\pi e^{-2u}/4}$. With $x=e^{-2u}$ the
entire integral over $u\le0$ is at most
$(12/\pi)\int_1^\infty x e^{-\pi x/4}\dd x<7$ in modulus.
The distance from $13/12$ to the Cauchy circle in $z$ is
$5/84$. The modulus of the generating function after subtracting
its principal part is therefore at most $2(150)(1700)+2(7)+1
<600{,}000$, as asserted. It follows that
\[
 \left|\sum_{j=0}^k(-1)^j\binom{k}{j}e_{n+j}\right|
       \le 600{,}000(7/8)^n(15/8)^k .
\]
This gives an alternative reduction that uses no numerical
verification of zeta zeros and has exponentially small absolute
error. Set $\alpha=1/32$ and
$\widetilde S_n=\sum_{\alpha N<k\le N^3}|c_{n,k}|^2$. Then
\[
 0\le B_n^{\mathbb P}-\widetilde S_n\le 10^{12}(23/25)^N
                     \qquad(n\ge0).
\]
With $\widetilde Y_n=(N/2)\widetilde S_{n+1}
-u_0\widetilde S_n-J_n^D$, the original remainder consequently satisfies
\[
 \widetilde Y_n-u_0\,10^{12}(23/25)^N
 \le I_n\le \widetilde Y_n+(N/2)\,10^{12}(23/25)^{N+1}.
\]
Indeed $\sum_{k\le l}\binom{n+k}{k}=\binom{N+l}{l}$.
At $l=\lfloor\alpha N\rfloor$ this is at most
$[(1+\alpha)^{1+\alpha}\alpha^{-\alpha}]^N$, by the
binomial theorem. The squared norm of the projected remainder
is consequently at most $(600{,}000)^2(8/7)^2q_\alpha^N$, where
\[
 q_\alpha=(7/8)^2(15/8)^{2\alpha}
                    (1+\alpha)^{1+\alpha}\alpha^{-\alpha}<0.917 .
\]
The cube term has full squared norm at most $(0.004)^2(6/7)^N$.
The squared norm of the sum is at most twice the sum of these
bounds. Finally $N[(16/e^3)/(23/25)]^N\le\sum_{j\ge1}j(0.87)^j<60$.
Adding the order bound above $N^3$ gives the stated constant,
since $2(600{,}000)^2(8/7)^2+2(0.004)^2+5100(60)<10^{12}$.
At larger $k$ the factor $(15/8)^k$ no longer controls the
binomial differences uniformly at the scale of the diagonal.

\head{Upper bounds for the intermediate Laguerre energy}
Direct partial summation against $\vartheta(t)-t$, using
$\int K=0$ and the Vinogradov--Korobov bound \cite{johnstonyang},
gives $|P(u)|\ll \exp[u/2-c\omega(u)]$ for large $u$,
where $\omega(u)=u^{3/5}(\log u)^{-1/5}$. This follows by the
split at $v=e^{-u/2}$ in Sec.~\ref{sec:primelossbudget}.
Since $e^u w_{n,2}=2^Nw_{n,1}$, concentration of the latter
gamma law gives $S_n^{\mathrm L}\le B_n^{\mathbb P}
\ll 2^N e^{-c'\omega(N)}$ for some $c'>0$.
The Gaussian kernel gives a stronger unconditional estimate.
For large $N$, let
\[
 \delta_N=\frac{1}{55.241(\log N)^{2/3}(\log\log N)^{1/3}}.
\]
The Vinogradov--Korobov zero-free region of Mossinghoff, Trudgian
and Yang \cite{mtyregion} gives $\re(\rho)\le1-\delta_N$ at
$3\le|\Im\rho|\le N$. In the full gamma norm the single
exponential formula therefore costs at most
$[2/(1+2\delta_N)]^{N/2}$ for each such mode. Their total
absolute amplitude is bounded independently of $N$. Above $N$
the same norm is bounded by $2^{N/2}$ per mode and the sum
of amplitudes is at most $N^{7/2}e^{-\pi N/4}$.
The contour and the proper powers have uniformly bounded norms.
By the triangle inequality one obtains
\[
           S_n^{\mathrm L}\le B_n^{\mathbb P}
                    \ll \left(\frac{2}{1+2\delta_N}\right)^N .
\]
The right-hand side has limiting $N$-th root $2$.
For one spectral term $e^{\lambda u}$, the normalized squared
Laguerre coefficients have a negative binomial distribution with
mean $N|\lambda|^2/[4(1-\Re\lambda)]$ and variance proportional
to $N$ for fixed $\lambda$. A fixed right-hand zero above $T$
places this mean strictly inside the intermediate range for sufficiently
large $n$, with single-mode squared norm growing as
$(1-\Re\lambda)^{-N}$. In the full signal the modes are summed before
squaring. The lower bound on the growth of that sum follows from
Theorem~\ref{th:distinctprimeenergy}. A proof directly from the primes
requires an estimate on $S_n^{\mathrm L}$ at unbounded degrees
sufficient to determine the sign of $I_n$.

\subsection{Curvature in the degree of the prime moments}
\label{sec:primecurvature}
Two different remainders occur in the induction. The truncation of
the Laguerre expansion discards a positive sum of squares. The
induction itself uses a signed integral over pairs of distinct primes.
We compute differences for the first remainder and derivatives with
respect to the real degree for the second.

\begin{proposition}[Differences of the discarded Laguerre energy]
\label{pr:truncationcurvature}
Let $N=n+1$, $E_n=B_n^{\mathbb P}-\widetilde S_n$,
$\delta E_n=E_{n+1}-E_n$, $q=6/7$ and $\kappa=0.853$.
With $A=[\widehat K(1/3)/3]^2=0.00000963843888823543\ldots$,
for every integer $n\ge0$ and $j\ge0$,
\[
 \left|\delta^jE_n-A(q-1)^jq^N\right|
       \le 10^{12}(1+\kappa)^j\kappa^N.
\]
In particular $E_n>0,\ \delta E_n<0,\ \delta^2 E_n>0$ for
$N\ge10^4$. At the initial degrees one has
\[
 E_2<0.000095,\qquad E_3>0.000100,\qquad E_4<0.000075 .
\]
Thus the discrete slope changes sign at degree three and
$\delta^2E_2<-0.000030$.
\end{proposition}
\begin{proof}
Enlarge the Cauchy circle used in Sec.~\ref{sec:primelaguerre} to
$|z|=19/16$, so that $|s-5/2|=19/8$ and $\Re s\ge1/8$.
The terms of M\"obius index $j\ge9$ in $\mathcal P(s)$
converge normally on a neighborhood of the disk, and
$\mu(4)=\mu(8)=0$. For $0<\Re(js)\le1$ in the disk one has
\[
 |\Im(js)|^2\le-39j^2/64+5j-1\le361/39<16 .
\]
The Hadamard estimate in the preceding subsection excludes zeros
at these heights. After cancellation at $s=1,1/2$, the poles
that remain arise from $j\in\mathcal J=\{3,5,6,7\}$. Put
$a_j=\mu(j)\widehat K(1/j)/j$, $r_j=2/(5/2-1/j)$,
and $r=16/19$. Subtraction of the four principal parts gives
\[
 m_n^{\mathbb P}=\sum_{j\in\mathcal J}a_jr_j^N+e_n,
                 \qquad |e_n|\le500{,}000\,r^n .
\]
The uniform constant follows from estimates on the circle, where
$|js-1|\ge1-j/8$ for $1\le j\le7$. For $w=js$ with
$\Re w\le3/2$, $|w|^2\le(15/2)j-39j^2/64\le300/13<25$.
Thus $|w-\rho|>|\rho|/6$, and the preceding Hadamard estimate
gives $|\xi'/\xi(w)|<0.024+6(5)(0.05)<2$. The digamma series
used there, $|w|\ge1/8$, and $|w-1|\ge1/8$ give
$|\zeta'/\zeta(w)|<30$. For $\Re w\ge3/2$ the absolutely
convergent Dirichlet series gives the same bound. The remaining
indices contribute in modulus at most
$(1-2^{-1/8})^{-1}\sum_{m\ge2}(\log m)m^{-9/8}
 <13(2+8^2)$, by integral comparison. Thus $|\mathcal P|<1100$.
The Gaussian envelope and the negative-time integral in the
preceding subsection give $|\widehat K|<200$ and an integral of
modulus below seven. Each $|a_j|<0.005$, and the four poles
in the $z$-plane are at least $1/112$ from the circle.
The resulting bound is $2(200)(1100)+14+112\sum|a_j|<500{,}000$.

Write $\alpha=1/32$. The squared norm of the coefficients from
$e_n$ at orders $k\le\alpha N$ is at most
$(500{,}000/r)^2p^N$, where
\[
 p=r^2(1+r)^{2\alpha}(1+\alpha)^{1+\alpha}\alpha^{-\alpha}
      <0.848.
\]
This follows from $|\sum_{i=0}^k(-1)^i\binom ki e_{n+i}|
\le500{,}000\,r^n(1+r)^k$ and the binomial sum in
Sec.~\ref{sec:primelaguerre}. The full squared norm of the
$j$-term is $|a_j|^2(3/2-1/j)^{-N}$, since its coefficients
are those of $a_j\exp[(1/j-1/2)u]$. For $j=5,6,7$ the
exponential base is at most $10/13$. The norm of all these
terms and the remainder together at $k\le\alpha N$ is less
than $600{,}000(0.848)^{N/2}$.
For the cube term the full squared norm is $Aq^N$.
The relative energies have probability generating function
$[(1-1/169)/(1-t/169)]^N$. At $t=2$ the energy above
$\alpha N$ is therefore at most
$A[q(168/167)2^{-1/32}]^N$. Both this base and the base
$\sqrt{0.848q}$ of the cross term are below $\kappa$.
Finally the squared norm of the actual signal above $N^3$ is at most
$5100N(16/e^3)^N$, with $\sum_{l\ge1}l[(16/e^3)/\kappa]^l<250$.
Adding the four error bounds gives the stated constant at $j=0$.
Taking differences proves the general estimate. At $N=10^4$ the
error for $j=2$ is less than $0.016$ of the leading term.
For $j=0,1$ that fraction is smaller, and it decreases with $N$.

For the small degrees let $P_8=P\mathbf1_{[0,8]}$ and
$h_n(t)=e^{-t}\sum_{i=0}^n t^i/i!$. From
Proposition~\ref{pr:primewindowprefix} and its proof,
$|P|<0.057$ on $[8,18]$, $|P|<0.3$ on $[18,64]$, and
$P^2\le576e^u$ globally for $u\ge0$. Thus
\begin{multline*}
 \|P-P_8\|_{L^2(w_{n,2})}^2\le b_n\\
 =0.057^2h_n(16)+0.3^2h_n(36)+576\,2^Nh_n(64).
\end{multline*}
For $n=2,4$, project $P_8$ on orders $1$ through $27,64$,
respectively. Both spaces lie in the retained orders. Add
$\sqrt{b_n}$ to the distance to each subspace before squaring.
For $n=3$ the omitted order zero gives
$E_3\ge(m_3^{\mathbb P})^2$, and the signed integral outside
$[0,8]$ has modulus at most $\sqrt{b_3h_3(16)}$.
On each half-unit interval in $[0,8]$ use 96-point Gauss--Legendre
quadrature with outward rounding. Truncating at $p\le\lceil16e^u\rceil$
has pointwise error below $(36/\pi)e^4e^{-16\pi}<10^{-18}$.
On the complex circle of radius $1/2$ about each midpoint,
$|P|<2000$ by the strip estimate, $|w_{n,2}|<500{,}000$, and
$|L_k^{(n)}(2u)|\le2^{68}e^{18}$ by its finite series.
Each integrand has modulus below $10^{38}$. Summing the Taylor
tail from degree 192 on the real intervals and using exactness
through degree 191 and positivity of the quadrature weights gives
a total error below $32\cdot10^{38}2^{-192}<10^{-17}$ per integral.
The tests in \texttt{anc/prime\_curvature\_check.py} then give
upper bounds $0.00009406,0.00007316$ and the lower bound
$0.00010097$, respectively.
\end{proof}
In particular $\delta E_n/E_n\to-1/7$ and
$\delta^2E_n/E_n\to1/49$. The signed truncation correction to
the induction itself is
\[
\begin{aligned}
 I_n-\widetilde Y_n&=(N/2)E_{n+1}-u_0E_n\\
        &=Aq^N(qN/2-u_0)+O(N\kappa^N).
\end{aligned}
\]
It is eventually positive, decreasing and discretely convex. Its
size leaves the contribution of the retained orders to $I_n$ open.

\begin{theorem}[The real-degree induction remainder]
\label{th:realdegreecurvature}
For $x>-1$ write $N=x+1$ and use the gamma density
$w_{x,2}(u)=2^Nu^{N-1}e^{-2u}/\Gamma(N)$. Set
\[
 I(x)=\int_0^\infty w_{x,2}(u)(u-u_0)[P(u)^2-D(u)]\,\dd u .
\]
It extends $I_n$ at the integers. On RH, as $x\to\infty$,
\[
\begin{gathered}
 I(x)\sim-C_0N^2/4,\qquad I'(x)\sim-C_0N/2,\\
 I''(x)\longrightarrow-C_0/2 .
\end{gathered}
\]
If a zero lies to the right, let $L=(3/2-\Theta)^{-1}>1$.
Then $I(x)>0$ eventually, $I(x)^{1/N}\to L$, and
\[
 I'(x)/I(x)\longrightarrow\log L,\qquad
 I''(x)/I(x)\longrightarrow(\log L)^2 .
\]
Thus boundedness from above of $I''$ along any sequence tending
to infinity is equivalent to RH.
\end{theorem}
\begin{proof}
Write $M_f(x)=\int w_{x,2}f$, $B(x)=M_{P^2}(x)$,
$\ell=\log(2u)-\psi(N)$, where $\psi=\Gamma'/\Gamma$.
Differentiation under the integral gives
\[
 M_f'(x)=M_{f\ell}(x),\qquad
 M_f''(x)=M_{f(\ell^2-\psi'(N))}(x).
\]
On compact subsets of $x>-1$ the exponential tail and local
integrability of $u^{x}(1+|\log u|^2)$ justify both operations
for the functions here. The scores have absolute means at most
$\sqrt{\psi'(N)}$ and $2\psi'(N)$, since
$M_\ell=0$ and $M_{\ell^2}=\psi'(N)$ \cite[Chap.~5]{Olver}.
With $C_0=\int K^2$ and $C_{\log}=\int K^2\log v\,\dd v$,
the polynomial part of the diagonal in $I$ is
\[
 J_D(N)=C_0 N(N+1)/4+(C_{\log}-u_0C_0)N/2-u_0 C_{\log}.
\]
Partial summation, as in Sec.~\ref{sec:positivepowers}, gives
$D(u)-C_0u-C_{\log}=O_a((1+u)^{-a})$ for every fixed $a>0$.
After gamma averaging this remains $O_a(N^{-a})$, also after
multiplication by $u-u_0$ and up to two degree derivatives.
On $u\ge N/4$, increasing the exponent in the partial summation
bound absorbs the score factors by $u w_{x,2}=(N/2)w_{x+1,2}$.
On the complementary interval the gamma probability is exponentially
small, also after polynomial factors in $u,\log(2u),\log N$,
by the exponential Markov bound and finite inverse gamma moments.

On RH the explicit formula and the proper-power estimate give
uniform approximation of $P^2$ on sufficiently late tails by
trigonometric polynomials (with bounded error on compact intervals).
A mode $e^{itu}$ has mean $z_t^N$ for
$z_t=2/(2-it)$. Its derivatives have multipliers
$(\log z_t)^j z_t^N$, using the principal logarithm.
For $t\ne0$ each fixed mode decays exponentially, and the
constant has vanishing derivatives. For the uniform approximation
error use the two score bounds and $\psi'(N)\sim1/N$.
The integral on a fixed initial interval is negligible, also with the
scores. Letting the uniform error tend to zero gives
$B=O(1)$, $B'=o(N^{-1/2})$, and $B''=o(N^{-1})$.
Consequently the second derivative of
$(N/2)B(x+1)-u_0 B(x)$ tends to zero, and its first derivative
is $o(N)$. Subtracting $J_D$ and its derivatives proves
the three assertions on RH.

If $L>1$, Corollary~\ref{co:primewindowinduction} gives
$B(n+1)/B(n)\to L$. The function
$\Phi(x)=\log\int_0^\infty u^x e^{-2u}P(u)^2\,\dd u$
is convex by H\"older's inequality. The secant slope on
$[n,n+1]$ is $\log((n+1)/2)+\log L+o(1)$.
Bounding $\Phi'$ there by the slopes on its two neighboring
integer intervals and using $\psi(N)=\log N+o(1)$ yields
$B'/B\to\log L$. Hence $B(x+j)/B(x)\to L^j$ for every
fixed real $j$. Under the probability measure with density
$w_{x,2}P^2/B(x)$, let $Y=2u/N$. Its $j$-th moment is
$\Gamma(N+j)B(x+j)/(N^j\Gamma(N)B(x))$. At $j=1,2,-2$
these moments tend to $L,L^2,L^{-2}$. Thus $Y\to L$ in mean
square, and $|\log Y|^4\ll Y^2+Y^{-2}$ gives uniform integrability
of $(\log Y)^2$. In the second-derivative identity this gives
$B''/B\to(\log L)^2$. The term $(N/2)B(x+1)$ now dominates
$u_0B(x)$ and the polynomial diagonal, with the same dominance
for the first two derivatives. This proves the alternative.
\end{proof}

\begin{proposition}[An unconditional interval of strict concavity]
\label{pr:finitecurvature}
Throughout the real interval $10^7\le x+1\le10^{12}$,
\[
 I'(x)<0,\qquad \left| I''(x)+C_0/2\right|<0.00021 .
\]
In particular $I(x),I'(x),I''(x)$ are all strictly negative there.
\end{proposition}
\begin{proof}
First, $|D(u)-C_0u-C_{\log}|<100/u$ for $u\ge100$.
Integration of $\log(yv)K(v)^2$ against
$\dd[\vartheta(yv)-yv]/y$, with $y=e^u$, permits integration by
parts with vanishing boundary terms. The three Gaussians give
$|K(v)|<2e^{-\pi v^2/8}$ and
$|K'(v)|\le40v e^{-\pi v^2/8}$, using
$(5/2+2t)e^{-t/2}\le3$. Consequently
$\int v|KK'|\le160/\pi<51$ and
$\int v|\log v|\,|KK'|\le80(2/\pi+8/\pi^2)<120$,
using $|\log v|\le v+1/v$.
Below $w=e^{-u/2}$, the elementary counting error $2yv$
costs at most $8w+320w^3(u/2+1/9)$.
Above $w$, use Dusart's bound $0.2yv/\log^2(yv)$ at
$yv\ge e^{50}>3{,}594{,}641$ \cite{dusarttheta}. Since
$C_0<1$, this costs at most $0.8(1+102u+240)/u^2$.
The sum is less than $100/u$.

Write $U_m=10^6$, $U_*=4.77\cdot10^{12}$, and
$J=P^2-D+C_0u+C_{\log}$.
On $[U_m,U_*]$ the explicit bound on $P$ from
Proposition~\ref{pr:primewindowprefix} gives $|P|<0.000224$.
Thus $|J|<0.000101$ there. Globally on $u\ge0$ the earlier
bounds $P^2\le576e^u$, $D\le8(u+2)$ give $|J|<610e^u$.
The digamma and trigamma recurrences give the two expectations
\[
\begin{gathered}
 M_{(u+u_0)(\ell^2+\psi'(N))}=(N+2u_0)\psi'(N)<2,\\
 M_{(u+u_0)|\ell|}\le(N/2+u_0)\sqrt{\psi'(N)}<\sqrt N .
\end{gathered}
\]
Indeed under the adjacent shape $N+1$ the second moment of
$\ell$ is still $\psi'(N)$. Here
$\psi'(N)\le1/N+1/N^2$ by its defining series.

To bound both tails, set $v=2u/N$. Both absolute scores are
bounded by $8(v^2+v^{-2})$ since
$0<\log N-\psi(N)<1/N$ (by the integral bounds on the
digamma series). On $u<U_m$ multiply the global estimate by
$e^{3(U_m-u)}$, and on $u>U_*$ multiply by $e^{(u-U_*)/2}$.
Completing the gamma densities leaves rates $4$ and $1/2$,
respectively. At each rate $a$, the exact moments
$\Gamma(N+j)/(a^j\Gamma(N))$ imply
$\mathbb E[(1+u)(v^2+v^{-2})]<100(N+1)$ for $N\ge10$.
The combined tail for either derivative is therefore at most
\[
 8(610)(100)(N+1)
 \left[2^{-N}e^{3U_m}+4^Ne^{-U_*/2}\right]<10^{-15}
\]
throughout the interval in the statement. We have
$I+J_D=M_{(u-u_0)J}$ exactly. The central bound now gives
$|I''+C_0/2|<2(0.000101)+10^{-15}$, and
$|I'+J_D'|<0.000101\sqrt N+10^{-15}$.
The linear expression for $J_D'$ proves the strict sign of
$I'$. Corollary~\ref{co:primewindowinduction} gives the sign of
$I$ at the left endpoint. By Theorem~\ref{th:realdegreecurvature},
an upper bound for $I''$ along any sequence tending to infinity
would imply RH.
\end{proof}

\head{Differentiation in the prime cutoff}
The independent variable matters. Consider the centered cutoff slope
\[
 d(X)=\sum_{m\le X}\frac{\Lambda(m)}m-\log X+\gamma_{\rm E},
       \qquad X>1 .
\]
Its distributional derivative is the signed measure
\[
 \mathrm d d(X)=-\frac{\mathrm dX}{X}
          +\sum_{p,\,j\ge1}\frac{\log p}{p^j}\,\delta_{p^j}.
\]
In particular every late tail contains intervals of negative density
and positive atoms. The slope itself takes both signs beyond every
threshold. Indeed, for each nontrivial zero $\rho=\beta+i\tau$
of multiplicity $m_\rho$, one has
\[
\begin{aligned}
 \limsup_{X\to\infty} X^{1-\beta}d(X)&\ge
       \frac{m_\rho}{|1-\rho|},\\
 \liminf_{X\to\infty} X^{1-\beta}d(X)&\le
       -\frac{m_\rho}{|1-\rho|}.
\end{aligned}
\]
For $\Re s>0$, direct integration of the prime-power series gives
\[
 F(s)=\int_1^\infty d(X)X^{-s-1}\,\dd X
    =-\frac{\zeta'(1+s)}{s\zeta(1+s)}-\frac1{s^2}
                                      +\frac{\gamma_{\rm E}}s .
\]
The Laurent expansion cancels the pole at zero. The function is
regular on the real interval $(-1,\infty)$, since the alternating
eta series excludes real zeros of $\zeta$ in $(0,1)$.
Its pole at $s_0=\rho-1$ has residue $-m_\rho/s_0$.
Let $a=\beta-1$ and suppose that $H(X)=cX^a\pm d(X)$
is nonnegative for all $X\ge X_0$, with $c>0$. Its tail
Mellin transform is $c/(s-a)\pm F(s)$ up to an entire function.
It converges absolutely at all real $\sigma>a$ by Landau's
theorem for Mellin integrals. A nonnegative tail with finite
abscissa $b$ has a singularity at $b$. Otherwise a Taylor
disk centered a little to the right of $b$ would extend to its
left. At a real point to the left of $b$ the Taylor series sums
the nonnegative log moments by monotone convergence, which gives
convergence beyond the abscissa. Here the real transform is regular
everywhere above $a$. Write $G$ for the tail transform. Positivity
gives $|G(\sigma+i\tau)|\le G(\sigma)$. Multiply by $\sigma-a$
and let $\sigma\downarrow a$. The residue at $s_0$ gives
$m_\rho/|s_0|\le c$. Applying this to the two signs gives the
displayed inequalities. Thus oscillation is already present at the
cutoff slope even for a zero on the critical line.

\subsection{Extracting a degree of negative curvature from prime sums}
\label{sec:twodegreecurvature}
Set $B(x)=M_{P^2}(x)$ and $B_n=B_n^{\mathbb P}$. An estimate
at one integer moment forces negative curvature at some degree beyond
a prescribed threshold. Two neighboring moments instead bound the
curvature at a specified degree.

\begin{proposition}[A finite quadratic form and a one-moment induction]
\label{pr:curvatureinduct}
For $n\ge1$, put $N=n+1$, $Y=e^{9N}$ and define
\[
\begin{aligned}
 P_Y(u)&=e^{-u/2}\sum_{p\le Y}(\log p)K(pe^{-u}),\\
 \mathcal B_n&=\int_0^{8N} w_{n,2}(u) P_Y(u)^2\,\dd u .
\end{aligned}
\]
For every such $n$,
\[
 |\mathcal B_n-B_n|\le1000(4e^{-4})^N .
\]
If $\liminf_{n\to\infty}[\mathcal B_n-C_0(n+1)/2]=-\infty$,
there is a sequence $x_j\to\infty$ with $I''(x_j)<0$.
\end{proposition}
\begin{proof}
Both prime sums are bounded by $24e^{u/2}$ in modulus using the
Gaussian majorant for $|K|$ and $\vartheta(t)\le2t$.
For $0\le u\le8N$, an omitted prime satisfies $pe^{-u}>e^N$.
Extracting half the exponent in the bound
$(6/\pi)e^{-\pi v^2/8}$ and applying partial summation to the
remaining Gaussian gives
\[
 |P-P_Y|\le (24/\pi)e^{u/2}\exp(-\pi e^{2N}/16).
\]
Since $M_{e^u}(n)=2^N$, the resulting error in the integral of
the square is at most $384\,2^N\exp(-\pi e^{2N}/16)$. The tail
above $8N$ is bounded by $576\,4^Ne^{-4N}$, by multiplying
$576e^u$ by $\exp[(u-8N)/2]$ and integrating.
For $N\ge2$ one has $\pi e^{2N}/16>(4-\log2)N$, proving the
stated error.

All the weights in this finite sum over pairs can also be evaluated
without integrating the prime sum. For each pair of Gaussian terms
let $a=\pi c p^2/4$, $b=\pi d q^2/4$, $c,d\in\{1,4,16\}$,
and $L=a+b$. Put
\[
 H_n(z,L)=(-\partial_z)^n\!
 \left[L^{-z}\{\gamma(z,L)-\gamma(z,Le^{-16N})\}\right],
\]
where $\gamma(z,t)=\int_0^t r^{z-1}e^{-r}\,\dd r$.
After the change $t=e^{-2u}$, the weight of
$d_c d_d(\log p)(\log q)/(4\pi^2n!)$ is
$4ab H_n(7/2,L)-L H_n(5/2,L)+H_n(3/2,L)/4$,
with $(d_1,d_4,d_{16})=(1,-3,2)$.
For example, putting $q_n(v)=e^{-v}\sum_{j=0}^n v^j/j!$,
direct expansion of $e^{-Lt}$ gives the absolutely convergent series
\[
\begin{gathered}
 \frac{H_n(z,L)}{n!}
 =\sum_{l\ge0}\frac{(-L)^l\,[1-q_n(16N(z+l))]}{l!\,(z+l)^N}\\
              (z,L>0).
\end{gathered}
\]
The absolute tail after $l=M$, for $M+2>L$, is at most
$L^{M+1}/[(M+1)!(z+M+1)^N(1-L/(M+2))]$.

For the induction let $d_n=B_n-C_0(n+1)/2$ and
$d_*=C_{\log}-u_0 C_0$. The polynomial diagonal computed in
Theorem~\ref{th:realdegreecurvature} gives the exact asymptotic
identity, for every fixed $A>0$,
\[
 I(n-1)=\tfrac n2(d_n-d_*)-u_0B_{n-1}
                 +u_0 C_{\log}+O_A(n^{-A}).
\]
Fix any real $t>-1$. The bound $B_{n-1}\ge0$ and a subsequence
with $d_n\to-\infty$ imply that for some $n>t+1$,
\[
 I(n-1)<I(t)+(n-1-t)I'(t).
\]
Taylor's theorem supplies $\xi\in(t,n-1)$ with
$I''(\xi)=2[I(n-1)-I(t)-(n-1-t)I'(t)]/(n-1-t)^2<0$.
Repeat with $t=\xi+1$. This constructs the required sequence
under the stated arithmetic hypothesis. The finite truncation
changes $d_n$ by $o(1)$. Conversely, on RH $B_n=O(1)$, so the
deficit tends to minus infinity. Together with
Theorem~\ref{th:realdegreecurvature}, the arithmetic hypothesis is
therefore equivalent to RH.
\end{proof}

\begin{proposition}[Curvature from two prime moments]
\label{pr:twodegreecurvature}
For real $N\ge25$, put $h=\sqrt N$ and
\[
\begin{gathered}
 T(N)=\tfrac12(1-h^{-1})B(N-h)+\tfrac12(1+h^{-1})B(N+h),\\
 C_N=\exp\{N\psi'(N-h)/2\},\qquad
 A_N=\tfrac12N\psi'(N) C_N .
\end{gathered}
\]
Then $C_N<2$, $A_N\to\sqrt e/2$, and
\begin{multline*}
 \left| I''(N-1)+C_0/2\right|
 \le A_N T(N)+100 C_N\psi'(N)\\
       +\frac{340000}{N}e^{200}2^{-(N-h)} .
\end{multline*}
At $N=k^2$ the two moments are precisely those at degrees
$k(k-1)$ and $k(k+1)$. In particular,
$\liminf_{k\to\infty}T(k^2)<C_0/\sqrt e$ suffices to
obtain $I''(k^2-1)<0$ for an unbounded sequence of integers
$k$. Replacing the two moments by the finite sums
$\mathcal B_{k(k\pm1)}$ gives the same criterion.

Unconditionally one has
\[
 I''(x)<-0.00059\qquad
       (14{,}000\le x+1\le6.91899\cdot10^{12}).
\]
On the part with $x+1\ge10^{12}$,
$\left|I''(x)+C_0/2\right|<0.0000075$.
\end{proposition}
\begin{proof}
Let $w=w_{N-1,2}$, $\ell=\log(2u)-\psi(N)$,
$s_N=\ell^2-\psi'(N)$. The trigamma series gives
$1\le N\psi'(N)\le1+1/N\le2$ by comparison with its integral.
For real $z$, $z^2\le1+\cosh z$, since
$\cosh(z/2)\ge1+z^2/8\ge |z|/\sqrt2$.
Thus $|z^2-a|\le a\cosh z$ whenever $a\ge1$.
Taking $z=h\ell$, $a=N\psi'(N)$ yields
$|s_N|\le\psi'(N)\cosh(h\ell)$. Moreover
\[
\begin{gathered}
 w e^{\pm h\ell}=c_\pm w_{N-1\pm h,2},\\
 c_\pm=\frac{\Gamma(N\pm h)}{\Gamma(N)}e^{\mp h\psi(N)}\le C_N.
\end{gathered}
\]
Taylor's formula for $\log\Gamma$ and the decrease of $\psi'$
give the inequality. Indeed
$N\psi'(N-h)\le(1-1/h)^{-1}+(h-1)^{-2}\le21/16$ for
$h\ge5$, and $N\psi'(N-h)\to1$.
Multiplication of a gamma density by $u$ increases its shape
by one. It follows that
\[
 \int_0^\infty u P^2 w |s_N|\,\dd u\le A_N T(N).
\]
Also $M_{|s_N|}\le C_N\psi'(N)$.
For $u\ge100$, Proposition~\ref{pr:finitecurvature} bounds
$|(u-u_0)(D-C_0u-C_{\log})|$ by 100. For $0\le u\le100$
the global bound $D\le8(u+2)$ bounds this expression by 82000.
The extra term needed when replacing $|u-u_0|P^2$ by $uP^2$
is bounded by $1600\mathbf1_{[0,1]}$, since $0<u_0<1$ and
$P^2\le576e^u$. For either neighboring shape $N\pm h$,
its gamma probability below 100 is at most
$e^{200}2^{-(N-h)}$, by the exponential Markov inequality.
The displayed comparison with $e^{\pm h\ell}$ therefore gives
\[
 \int_0^{100} w|s_N|\,\dd u
            \le (4/N)e^{200}2^{-(N-h)} .
\]
Combining these bounds with $J_D''=C_0/2$ proves the first
assertion, since $4(82000+1600)<340000$.

For the finite range, the proof of
Proposition~\ref{pr:primewindowprefix} gives $B_m<9\cdot10^{-6}$
at every integer $250\le m\le n_*=6{,}919{,}000{,}000{,}000$,
using the verified zero height of Platt and Trudgian \cite{pt}.
To pass to real degrees $m\le y\le m+1$, H\"older's inequality
makes $\log(2^{-(y+1)}\Gamma(y+1)B(y))$ convex.
The error between the chord of $\log\Gamma(y+1)$ and its
value at $y$ is at most $\psi'(m+1)/8$ by the interpolation
remainder for a unit interval. Consequently
\[
 B(y)\le e^{\psi'(m+1)/8}\max\{B_m,B_{m+1}\}.
\]
Throughout the stated range both shifted degrees lie between
13880 and $n_*-1$. The last bound gives $T(N)<9.01\cdot10^{-6}$.
Using the decrease of $(1-1/h)^{-1}+(h-1)^{-2}$ gives
$C_N<1.657$ and $A_N<0.829$. The exponential error is below
$10^{-15}$, so the right side is at most
$0.829(9.01\cdot10^{-6})+165.7(1/N+1/N^2)+10^{-15}$.
At $N=14000$ and at $N=10^{12}$ this gives the two stated bounds.
The evaluations are reproduced
in \texttt{anc/prime\_curvature\_check.py}.
\end{proof}

\head{Summability of the prime moments at fixed orders}
Partial summation against $\vartheta(t)-t$, using $\int K=0$,
gives $|P(u)|\ll_A e^{u/2}(1+u)^{-A}$ for every fixed $A>0$
by the prime number theorem \cite{johnstonyang}. For each integer
$r\ge0$, Tonelli's theorem and the exponential series give
the exact identity
\[
 \sum_{m\ge r}m^{\underline r}\frac{B_m}{2^{m+1}}
       =\int_0^\infty u^r e^{-u}P(u)^2\,\dd u<\infty ,
\]
where $m^{\underline r}=m!/(m-r)!$. Summability at every fixed
order $r$ permits exponential growth of $B_m$. For example, the
gamma means of $e^{u/2}$ equal $(4/3)^{m+1}$ and satisfy all the
displayed summability bounds. For the prime weights, the unbounded
additive deficit of Proposition~\ref{pr:curvatureinduct} and the
two-moment bound along unbounded square degrees require stronger
arithmetic control.

\section{Linear prime sums and the signed error in the prime number theorem}\label{app:primelinear}
\subsection{A linear prime sum at a varying Mellin frequency}
\label{sec:movingprime}
The prime-return times $u=\log p$ determine the Dirichlet series
$\mathcal P(s)=\sum_p(\log p)p^{-s}$. In the bilateral transform
of $P(u)$, the substitution $v=pe^{-u}$ gives
$\widehat K(s)\mathcal P(s)$ for $\Re s>1$. The linear
coefficient below uses the prime weights before multiplication by
the detector response $\widehat K$. Disk transforms also occur in
Li's criterion and its arithmetic interpretation
\cite{bombierilagarias}. Varying the sampling frequency with the
degree separates poles of equal modulus.
For $n\ge2$, let $\tau_n=n^{-1/2}$, $Y_n=4096^n$, and put
\[
\begin{aligned}
 s_\tau(z)&=\tfrac12+i\tau+\frac{1+z}{1-z},\\
 g_\tau(z)&=(1-z)(s_\tau(z)-1)
          =\tfrac12+i\tau+(\tfrac32-i\tau)z .
\end{aligned}
\]
Define the finite arithmetic coefficient
\begin{equation}
 {\cal A}_n=[z^n]\,g_{\tau_n}(z)
       \sum_{p\le Y_n}(\log p)\,p^{-s_{\tau_n}(z)} .
 \label{eq:movingprime}
\end{equation}
Here $[z^n]$ takes the Taylor coefficient at zero. With
$\ell_j=L_j^{(-1)}$ the weight multiplying
$(\log p)p^{-3/2-i\tau_n}$ is
\[
 (\tfrac12+i\tau_n)\ell_n(2\log p)
       +(\tfrac32-i\tau_n)\ell_{n-1}(2\log p).
\]
This follows from $\sum_{j\ge0}\ell_j(v)z^j=\exp[-vz/(1-z)]$
\cite{Olver}. Replacing the infinite prime measure by $\dd t$
on $[1,\infty)$ gives $g_\tau(z)/(s_\tau(z)-1)=1-z$.
Thus its coefficient vanishes exactly for $n\ge2$.

\begin{theorem}[An unbounded sequence of finite linear prime sums]
\label{th:movingprime}
The Riemann hypothesis holds if and only if
$\lvert{\cal A}_n\rvert\le2$ for infinitely many integers $n$.
On RH, with $Q=2+\gamma_{\rm E}-\log(4\pi)$,
\begin{equation}
\begin{aligned}
 1-4Q&\le\liminf_{n\to\infty}|{\cal A}_n|
      \le\limsup_{n\to\infty}|{\cal A}_n|\\
     &\le1+4Q<1.185.
\end{aligned}
 \label{eq:movingband}
\end{equation}
If there is a zero with $\re(\rho)>1/2$, put
\[
 r_*=\min_{\re(\rho)>1/2}\left|\frac{\rho-3/2}{\rho+1/2}\right| .
\]
Then $0<r_*<1$ and $n^{-1}\log(1+|{\cal A}_n|)\to-\log r_*>0$.
The assertions hold for every $0<\tau_n\le1$ with
$\tau_n\to0$ and $n\tau_n\to\infty$.
\end{theorem}

\begin{proof}
Write $a_n(\tau)=[z^n]g_\tau(z)\mathcal P(s_\tau(z))$.
On $|z|=1/7$ one has $\Re s_\tau(z)\ge5/4$ and
$|g_\tau(z)|<2$ for $0\le\tau\le1$. Partial summation using
$\vartheta(t)\le2t$ gives
$\sum_{p>Y}(\log p)p^{-5/4}\le10Y^{-1/4}$.
Cauchy's estimate therefore gives the uniform truncation bound
\begin{equation}
 |{\cal A}_n-a_n(\tau_n)|\le20(7/8)^n .
 \label{eq:movingtruncation}
\end{equation}

Set $L(s)=-\zeta'(s)/\zeta(s)$. Separating the odd prime powers
in the Euler product gives
\[
\begin{gathered}
 \mathcal P(s)=L(s)-L(2s)-R_3(s),\\
 R_3(s)=\sum_p(\log p)\frac{p^{-3s}}{1-p^{-2s}} .
\end{gathered}
\]
The series for $R_3$ and its first derivative converge absolutely
and uniformly on $\Re s\ge1/2$. The function
$g_0(z)R_3(s_0(z))$ is bounded and analytic in the unit disk.
Its coefficients are square summable by Parseval on concentric
circles and hence tend to zero. Varying $\tau$ from zero to
$\tau_n$ changes the function in the disk supremum norm by
$O(\tau_n)$. The corresponding coefficient of $g_{\tau_n}R_3$
is consequently $o(1)$.

For a pole at $s=\lambda$ write
\[
\begin{gathered}
 w_\lambda(\tau)=\frac{\lambda-3/2-i\tau}{\lambda+1/2-i\tau},\\
 D_\lambda(\tau)=\frac{(\lambda-1)(1-w_\lambda(\tau))^3}
                                  {2w_\lambda(\tau)} .
\end{gathered}
\]
A simple fraction $b/(s-\lambda)$, after composition and
multiplication by $g_\tau$, contributes
$-bD_\lambda(\tau)w_\lambda(\tau)^{-n}$ for $n\ge2$.
Indeed, $s_\tau'(w)=2/(1-w)^2$, and removing that pole leaves
a polynomial of degree at most one. The pole of $L(s)$ at
$s=1$ contributes only $1-z$.

Suppose first that RH holds. Apply the genus-one Hadamard product
for the entire function $(s-1)\zeta(s)$ to $L(s)$ and $L(2s)$
\cite{Titchmarsh}. Its logarithmic derivative contributes the
simple fractions at the zeros, including the trivial zeros, with
their multiplicities. The regularizing constants multiply the
affine polynomial $g_\tau$ and drop out at degree $n\ge2$.
Termwise coefficient extraction follows from local uniform convergence
near $z=0$. For the zeros of either factor,
$|D_\lambda(\tau)|\ll(1+|\lambda|^2)^{-1}$ uniformly for small
$\tau$, and this bound is summable with multiplicities. Every
fixed pole strictly left of the line $\Re s=1/2$ has
$|w_\lambda(\tau_n)|>1$ tending to a number greater than one.
Dominated convergence disposes of all such poles, including those at
$s=\rho/2$ arising from the zeros of $\zeta(2s)$.
The pole of $-L(2s)$ at $s=1/2$ has residue $-1/2$
and supplies the explicit coefficient
\[
\begin{gathered}
 H_n=\frac{1}{(1+\tau_n^2)(1-i\tau_n)}
      \left(-\frac{1-i\tau_n}{1+i\tau_n}\right)^n,\\
 |H_n|=(1+\tau_n^2)^{-3/2}.
\end{gathered}
\]
Every remaining pole is a zero $\rho=1/2+i\gamma$ of $\zeta$
and has $|w_\rho(\tau)|=1$. At $\tau=0$ its weight satisfies
\[
 |D_\rho(0)|=\frac{4\sqrt{\gamma^2+1/4}}{(1+\gamma^2)^{3/2}}
       \le\frac{4}{\gamma^2+1/4}.
\]
The Hadamard product for $\xi$ and its functional equation give
$\sum_\rho[\rho(1-\rho)]^{-1}=2\xi'(1)/\xi(1)=Q$, counting
multiplicities. On RH this is a sum of positive real terms.
Thus $\limsup|{\cal A}_n-H_n|\le4Q$, which proves
\eqref{eq:movingband}.

Suppose now that a right-hand zero exists. For
$\rho=1/2+d+i\gamma$, $d>0$, its image has squared modulus
\[
 |w_\rho(\tau)|^2=1-\frac{4}{q_\rho(\tau)},\qquad
 q_\rho(\tau)=\frac{(1+d)^2+(\gamma-\tau)^2}{d}.
\]
These moduli tend to one as $|\gamma|\to\infty$, uniformly for
$0\le\tau\le1$, since $d<1/2$. No pole maps to the origin,
where $\Re s_\tau(0)=3/2$. Thus the minimum $r_*$ is attained
by a finite set $S$ of distinct zeros. For some fixed
$r_*<r_2<1$ and all sufficiently small $\tau$, precisely $S$
lies inside the circle of radius $r_2$ after transformation.
Subtracting these poles and applying Cauchy's estimate on that
circle gives, uniformly in $\tau$,
\[
 a_n(\tau)=\sum_{\rho\in S}m_\rho D_\rho(\tau)w_\rho(\tau)^{-n}
                  +O(r_2^{-n}),
\]
where $m_\rho$ is the multiplicity. The error is uniform by
compactness and continuity on the pole-free circle.

The common value $q_\rho(0)$ of the members of $S$ splits at
first order. Indeed, two equal slopes would give
$\gamma_1/d_1=\gamma_2/d_2=v$. If $d_1=d_2$ these are the
same zero. Otherwise equality of the constant coefficients gives
\[
 \frac1{d_1}+(1+v^2)d_1=\frac1{d_2}+(1+v^2)d_2,\qquad
 \gamma_1\gamma_2=1-d_1d_2<1 .
\]
Every zero with $\re(\rho)\ge1/2$ has $|\gamma|>1$. The
Hadamard identity and reflection in the critical line give
$\sum_\rho\Re(1/\rho)=Q/2<0.024$, with positive summands.
A zero with $\re(\rho)\ge1/2$ and $|\gamma|\le1$ would itself
contribute at least $1/4$. Since the two ordinates above have
the same sign, their product exceeds one, a contradiction.
Thus a unique member $\rho_*$ is smallest for all small
positive $\tau$, and every gap to another member is at least
$c\tau$ for some $c>0$. Since $r_*>0$, the derivative of
$\tfrac12\log(1-4/q)$ is finite and strictly positive at
$q=4/(1-r_*^2)$. The other terms are smaller in relative modulus
by $O(\exp[-c' n\tau_n])$, which tends to zero. In particular,
the square-root choice gives $O(\exp[-c'\sqrt n])$.
The leading amplitude tends to
$m_{\rho_*}D_{\rho_*}(0)\ne0$. This proves the asserted
exponential rate, also after the truncation \eqref{eq:movingtruncation}.
The functional equation supplies a right-hand zero whenever RH fails.
\end{proof}

\begin{corollary}[Faster frequency shifts and real prime sums]
\label{co:movingblocks}
Let $\tau_n=n^{-\alpha}$ with any fixed $\alpha>0$, or take
$\tau_n=0$ so that the sums are real. Put $h_n=\lceil\log(n+1)\rceil$.
Then RH holds if and only if
$\max_{0\le k<h_n}|{\cal A}_{n+k}|\le2$ for infinitely many $n$.
In particular the single-degree theorem applies to $\alpha=1/2$,
and the consecutive-degree criterion applies also to
$\alpha=3/2,5/2,7/2$.
\end{corollary}
\begin{proof}
On RH the proof of \eqref{eq:movingband} uses only $\tau_n\to0$.
Otherwise let $J=|S|$. At degrees $n+k$, $0\le k<J$, factor
$m_\rho D_\rho(\tau_n)w_\rho(\tau_n)^{-n}$ from the contribution
of each $\rho\in S$ in the residue formula. The resulting
$J$-by-$J$ matrix tends entrywise to
$(w_\rho(0)^{-k})_{k,\rho}$ because
$n|\tau_{n+k}-\tau_n|\to0$ for fixed $k$.
The limit is an invertible Vandermonde matrix since distinct
zeros have distinct images. Its approximants have uniformly bounded
inverses. The largest modulus of the factored amplitudes is
$\exp[(-\log r_*+o(1))n]$, whereas the remainder has size
$O(r_2^{-n})$. Hence some member of each such block exceeds two
for all large $n$, when $h_n\ge J$.
\end{proof}

The criterion has a direct form in terms of the prime counting error.
Put $f_n(t)=[z^n]g_{\tau_n}(z)t^{-s_{\tau_n}(z)}$ and
$E(t)=\vartheta(t)-t$. Abel summation gives, with $Y=Y_n$,
\begin{equation}
 \left|{\cal A}_n+\int_1^Y E(t)f_n'(t)\,\dd t\right|
       \le \tfrac{22}{7}(7/8)^n .
 \label{eq:movingabel}
\end{equation}
Indeed $f_n(1)=0$, and the left-hand expression inside the
modulus equals $E(Y)f_n(Y)+\int_1^Y f_n(t)\,\dd t$.
The last integral is $-[z^n](1-z)Y^{1-s_{\tau_n}(z)}$.
On $|z|=1/7$ one has $|1-z|\le8/7$, $|g_{\tau_n}|<2$,
and $|E(Y)|\le Y$. Cauchy's estimate gives the stated bound.
Hence a modulus at most two for this integral at infinitely many
degrees is also equivalent to RH.

One range of primes can be bounded unconditionally. Write
${\cal A}_n(p\le Z)$ for the contribution of these primes to
\eqref{eq:movingprime}. Cauchy's bound on circles tending to the unit
circle, where $\Re s_\tau>1/2$ in the open disk and $|g_\tau|<3$,
gives
\[
 |{\cal A}_n(p\le Z)|\le3\sum_{p\le Z}\frac{\log p}{\sqrt p}
                      \le12\sqrt Z
\]
by partial summation. In particular the primes through $e^{\sqrt n}$
contribute at most $12e^{\sqrt n/2}$ at every degree. An equivalent
condition for RH is therefore
\begin{equation}
 \liminf_{n\to\infty}\frac1n
 \log\bigl(1+|{\cal A}_n(e^{\sqrt n}<p\le4096^n)|\bigr)=0 .
 \label{eq:movingrange}
\end{equation}
Indeed, on RH the bounded full sum minus this initial range has
subexponential growth. With a right-hand zero the initial range
has strictly smaller exponential rate than the full sum. Moreover,
the verified height $T=3\cdot10^{12}$ \cite{pt} gives the
unconditional upper bound $-\tfrac12\log(1-2/T^2)<1.12\cdot10^{-25}$
for the exponential rate of the full sum, by the formula for
$|w_\rho(0)|^2$. An arithmetic proof of \eqref{eq:movingrange}
would give $I''(x)\to-C_0/2<0$ by
Theorem~\ref{th:realdegreecurvature}. A derivation from the retained
prime range alone requires further cancellation. At degree two a direct sum over the
$1\,077\,871$ primes through $4096^2$ gives
$0.64895459931072<|{\cal A}_2|<0.64895459931074$.
The interval calculation and independent Taylor checks of the
weights and residues are in \texttt{anc/prime\_linear\_check.py}.

\head{Point-count weights in the linear coefficient}
Let $w_{p,r}=N_{p,r}/(p^r+1)$ for odd $p$ and $w_{2,r}=1$ be
the weights of Sec.~\ref{sec:stepprimes} for the curve $y^2=x^3-x$.
Define
\[
 {\cal A}_{n,E,r}=[z^n]\,g_{\tau_n}(z)
       \sum_{p\le Y_n}w_{p,r}\,(\log p)\,p^{-s_{\tau_n}(z)} .
\]
\begin{proposition}[The linear coefficient with point-count weights]
\label{pr:movingweights}
For every fixed integer $r\ge1$,
${\cal A}_{n,E,r}-{\cal A}_n\to0$ as $n\to\infty$. Consequently
RH holds if and only if $|{\cal A}_{n,E,r}|\le2$ for infinitely many
$n$. For $r\ge2$ the convergence is uniform over all weights with
$|w_p-1|\le2/p$.
\end{proposition}
\begin{proof}
For $|z|<1$ and $0\le\tau\le1$ one has $\Re s_\tau(z)>1/2$ and
$|g_\tau(z)|<3$. Cauchy's bound on circles of radius tending to one
gives $|[z^n]g_\tau(z)p^{-s_\tau(z)}|\le3p^{-1/2}$. Moreover
$g_\tau(z)p^{-s_\tau(z)}=p^{-i\tau}[g_0(z)+i\tau(1-z)]p^{-s_0(z)}$.
The coefficients of the bounded analytic functions
$g_0p^{-s_0}$ and $(1-z)p^{-s_0}$ are square summable. Hence the
coefficient at $\tau_n$ tends to zero for each fixed $p$.
For $r\ge2$, $|w_{p,r}-1|\le2p^{-r/2}\le2/p$. Dominated convergence
with the majorant $6(\log p)p^{-3/2}$ proves the first assertion,
uniformly in the weights and with the moving truncation.

Let $r=1$, so that $1-w_{p,1}=a_p/(p+1)$. Put
$\mathcal D_E(s)=\sum_{p>2}(\log p)a_p\,p^{-s}/(p+1)$ for
$\Re s>1/2$. The $L$-function
$L_E(w)=\prod_{p>2}(1-a_pp^{-w}+p^{1-2w})^{-1}$ is that of a Hecke
character of $\mathbb Q(i)$ with conductor $(2+2i)$. It is entire
of order one, satisfies a functional equation under $w\mapsto2-w$
with gamma factor $\Gamma(w)$, and does not vanish for
$\Re w\ge3/2$ \cite{irelandrosen,hecke1920,iwanieckowalski}. Then
\[
 \mathcal D_E(s)=-\frac{L_E'(s+1)}{L_E(s+1)}+R_E(s),
\]
where $R_E$ and $R_E'$ are bounded on $\Re s\ge1/2$. Indeed
$R_E$ consists of the difference between $1/(p+1)$ and $1/p$
and of the Euler terms of degree $j\ge2$, with the summable
majorant $4\sum_p(\log p)\sum_{j\ge2}p^{-j}$ and its analogue
after differentiation. As for $R_3$ in the proof of
Theorem~\ref{th:movingprime}, the coefficient of
$g_{\tau_n}R_E(s_{\tau_n})$ tends to zero. Every zero $\lambda$ of
$L_E(s+1)$, including the trivial zeros $-1,-2,\ldots$, has
$\Re\lambda<1/2$, and $\sum_\lambda(1+|\lambda|^2)^{-1}<\infty$.
A simple fraction $1/(s-\lambda)$ contributes
$-D_\lambda(\tau)w_\lambda(\tau)^{-n}$ in the notation of that
proof. Explicitly
\[
 -D_\lambda(\tau)w_\lambda(\tau)^{-n}
 =\frac{4(\lambda-1)}{(\tfrac32+i\tau-\lambda)^3}
  \left(\frac{\lambda+\tfrac12-i\tau}{\lambda-\tfrac32-i\tau}\right)^{n-2}.
\]
The last ratio has modulus less than one because $\Re\lambda<1/2$,
and its limit at $\tau=0$ has modulus less than one. The prefactor
is $O((1+|\lambda|^2)^{-1})$ uniformly for $0\le\tau\le1$. The
regularizing constants of the genus-one Hadamard product multiply the
affine function $g_\tau$ and drop out at degree $n\ge2$. Dominated
convergence gives $[z^n]g_{\tau_n}\mathcal D_E(s_{\tau_n})\to0$.
Finally on $|z|=1/7$ the tail of $\mathcal D_E$ beyond $Y$ is
$O(Y^{-3/4})$, since $|a_p|\le2\sqrt p$ and $\Re s\ge5/4$.
Cauchy's estimate bounds the truncation error by
$O(7^nY_n^{-3/4})=O((7/512)^n)$. The equivalence follows from
Theorem~\ref{th:movingprime}. On RH the upper limit of
$|{\cal A}_n|$ is less than $1.185$, and otherwise
$|{\cal A}_n|\to\infty$.
\end{proof}
The convergence is not uniform at small degrees. At $n=2$ the direct
sums give ${\cal A}_{2,E,1}-{\cal A}_2=-0.1469687\ldots+0.2172274\ldots i$
and ${\cal A}_{2,E,2}-{\cal A}_2=0.2518146\ldots+0.6685856\ldots i$.
The small primes dominate: $w_{3,2}=8/5$. The comparison at large
$n$ is unconditional. The full coefficient on either side still
contains the ordinary prime sum.

\subsection{The Abel kernel on power errors and a divisor decomposition}
\label{sec:powerbilinear}
Let $Y=4096^n$, $\tau=n^{-1/2}$, $s=s_\tau(z)$,
$g=g_\tau(z)$, and $f_n(t)=[z^n]g\,t^{-s}$, as in
Eq.~\eqref{eq:movingabel}. The response of the Abel kernel to a single
power quantifies the contribution permitted by a pointwise bound on
the prime error. Estimating the coefficient more sharply requires
control of its phase in the integral.

\begin{proposition}[Response to a power error]\label{pr:powerprofile}
For $1/2<\beta<1$ and integer $n\ge2$,
\[
 \lim_{n\to\infty}
 \left|\int_1^{4096^n}t^\beta f_n'(t)\,\dd t\right|^{1/n}
     =\frac{\beta+1/2}{3/2-\beta}.
\]
For $\beta=1$ the integral on the infinite interval is zero
at every such degree and the finite-interval integral is
$O((7/8)^n)$.
\end{proposition}
\begin{proof}
On $|z|=1/7$, $\Re s\ge5/4$, $|s|<3$, and $|g|<2$.
Differentiation, integration and coefficient extraction commute there
by absolute convergence. For $1/2<\beta\le1$,
\[
\begin{gathered}
 \int_1^\infty t^\beta f_n'(t)\,\dd t
   =[z^n]\frac{-gs}{s-\beta},\\
 \left|\int_Y^\infty t^\beta f_n'(t)\,\dd t\right|
   \le 24\,7^n Y^{\beta-5/4}\le24(7/8)^n .
\end{gathered}
\]
At $\beta=1$ the rational function reduces to $-(1-z)s$,
which is affine. For $\beta<1$ the term $-g$ is affine and
the remaining fraction is $-\beta g/(s-\beta)$.
The pole formula in the proof of Theorem~\ref{th:movingprime}
gives the coefficient $\beta D_\beta(\tau)w_\beta(\tau)^{-n}$.
As $\tau\to0$, $D_\beta(\tau)$ tends to a nonzero limit and
$w_\beta(0)=(\beta-3/2)/(\beta+1/2)$. Its inverse modulus
exceeds one. Taking $n$-th roots proves the assertion.
\end{proof}
For $\beta=3/4$ the rate is $5/3$. More generally, write
$B(t)=t(\log t)^{1.801}\exp[-0.1853(\log t)^{3/5}
(\log\log t)^{-1/5}]$ for $t\ge23$.
Theorem~1.4 of Johnston and Yang gives
$|\vartheta(t)-t|\le0.027 B(t)$ on this range
\cite{johnstonyang}. For fixed $\beta<1$, $t^\beta/B(t)\to0$.
Thus this size envelope alone allows a test error whose integral
grows at the rate just computed. Cutting off a fixed initial interval
does not change that rate. Indeed, on every $|z|=r<1$ the transform
of that interval is bounded uniformly for $0\le\tau\le1$. Cauchy's
bound and $r\uparrow1$ give upper $n$-th-root limit at most one.
This is a test of the pointwise envelope. It makes no assertion
that a power error is realized by the primes.

\begin{proposition}[Signed divisor sums for the actual weight]
\label{pr:bilinearweights}
Let $S_{\Lambda,n}=\sum_{k\le Y}\Lambda(k)f_n(k)$.
Unconditionally,
\begin{equation}
\begin{aligned}
 |S_{\Lambda,n}-\mathcal A_n|
 &\le 3\log Y+6+3\sum_{j\ge2}\frac{\log j}{j(\sqrt j-1)}\\
 &=O(n+1).
\end{aligned}
 \label{eq:properpowerslinear}
\end{equation}
For integers $U,V\ge1$, set
$\mu_{\le U}(d)=\mu(d)\mathbf1_{d\le U}$ and
$\Lambda_{\le V}(d)=\Lambda(d)\mathbf1_{d\le V}$, with
$>U,>V$ denoting the complementary functions.
If $*$ denotes Dirichlet convolution and $\mathbf1(k)=1$, then
\[
 \Lambda=\Lambda_{\le V}+\mu_{\le U}*\log
       -\mu_{\le U}*\mathbf1*\Lambda_{\le V}
       +\mu_{>U}*\mathbf1*\Lambda_{>V}.
\]
In particular, the last contribution to $S_{\Lambda,n}$ is exactly
\begin{equation}
 \sum_{\substack{d>U,\ m>V\\dm\le Y}}\mu(d)
   \left(\sum_{\substack{e\mid m\\e>V}}\Lambda(e)\right)f_n(dm).
 \label{eq:signedtypeII}
\end{equation}
\end{proposition}
\begin{proof}
For $0\le\tau\le1$ the unit disk has $\Re s>1/2$ and
$|g|<3$. Cauchy's coefficient bound as the radius tends to one
gives $|f_n(k)|\le3/\sqrt k$ for each fixed $k$.
For prime squares use partial summation and
$\vartheta(t)\le2t$ to get
$\sum_{p\le y}(\log p)/p\le 2+2\log y$ for $y\ge1$.
The contribution of squares is at most $6+3\log Y$.
For higher powers summation of
the geometric series at each prime is bounded by the convergent
series over $j$ in \eqref{eq:properpowerslinear}.
The convolution identity follows by substituting
$\log=\mathbf1*\Lambda$ and $\mu*\mathbf1=\delta_1$ on its
right. Grouping the two remaining factors in the last convolution
gives \eqref{eq:signedtypeII}. All sums at a fixed $Y$ are finite.
\end{proof}
The decomposition is Vaughan's identity
\cite{vaughanprime,heathbrownprime}. The small-divisor terms yield
Type~I sums and \eqref{eq:signedtypeII} is a Type~II sum.
At the scale $n^{-1}\log(1+|\cdot|)$, the deterministic
$O(n+1)$ discrepancy in \eqref{eq:properpowerslinear} tends
to zero in either direction of comparison. At the fixed threshold
$|\mathcal A_n|\le2$ the proper-power correction must be retained.
The identity permits one to test cancellation before any triangle
inequality is applied. A bound for the full signed combination
with this weight remains the arithmetic input to
Theorem~\ref{th:movingprime}.

\section{Divisor recursions and Fourier bounds at the ends of a window}\label{app:fourierwindows}
\subsection{A divisor identity using only earlier prime powers}
\label{sec:factorialdescent}
The integer in Proposition~\ref{th:cycunit} yields a recursion
whose remainder can be estimated separately. Chebyshev's factorial
construction uses the divisor identity
$\log n=\sum_{b\mid n}\Lambda(b)$ \cite{diamonderdos}. Define
\[
 L(x)=\sum_{n\ge1}(\log n)\,k(n/x)=\sum_{j\ge1} V(x/j).
\]
Both sums are finite. The first can be computed from ordinary
integers alone. Partial summation expresses it in terms of
$\log(\lfloor u\rfloor!)$. Fix $w=210=2\cdot3\cdot5\cdot7$,
$D=2310$, and put
\[
\begin{gathered}
 a_j=\mathbf 1_{(j,w)=1}(1-11\mathbf 1_{11\mid j}),\\
 B(x)=\sum_{d\mid w}\mu(d)\{L(x/d)-11L(x/(11d))\}.
\end{gathered}
\]
This choice removes the terms $2,\ldots,10$. Since $11$ is
coprime to $w$, ordinary finite inclusion--exclusion gives
$B(x)=\sum_{j\ge1}a_j V(x/j)$. The $a_j$ are even as a
function of the residue class modulo $D$, their average is zero,
and $a_D=0$. In particular
$\sum_{r=1}^D a_r=\sum_{r=1}^D r a_r=0$.

\begin{theorem}\label{th:cycdescent}
With $M=D\lceil x^{1/3}/D\rceil$ and $A=6\log2$, for every
sufficiently large real $x$ one has, unconditionally,
\begin{equation}
 E(x)=-\sum_{11\le j\le M} a_j E(x/j)+O(x^{1/3}).
 \label{eq:cycdescent}
\end{equation}
All prime powers in this finite sum are at most $8x/11<x$.
The constant is independent of $x$.
\end{theorem}
\begin{proof}
Write $b=\int_1^8 k(v)\log v\,\dd v$. The trapezoidal rule on
the full integer grid gives
\[
 L(y)=Ay\log y+by+O\big(\log(2y)/y\big)\qquad(y\ge2).
\]
Indeed the compactly supported continuous function
$f(u)=\log(u)k(u/y)$, extended by zero, has derivative of total
variation $O(\log(2y)/y)$. This follows on the three smooth
pieces by integration of $|f''|$, including the four jumps of
$f'$. The integral error on each unit cell is bounded by the
variation of $f'$ there, by linear interpolation. The grid sum
equals the composite trapezoidal sum. For the fixed $w$ it follows
that $B(x)=Ax(\varphi(w)/w)\log11+O(\log(2x)/x)$.
Also $\sum_{j\ge1}a_j/j=(\varphi(w)/w)\log11$, by taking
$s\downarrow1$ in
\[
 \sum_{j\ge1}a_j j^{-s}=(1-11^{1-s})
             \prod_{p\mid w}(1-p^{-s})\,\zeta(s).
\]
Bounded partial sums ensure convergence at 1 and justify this
limit. Hence $\sum_{j\ge1}a_j E(x/j)=O(\log(2x)/x)$.

For the truncation put $f_x(u)=E(x/u)$. Chebyshev's estimate gives
$\sup|E'_+|\ll1$. The jumps in $E'$ occur at $n/8,n/4,n/2,n$,
each of magnitude at most $8\Lambda(n)/n$ before coincident
jumps are added. Consequently, for all $Y>0$,
\[
 \int_{(0,Y]} v^2\,|dE'(v)|
   \ll \sum_{n\le8Y}n\Lambda(n)\ll Y^2 .
\]
Since $f'_x(u)=-xE'(x/u)/u^2$ away from jumps, the variation of
$f'_x$ on $[M,\infty)$ is $O(x/M^2)$. The variation of its
smooth factor contributes $O(x/M^2)$ and the weighted jumps
contribute $O((x/M)^2/x)$ by the last display. On each full
period $j=N+r,\ 1\le r\le D$, the constant and linear terms in
the expansion about $N$ sum to zero. The remainder of the
weighted sum on that period is bounded by
$\sum_{r=1}^D r|a_r|$ times the variation of $f'_x$ there.
Summing over periods gives
$\left|\sum_{j>M}a_jE(x/j)\right|\ll_D x/M^2$, which proves
\eqref{eq:cycdescent}.
\end{proof}

Together with \eqref{eq:cycround}, this reduces an upper bound in
\eqref{eq:cyctarget} to the signed sum in \eqref{eq:cycdescent}
over earlier scales. A termwise estimate
$|E(y)|\le C\sqrt y$ would give the normalized multiplier
\[
 \sum_{11\le j\le M}\frac{|a_j|}{\sqrt j}
       =\frac{64}{77}\sqrt M+O_D(1).
\]
Indeed $|a_j|$ has mean
$(48/210)(10/11+10/11)=32/77$.
The loss is in correlations among those terms. In the critical
strip the Dirichlet series of the full arithmetic recursion
vanishes exactly at the zeros of $\zeta$, since its other factors
have no zeros there. In particular a hypothetical mode $x^\rho$
satisfies the homogeneous identity
$\sum_{j\ge1}a_j(x/j)^\rho=0$. The independent $O(x^{1/3})$
bound controls the truncation. A proof by this recursion requires a uniform
upper bound for the retained sum on the windows of
Theorem~\ref{th:cycwindow}.

\subsection{Fourier analysis on the interval of length \texorpdfstring{$X^{3/4}$}{X to the power 3/4}}
\label{sec:fourierwindow}
Write
\[
\begin{gathered}
 D(x)=\frac{Ax-V(x)}{\sqrt x},\\
 D_{\mathrm{cyc}}(x)=
 \frac{\log|\mathcal I(x)|-\log|\mathcal A(x)|+mAx}{m\sqrt x}.
\end{gathered}
\]
Here $A=6\log 2$, $k$ is given by \eqref{eq:stepweight},
and $0\le D_{\mathrm{cyc}}(x)-D(x)\le(\log2)/\sqrt x$
by \eqref{eq:cycround}. The exact expansion for $x>1$ is
\begin{equation}
\begin{gathered}
 D(x)=\sum_\rho K(\rho)x^{\rho-1/2}+r(x),\\
 0\le r(x)=\sum_{j\ge1}K(-2j)x^{-2j-1/2}
          \le \frac{x^{-5/2}}{3(1-x^{-2})}.
\end{gathered}
 \label{eq:fourierdeficit}
\end{equation}
It follows from \eqref{eq:explicit} and \eqref{eq:steptransfer}.
Both series are absolutely convergent. The sum includes conjugate zeros.

\begin{theorem}\label{th:fourierendpoints}
Let $h=X^{3/4}$, $w=(\log X/\log\log X)^{1/3}$ and
$T_X=X^{1/4}\exp(-w/2000)$. Unconditionally, as $X\to\infty$,
\begin{gather}
 \sup_{X\le t\le X+h}\ 
 \sum_{|\im\rho|>T_X}|K(\rho)|t^{\re(\rho)-1/2}
 \ll e^{-w/500},\label{eq:fouriertail}\\
 \begin{gathered}
 \sup_{0\le v\le1}
 |D(X+vh)-(1-v)D(X)-vD(X+h)|\\
 \ll e^{-w/500}.
 \end{gathered}\label{eq:fourierchord}
\end{gather}
In particular the minimum of $D$ on the interval differs from the
smaller endpoint value by $O(e^{-w/500})$. The same estimate for
the infimum holds for $D_{\mathrm{cyc}}$.
\end{theorem}
\begin{proof}
Put $L=\log X$. Write $N(\sigma,U)$ for the number of zeros with
$\beta\ge\sigma,\ 0<\gamma\le U$, counted with multiplicity.
Ingham's density theorem and the Vinogradov--Korobov zero-free
region \cite{ingham40,mtyregion} give the bound
\[
 N(\sigma,U)\ll U^{3(1-\sigma)}(\log U)^5,\qquad
 \tfrac12\le\sigma\le1,\quad U\ge3 ,
\]
and $\beta\le1-\lambda$ for $0<|\gamma|\le2X$, where
$\lambda=w/(60L)$ for all sufficiently large $X$. Indeed the
published zero-free constant is $55.241$ and replacing $X$ by
$2X$ in its logarithms changes their product by $1+o(1)$.
The finitely many small ordinates are absorbed by increasing $X$.

For a dyadic block $U<\gamma\le2U$, $T_X\le U<X$, the contribution
to \eqref{eq:fouriertail} is bounded, up to an absolute constant, by
\[
 \frac{\log(2U)}{U}
   +L^6\int_{1/2}^{1-\lambda}X^{s-1/2}U^{-2+3(1-s)}\,\dd s .
\]
Here $|K(\rho)|\ll\gamma^{-2}$ and $X\le t\le2X$. Count the
baseline weight at $s=1/2$ by $N(2U)\ll U\log(2U)$, then
integrate the derivative of $X^{s-1/2}$ against the tail count.
The exponent of $U$ is negative throughout the integral, so the
integrand is at most $e^{w/1000}X^{-(1-s)/4}$.
There are $O(L)$ blocks. Their total is
$O(L\log(2T_X)/T_X+L^7e^{w/1000-w/240})$.
For ordinates above $X$, direct summation using $N(U)\ll U\log U$
costs $O(L/\sqrt X)$. These estimates imply \eqref{eq:fouriertail}.

For $|\gamma|\le T_X$, interpolate $K(\rho)t^{\rho-1/2}$
linearly between $X$ and $X+h$. The elementary interpolation
remainder is bounded by $h^2/8$ times the supremum of the second
derivative in modulus. Since
$|K(\rho)(\rho-1/2)(\rho-3/2)|\ll1$, this is $O(X^{\beta-1})$
per zero. Another integration of the tail count gives
\[
\begin{aligned}
 \sum_{|\gamma|\le T_X}X^{\beta-1}
 &\ll X^{-1/2}T_X\log(2T_X)\\
 &\quad+L^6\int_{1/2}^{1-\lambda}X^{-(1-s)/4}\,\dd s\\
 &\ll e^{-w/500}.
\end{aligned}
\]
The contribution of the other terms, both in the interior and at
the endpoints, follows from \eqref{eq:fouriertail} and
\eqref{eq:fourierdeficit}. This proves \eqref{eq:fourierchord}.
\end{proof}

The corresponding integral on the prime side can also be evaluated
locally. Put $(d_1,d_2,d_4,d_8)=(2,-2,-8,8)$. For $0\le v\le h$,
\begin{multline}
 V(X+v)-V(X)-vV'_+(X)\\
 =\int_0^v\ \sum_{\ell\in\{1,2,4,8\}}d_\ell
     \sum_{\ell X<n\le\ell(X+s)}\frac{\Lambda(n)}{n}\,\dd s .
 \label{eq:fourierprime}
\end{multline}
This follows by integrating the four jumps in the slope of $k(n/x)$.
The prime number theorem in short intervals \cite{vaughanshort}
gives each inner sum as $s/X+o(s/X)$, uniformly for
$X^{7/10}\le s\le h$. Replacing $1/n$ by $1/(\ell X)$ costs
$O((s/X)^2)$ there. The four main terms cancel since
$\sum_\ell d_\ell=0$. On shorter intervals monotonicity bounds
each inner sum by $O(X^{-3/10})$. The integral is therefore
$o(\sqrt X)$, uniformly in $v$. This calculation controls the
change of slope. Its two integration constants are $V(X)$ and
$V'_+(X)$.

\head{Three Fourier moments from prime powers}
Let $H=3\cdot10^{12}$ and $e_H<1.48\cdot10^{-12}$ be the
positive-ordinate tail bound of Appendix~\ref{app:check}
\cite{pt,bpt}. For $|\gamma|\le H$ the zeros are on the critical line.
Writing $P(s)=(4^s-1)(2^s-2)$, $\theta=\gamma\log2$, one obtains
\[
 |P(1/2+i\gamma)|^2
 =30-12\sqrt2\cos\theta-24\cos(2\theta)+8\sqrt2\cos(3\theta).
\]
Set $C_j=C(j\log2)$ as in \eqref{eq:C}, and
$Q=2+\eul-\log(4\pi)$. The contribution of unverified zeros to
$C(\log q)=\sum_\rho c_\rho q^{\rho-1/2}$ has modulus at most
$2\sqrt q\,e_H$ for $q\ge1$. The verified weights $c_\rho$
are positive and sum to at most $Q$. Cauchy--Schwarz therefore gives
\begin{multline}
 \sum_{|\gamma|\le H}|K(\rho)|\\
 \le \{Q(30Q-12\sqrt2 C_1-24C_2+8\sqrt2 C_3+208e_H)\}^{1/2}\\
 <0.280827 .\label{eq:fourierlow}
\end{multline}
The constants $C_j$ require prime powers only up to $8$.
In particular $T(8)=4\log2+\tfrac53\log3+\tfrac35\log5+\tfrac17\log7$.
The bound on the unverified mass uses the zero count with multiplicity.

For every zero, $|K(\rho)|\le16/\gamma^2$. Indeed, with
$z=2^\rho$ the maximum of $|(z^2-1)(z-2)|$ on $|z|\le2$
is attained on the boundary $|z|=2$.
Its square there is $8(25-16v^2)(1-v)$ with $v=\cos(\arg z)$.
Its maximum is $(3344+364\sqrt{91})/27<256$.
Divide by $|\rho(1-\rho)|\ge\gamma^2$.
Consequently all zeros with $\beta\le1/2$ together contribute
at most $0.280827+32e_H<0.280828$ in absolute value to
\eqref{eq:fourierdeficit}, for every $x>1$. On RH this improves
the constant in the eventual bound \eqref{eq:cyctarget} to $0.281$.

\head{A positive integral with a two-component transfer matrix}
Let $\sigma_1,\sigma_3$ be the Pauli matrices and $c,\eta>0$.
Define
\begin{equation}
\begin{aligned}
 M(z)&=\frac{\operatorname{Tr}\exp(c\sigma_1+\eta z\sigma_3)}
                  {2\cosh c}
     =\frac{\cosh\sqrt{c^2+(\eta z)^2}}{\cosh c}\\
     &=\int_{-\eta}^{\eta} e^{zu}\,\dd\mu(u).
\end{aligned}
 \label{eq:fourierfilter}
\end{equation}
The square-root expression is entire, by its even power series.
The probability $\mu$ has an atom of mass $1/(2\cosh c)$ at each
endpoint. In the interior, with $v=u/\eta$, its density is
\[
 \frac{c}{2\eta\cosh c}\,
    \frac{I_1(c\sqrt{1-v^2})}{\sqrt{1-v^2}} .
\]
For a proof average $\exp(r(c n_1+\eta z n_3))$ over the unit
sphere. Rotational symmetry for real $z$, followed by analytic
continuation, gives
\[
 \frac12\int_{-r}^{r}I_0(c\sqrt{r^2-v^2})e^{\eta zv}\,\dd v
 =\frac{\sinh(r\sqrt{c^2+(\eta z)^2})}{\sqrt{c^2+(\eta z)^2}} .
\]
Differentiate at $r=1$ and divide by $\cosh c$. Positivity of
$I_1$ and $M(0)=1$ give a probability. In the basis of
$\sigma_3$, the two signs represent opposite shifts of the
logarithmic clock. Equivalently, the Lie product formula expands
the trace in closed paths through these two states with nonnegative
weights, since every entry of $e^{c\sigma_1/N}$ is positive.
This gives a two-channel time average of the prime signal. Its
transfer matrix is specified independently of the Majorana
Hamiltonian. For comparison, the compact Logan filter used by
B\"uthe has a $\sinh(z)/z$ multiplier \cite{buthefilter}.

\begin{theorem}\label{th:fourierfinite}
For every real $64\le X\le8\cdot10^{42}$ there exists a real
$t\in[X,X+X^{3/4}]$ with $D_{\mathrm{cyc}}(t)<0.357$.
This uses the published verification to height $H$ and an
unconditional upper bound for all higher zeros.
\end{theorem}
\begin{proof}
Take $c=28.2,\ \eta=c/H,\ x_0=Xe^\eta$. At the upper endpoint of
the stated range $(H/2)\log(1+X^{-1/4})>28.2045$.
The left side decreases with $X$, so
$[x_0e^{-\eta},x_0e^\eta]\subseteq[X,X+X^{3/4}]$ on the full range.
The integral of $D(x_0e^u)$ has the prime-power expansion
\begin{multline*}
 A\sqrt{x_0}\,M(1/2)\\
 -\frac1{\sqrt{x_0}}\sum_n\Lambda(n)
    \int e^{-u/2}k\big(n/(x_0e^u)\big)\,\dd\mu(u).
\end{multline*}
Equivalently integrate \eqref{eq:fourierdeficit} term by term.
For the verified part $|M(i\gamma)|\le1$, giving
\eqref{eq:fourierlow}. For an unverified zero put $a=\eta(\beta-1/2)$,
$b=\eta\gamma$, so $|b|\ge c$. For $W^2=c^2+(a+ib)^2$,
\[
\begin{gathered}
 (\re W)^2=a^2+\tfrac12(\sqrt{d^2+4a^2c^2}-d)
 \le a^2+|a|c,\\
 d=b^2+a^2-c^2\ge0.
\end{gathered}
\]
Hence $|\re W|<10^{-4}$ and
$|M(\rho-1/2)|\le\cosh(10^{-4})/\cosh c$.
Pair $\beta+i\gamma$ with $1-\beta+i\gamma$ and also include
their conjugates. Since
\[
 x_0^{\beta-1/2}+x_0^{1/2-\beta}\le\sqrt{x_0}+1/\sqrt{x_0},
\]
the entire unverified contribution is at most
\[
 b_H(\sqrt{x_0}+1/\sqrt{x_0}),\qquad
 b_H=16e_H\,\frac{\cosh(10^{-4})}{\cosh c}.
\]
This pairing bound includes zeros on the line by symmetrization.
There is a point where $D$ is at most its average. Write $y=\sqrt X$.
Including the group rounding at that point bounds $D_{\mathrm{cyc}}$ by
\[
 0.280827+b_H(e^{\eta/2}y+e^{-\eta/2}/y)
     +\frac{\log2}{y}+\frac{X^{-5/2}}{3(1-X^{-2})}.
\]
For $65536\le X\le8\cdot10^{42}$, bound the last term at $65536$.
The remaining expression is convex in $y>0$ and its upper bounds at
the two endpoints are $0.28354$ and $0.35666$, even after adding
the trivial-zero term. On $64\le X\le65536$ the direct prime-power
bound following \eqref{eq:cycpnt} already gives $0.211$ at $t=X$.
All decimal bounds are enclosed outwards in \texttt{anc/verify.py}.
\end{proof}
Without the filter, the same symmetry and \eqref{eq:fourierlow} give
at every $t>1$, with $y=\sqrt t$,
\begin{multline*}
 D_{\mathrm{cyc}}(t)<
 0.280827+16e_H(y+1/y)+(\log 2)/y\\
             +\frac{1}{3y(t^2-1)} .
\end{multline*}
On $64\le t\le1.4\cdot10^{19}$ the last term can be bounded at
$64$. The convex remainder evaluated at the two endpoints gives
total upper bounds $0.367481$ and $0.369440$, respectively.
Thus the threshold $0.370$ holds pointwise throughout this range.

\begin{theorem}\label{th:fixedfourierwindow}
There are fixed constants $R>1$ and $b>0$ such that every
sufficiently large interval $[X,RX]$ contains three points
$t_+,t_-,t_0$ satisfying
\[
 D(t_+)>b,\qquad D(t_-)<-b,\qquad V(t_0)=At_0 .
\]
In particular $D_{\mathrm{cyc}}(t_-)<-b/2$ eventually.
If RH fails, the first two inequalities hold with right sides
$b\,t_+^{\beta_*-1/2}$ and $-b\,t_-^{\beta_*-1/2}$,
respectively, for some fixed $\beta_*>1/2$.
\end{theorem}
\begin{proof}
First suppose an off-line zero exists to the right of $1/2$.
Choose $\sigma>1/2$ below its real part and distinct from the real
parts of all zeros. The full exponent in Ingham's theorem is
$a_\sigma=3(1-\sigma)/(2-\sigma)<1$ \cite{ingham40}. Thus
\[
 \sum_{\beta\ge\sigma,\ \gamma>0}\frac1\gamma<\infty ,
\]
by partial summation of $N(\sigma,U)\ll_\sigma U^{a_\sigma}\log^5 U$.
Let $\gamma_*$ be the least positive ordinate in this set, and
$\beta_*>\sigma$ the largest real part at that ordinate.
For each conjugate pair in the set with $\gamma>\gamma_*$, take
the probability on $0,d_\gamma=\pi/\gamma$ with masses proportional
to $e^{(\beta-\beta_*)d_\gamma},1$. Its Laplace multiplier vanishes
at $\rho-\beta_*$ and its conjugate. The infinite convolution
$\nu$ of these measures has support in $[0,L_*]$ with
$L_*=\sum d_\gamma<\infty$. Indeed the partial random sums of
independent variables at $0,d_\gamma$ converge uniformly in the
choices of their values. Their transforms therefore converge at
every fixed complex argument.

At $i\gamma_*$ the projected length of the chord joining the two
phases gives the lower bound
\[
 |\widehat\nu(i\gamma_*)|
 \ge\prod_{\beta\ge\sigma,\ \gamma>\gamma_*}
       \cos\big(\pi\gamma_*/(2\gamma)\big)>0 .
\]
Finitely many initial factors are positive, and the logarithms of
the remaining factors are $O_{\gamma_*}(\gamma^{-2})$.
Average $D(t)/t^{\beta_*-1/2}$ at $t=Xe^{v+u}$ over $u\sim\nu$.
The zeros with $\beta<\sigma$ contribute $O(X^{\sigma-\beta_*})$
by absolute summability. The finitely many terms with ordinate
$\pm\gamma_*$ and smaller real part contribute $o(1)$.
Termwise integration of all terms is valid on the fixed compact
logarithmic interval. The mean is consequently
\[
 2m_*\re\!\left(K(\beta_*+i\gamma_*)\,\widehat\nu(i\gamma_*)
               e^{i\gamma_*(\log X+v)}\right)+o(1)
\]
uniformly for $0\le v\le2\pi/\gamma_*$. The amplitude is positive,
so the maximum and minimum over a period have opposite signs.
Positivity of $\nu$ yields both bounds in the theorem with
$R=\exp(L_*+2\pi/\gamma_*)$. Continuity gives $t_0$.

If RH holds, take $\beta_*=1/2$ and the least positive ordinate
$\gamma_*$. Convolve just the two-point probabilities for
$\gamma_*<\gamma\le T$, now with equal masses. The same chord
bound, with the product extended over all $\gamma>\gamma_*$,
is a positive lower bound independent of $T$.
Fix $T$ so large that $\sum_{|\gamma|>T}|K(\rho)|$ is less than
one quarter of the resulting lower bound on the sinusoid's amplitude.
Those terms have multipliers of modulus at most one. Sliding through
one period again supplies both signs, now using only a finite
convolution of fixed support length. These two cases prove existence
with a fixed $R$ whose numerical value is left undetermined.
\end{proof}

\begin{corollary}\label{th:fourierexception}
For every $0<\alpha<1$ there are $R_\alpha>1,\ b_\alpha>0$
such that each sufficiently large $[X,R_\alpha X]$ contains an
interval of length $X^\alpha$ on which $D(t)<-b_\alpha$ throughout,
and another of that length on which $D(t)>b_\alpha$ throughout.
Moreover RH is equivalent to the existence of a fixed $\delta>0$
such that
\begin{equation}
\begin{gathered}
 \int_X^{2X}\mathbf1_{\{V(t)<At-0.370\sqrt t\}}\,\dd t
                  =O(X^{1-\delta})\\
 (X\to\infty).
\end{gathered}
 \label{eq:fourierexception}
\end{equation}
The same equivalence holds for the exceptional set
$\{t:D_{\mathrm{cyc}}(t)>0.370\}$.
\end{corollary}
\begin{proof}
Let $\Theta=\sup_\rho\re(\rho)$. The truncated explicit formula and
partial summation as in Theorem~\ref{th:cycwindow} give
$|V'(x)-A|\ll x^{\Theta-1}\log^2(2x)$ off the corners if
$\Theta<1$. When $\Theta=1$ Chebyshev's bound gives $O(1)$.
In the off-line case of Theorem~\ref{th:fixedfourierwindow} choose
$\max(1/2,\Theta+\alpha-1)<\sigma<\Theta$, avoiding zero real
parts. Its selected zero satisfies $\beta_*>\Theta+\alpha-1$.
In the RH case take $\beta_*=\Theta=1/2$, which has the same
strict inequality. The two raw deviations $At_\pm-V(t_\pm)$
have opposite signs and absolute values $\gg X^{\beta_*}$.
Over distance $X^\alpha$ the change in the raw deviation is
$O(X^{\Theta+\alpha-1}\log^2 X)=o(X^{\beta_*})$.
Enlarge the fixed multiplicative endpoint by one to include
the two full intervals. In the off-line case $D$ on the positive
interval is in fact $\gg X^{\beta_*-1/2}$.

On RH both exceptional sets in the statement are empty on a tail,
by \eqref{eq:fourierlow} and the rounding estimate. Conversely,
if RH fails, fix $\alpha$ with $\max(0,1-\delta)<\alpha<1$.
The positive interval just constructed contributes at least
$X^\alpha$ to either exceptional set inside $[X,R_\alpha X]$.
This last interval is covered by a fixed finite list of dyadic
intervals $[2^jX,2^{j+1}X]$. The proposed estimate gives a total
exceptional length $O_{\alpha,\delta}(X^{1-\delta})$ there,
contradicting $X^\alpha$. This proves the equivalence.
\end{proof}
For $64\le X\le7\cdot10^{18}$ the integral
\eqref{eq:fourierexception} is exactly zero by the pointwise bound.
For arbitrary finite $X$, it also has an exact evaluation using
only prime powers up to $16X$. Partition $[X,2X]$ at the numbers
$p^j/\ell$ it contains for $\ell=1,2,4,8$. On each resulting
interval $[l,r]$, the weight \eqref{eq:stepweight} gives
$V(t)=a t+d$. Set $y=\sqrt t$ and take the intervals $(u,v)$
where $(A-a)y^2-0.370y-d>0$ inside $[\sqrt l,\sqrt r]$.
The contribution to \eqref{eq:fourierexception} is exactly
$\sum(v^2-u^2)$. Endpoints have measure zero. The open estimate
concerns the sum of these lengths as $X$ tends to infinity.

\head{Frequency resolution on short logarithmic intervals}
For any probability on $[-\delta/2,\delta/2]$ its Laplace
multiplier $M_\delta$ satisfies, when $|\gamma|\delta\le\pi$,
\[
 |M_\delta(a+i\gamma)|
 \ge e^{-|a|\delta/2}\cos(|\gamma|\delta/2).
\]
The bound follows by taking the real part under the integral.
Annihilating a frequency $\gamma\ne0$ by positive weights requires
logarithmic length at least $\pi/|\gamma|$. Two endpoint masses
of equal exponentially weighted size and opposite phase attain this
length. For the window $\delta=\log(1+X^{-1/4})$, the phase
excursion is at most $e^{-w/2000}$ for every $|\gamma|\le T_X$.

In particular let
\[
 B_X(t)=2\re\!\sum_{\substack{\rho=\beta+i\gamma,\ \beta>1/2\\
                               H<\gamma\le T_X}}
                          K(\rho)t^{\rho-1/2}.
\]
By \eqref{eq:fouriertail} and \eqref{eq:fourierlow}, for
$t\in[X,X+h]$ one has
$D_{\mathrm{cyc}}(t)\le B_X(t)+0.280828+O(e^{-w/500})$.
Thus the purely spectral statement
\begin{equation}
\begin{gathered}
 \min\{B_X(X),B_X(X+X^{3/4})\}\le0.08916\\
 \hbox{for every sufficiently large }X
\end{gathered}
 \label{eq:fourieropen}
\end{equation}
would imply \eqref{eq:cyctarget} with $C=0.370$, hence RH. On RH
this sum is empty. A fixed off-line mode has size proportional to
$t^{\beta-1/2}$ and its interpolation error on the window is only
$O_\rho(X^{\beta-1})$. It therefore obeys the curvature estimate
even as its maxima grow. An independent use of this argument would require
a bound for the signed sum in \eqref{eq:fourieropen} at one of the
two endpoints.

\begin{theorem}\label{th:fourierobstruction}
Suppose $\rho_0=\beta+i\gamma_0$ is a zero with
$\beta>1/2,\ \gamma_0>0$, of multiplicity $m_0$. Then
\begin{equation}
 \limsup_{X\to\infty}\min_{X\le t\le X+X^{3/4}}
       \frac{\pm D(t)}{t^{\beta-1/2}}\ \ge\ m_0|K(\rho_0)|>0
 \label{eq:fourieromega}
\end{equation}
for each sign separately. If there is a zero farther to the right,
both upper limits in \eqref{eq:fourieromega} are infinite.
\end{theorem}
\begin{proof}
Fix either sign and put $f(t)=\pm(At-V(t))$. If its upper limit
in \eqref{eq:fourieromega} is finite or equals $-\infty$, then for
every $b$ larger than that limit and each large $X$, some point
in the window satisfies $f(t)\le b\,t^\beta$.
Since $f'=O(1)$, this gives an eventual one-sided bound
$f(X)\le C X^{\max(\beta,3/4)}$. The nonnegative-tail Mellin
argument of Theorem~\ref{th:cycwindow} works for either sign and
excludes zeros to the right of $\max(\beta,3/4)$.
If $\beta<3/4$, the resulting bound
$f'(X)=O(X^{-1/4}\log^2(2X))$ away from corners gives
$f(X)\le C' X^\beta$, because $X^{1/2}\log^2(2X)=o(X^\beta)$.
A second application excludes all zeros to the right of $\beta$.
This proves the assertion when a farther zero exists.

In the remaining case \eqref{eq:fourierdeficit} gives
\[
\begin{gathered}
 e^{-\beta u}f(e^u)=S_\beta(u)+o(1),\\
 S_\beta(u)=\pm\sum_{\re(\rho)=\beta}K(\rho)e^{i(\im\rho)u}.
\end{gathered}
\]
The terms strictly to the left tend to zero by domination by
$\sum|K(\rho)|$. Absolute uniform convergence makes $S_\beta$
real and uniformly continuous. The logarithmic length of the
window tends to zero, hence $S_\beta(u)\le b+o(1)$ on a tail.
Its time mean is zero and its coefficient at $\gamma_0$ is
$\pm m_0K(\rho_0)$, by termwise integration of the series against
$e^{-i\gamma_0 u}$ over long intervals. For each $\varepsilon>0$
the function $b+\varepsilon-S_\beta$ is eventually nonnegative.
The absolute value of any of its Fourier coefficients is at most
its mean, by the triangle inequality under the time integral.
Thus $m_0|K(\rho_0)|\le b+\varepsilon$. Let $\varepsilon\downarrow0$
and take the infimum over $b$. The numerator of $K$ has no zero
in $0<\re s<1$.
\end{proof}
Thus the prime sum determines the supremum of the real parts of the
zeros. With $[u]_+=\max(u,0)$, put
$M_\pm(X)=[\min_{X\le t\le X+X^{3/4}}(\pm D(t))]_+$. Then
\[
 \Theta=\tfrac12+
   \limsup_{X\to\infty}\frac{\log(1+M_\pm(X))}{\log X}
\]
for either sign. The upper bound follows by taking absolute values in
\eqref{eq:fourierdeficit} and using $\beta\le\Theta$. For
$\Theta>1/2$ the reverse bound follows from \eqref{eq:fourieromega}
at any zero to the right of $1/2$, then by taking the supremum of
their real parts. When $\Theta=1/2$ the logarithm is nonnegative
and bounded.

There is also a prime-power recurrence for integer indices. Set
$d_n=D(n^4)$ and $\delta_n=d_{n+2}-2d_{n+1}+d_n$. Applying
\eqref{eq:fourierchord} at $t=(n+v)^4$, $0\le v\le2$, gives
$\delta_n=O(e^{-w(n^4)/500})$ and bounds the interpolation error
on $0\le v\le1$ by the same quantity. At low
ordinates the second derivative in $v$ of each term is
$O((n^4)^{\beta-1})$, and the high-ordinate bound in the proof is
uniform on $[n^4,2n^4]$. Consequently RH is equivalent to eventual
boundedness from above of the explicit prime-power sequence $d_n$.
Under RH this follows from \eqref{eq:fourierlow}. In the other
direction interpolation bounds $D$ above everywhere on a tail,
and Theorem~\ref{th:cycwindow} applies. The exact induction is
\[
\begin{gathered}
 d_N=d_m+(N-m)(d_{m+1}-d_m)
           +\sum_{j=m}^{N-2}(N-1-j)\delta_j,\\
 N\ge m+1 .
\end{gathered}
\]
Its step follows by subtraction of consecutive identities. The
absolute majorant for $\delta_j$ above decays more slowly than
any inverse power of $j$. The sign of the accumulated second
differences therefore matters, even after controlling all local
interpolation errors.

\end{document}